\documentclass{amsart}

\usepackage[hmarginratio={1:1}]{geometry}
\usepackage[shortlabels]{enumitem}
\usepackage[colorlinks,allcolors=blue!50!black,bookmarksdepth=2]{hyperref}
\usepackage[style=alphabetic,
urldate=iso,
seconds=true, 
doi=true,
isbn=false,
maxnames=99,
maxalphanames=99,
backref=true]{biblatex}

\AtEveryBibitem{%
  \clearfield{eprintclass}%
  \iffieldundef{eprint}
    {}
    {\clearfield{url}}%
  \iffieldundef{doi}
    {}
    {\clearfield{url}}%
}

\usepackage{amssymb}
\usepackage{tikz-cd}
\usepackage{makecell}
\usepackage{xpatch}
\usepackage{graphicx}
\usepackage{eurosym}
\usepackage{mdframed}
\usepackage[most]{tcolorbox}
\usepackage{etoc}
\usepackage{xstring}
\usepackage{xifthen}
\usepackage{mathtools}

\newcommand\abs[2][]{#1\lvert{#2}#1\rvert}
\newcommand\alttext{\texorpdfstring}
\newcommand\AND{\text{ and }}
\newcommand\ari[1][]{\ar["\sim"', sloped, #1]}
\newcommand\Ar{\operatorname{Ar}}
\newcommand\bC{\mathbf{C}}
\newcommand\bD{\mathbf{D}}
\newcommand\bF{\mathbf{F}}
\newcommand\br[1]{\langle{}#1\rangle{}}
\newcommand\Cat{\mathcal{C}at}
\newcommand\cC{\mathcal{C}}
\newcommand\cF{\mathcal{F}}
\newcommand\clo{\operatorname{clo}}
\newcommand\Cls{\mathrm{Cls}}
\newcommand\cL{\mathcal{L}}

\newcommand\cod{\operatorname{cod}}
\newcommand\colim{\operatorname{co{\lim}}}
\newcommand\comp{\operatorname{comp}}
\newcommand\cP{\mathcal{P}}
\newcommand\cT{\mathcal{T}}
\newcommand\cU{\mathcal{U}}
\newcommand\cV{\mathcal{V}}
\newcommand\defeq{\mathrel{\vcenter{\baselineskip0.5ex \lineskiplimit0pt \hbox{.}\hbox{.}}} =}
\newcommand\defiff{\mathrel{\vcenter{\baselineskip0.7ex \lineskiplimit0pt \hbox{.}\hbox{.}}} \Leftrightarrow}
\newcommand\defword[1]{{\boldmath \textbf{#1}}}
\newcommand\DOF{\operatorname{DOF}}
\newcommand\dom{\mathrm{dom}}
\newcommand\dwnfn{{\downarrow}}
\newcommand\dwn{\mathpunct{\downarrow}}
\newcommand\Dwn{{\Downarrow}}
\newcommand\Ein{\mathrel{\text{\textup{\euro}}}}
\newcommand\ein{\mathrel{\varepsilon}}
\newcommand\El{\operatorname{El}}
\newcommand\enumbelow{\ }
\newcommand\ev{\mathrm{ev}}
\newcommand\floor[2][]{#1\lfloor #2 #1\rfloor}
\newcommand\FORALL{\text{ for all }}
\newcommand\fP{\mathfrak{P}}
\newcommand\fS{\mathfrak{S}}
\newcommand\fund{\mathrm{fund}}
\newcommand\Fun{\operatorname{Fun}}
\newcommand\gpd{\mathrm{gpd}}
\newcommand\hsSet{\mathrm{hsSet}}
\newcommand\hs{\mathrm{hs}}
\newcommand\hto{\hookrightarrow}
\newcommand\ibar{\bar{\imath}}
\newcommand\id{\mathrm{id}}
\newcommand\IFF{\text{ if and only if }}
\newcommand\IF{\text{ if }}
\newcommand\IMPLIES{\text{ implies }}
\newcommand\im{\operatorname{im}}
\newcommand\induc{\operatorname{induc}}
\newcommand\iso{\mathrm{iso}}
\newcommand\I{^{-1}}
\newcommand\len{\operatorname{len}}
\newcommand\loccitnodot{\textit{loc.~cit}}
\newcommand\loccitspace{\textit{loc.~cit.}\ }
\newcommand\loccit{\textit{loc.~cit.}}

\newcommand\MKo{\mathrm{MK}_{\omega}}
\newcommand\Mon{\operatorname{Mon}}
\newcommand\Nv{\operatorname{Nv}}
\newcommand\N{\mathbb{N}}
\newcommand\Ob{\operatorname{Ob}}
\newcommand\ol[1]{\overline{#1}}
\newcommand\oot{\twoheadleftarrow}
\newcommand\opcitspace{\textit{op.~cit.}\ }
\newcommand\opcit{\textit{op.~cit.}}
\newcommand\op{\mathrm{op}}
\newcommand\oring[1]{\mathring{#1}}
\newcommand\ot{\leftarrow}
\newcommand\pb[1][rd]{\ar[#1, phantom, "\lrcorner", pos=0]}
\newcommand\pbig[1]{\big({#1}\big)}
\newcommand\pBig[1]{\Big({#1}\Big)}
\newcommand\pbr{\ar[dd, phantom, "\lrcorner" {pos=0, rotate=-45}]}
\renewcommand\phi{\varphi}
\newcommand\po[1][lu]{\ar[#1, phantom, "\ulcorner", pos=0]}
\newcommand\pow{\mathcal{P}}
\newcommand\qefed{=\mathrel{\vcenter{\baselineskip0.5ex \lineskiplimit0pt \hbox{.}\hbox{.}}}}
\newcommand\Rel{\operatorname{Rel}}
\newcommand\Rrel{\mathrel{R}}
\newcommand\rstr[2]{{\left.#1\right|_{#2}}}
\newcommand\rS{\mathrm{S}}

\newcommand\rt{\mathrm{t}}
\newcommand\rV{\mathrm{V}}
\newcommand\set[2][]{#1\{ #2 #1\}}
\newcommand\Set{\mathbf{Set}}
\newcommand\SOF{\operatorname{SOF}}
\newcommand\sqbig[1]{\big[ #1 \big]}
\newcommand\SSOF{\mathcal{SOF}}
\newcommand\St{\mathrm{St}}
\newcommand\Sub{\operatorname{Sub}}
\newcommand\suc{\mathrm{suc}}
\newcommand\tc{\operatorname{tc}}
\newcommand\THEN{\text{ then }}
\newcommand\tm{\mathbf{1}}
\newcommand\tocell{\hspace{-5pt}\begin{tikzcd}[column sep=15pt, ampersand replacement=\&]{}\ar[r, shift left=3pt, ""{name=top}]\ar[r, shift right=3pt, ""'{name=bot}, start anchor=east, end anchor=west]\&{}\ar[Rightarrow, from=top, to=bot, shorten=0.5pt]\end{tikzcd}\hspace{-5pt}}
\newcommand\toix[1]{\xrightarrow[{\raisebox{3.5pt}[0pt][0pt]{\ensuremath{\scriptstyle{\sim}}}}]{#1}}
\newcommand\Toi{\xRightarrow{\raisebox{-1.5pt}[0pt][0pt]{\ensuremath{\scriptstyle{\sim}}}}}
\newcommand\toi{\xrightarrow{\raisebox{-2.5pt}[0pt][0pt]{\ensuremath{\scriptstyle{\sim}}}}}
\newcommand\tolon{{\,:\,}}
\newcommand\toot{\leftrightarrows}
\newcommand{\toox}[2][]{\xrightarrow[#1]{#2}\mathrel{\mkern-14mu}\rightarrow}
\newcommand\too{\twoheadrightarrow}
\newcommand\topp{\mathrm{top}}
\newcommand\toto{\rightrightarrows}
\newcommand\tot{\leftrightarrow}
\newcommand\ToT{\Leftrightarrow}
\newcommand\tox[2][]{\xrightarrow[#1]{#2}}
\newcommand\To{\Rightarrow}
\newcommand\tprod{{\textstyle \prod}}
\newcommand\tto{\rightarrowtail}
\newcommand\und[1]{\underline{#1}}
\newcommand\unein{{\ein}}
\newcommand\unEin{{\Ein}}
\newcommand\unex{{!}}
\newcommand\unin{{\in}}
\newcommand\univ{\mathrm{univ}}
\newcommand\unsim{{\sim}}
\newcommand\up{\uparrow}
\newcommand\vA{\vek{A}}
\newcommand\va{\vek{a}}
\newcommand\vb{\vek{b}}
\newcommand\vB{\vek{B}}
\newcommand\vek[1]{{\vec{#1}}}
\newcommand\vt{\vek{t}}
\newcommand\vx{\vek{x}}
\newcommand\vX{\vek{X}}
\newcommand\vy{\vek{y}}
\newcommand\wh[1]{\widehat{#1}}
\newcommand\WITH{\text{ with }}
\newcommand\wt[1]{\widetilde{#1}}
\newcommand\xot[2][]{\xleftarrow[#1]{#2}}
\newcommand\zer{\mathrm{zer}}
\newcommand\ZFo{\mathrm{ZF}_\omega}
\newcommand\zset{\mathbf{0}}

\DeclareFontFamily{U}{mathx}{\hyphenchar\font45}
\DeclareFontShape{U}{mathx}{m}{n}{
      <5> <6> <7> <8> <9> <10>
      <10.95> <12> <14.4> <17.28> <20.74> <24.88>
      mathx10
      }{}
\DeclareSymbolFont{mathx}{U}{mathx}{m}{n}
\DeclareFontSubstitution{U}{mathx}{m}{n}
\DeclareMathAccent{\widecheck}{0}{mathx}{"71} 

\makeatletter
\newcommand\prooflabel[1]{\gdef\proof@label{#1}}
\def\proof@label{}
\newtheoremstyle{plainwithref}%
  {}{}
  {\itshape}{}
  {\bfseries}{.}
  {.5em}{%
    \expandafter\IfBlankTF\expandafter{\proof@label}%
    {\thmname{#1}\thmnumber{ #2}}%
    {\hyperref[\proof@label]{\thmname{#1}\thmnumber{ #2}}\gdef\proof@label{}}\thmnote{ (#3)}%
  }
\makeatother

\makeatletter
\newcommand\deftag[2]{%
  \expandafter\newcommand\csname tag@#1\endcsname{#2}%
}
\newcommand\tagword[1]{\csname tag@#1\endcsname}
\makeatother

\makeatletter
\def\label@tag#1:#2\@nil{#1}
\newcommand\labeltag[1]{\label@tag#1:\@nil}
\makeatother

\newcommand\autothmwithsec[2]{\newtheorem{#1}{\tagword{#1}}[#2]}
\newcommand\autothmwiththm[2]{\newtheorem{#1}[#2]{\tagword{#1}}}

\newcommand\jref[1]{%
  \texorpdfstring%
    {\hyperref[#1]{\tagword{\labeltag{#1}}~\ref*{#1}}}%
    {\tagword{\labeltag{#1}}~\ref*{#1}}%
}
\newcommand\jrefs[2]{%
  \ifthenelse{\equal{\labeltag{#1}}{\labeltag{#2}}}%
  {%
    \texorpdfstring%
      {\tagword{\labeltag{#1}}s~\ref{#1}~and~\ref{#2}}%
      {\tagword{\labeltag{#1}}s~\ref*{#1}~and~\ref*{#2}}%
  }{\jref{#1} and \jref{#2}}%
}

\newcommand\jrefss[3]{%
  \ifthenelse{\equal{\labeltag{#1}}{\labeltag{#2}}}%
  {%
    \ifthenelse{\equal{\labeltag{#1}}{\labeltag{#3}}}%
    {%
      \texorpdfstring%
      {\tagword{\labeltag{#1}}s~\ref{#1},~\ref{#2},~and~\ref{#3}}%
      {\tagword{\labeltag{#1}}s~\ref*{#1},~\ref{#2}~,and~\ref*{#3}}%
    }{\jrefs{#1}{#2} and \jref{#3}}%
  }{%
    \ifthenelse{\equal{\labeltag{#2}}{\labeltag{#3}}}%
    {\jref{#1} and \jrefs{#2}{#3}}%
    {\jref{#1}, \jref{#2}, and \jref{#3}}%
  }%
}

\newcommand\prerefpost[3]{%
  \texorpdfstring%
    {\hyperref[#2]{#1\ref*{#2}#3}}%
    {#1\ref*{#2}#3}%
}
\newcommand\sref[1]{\prerefpost{\S}{#1}{}}
\newcommand\aref[1]{\prerefpost{Appendix~}{#1}{}}
\newcommand\pref[1]{\hyperref[#1]{p.~\pageref*{#1}}}

\deftag{defn}{Definition}
\deftag{rmk}{Remark}
\deftag{exm}{Example}
\deftag{thm}{Theorem}
\deftag{cor}{Corollary}
\deftag{propn}{Proposition}
\deftag{lem}{Lemma}
\deftag{question}{Question}
\deftag{claim}{Claim}

\theoremstyle{definition}

\autothmwithsec{defn}{subsection}

\autothmwiththm{rmk}{defn}
\autothmwiththm{exm}{defn}

\theoremstyle{remark}

\theoremstyle{plainwithref}

\autothmwiththm{thm}{defn}
\autothmwiththm{cor}{defn}
\autothmwiththm{propn}{defn}
\autothmwiththm{lem}{defn}
\autothmwiththm{question}{defn}

\autothmwithsec{claim}{subsection}

\newenvironment{settybox}{\begin{tcolorbox}[colback=blue!20,breakable%
    ,enhanced%
    ,sharp corners%
    ,frame hidden%
    ,boxsep=0pt
    ,left=5pt,right=5pt,top=5pt,bottom=5pt
    ,enlarge left by=-5pt,width=\linewidth+10pt
    ,enlarge top by=5pt,enlarge bottom by=5pt
    ,before skip balanced=0pt,after skip balanced=0pt
    ]}{\end{tcolorbox}}

\xpatchcmd{\proof}{\itshape}{\phantomsection\bfseries}{}{}

\title{Set-theoretic universes and paradoxes in 2-topoi}
\author{Joseph Helfer}
\date{}

\begin{document}

\begin{abstract}
  This paper continues the development of (elementary) 2-topos theory, a foundational theory based on an axiomatization of the 2-category of categories.
  We prove that any 2-topos contains a model of intuitionistic ZF set theory, and we use this to show that, if appropriate size restraints are not imposed, then a Burali-Forti type paradox can be deduced from the 2-topos axioms.
  The set-theoretic universe is produced as a special case of a general construction giving an internal category of models in a 2-topos of an arbitrary finite higher-order theory, which in turn is carried out using a notion of ``topos sketch''.
  As a byproduct of our construction, we also make contact with the subject of ``algebraic set theory'', introducing a novel approach to the construction of set-theoretic universes from a ``category of classes''.
\end{abstract}

\maketitle

\etocsettocstyle{\hfill\textsc{\contentsname}\hfill\,}{}

\etocsettagdepth{proofs}{section}
\tableofcontents

\section{Introduction}
The notion of (elementary) 2-topos was introduced in \cite{weber-2-toposes} as an axiomatization of the key elementary properties of \( \Cat \), the 2-category of categories, just as the concept of 1-topos does for \( \Set \), the 1-category of sets.
In \cite{helfer-topoi-in-topoi}, we further developed the concept of 2-topos, in particular proposing certain natural further axioms, though we did not suggest there that we thus arrive at a satisfactory final definition of the concept, which we regard as still being in development.

The central player in a 2-topos \( \cC \) is the \emph{discrete opfibration classifier}, or \emph{DOF classifier}, \( \rS \in \cC \), which plays an analogous role to that of the subobject classifier \( \Omega \) in a 1-topos, and which in the archetypical 2-topos \( \Cat \) is given by \( \Set \in \Cat \).
The main result of \cite{helfer-topoi-in-topoi} is that, as a consequence of our proposed axioms, the DOF classifier \( \rS \) is itself an \emph{internal 1-topos}, in a sense we make precise there.

A marked difference between \( \Omega \) and \( \rS \) is that, while the universal property of \( \Omega \) determines it uniquely up to isomorphism, \( \rS \) is not determined, even up to equivalence, by its universal property.
Intuitively, the reason is that there is, in fact, not a single ``category of sets'', but many of them, each depending on a choice of Grothendieck universe.
This circumstance manifests itself in the definition of DOF classifier in the following way.
The DOF classifier \( \rS \) is equipped with a \emph{universal DOF} \( p \colon \rS_* \to \rS \), and for each \( A \in \cC \), there is a natural functor \( \El^A \colon \cC(A,\rS) \to \DOF(A) \)---where \( \DOF(A) \) is the category of DOFs in \( \cC \) over \( A \)---obtained by pulling back \( p \); in the archetypical case \( \Set \in \Cat \), the functor \( \El^A \) is precisely the ``category of elements'' construction.
This functor is not required to be an \emph{equivalence}, as one might expect from the name ``DOF classifier'', but rather only \emph{fully faithful}.
Indeed, in the case of \( \Set \in \Cat \), the essential image of \( \El^A \) does not consist of \emph{all} DOFs over \( A \), but only those with \emph{small} fibres.
In a general 2-topos, we can thus \emph{define} a DOF over \( A \) to be ``small''---or more precisely, ``small relative to \( \rS \)''---if it is in the essential image of \( \El^A \).

This naturally raises the question of what happens if we \emph{do} demand that \( \El^A \) is an equivalence for every \( A \).
As we point out in \cite{helfer-topoi-in-topoi}, this is akin to the assumption in class-set theory that every class is a set, and we would therefore expect it to lead to inconsistencies.
The first goal of this paper is to prove that this is so:
{
  \renewcommand\thedefn{A}
  \begin{thm}[paraphrased]
    If \( \cC \) is a 2-topos, and there exists a DOF classifier \( \rS \in \cC \) such that every DOF in \( \cC \) is \( \rS \)-small, then \( \cC \) is equivalent to the trivial 2-category.
  \end{thm}
}

To make this statement more precise, we must say what we mean by ``2-topos''.
We recall from \cite{helfer-topoi-in-topoi} that our ``minimal'' proposed set of 2-topos axioms---those needed to prove that the DOF classifier \( \rS \) is an internal 1-topos---are that the 2-category \( \cC \) is \emph{groupoidant}, meaning roughly that it has finite limits and that objects in \( \cC \) have well-behaved \emph{cores} (or \emph{maximal subgroupoids}), and that \( \rS \) is \emph{plentiful}, meaning that the \( \rS \)-small DOFs satisfy certain closure conditions (see \sref{subsec:2-topoi-prelims} for the precise definitions).

We then have:
{
  \renewcommand\thedefn{A}
  \begin{thm}[see \jref{cor:2-topos-world-explodes}] \label{thm:paradox}
    If \( \cC \) is a groupoidant 2-category, and there exists a DOF classifier \( \rS \in \cC \) such that every SOF in \( \cC \) is \( \rS \)-small, then \( \cC \) is equivalent to the trivial 2-category.
  \end{thm}
}

\noindent
The appearance of \emph{SOFs} (``setoidal opfibrations'') in the statement is a technicality related to the fact that we are working, for the sake of convenience, in a strict 2-categorical setting; see \sref{subsec:2-topoi-as-foundation} below, and also the discussion in \cite[\S\S1.3~and~4]{helfer-topoi-in-topoi}.

The second goal of this paper is to show that we can construct a set-theoretic universe inside of any 2-topos:
{
  \renewcommand\thedefn{B}
  \begin{thm}[paraphrased]
    Any 2-topos \( \cC \) contains an internal model for intuitionistic ZF set theory.
  \end{thm}
}

Again, we should specify what we mean by ``2-topos''.
Here, we take the same ``minimal'' set of axioms as before, but now we also need the ``axioms of universes'' (UA) from \cite[\S4]{helfer-topoi-in-topoi} (see also \sref{subsec:2-topoi-prelims}), which roughly says that \emph{every} DOF is small with respect to \emph{some} DOF classifier.
We should also say what is meant by ``internal model''.
We mean that \( \cC \) contains a plentiful DOF classifier \( \rS \) such that the 1-topos \( \cC(\tm,\rS) \) contains an internal model of ZF: an object \( \cV \) equipped with a binary relation \( \unein \tto \cV \times \cV \), such that \( (\cV,\unein) \) satisfies the axioms of ZF in the usual sense of categorical logic.
By \cite[Theorem~5.4.7]{helfer-topoi-in-topoi}, this is equivalent to saying that the 1-topos \( \DOF(\tm) \) of DOFs over the terminal object \( \tm \) contains such an object \( \cV \).

The conclusion is also slightly stronger than what we wrote above: we in fact obtain a model of intuitionistic Morse-Kelley class-set theory (see \S\ref{subsubsec:relation-to-zf}).
We thus have:
{
  \renewcommand\thedefn{B}
  \begin{thm}[see \jrefs{thm:model-of-mk}{thm:infinity}] \label{thm:universe}
    Let \( \cC \) be a groupoidant 2-category satisfying (UA).
    Then the 1-topos \( \DOF(\tm) \) contains an internal model for MK.
  \end{thm}
}

\noindent
The relationship between these two goals is that we will prove \jref{thm:paradox} by reproducing the \emph{Burali-Forti paradox} in the set-theoretic universe \( \cV \) provided by \jref{thm:universe}.
As an easy consequence of \jref{thm:universe}, we also settle a matter which was not discussed in \cite{helfer-topoi-in-topoi}: whether the internal topos \( \rS \) has a natural numbers object (NNO).

{
  \renewcommand\thedefn{C}
  \begin{cor}[see \jref{thm:infinity}]
    Let \( \cC \) be a groupoidant 2-category satisfying (UA).
    Then the 1-topos \( \DOF(\tm) \) has an NNO.
  \end{cor}
}

\subsection{Abstract and concrete sets} \label{subsec:abstract-concrete}
The subject of producing models of set theory from topoi is an old one, starting with the papers \cite{cole-categories_of_sets,mitchell-boolean-topoi}, and the proof of \jref{thm:universe} is based on the same basic construction that appears in those papers, and which has been further developed in various directions since then.
Let us call this general construction the ``structural approach''; we will review it now.
But before we do, let us mention that there is a \emph{second} basic approach, beginning in \cite{fourman-sheaf-models,hayashi-set-theories-in-toposes}, which we call the ``universe approach''.
This has likewise been subject to further developments, notably including the notion of \emph{categories of classes} of \cite{joyal-moerdijk-ast}, which, as we will see below, also plays an important, though in a sense secondary, role in our construction.

As is emphasized in \cite{shulman-comparing,makkai-abstract-sets}, the construction underlying the ``structural approach'' is at root really more generally about comparing the notions of \emph{concrete} and \emph{abstract set}\footnote{ %
  In \cite{shulman-comparing}, the terms ``material set theory'' and ``structural set theory'' are used for these concepts. %
}---the former being the ``usual'' notion of set, axiomatized in ZF, and the latter being the notion which, upon reflection, is the one being axiomatized in topos theory, as was made explicit in \cite{lawvere-variable} and further elaborated upon in \cite{makkai-abstract-sets}.
We refer the reader to those texts for a thorough discussion of abstract sets, but briefly: whereas a concrete set is \emph{formed} out of other, pre-existing entities, an abstract set is simply \emph{given} to us, so to speak, together with its elements---which have no independent existence apart from \emph{being} elements of that set---and its equality relation.

It is clear how we can model abstract sets using concrete sets, since abstract sets are indeed an \emph{abstraction} of concrete sets.
That is, we can simply \emph{forget}, or ignore, the concrete nature of the elements of concrete sets; this is more or less precisely what we are doing when we form the category of sets, and proceed to only reason about it using categorical constructions---or, somewhat more precisely, when we regard the category of sets as a model for the theory of elementary topoi.
It is much less clear how to model concrete sets using abstract sets, and indeed, it seems at first glance that this might not be possible---but it is possible, based on the notion of \emph{transitive closure}, which serves as the basis for the ``structural approach'' mentioned above.

Before describing it, we must distinguish a particular important, and familiar, class of concrete sets, namely the \emph{pure}, or \emph{hereditary}, sets.
These are those concrete sets all of whose elements are sets---and in fact, all of whose elements are \emph{hereditary} sets, meaning that \emph{their} elements are (hereditary) sets, and so on.
Ordinarily, in ZF set theory or (our preferred setting) MK class-set theory (see \sref{subsec:axioms-of-mk}), one assumes---in the form of the axiom of foundation---that \emph{every} set is hereditary.
However, if we omit this axiom, we can regard the theory as being concerned more generally with not-necessarily-hereditary concrete sets, in which case the hereditary sets can be \emph{defined} recursively by saying that the class \( \rV \) of hereditary sets is the least class \( X \) which is closed under set-formation: if \( A \subset X \) is a set, then \( A \in X \).

At any rate, it is the \emph{hereditary} sets which can be modelled, via the transitive closure construction, using abstract sets.
Still in the context of ZF or MK set theory, given a hereditary set \( x \), the \emph{transitive closure} \( \tc(\set{x}) \) is defined to be the least transitive hereditary set \( y \) with \( x \in y \), where \( y \) being transitive means that \( a \in b \in y \) implies \( a \in y \) for all sets \( a \) and \( b \).
We equip \( \tc(\set{x}) \) with the binary relation \( \unein_x \) given by the restriction to \( \tc(x) \) of the element relation: for \( a,b \in x \), we have \( a \ein_x b \) iff \( a \in b \).
We thus obtain a structure \( \pbig{\tc(\set{x}), \unein_x} \) consisting of a set with a binary relation.
The first fundamental fact about the transitive closure operation is:
\[
  \text{given hereditary sets } x \AND y,
  \IF \pbig{\tc(\set{x}), \unein_x} \cong \pbig{\tc(\set{y}), \unein_y}
  \THEN x = y.
\]
In other words, a hereditary set is uniquely determined by the isomorphism type of its transitive closure.
The second fundamental fact is \emph{Mostowski's collapsing lemma} (\jref{propn:mostowski-collapse}), which gives sufficient and necessary conditions for an arbitrary set equipped with a binary relation to be isomorphic to one of the form \( \pbig{\tc(\set{x}), \unein_x} \) for a (by the first fundamental fact unique) hereditary set \( x \); the conditions in question are \emph{extensionality}, \emph{well-foundedness}, and \emph{existence of a top element} (see \sref{subsec:ein-sets}).
In this paper, we call structures of this kind \emph{hs-sets}.\footnote{%
  They are called ``well-founded extensional accessible pointed graphs'' in \cite{aczel-non-well-founded-sets,shulman-comparing} and ``pure set structures'' in \cite{makkai-abstract-sets} (though there is an error in the definition given there: the condition that \( \unein_x \) is \emph{transitive} is also included); in \cite{cole-categories_of_sets,mitchell-boolean-topoi,osius-categorical-set-theory,johnstone-topos-theory}, the concepts serving the same purpose as hs-sets are called ``E-sets'', ``trees'', ``transitive-set objects'' and ``transitive objects'', respectively.

  ``hs-set'' stands for ``hereditary set-set''; the reason we don't use the more reasonable name ``hs-structure'' is so as to avoid awkward phrases such as ``endowing \( X \) with an hs-structure structure'', and also to distinguish this concept from the topos-theoretic analogue, which we call ``hs-object'' (see below).
}

As a result of these two fundamental facts, the transitive closure operation induces a bijection between the class \( \rV \) of hereditary sets, and the class \( \rV_\hs \) of isomorphism classes of hs-sets.
Moreover, it is easy to characterize directly what the relation on \( \rV_\hs \) is which corresponds to the relation \( \in \) on \( \rV \) under this bijection (see \jref{defn:ein-hs-set-mor}).
The point now is that hs-sets are \emph{abstract structures} in which, as in all abstract structures, the particular nature of the concrete sets they are based on is completely irrelevant.
Thus, within abstract set theory, we can completely reproduce the talk of hereditary sets by instead talking about hs-sets.

In particular, we can talk about hs-sets (or, as we call them, \emph{hs-objects}) in any elementary topos \( \bC \), and it is by taking the isomorphism classes of these that one, in the ``structural approach'', obtains from \( \bC \) a model of (some fragment) of ZF set theory.
We note that additional assumptions are usually placed on \( \bC \) in order to ensure that the resulting model is well-behaved; in particular, it seems that in all the versions of this construction that have been considered, \( \bC \) is assumed to be \emph{well-pointed}.
It is one of the features of our version of this construction that these assumptions are relaxed (though we still can't quite say that we set out from an \emph{arbitrary} elementary topos---see \jref{question:internalizing-topoi} below).

A related issue to this is that, in the structural approach, one must pass in the end \emph{out} of the topos \( \bC \) in which one began.
That is, one considers hs-objects in \( \bC \), but the \emph{totality} of these is a set only in the sense of the ambient set theory in which the topos \( \bC \) is being considered, rather than an object internal to \( \bC \) itself.
Again, this is something that is overcome in our approach: our construction begins not with a genuine elementary topos \( \bC \), but with an \emph{internal topos} \( \rS \) inside some \emph{2-topos} \( \cC \)---and the resulting universe of hs-objects in \( \rS \) then ends up being something which is again \emph{internal to \( \cC \)}.

\subsection{Universes of sets in 2-topoi and the paradox} \label{subsec:unis-and-paradox}
We are now ready to explain the basic idea behind the proofs of \jrefs{thm:paradox}{thm:universe}.

Fix a ``2-topos'', i.e., a groupoidant 2-category \( \cC \) which contains a plentiful DOF classifier \( \rS \in \cC \), and fix such an \( \rS \).
The main theorem of \cite{helfer-topoi-in-topoi} states that \( \rS \) is an \emph{internal topos}, which means, primarily, that for every groupoidal object \( A \in \cC \), the category \( \cC(A, \rS) \) is an elementary topos (see \sref{subsubsec:dof-class-are-topoi}).
Since by definition of DOF classifier, \( \cC(A, \rS) \) is equivalent to the category \( \DOF_\rS(A) \) of \( \rS \)-small DOFs over \( A \), this is equivalent to the latter category being a topos for each groupoid \( A \).
The categories \( \DOF_\rS(A) \) tend to be easier to work with than \( \cC(A, \rS) \).

Given any topos \( \bC \), we can construct the category \( \hs(\rS) \) of hs-objects in \( \bC \) and isomorphisms between them.
The idea is now to construct the ``category \( \hs(\rS) \) of hs-objects in \( \rS \)'', which should again be an object of \( \cC \).
It should have the universal property that there is a (natural) equivalence for each groupoid \( A \in \cC \), between the category \( \cC\pbig{A,\hs(\rS)} \) and the groupoid of hs-objects in the topos \( \cC(A, \rS) \simeq \DOF_\rS(A) \) (see \sref{subsec:hs-2-topos}).
The next main result of this paper is that this is always possible; in fact, we show a much more general result, that we can construct the groupoid of models of \emph{any} finite theory formulated in higher-order logic:
{
  \renewcommand{\thedefn}{\ref{thm:theory-cotensors}}
  \begin{thm}
    For any internal topos \( E \) in a groupoidant 2-category \( \cC \) and any finite higher-order theory \( \cT \), there exists a cotensor \( (E^\cT)^\iso \) of \( E \) by \( \cT \).
  \end{thm}
}
\noindent
(See \sref{subsec:thy-cotensors} for the meaning of ``cotensor'').

The proof of this statement proceeds in two stages: we translate the given higher-order theory \( \cT \) into a \emph{topos sketch} \( J_\cT \) having the same categories of models, and we then show that cotensors \( (E^J)^\iso \) by arbitrary finite topos sketches \( J \) exist; see \sref{sec:univs-in-2-topoi} for more details.

Thus, returning to our context of a plentiful DOF classifier \( \rS \) in a groupoidant 2-category \( \cC \), we may now form the object \( \hs(\rS) \).
It turns out that this object is \emph{discrete}---or more precisely, setoidal (see \sref{subsec:hs-2-topos}).
This corresponds to the fact that the only automorphism of an hs-set, or of an hs-object in any topos, is the identity.
Thus, \emph{if we assume} that \( \rS \) is truly a \emph{universal} DOF-classifier, in the sense of \jref{thm:paradox}---or equivalently, that \( \DOF_\rS(A) = \DOF(A) \) for every object \( A \)---then \( \hs(\rS) \) is classified by a morphism \( \tm \to \rS \); that is, informally ``the class \( \hs(\rS) \) of hs-objects in \( \rS \) is itself an element of \( \rS \)''.
We see that we are close to a paradox.

To finish the argument of Burali-Forti, we really want \( \hs(\rS) \) to be an element \emph{of \( \hs(\rS) \)} (rather than just of \( \rS \)), which means that we need to endow it with a relation making it into an hs-object; intuitively, this should be the ``element relation'' between hs-sets.
Thus, we next apply \jref{thm:theory-cotensors} again to produce a second setoidal object \( \unEin_\rS \) and a monomorphism \( \unEin_\rS \to \hs(\rS) \times \hs(\rS) \) giving the desired ``element relation''; i.e., a morphism \( \br{x,y} \colon A \to \hs(\rS) \times \hs(\rS) \) from a groupoid \( A \) factors through \( \unEin_\rS \) if and only if the hs-object in \( \DOF_\rS(A) \) classified by \( x \) is an ``element'' of that classified by \( y \), in the sense of \jref{defn:ein-hs-mor} below.
The resulting structure \( \pbig{\hs(\rS), \Ein_\rS} \) in the topos \( \DOF(\tm) \) is \emph{nearly} an hs-object---it is well-founded and extensional---but it has no top element; intuitively, this is because there is no hs-set containing all other hs-sets.
We may, however, artificially add a top element, as explained in \sref{subsec:top-element}, to obtain an hs-object \( \pbig{\wh \hs(\rS), \wh{\unEin}_\rS} \).
With this hs-object in hand, we proceed with the usual Burali-Forti argument: we conclude that the element \( v \colon \tm \to \wh \hs(\rS) \) classifying the hs-object \( \pbig{\wh \hs(\rS), \wh \Ein_\rS} \) itself satisfies \( v\ \wh \Ein_{\rS}\ v \), contradicting the well-foundedness of \( \wh \unEin_{\rS} \).
The conclusion of this ``contradiction'' is that \( \mathsf{true} = \mathsf{false} \colon \tm \to \Omega \) in the topos \( \DOF_\rS(\tm) \), which implies that this topos is trivial; and from there, we argue that \( \cC \) itself is trivial.
This is done in \sref{subsec:burali-forti}.

Before moving on, let us comment on two other variants of the Burali-Forti paradox.
The first is Burali-Forti's original version (reproduced in \cite{heijenoort-frege-to-godel}), which deals not with hs-sets, but with ordinals: one shows that, assuming every class is a set, the set of all ordinals is itself an ordinal greater than every ordinal, hence greater than itself, which is a contradiction.
This version of the paradox easily goes through in the present context, once the work has already been done to establish the version in \jref{thm:paradox}, and using the appropriate intuitionistic notion of ``ordinal''.
Interestingly, though the abstract notion of ordinal (well-founded linear order) seems simpler than that of ``hs-set'', the correct intuitionistic definition of ordinal (``a transitive set of transitive sets''---see \cite{powell-extending-godels-negative}), when understood abstractly, is actually an hs-set satisfying an \emph{additional} condition.
Hence it is actually somewhat simpler to carry out the version of the argument in \jref{thm:paradox}---not to mention the important byproduct of constructing a model of hereditary set theory along the way (to which we will return shortly).

The second variant we want to mention is the one in appearing in ``Girard's paradox'' in dependent type theory, as explained in \cite{ml72}.
This version deals not with ordinals or hs-sets, but simply with the class of \emph{all} well-founded relations, and shows that this class itself admits a well-founded relation and thus contains itself.
This version is \emph{not} suited to the present context.
The reason is that it was essential for us that the universe \( \hs(\rS) \) is \emph{discrete}, so that it could be considered an element of itself.
The discreteness was due to the fact that an hs-set has no non-trivial automorphisms; however, the same is not true for a general well-founded relation.

\subsection{Topoculi} \label{subsec:intro-topoculi}
We have seen how we use the object \( \hs(\rS) \) of hs-objects in \( \rS \) to arrive at a contradiction under the too-strong assumption that every SOF is \( \rS \)-small.
We next want to explain how \( \hs(\rS) \) gives us a model of hereditary set theory when we don't make this assumption.

In the previous argument, we used that \( \hs(\rS) \) was \( \rS \)-small so that we could consider it as an object of the topos \( \DOF_{\rS}(\tm) \).
Now, instead, we assume the existence of a \emph{second} DOF classifier \( \rS' \) with respect to which \( \hs(\rS) \) is small; such an \( \rS' \) exists, for instance, if we assume axiom (UA) as in \jref{thm:universe} (and in fact, all we need for that theorem is the existence the two DOF classifiers \( \rS \) and \( \rS' \), which is a much weaker assumption than (UA)).

We thus find ourselves faced with the topos \( \DOF_{\rS'}(\tm) \), in which we have a distinguished set of morphisms, namely those which are not only \( \rS' \)-small but \( \rS \)-small.
It turns out to be convenient, and to render the arguments more transparent---both for the present purpose of producing a model of hereditary set theory, and for the Burali-Forti argument above---to axiomatize this situation, and work in a general topos together with a distinguished set of morphisms satisfying appropriate conditions.
We call such a set of morphisms a \emph{topoculus} (little topos); see \sref{subsec:topoculi}.
This is where our approach makes contact with the ``universe approach'' mentioned above, and in particular with the ``categories of classes'' of \cite{joyal-moerdijk-ast}, the idea of which we will now recall.

As mentioned above, the starting point of this approach was \cite{fourman-sheaf-models,hayashi-set-theories-in-toposes}, in which a model for hereditary set theory is constructed from a \emph{cocomplete} topos \( \bC \) (for example, a Grothendieck topos), by imitating the construction of the cumulative hierarchy: we being with an initial object \( V_0 = \zset \), iterate the power object construction to obtain objects \( V_{\alpha + 1} = \pow V_\alpha \) for successor ordinals \( \alpha \), and take colimits \( V_\lambda = \colim_{\alpha < \lambda} V_\alpha \) for limit ordinals \( \lambda \).
Unlike in the structural approach, which, as we mentioned above, produces a model which is \emph{external} to the topos \( \bC \)---viz.\ the (external) \emph{set} of hs-objects in \( \bC \)---this approach produces an object \( V \) \emph{in} \( \bC \), with a binary relation \( \unin \tto V \times V \), which is a model for ZF in the usual sense of categorical logic.

Also unlike the structural approach, however, the above construction begins with the \emph{non-elementary} input of a cocomplete topos; that is to say, in order to talk about cocompleteness, we need to have a notion of ``set'' prior to the one provided by the topos under consideration.
As a result, this construction is not well-suited to comparing the logical strength of topos theory and ZF as foundational theories---a matter we return to in \sref{subsec:2-topoi-as-foundation}.
Thus, the question remained if something like the construction of \opcitspace could be produced in an elementary context.

The ``category of classes'' framework of \cite{joyal-moerdijk-ast} approached the problem as follows.
Instead of axiomatizing just the category of sets, we instead consider the usual situation in which we make a distinction between ``large'' and ``small'' sets (i.e., classes and sets, or subsets and elements of some fixed Grothendieck universe), and axiomatize the \emph{category of classes} \( \bC \).
Thus a general object of \( \bC \) is to be considered a ``class'', and there is a distinguished collection of objects, namely the ``small'' classes, or sets---and more generally, a distinguished collection of morphisms, namely those whose fibres are small.
The universe \( V \) \emph{of sets} should thus be an object of \( \bC \), and the main purpose of \opcitspace is to describe universal properties describing \( V \) (one of which we describe in \jref{rmk:initial-algebra}), as well as to give axioms on the category \( \bC \) together with its collection of ``small maps'' which guarantee the existence of such a \( V \).

Thus, we may say that from the perspective of the ``universe approach'' to categorically constructing models of hereditary set theory, the role played by 2-topoi is to give a particularly natural source of categories of classes (namely, the category \( \DOF_{\rS'}(\tm) \) equipped with the set of \( \rS \)-small morphisms).
Since our axioms for the collection of small maps are different from the ones appearing in the literature, we chose to use a different name.
The reason for the name ``topoculus'' is that in our situation, the small objects do form a subtopos of the given topos \( \bC \), and similarly the small morphisms with each codomain \( U \in \bC \) form a subtopos of \( \bC / U \).

Thus, the main construction of \jref{thm:universe} is factored into two parts: first we show how to construct a topoculus from a 2-topos, and thereafter, most of the paper---in particular, the construction of a model of hereditary set theory, and the derivation of the paradox---takes place in an \emph{arbitrary} topoculus.
We note that, while much of the work on categories of classes has been devoted to finding the weakest possible setting in which a universe can be constructed (for example, one usually works in a general Heyting pretopos, rather than a topos), we do not take any care to make the topoculus axioms as weak as possible; we simply note those properties that we get from our construction of the topos \( \DOF_{\rS'}(\tm) \) with its collection of \( \rS \)-small morphisms, and that we need in order to carry out our arguments.
However, it would be interesting to see how far the topoculus axioms could be weakened, and more generally, to compare them to the existing axiom systems for collections of small maps; we return to this question at the end of \sref{subsec:2-topoi-as-foundation} below.

Now, we recall that our starting point was the ``structural approach'', in which we produced the universe \( \hs(\rS) \) of hs-objects in our internal topos \( \rS \in \cC \).
The way that this feeds into the ``universe approach''---which is to say, into our topoculus on the topos \( \bC = \DOF_{\rS'}(\tm) \), is that it gives us an object \( \cV \in \bC \) with a particular universal property: morphisms \( U \to \cV \) from objects \( U \in \bC \) classify \emph{small} (with respect to the topoculus) hs-objects in the slice \( \bC / U \).
We call such an object \( \cV \) an \emph{hs-classifier}.
It is a very natural assumption to make about a category of classes that it contain such an hs-classifier, and it seems not to have been considered thus far in the literature.
Our main theorem regarding topoculi is then:
{
  \prooflabel{thm:model-of-mk}
  \renewcommand{\thedefn}{~\ref*{thm:model-of-mk}}
  \begin{thm}
    Any hs-classifier \( \cV \) for a topoculus \( \fS \) on a topos \( \bC \) gives a model of MK in \( \bC \).
  \end{thm}
}

The paper \cite{shulman-stack-semantics} of Shulman also has as one of its aims to blend the structural and universe approaches to constructing models of hereditary set theory; and as we mentioned in \cite{helfer-topoi-in-topoi}, Shulman's other writings \cite{shulman-nlab-page} on 2-topoi are related in various ways to the that paper and the present one.
In particular, many of the arguments in this paper are clearly related to the \emph{stack semantics} introduced in \cite{shulman-stack-semantics}, though not in an explicit way; we will return to this point in \sref{subsec:sets-and-topoi} below.

To end the discussion of our approach, let us return to the comment made in \sref{subsec:abstract-concrete} that the structural approach normally begins from the assumption that the given topos \( \bC \) out of which we construct a universe of hs-objects is \emph{well-pointed}, whereas we do not need to impose this, or any other, assumption on \( \bC \).
The reason this is not quite true is that we \emph{are} assuming that \( \bC \) arises as the category \( \cC(\tm, \rS) \) for some plentiful DOF classifier \( \rS \) in a 2-topos \( \cC \).
This naturally gives rise to the question---a question which was already a very reasonable one to ask in the context of \cite{helfer-topoi-in-topoi}, but which we did not ask there---of whether \emph{every} elementary topos arises in this way.
This is a ``categorification'' of the---as far as we know still open---basic question in topos theory, of whether every Heyting algebra is of the form \( \bC(\tm, \Omega) \) for some topos \( \bC \).
\begin{question} \label{question:internalizing-topoi}
  Is every elementary topos \( \bC \) equivalent to the topos \( \cC(\tm, \rS) \) for some plentiful DOF classifier \( \rS \) in some groupoidant 2-category \( \cC \)?
\end{question}
\noindent
Similar questions have been considered in the context of categories of classes; see \cite{awodey-butz-simpsons-streicher-relating}.

\subsection{Set theoretic arguments in topos theory} \label{subsec:sets-and-topoi}
We wish to make some remarks on the methods of proof employed in this paper.
Once we make the passage from 2-topoi to topoculi described in the previous section, we are no longer doing 2-topos theory, but rather 1-topos theory.
In particular, we end up carrying out many of the same kinds (and simply the same) arguments as have been made in previous work on the structural approach to producing models of hereditary set theory from topoi (cf.\ the remarks in \cite[\S8]{shulman-comparing}).
However, we give all these arguments in full, for one, since we are working in a more general context than has previously been considered (e.g., no assumption of well-pointedness).

Now, seeing as topos theory is a formalization of intuitionistic abstract set theory, it has the general feature that the statements which are true in a general topos tend to coincide with those which can be proven in ``naive'' intuitionistic set theory.
Moreover, the general experience---and certainly the author's experience in writing this paper---is that, most of the time, the only reasonable way to proceed is to \emph{first} make the argument in naive (intuitionistic) set theory, and then ``translate'' the argument into topos theory; to try to make the topos-theoretic argument directly involves one in too many particular category-theoretic details and distracts from the substance of the argument.
(Of course, this is a familiar feature of \emph{any} sufficiently abstract mathematical theory; for example, one still draws pictures of arrows and parallelograms even when reasoning about a vector space of arbitrary dimension over an arbitrary field.)

Because the set-theoretic arguments are easier to understand, we often present proofs in this paper by first giving that argument, and then proceed with the topos-theoretic ``translation'' (we do not bother with this when the translation is very faithful so that the two proofs are practically identical).
However, since the set-theoretic argument is not part of the actual proof, we usually visually distinguish it by putting it in a coloured (purple) box as follows:
\begin{settybox}
  A set-theoretic argument will usually appear in a coloured box like this.
\end{settybox}
For the reader interested in getting the main gist of the proofs, without technical details, the arguments in the coloured boxes will likely suffice.
We likewise put definitions and propositions specifically related to the set-theoretic counterparts in such boxes.

The translation from the set-theoretic to the topos-theoretic argument is usually more or less automatic.
The general reason for is the well-known and fundamental interpretation of intuitionistic higher-order logic (AKA type theory; see \cite{lambek-scott}) into topoi (which we review in \sref{subsec:logic-in-topoi}).
Because of this, whenever the naive set-theoretic argument can be formalized as a proof of a certain sentence in intuitionistic higher-order logic, the interpretation of that sentence in any topos will automatically be true.

However, the translation is not \emph{completely} automatic, which is why we still bother to give the topos-theoretic proof after the set-theoretic one.
One reason for this is that even when a given argument \emph{can} be formalized in intuitionistic higher-order logic, one must still \emph{produce} this formalization, and this tends to be no easier than simply carrying out the topos-theoretic argument directly---especially in light of the \emph{Kripke-Joyal} (KJ) semantics (which we review in \sref{subsubsec:kj-semantics}), since this allows one systematically to reason with \emph{generalized elements}, and thus brings the topos-theoretic reasoning very close to naive set-theoretic reasoning (with a few twists).
In addition, it is often not the most efficient or clear approach to literally take over the higher-order logic proof step by step into the topos; there are often useful specifically topos-theoretic constructs that one can introduce, which allow one to structure the proof in a different way, and which bring out interesting features of the topos-theoretic situation.

The more fundamental reason that the translation is not always automatic, however, is that the set-theoretic argument can \emph{not} always be formalized in higher-order logic.
This happens primarily when some kind of ``quantification over sets'' is required, as this is not available in higher-order logic, where we can only quantify over the elements of one set at a time; this happens, for instance, when we are trying to prove properties about the universe \( \hs(\rS) \in \cC \) of hs-objects in an internal topos \( \rS \in \cC \), as this involves proving properties about \emph{all} hs-objects in the topoi \( \bC(U, \rS) \).
(We note that in \cite{shulman-stack-semantics}, the passage to slice categories in order to interpret statements involving quantification over sets is systematized using an extension of the KJ-semantics.
As mentioned above, we do not explicitly use the theory developed in that paper, but several of our constructions are clearly in a similar spirit.)

\subsection{2-topoi as a foundational theory} \label{subsec:2-topoi-as-foundation}
Although the role of 2-topos theory as a foundational theory is very much the primary motivation for both this paper and \cite{helfer-topoi-in-topoi}, we did not discuss it explicitly there.
Let us take a moment to do so here.

Lawvere's original paper \cite{lawvere-etcs} which led to the development of topos theory is titled \emph{An elementary theory of the category of sets}; later, he also considered an \emph{elementary theory of the category of categories} \cite{lawvere-cat-of-cats}.
The simplest and clearest description of 2-topos theory as a foundational theory is that it is an \emph{elementary theory of the 2-category of categories}.
That it really does serve as a foundational theory---in the sense that all of the usual mathematical (set-theoretic) constructions can be formulated in it is demonstrated by the fact, proven in \cite{helfer-topoi-in-topoi}, that it subsumes elementary topos theory, and the fact, proven in the present paper, that it subsumes ZF set theory.

The word \emph{elementary} here refers to the fact that the theory posits certain \emph{undefined/primitive notions}, each of which is either a \emph{type}/\emph{sort} of entity (like ``point'' or ``line'' in Euclidean geometry) or a \emph{relation} between these (like the incidence relation between points and lines).
In the present case, the sorts are the \emph{0-cells}, \emph{1-cells}, and \emph{2-cells}, and the relations are: that of a 0-cell or 1-cell being the \emph{domain/codomain} of a 1-cell or 2-cell, respectively; \emph{horizontal composition} of 1-cells; \emph{vertical composition} of 2-cells; and \emph{horizontal composition} of 2-cells (or alternatively, \emph{whiskering} of 1-cells with 2-cells on the left or the right).
No further mathematical objects of any kind need to be presupposed for the development of the theory---this is in contrast, for example, to the concept of \emph{Grothendieck topos} which, as mentioned in \sref{subsec:intro-topoculi} above, requires the concept of \emph{small cocompleteness}, and hence, the notion of \emph{set}.
In 2-topos theory, the notion of ``set'' is \emph{provided} by the theory, namely as ``1-cell \( X \to \rS \) where \( X \) is terminal and \( \rS \) is a plentiful DOF classifier''.

In addition to these primitive notions, there are then the \emph{axioms}, which are purely logical statements describing our assumptions regarding these various primitive relations.
It is an ``empirical'' fact (sometimes called ``Hilbert's thesis''; see e.g.\ \cite[\S A.1.5]{barwise-handbook}) that such purely logical statements can always be formulated in \emph{first-order logic}---and hence, ``elementary theory'' is often taken to be synonymous with ``first-order theory''---and 2-topos theory is no exception.

We will not spell out all of the 2-topos axioms (namely, the 2-category axioms, together with the assumptions of being groupoidal and having a plentiful DOF classifier) as first-order sentences, but readers habituated to formulating statements precisely with first-order logic will easily see that this can be done in each case.
As a (trivial) explicit example, let us consider the very ``first'' axiom---existence of composites of 1-cells:
\[
  \forall f,g \tolon \Ar.\ \cod(f) = \dom(g) \to
  \exists h \tolon \Ar.\ \comp(f,g,h)
\]
Here, we write \( x \tolon X \) under a quantifier to indicate that the variable \( x \) has sort \( X \), we write \( \Ar \) for the sort of 1-cells, and we write \( \comp(f,g,h) \) for the (primitive) relation ``\( h \) is the composite of \( f \) and \( g \)''.

We must now make an important comment regarding the above formalization of the axioms, regarding the use of \emph{equality}, which we, as usual, took to be a part of the language of first-order logic that we are free to use in formulating the axioms, and which we used above to specify that the objects given by the codomain of \( f \) and the domain of \( g \) are \emph{equal}.
We note that doing this violates the general principle in category theory that one should never speak of \emph{equality} of objects, but only \emph{isomorphism}---and the analogous principle regarding morphisms in 2-category theory.
This principle is in keeping with the \emph{abstract} view of sets---the natural one from the point of view of category theory---since to assert that two sets are \emph{equal} is to say that they have the \emph{same} elements---something which only makes sense with concrete sets (this point is made explicitly in \cite{makkai-categorical-foundation,makkai-fom-both,makkai-abstract-sets}).

There is a certain constrained version of first-order logic called \emph{first-order logic with dependent sorts} (FOLDS), which was introduced in \cite{makkai-folds}, which is closely related with the dependent theory of types of \cite{ml72}, and which is precisely suited to this problem of avoiding equality when axiomatizing category theory.
In FOLDS, the sort of each variable is made to to \emph{depend} in a prescribed way on other variables; thus, rather than a variable simply denoting a 1-cell, it would denote a 1-cell \emph{with prescribed domain and codomain}, each of which are variables denoting 0-cells.
In this way, we avoid needing to speak of the domain of a given arrow being \emph{equal} to a given object, since the domain of the arrow must be given before the arrow can even be introduced.
In FOLDS, the above sentence would look as follows, where we write \( \Ob \) for the sort of objects and \( \Ar(X, Y) \) for the sort of 1-cells with domain \( X \) and codomain \( Y \):
\[
  \forall X,Y,Z \tolon \Ob.\ %
  \forall f \tolon \Ar(X,Y).\ %
  \forall g \tolon \Ar(Y,Z).\ %
  \exists h \tolon \Ar(X,Z).\ %
  \comp(f,g,h).
\]
A crucial consequence of expressibility in FOLDS is the ``General Invariance Theorem'' of \cite{makkai-folds}: the sentences which are expressible in the FOLDS language for 2-categories (for which see \cite[\S7]{makkai-folds}) are exactly those which are \emph{invariant under equivalences of 2-categories}.
Thus, \emph{if} all of the 2-topos axioms were FOLDS-expressible, it would follows that any 2-category equivalent to a 2-topos is again a 2-topos.

However, it is \emph{not} the case that the 2-topos axioms, as we have given them, are so expressible.
The reason is that we are, for the sake of convenience, working in a \emph{strict} context, i.e., with \emph{strict} 2-categories, and with the strict variants of the various relevant universal properties (limits, colimits, etc.).
As we explain in \cite{helfer-topoi-in-topoi}, we are ``permitted'' to work with strict 2-categories because of the well-known theorem that every weak 2-category (AKA bicategory) is equivalent to a strict one; the ``legitimacy'' of also working with the other strict notion depends on similar results (some of which have been proven, and some not, but all of which we \emph{expect} to be true), for example to the effect that any bicategory with finite bilimits is equivalent to a strict 2-category with strict finite 2-limits.

As we also explain in \opcit, however, the use of the strict context really is only an expedient, rather than the final or ``correct'' version of the theory, which should be the ``fully weak'' one---or, as we may now say, the one in which all of the axioms are completely expressible in FOLDS.
This position is plausible from experience with category theory, in which rigid adherence to isomorphism-invariance is a familiar practice.
However, it might be objected that the intended model for the theory, namely the bicategory of categories, \emph{is} a strict 2-category, and thus, it is only reasonable to include this in the axioms.

The response to this is that, for reasons explained in \cite{makkai-categorical-foundation}, the 2-category of categories is \emph{not} in fact the ultimate intended model of the theory; it is rather the \emph{anabicategory of categories, saturated anafunctors, and natural transformations}.
Without getting into the details of this here, the point is that, if we really imagine a world in which there are only \emph{abstract} sets, the usual notion of functor---say, on the category of sets---which, for each input, \emph{chooses} a particular set, among many possible alternative isomorphic sets---for instance, in the case of the power set functor, because the chosen set is ``the'' set of subsets of the given input set---no longer makes sense.
Rather, the only sensible notion is an analogue of a functor which does not discern between the different possible isomorphic choices in the target category, but ``chooses them all''; this idea is precisely captured by the concept of saturated anafunctor.
And it is easy to see that composition of saturated anafunctors is not only not strictly associative, but is in fact only well-defined up to isomorphism (hence the ``ana'' in ``anabicategory'').

We will make one final comment regarding foundations.
\jref{thm:universe} says that we can construct a model of ZF or MK set theory ``inside'' of a 2-topos, if we assume the 2-topos satisfies axiom (UA).
Though we do not do this, it is not hard to see that this can be reinterpreted as giving an \emph{interpretation} of MK into the theory (GT)+(UA), where we write (GT) for the theory of groupoidant 2-categories---that is, a translation of the language of MK into that of 2-categories which transforms each axiom of MK into a theorem of (GT)+(UA).
This shows in a direct way that (GT)+(UA) is at least as ``strong'' as MK, i.e., the consistency of the former implies that of the latter.

As we explained above, the full axiom (UA) isn't really necessary for this interpretation, since we only need a certain weaker axiom (UU) asserting the existence of \emph{two} DOF classifiers \( \rS \) and \( \rS' \), rather than infinitely many, as is implied by (UA).
However, this is in a sense still ``more'' than we need, since the resulting model \( \cV \) of MK is an \emph{object} of the topos \( \cC(\tm, \rS') \), which thus contains all the iterated power objects of \( \cV \), and hence is a model for an extension of MK in which there are not only classes, but classes-of-classes, and so on.
The result is that there is \emph{not} an interpretation of (GT)+(UU) into MK, but only into MK plus the assumption of an inaccessible cardinal.

It would be interesting to identify precisely what the appropriate weakening of (GT)+(UU) is to obtain a theory which is bi-interpretable with MK.
In this theory, it would seem, the category \( \DOF(\tm) \) would have to be, not a topos, but only a pretopos; and the use of topoculi would perhaps be substituted for something more like the ``category of classes'' of \cite{joyal-moerdijk-ast}.

\subsection{Acknowledgments}
This paper is a continuation of the project begun in \cite{helfer-topoi-in-topoi}.
As I mentioned there, it was first suggested to be by M.~Makkai, who has maintained an active interest in the project, for which I am grateful.
I would also like to thank the NY Category Theory Seminar at CUNY Graduate Center, and the organizers of CT26, for allowing me to present on this material; a recording of and the slides for the former talk are available (at the time of writing) at \url{https://www.youtube.com/watch?v=eXJaRsE_i68}.

I would also like to thank F.~Abellán for telling me about his related work \cite{abellan-martini-2-topoi} with L.~Martini on \( (\infty,2) \)-topoi, which in particular lead me to notice a mistake in a conjecture made in \cite[\S7]{helfer-topoi-in-topoi} (this issue was also brought to my attention independently by Colin Zwanziger).
There, it is shown that the 2-category of \( \Cat \)-valued functors on any small 1-category is a 2-topos, and it was conjectured that this should be true more generally for any Grothendieck 2-topos in the sense of Street, and in particular, for the 2-category of \( \Cat \)-valued functors on any small 2-category \( \cC \).
However, this fails already for very simple 2-categories \( \cC \), such as the ``walking adjunction'' or ``walking monad'', since the resulting functor 2-categories have ``very few'' groupoids, and hence are not groupoidant.
(By contrast, such 2-categories---or rather, the analogous \( (\infty,2) \)-categories in which \( \Cat \) is replaced by \( \Cat_\infty \)---\emph{are} \( (\infty,2) \)-topoi in the sense of \cite{abellan-martini-2-topoi}.)
It would be very interesting to see if the main theorems of \cite{helfer-topoi-in-topoi} and the present paper could be proven under appropriate weaker assumptions which would allow for such examples.

\section{Preliminaries}
\subsection{Categories} \label{subsec:cat-prelims}
We discuss our notational conventions regarding categories.

We write \( \bC(X, Y) \) for the set of morphisms \( X \to Y \) between objects \( X,Y \in \bC \) in a category \( \bC \), and similarly \( \cC(X,Y) \) for the category of morphisms in a 2-category \( \cC \).
We write \( f_* \colon \bC(A, X) \to \bC(A, Y) \) and \( f^* \colon \bC(X, A) \to \bC(Y, A) \) for the maps given by composition with a morphism \( f \colon X \to Y \) in a category \( \bC \), and likewise for the composition-functors in a 2-category.
As in \cite{helfer-topoi-in-topoi}, we use the ``geometric order'' of composition, writing \( X \tox{f g} Z \) (or \( f \cdot g \) or \( f \circ g \)) for the composite of morphisms \( X \tox{f} Y \tox{g} Z \), but the conventional order \( f(x) \) for \emph{application} of a function \( f \colon X \to Y \) to an element \( x \in X \), including the application of function symbols to terms in first-order formulas (see \jref{rmk:kj-bookkeeping}).
We also use the geometric order \( f \alpha \colon W \tocell Y \) and \( \alpha g \colon X \tocell Z \) for the whiskering of a 2-cell \( \alpha \colon X \tocell Y \) with a morphism \( f \colon W \to X \) or \( g \colon Y \to Z \).

For a category \( \bC \), we write \( \bC^\iso \) for the core (or maximal sub-groupoid) of \( \bC \).

We say that a category \( \bC \) is \defword{of finite presentation} if there is a category \( \bF \) which is free on a finite graph, and a finite set \( R \subset \Ar \bF \times_{\Ob \bF \times \Ob \bF} \Ar \bF \) of pairs of parallel morphisms in \( \bF \), such that \( \bC \) is isomorphic to the quotient of \( \bF \) by the congruence relation generated by \( R \).

\subsubsection{Categorical discipline} \label{subsubsec:discipline}
When we speak of, for example, the product \( X \times Y \) of two sets or of two categories \( \bC \times \bD \), we have in mind the usual set-theoretic constructions.
However, when we speak of a product \( X \times Y \) of two objects in an \emph{arbitrary} category \( \bC \), we never assume, unless stated explicitly, that \( \bC \) is endowed with a preferred product operation; rather \( X \times Y \) denotes an \emph{arbitrary} object in \( \bC \), equipped with a product diagram, which we implicitly fix at the time when we introduce the product \( X \times Y \).
We also often use the word ``the'' when we should really say ``a'' for ease of reading.
Thus, if we write something like ``let \( f \colon A \to X \times Y \) be a morphism into the product \( X \times Y \)'', we really mean ``let \( X \xot{\pi_X} X \times Y \tox{\pi_Y} Y \) be a product diagram, and let \( f \colon A \to X \times Y \) be a morphism''.
Similar comments apply to other objects determined by universal properties (limits, colimits, power objects, etc.).
In places where we would want to make use of a \emph{functor} determined by making choices of objects determined by universal properties (such as product functors or pullback functors), we instead use the corresponding canonically defined \emph{anafunctor} (see \cite{makkai-avoiding-choice}).
We generally assume familiarity with anafunctors, and we will take for granted the generalization to anafunctors of various familiar facts regarding functors.

In connection with this convention, we allow ourselves a certain abuse of notation regarding the equality sign, of the following kind: if we are given an object \( P \) and a diagram \( X \ot P \to Y \), to express that this diagram exhibits \( P \) as a product of \( X \) and \( Y \), we may write \( P = X \times Y \).

As indicated above, we usually denote product projections out of a product \( X \times Y \) by \( \pi_X \) and \( \pi_Y \); when this would be ambiguous (for example, when \( X = Y \)), we write instead \( \pi_0 \) and \( \pi_1 \).
We sometimes also use the notation \( \pi_{ij} = \br{\pi_i,\pi_j} \colon X^n \to X^2 \) and \( \pi_{ijk} = \br{\pi_i,\pi_j,\pi_k} \colon X^n \to X^3 \).
We write \( \unex_X \) or \( \unex \) for the unique morphism \( X \to \tm \) from an object \( X \) to a given terminal object \( \tm \).
Given a morphism \( f \colon A \to B \) in a category \( \bC \) having pullbacks, we write \( f^* \colon \bC / B \to \bC / A \) for the pullback anafunctor.

\subsubsection{Underobjects and slice categories} \label{subsubsec:under-slice-prelims}
By an \defword{underobject} of an object \( X \), we mean a pair \( (A, i_A) \) consisting of an object \( A \) and a monomorphism \( i_A \colon A \tto X \); we would perhaps prefer to call this a \emph{subobject}, but to avoid confusion, we retain the usual meaning of that term, namely an equivalence class of underobjects.
We will often conflate an underobject \( (A, i_A) \) with its underlying object \( A \); by default, we use the notation \( i_A \) for the monomorphism \( i_A \colon A \tto X \) underlying an underobject \( A \).
Note that, given two monomorphisms \( P \tox{i_P} X \tox{f} Y \), if we write ``the underobject \( P \)'', this could refer to either the underobject \( (P, i_P) \) of \( X \) or \( (P, i_P f) \) of \( Y \), but it should always be clear from context what is meant.
We write \( A \le B \) and \( A \approx B \) for the usual preorder and equivalence relation on underobjects.

Given an object \( U \in \bC \) in a category \( \bC \), we denote an object \( p \colon X \to U \) in the slice \( \bC / U \) by \( (X, p) \).
When it is clear from context what the morphism \( p \) is, we may omit it from the notation, but in such cases, we always underline the object and write \( \und X \in \bC / U \) to distinguish it from \( X \) regarded as an object of \( \bC \).
When we introduce an object \( \und X \in \bC / U \) in a slice without introducing a notation for the underlying morphism \( X \to U \), we will by default denote it by \( \delta_X \colon X \to U \).
Given a morphism \( f \colon U \to V \) and objects \( \und X \in \bC / U \) and \( \und Y \in \bC / V \), we say a morphism \( g \colon X \to Y \) is \emph{over} \( f \) if \( g \delta_Y = \delta_X f \).
Note that, to denote a product in a slice category \( \bC / U \), we may write either \( \und X \times \und Y \) or \( \und{X \times_U Y} \).

We also recall the \emph{absoluteness} of certain notions between the slice category \( \bC / U \) and the underlying category \( \bC \): the forgetful functor \( \Sigma \colon \bC / U \to \bC \) preserves and reflects monos and pullbacks, and reflects colimits and epis; and if \( \bC \) has finite products, then \( \Sigma \) preserves colimits; and if moreover \( \bC \) has pushouts, then \( \Sigma \) preserves epis.
Finally, we recall the isomorphism \( (\bC / U) / \und A \toi \bC / A \), with respect to which the forgetful functor \( (\bC / U) / \und A \to \bC / U \) corresponds to the composition-with-\( \delta_A \) functor \( \Sigma_{\delta_A} \colon \bC / A \to \bC / U \).

\subsubsection{Topoi} \label{subsubsec:topoi-prelims}
By \emph{topos}, we always mean elementary topos.
We assume familiarity with the basic notions and results of topos theory, for which we refer to \cite{johnstone-topos-theory,johnstone-elephant-vols,mac-lane-moerdijk}.

Given an object \( X \) in a topos \( \bC \), a power object \( \pow X \) for \( X \) by definition comes with a universal relation \( \unin_X \tto X \times \pow X \), so that any other binary relation \( R \tto X \times U \) in \( \bC \) is pulled back along \( \id_X \times f \colon X \times U \to X \times \pow X \) for some unique \( f \colon U \to \pow X \).
In this case, we say that \( f \) \emph{classifies} the underobject \( R \).
It will often be more natural to have the product the other way around, so we may also speak of \( f \colon U \to \pow X \) classifying an underobject \( R \tto U \times X \) (meaning that \( R \) is the pullback of \( \unin_X \) under \( \br{\pi_1, \pi_0 f} \colon U \times X \to X \times \pow X \)).
In the case \( U = X \), where ambiguity may arise, we will say explicitly which of these conventions we are using.

We write \( {\le} = {\le_X} \tto \pow X \times \pow X \) for the subset relation on a power object \( \pow X \) and we write \( \top_X \colon \tm \to \pow X \) for the maximal element.
We write \( P \wedge Q \tto X \) for the intersection (pullback) of underobjects \( P, Q \tto X \).
As a maximal underobject, we always simply use \( X = (X, \id_X ) \), and we also do not introduce special notation for the minimal underobject as this is the same thing as an initial object \( \zset \tto X \).

We write \( \im(f) \tto Y \) for the image of a morphism \( f \colon X \to Y \), i.e., an underobject of \( Y \) through which \( f \) factors via an epi \( X \to \im(f) \), or equivalently, the least underobject of \( Y \) through which \( f \) factors.
We write \( \exists_f \colon \pow X \to \pow Y \) for the image \emph{morphism} along a morphism \( f \colon X \to Y \) and \( f\I \colon \pow Y \to \pow X \) for the preimage morphism.
Thus for \( p \colon U \to X \) and \( q \colon U \to Y \), we have \( p \exists_f \le_X q \) iff \( p \le_X q f\I \).
In \jrefs{exm:im-as-setbuilder}{exm:im-preim-by-comp}, we recall the description of \( \im(f) \) and of \( \exists_f \) and \( f\I \), as well as that of the subset relation \( {\le_X} \), in terms of generalized elements.
Note that for an underobject \( i_P \colon P \tto X \), the global element of \( \pow X \) classifying it is \( \top_P \exists_{i_P} \colon \tm \to \pow X \).

We write \( \sigma = \sigma_X \colon X \to \pow X \) for the singleton morphism, i.e., the morphism classifying the diagonal \( \Delta_X \colon X \tto X \times X \).

We call a collection of morphisms \( U_0,\ldots,U_{n - 1} \to U \) with a common codomain a \emph{cover} if these morphisms are jointly epimorphic, or equivalently, if the union of their images is equal to \( U \).
Thus, for instance, the coprojections \( A \tox{i_A} A + B \xot{i_B} B \) into a coproduct always form a cover---and we recall moreover that, in a topos, the coprojections are always monic and \emph{disjoint}, in the sense that \( A \wedge B \approx \zset \) \cite[Corollary~1.57]{johnstone-topos-theory}.

In \aref{sec:topos-logic-background}, we recall our conventions concerning interpreting logic in topoi, and the KJ semantics.
In particular, we explain what we mean by the \emph{internal language} \( \cL_\bC \) of a topos \( \bC \), what we mean by the \emph{satisfaction} \( \bC \vDash \phi \) of a higher-order formula \( \phi \) over \( \cL_\bC \), and more generally the \emph{forcing relation} \( U \Vdash_\vx \phi \) of such a higher-order formula \( \phi \) with respect to an object \( U \in \bC \) and a sequence \( \vx \) of morphisms \( x_i \colon U \to X_i \in \bC \).
In \sref{subsubsec:kj-semantics}, we also introduce some special notational conventions, according to which we, for example, denote the composite of two morphisms \( U' \to U \tox{x} X \) by \( x' \colon U' \to X \); and according to which, given morphisms \( x,y \colon U \to X \) and an underobject \( R \tto X \times X \), we write \( x R y \) to express that \( \br{x,y} \colon U \to X \times X \) factors through \( R \) (which we in fact already did above with \( R = {\le_X} \)).
And in \sref{subsubsec:interpretations}, we explain our use of the set-builder notation \( \set{x \mid \phi} \) to specify underobjects in topoi.
See there for more details on all of these things.

\subsection{2-topoi} \label{subsec:2-topoi-prelims}
We will frequently be making use of the notions from \cite{helfer-topoi-in-topoi} of a \emph{groupoidant} 2-category and of a \emph{plentiful DOF classifier} in a 2-category; as mentioned in the introduction, our ``minimal'' notion of 2-topos is that of a groupoidant 2-category containing a plentiful DOF classifier.
In fact, the details of these definitions will for the most part not matter, but rather only the main results that were proven in \opcitspace concerning them, which we recall in \jref{thm:topoi-in-topoi-main-thm} below.
Nonetheless, for convenience, we briefly recall the relevant definitions.
We refer to \opcitspace for a more detailed explanation of and motivation behind these definitions.

\subsubsection{Fibrations}
A functor \( p \colon \bC \to \bD \) between categories \( \bC \) and \( \bD \) is an \defword{isofibration} if for each isomorphism \( f \colon X \to Y \) in \( \bD \) and each \( \wt X \in \bC \) with \( p(\wt X) = X \), there exists an isomorphism \( \tilde f \) in \( \bC \) with domain \( \wt X \) and \( p(\tilde f) = f \).
A \defword{discrete opfibration (or DOF)} is defined in the same way, with both instances of ``isomorphism'' replaced by ``morphism'', and where the lift \( \tilde f \) with domain \( \wt X \) and \( p(\tilde f) = f \) is required to be \emph{unique}.
A \defword{setoidal opfibration (or SOF)} is defined in the same way as a DOF, except that \( p \) is now required to be faithful, the lift \( \tilde f \) is no longer required to be unique, and given any two lifts \( \tilde f_j \colon \wt X \to \wt Y_j \) (\( j = 0,1 \)), there is required to be a unique morphism \( g \colon \wt Y_0 \to \wt Y_1 \) with \( \tilde f_0 g = \tilde f_1 \) and \( p(g) = \id_Y \).

By a \defword{setoid}, we mean a category which is both a groupoid and a preorder, or equivalently, a category which is equivalent to a discrete category; indeed, any setoid \( X \) is equivalent to its \emph{quotient set}, obtained by identifying isomorphic objects.
The \defword{fibre} of a functor \( p \colon \bC \to \bD \) over \( X \in \bD \) is the subcategory of \( \bC \) consisting of objects mapping to \( X \) and morphisms mapping to \( \id_X \).
The fibres of a DOF are discrete categories, and those of a SOF are setoids.
Every DOF is a SOF and every SOF is an isofibration.

A morphism \( p \colon X \to Y \) in a 2-category \( \cC \) is an \defword{isofibration}, a \defword{DOF}, or a \defword{SOF} if the functor \( p_* \colon \cC(A,X) \to \cC(A, Y) \) is one for each \( A \in \cC \).
For an object \( A \) in a 2-category \( \cC \), we write \( \DOF(A) \) for the category of DOFs over \( A \) in \( \cC \), in which a morphism \( (X, p) \to (Y, q) \) is a morphism \( f \colon X \to Y \) with \( f p = q \).
We will say more about this in category \sref{subsubsec:dof-sof-digression} below.

\subsubsection{Corepita 2-categories} \label{subsubsec:corepita-prelims}
Fix a 2-category \( \cC \).
A square in \( \cC \) in a \defword{strict pullback square} if its image under \( \cC(A,-) \colon \cC \to \Cat \) is a pullback square (in the underlying 1-category of \( \Cat \)) for each \( A \in \cC \); \defword{strict terminal objects} are defined similarly.
A \defword{strict arrow object} for \( X \in \cC \) is an object \( X^\to \) equipped with a 2-cell \( X^\to \tocell X \) such that the induced functor \( \cC(A, X^\to) \to \cC(A, X)^\to \) (the codomain being the category of arrows in \( \cC(A, X) \)) is an isomorphism for each \( A \).
More generally, a \defword{strict cotensor} of \( X \in \cC \) by a category \( J \) is an object \( X^J \in \cC \) equipped with a functor \( J \to \cC(X^J, X) \) inducing an isomorphism \( \cC(A, X^J) \to \cC(A, X)^J \) for each \( A \in \cC \).

An object \( X \in \cC \) is \defword{groupoidal} (or simply \defword{a groupoid}) if \( \cC(A, X) \) is a groupoid for each \( X \in \cC \).
A morphism \( f \colon X \to Y \) is an \defword{arrow-wise iso} if \( \alpha f \colon A \tocell Y \) is an invertible 2-cell for each 2-cell \( \alpha \colon A \tocell X \).
Every morphism out of a groupoid is an arrow-wise iso.
A (strict) \defword{core} for \( X \in \cC \) is a groupoidal object \( X^\iso \) equipped with an arrow-wise iso \( i \colon X^\iso \to X \) such that every arrow-wise iso \( f \colon A \to X \) factors uniquely through \( i \), as does every invertible 2-cell \( A \tocell X \) between arrow-wise isos.

A 2-category \defword{has pita limits}, or \defword{is pita}, if it has a strict terminal object, strict pullbacks of isofibrations,\footnote{%
  As mentioned in \cite{helfer-topoi-in-topoi}, we only ever really need pullbacks of \emph{normal} isofibrations in the sense of \cite{bourke-accessible-aspects}.
}
and strict arrow objects, and is \defword{corepita} if it also has cores.
Any pita 2-category has cotensors by arbitrary finite categories (and more generally by categories which are of finite presentation, in the sense of \sref{subsec:cat-prelims}---see \cite[Proposition~2.7]{power-robinson-pie-limits}).

\subsubsection{Groupoidant 2-categories}
A \defword{gbo-congruence} (``groupoidal bijective-on-objects congruence'') in a pita 2-category \( \cC \) is a a truncated simplicial diagram \( C \colon \Delta_{\le 2}^\op \to \cC \) such that (i) \( C \) is an internal category, in the usual sense, in the underlying 1-category of \( \cC \), (ii) the morphism \( C_1 \to C_0 \times C_0 \) is a DOF, and (iii) the objects \( C_0 \), \( C_1 \), and \( C_2 \) are all groupoidal.

Any gbo-congruence \( C \in \cC \) determines a category for each \( A \in \cC \) with object set \( \Ob \cC(A, C_0) \) and arrow set \( \Ob \cC(A, C_1) \), and composition determined by the internal category structure on \( C \).
A morphism \( A \to C_1 \) is \defword{\( C \)-invertible} if it is invertible when considered as a morphism in this category.
Any gbo-congruence in a pita 2-category admits an \defword{object of isomorphisms} \( C_{\cong} \to C_1 \), which is by definition a \( C \)-invertible morphism such that every \( C \)-invertible morphism \( A \to C_1 \) factors through it uniquely.
The identity morphism \( e \colon C_0 \to C_1 \) of the internal category \( C \) is always \( C \)-invertible and hence factors through a unique morphism \( \bar e \colon C_0 \to C_{\cong} \).
We say that the gbo-congruence \( C \) is a \defword{complete congruence} if \( \bar e \) is an equivalence.

A \defword{cocone with vertex} \( Y \in \cC \) under a gbo-congruence \( C \) is a pair \( (g, \gamma) \) with \( g \colon C_0 \to Y \) a morphism and \( \gamma \colon \pi_0 g \To \pi_1 g \) a 2-cell satisfying \( e \gamma = \id_g \colon C_0 \tocell Y \) and \( \pi_{02} \gamma = (\pi_{01} \gamma) (\pi_{12} \gamma) \colon C_2 \tocell Y \), where we write \( \pi_j \colon C_1 \to C_0 \) and \( \pi_{j k} \colon C_2 \to C_1 \) for the images in the diagram \( C \colon \Delta_{\le 2}^\op \to \cC \) of the injective maps \( [0] \to [1] \) and \( [1] \to [2] \) in \( \Delta_{\le 2}^\op \) with image \( \set{j} \) and \( \set{j,k} \), respectively.
A cocone \( (g, \gamma) \) with vertex \( Y \) is a (strict) \defword{quotient} if every cocone \( (h, \eta) \) under \( C \) with vertex \( Z \) factors through a unique morphism \( f \colon Y \to Z \) (in the sense that \( h = g f \) and \( \eta = \gamma f \)), and every morphism of cocones \( \alpha \colon (h, \eta) \to (k, \kappa) \) with vertex \( Z \) (a 2-cell \( \alpha \colon h \To k \) with \( (\pi_0 \alpha) \kappa = \eta (\pi_1 \alpha) \colon C_1 \tocell Z \)) factors uniquely through a 2-cell \( \beta \colon Y \tocell Z \) (in the sense that \( \alpha = g \beta \)).

From any object \( X \in \cC \) in a corepita 2-category, we can form a complete congruence \( \Nv(X) \) called its \defword{nerve}, with objects \( (X^{[0]})^\iso \), \( (X^{[1]})^\iso \), and \( (X^{[2]})^\iso \) (these are cores of cotensors of \( X \) by the categories \( [k] \in \Delta_{\le 2} \)), and where the arrows in the diagram \( \Nv(X) \colon \Delta_{\le 2}^\op \to \cC \) are the ones determined by the morphisms in \( \Delta_{\le 2} \) between \( [0] \), \( [1] \), and \( [2] \).
The nerve comes with a canonical cocone \( (g, \gamma) \) with vertex \( X \), where \( g \colon (X^\iso)^{[0]} \to X \) is the core inclusion, and \( \gamma \colon (X^{[1]})^\iso \tocell X \) is the core inclusion composed with the universal 2-cell from the arrow object \( X^{[1]} \).

A corepita 2-category is \defword{groupoidant} if the canonical cocone under any nerve is a quotient, and if every complete congruence is isomorphic to the nerve \( \Nv(X) \) of some object \( X \).\footnote{%
  Or rather, we might just say that every complete congruence \emph{is} a nerve of some object, since the nerve is only determined up to isomorphism anyway).
}

\subsubsection{Plentiful DOF classifiers} \label{subsubsec:plenitude}
Next, fix a DOF \( p \colon \rS_* \to \rS \) in a 2-category \( \cC \).
For any \( A \in \cC \), there is an induced anafunctor \( \El^A \colon \cC(A, \rS) \to \DOF(A) \) taking \( f \colon A \to \rS \) to the strict pullback \( f^* p \colon f^* \rS_* \to A \)---say, with pullback projection \( \pi_f \colon f^* \rS_* \to \rS \)---and taking a 2-cell \( \gamma \colon f \To h \) to the unique morphism \( g \colon (f^* \rS_*, f^* p) \to (h^* \rS_*, h^* p) \) in \( \DOF(A) \) for which there exists a 2-cell \( \ddot \gamma \colon \pi_f \To g \pi_h \) with \( \ddot \gamma p = (f^* p) \gamma \colon f^* \rS_* \tocell \rS \).
We say that \( p \) is a \defword{generic DOF}, and \( \rS \) a \defword{DOF classifier}, if \( \El^A \) is fully faithful for all \( A \in \cC \).
We also say that the morphism \( f \) above \defword{classifies} the DOF \( (f^* \rS_*, f^* p) \in \DOF(A) \).
The DOFs in the essential image of \( \El^A \) are called \defword{\( p \)-small} or \defword{\( \rS \)-small}; they form a full subcategory of \( \DOF(A) \) denoted \( \DOF_\rS(A) \).
Similarly, we call a SOF over \( A \) \( p \)-small or \( \rS \)-small if it is in the essential image of the composite \( \cC(A, \rS) \tox{\El^A} \DOF(A) \hto \SOF(A) \) (see below for the definition of the category \( \SOF(A) \)), and write \( \SOF_\rS(A) \subset \SOF(A) \) for the category of \( \rS \)-small SOFs over \( A \).

A DOF classifier \( \rS \) is \defword{pre-plentiful} if every DOF which is a monomorphism (in the underlying category of \( \cC \)) is \( \rS \)-small, and if \( \rS \)-small DOFs are closed under composition.
It follows from this that, given DOFs \( X \tox{f} Y \tox{g} Z \), if \( f g \) is \( \rS \)-small, then so is \( f \) (see \prerefpost{footnote~}{footnote:1-of-2-of-3}{} on \pref{footnote:1-of-2-of-3}); it also holds (for \emph{any} DOF classifier \( \rS \), in fact) that \( \rS \)-small DOFs are closed under strict pullbacks.

For any pre-plentiful DOF \( \rS \) in a pita 2-category \( \cC \), there exists a \defword{monomorphism classifier} \( \Mon(\rS) \); this is equipped with a stable monomorphism 2-cell \( \Mon(\rS) \tocell \rS \) such that, for each \( A \in \cC \), the induced functor \( \cC\pbig{A, \Mon(\rS)} \to \cC(A, \rS)^\to \) is an isomorphism onto the full subcategory of \( \cC(A, \rS)^\to \) spanned by the stable monomorphisms.
Here, a 2-cell \( \alpha \colon X \tocell Y \) is a \defword{stable monomorphism} if \( f \alpha \colon A \tocell Y \) is a monomorphism in \( \cC(A, Y) \) for all \( f \colon A \to X \).
A pre-plentiful DOF classifier \( \rS \) is \defword{plentiful} if it is pre-plentiful, and the morphism \( \partial_1^\iso \colon \Mon(\rS)^\iso \to \rS^\iso \) (which is always a SOF) is \( \rS \)-small, where \( \partial_1 \) is the codomain of the universal 2-cell \(  \Mon(\rS) \tocell \rS \).

We are now ready to state the main results we will need from \cite{helfer-topoi-in-topoi}.
But first, we make a brief aside regarding some generalities on DOFs and SOFs.

\subsubsection{Digression on DOFs and SOFs} \label{subsubsec:dof-sof-digression}
The 2-category \( \DOF(A) \) of DOFs over an object \( A \) in a 2-category \( \cC \), has a natural 2-category structure, in which a 2-cell \( \alpha \colon (X, p) \tocell  (Y, q) \) is a 2-cell \( \alpha \colon X \tocell Y \) in \( \cC \) such that \( \alpha q = \id_{p} \).
However, every such \( \alpha \) turns out to be an identity, so that \( \DOF(A) \) really is naturally a 1-category.

Let us now write \( \SSOF(A) \) for the corresponding 2-category of \emph{SOFs} over \( A \); now, the 2-cells may be non-trivial, so this really is a 2-category, albeit of a very trivial nature: it is a \emph{locally setoidal} 2-category (i.e., all hom-categories are setoids), and thus ``essentially'' a 1-category.
Namely, it is equivalent, as a 2-category, to the locally discrete 2-category (i.e., 1-category) obtained by passing from each hom-setoid to its quotient set.
Let us write \( \SOF(A) \) for the latter category, so that we have a quotient 2-functor \( \SSOF(A) \to \SOF(A) \), and full inclusions \( \DOF(A) \hto \SSOF(A) \) and \( \DOF(A) \hto \SOF(A) \) (the latter because the morphisms in \( \DOF(A) \subset \SSOF(A) \) are unaffected by the quotient procedure producing \( \SOF(A) \)).

A \defword{DOF collapse} of a SOF over \( A \) is a DOF over \( A \) which is isomorphic to it in \( \SOF(A) \); we say that \( \cC \) \defword{has DOF collapses} if every SOF in \( \cC \) has a DOF collapse, i.e., if \( \DOF(A) \hto \SOF(A) \) is an equivalence for all \( A \in \cC \).
Note that for a DOF classifier \( \rS \), the inclusion \( \DOF_\rS(A) \hto \SOF_\rS(A) \) is \emph{always} an equivalence, by the definition of ``\( \rS \)-small SOF''.
We say that a groupoidant \( \cC \) \defword{satisfies the axiom (UA)} if it has DOF collapses and if any two DOFs in \( \cC \) are \( \rS \)-small with respect to some plentiful DOF-classifier \( \rS \).
Equivalently, \( \cC \) satisfies (UA) if any two SOFs in \( \cC \) are \( \rS \)-small with respect to some plentiful DOF classifier \( \rS \).

\subsubsection{Plentiful DOF classifiers are internal topoi} \label{subsubsec:dof-class-are-topoi}
We now turn to the promised results from \cite{helfer-topoi-in-topoi}.
We recall from there that an \defword{internal topos} in a 2-category \( \cC \) is an object \( E \in \cC \) such that \( \cC(A, E) \) is a topos for all \emph{groupoids} \( A \in \cC \), and the functor \( f^* \colon \cC(A', E) \to \cC(A, E) \) is logical (see \sref{subsec:logical-functors}) for all morphisms \( f \colon A \to A' \) between groupoids.
(Much of \opcitspace is dedicated to ``justifying'' the fact that this definition only makes reference to \emph{groupoids} \( A \)---or more precisely, justifying this \emph{in the context of a groupoidant 2-category}.)

We now have:
\begin{thm}[{\cite[Theorems~5.4.6~and~5.4.7]{helfer-topoi-in-topoi}}] \label{thm:topoi-in-topoi-main-thm}
  Let \( \cC \) be a groupoidant (or more generally, corepita) 2-category.
  Then:
  \begin{enumerate}[(i)]
  \item Any plentiful DOF classifier \( \rS \in \cC \) is an internal topos; thus, \( \DOF_\rS(A) \) is a topos for each groupoid \( A \in \cC \), and (by \jref{propn:el-pb-compat} below) the pullback anafunctor \( f^* \colon \DOF_\rS(A) \to \DOF_\rS(A') \) is logical for each morphism \( f \colon A' \to A \) between groupoids.
  \item\label{item:topoi-in-topoi-main-thm-inclusion} Given plentiful DOF classifiers \( \rS \) and \( \rS' \) such that the generic DOF \( p \colon \rS_* \to \rS \) over \( \rS \) is \( \rS' \)-small, and letting \( f \colon \rS \to \rS' \) be a morphism classifying \( p \), the functor \( f_* \colon \cC(A, \rS) \to \cC(A, \rS') \) is fully faithful and logical for each \( A \in \cC \); hence, by \jref{propn:el-pb-compat}, the inclusion \( \DOF_{\rS}(A) \hto \DOF_{\rS'}(A) \) is also logical for each \( A \in \cC \).
  \item If \( \cC \) satisfies axiom (UA), then \( \DOF(A) \) is a topos for each \( A \in \cC \), and the inclusion \( \DOF_\rS(A) \hto \DOF(A) \) is logical for each plentiful DOF classifier \( \rS \).
  \end{enumerate}
\end{thm}

Note that the above theorem remains valid if the categories \( \DOF(A) \) (and \( \DOF_\rS(A) \), etc.) are all replaced by \( \SOF(A) \) (etc.), as does the following proposition.

\begin{propn}[{see \cite[Proposition~5.2.1~and~Theorem~5.4.7]{helfer-topoi-in-topoi}}] \label{propn:el-pb-compat}
  Let \( \cC \) be a corepita 2-category, let \( \rS \), \( \rS' \), and \( f \) be as in \jref{thm:topoi-in-topoi-main-thm}~\ref{item:topoi-in-topoi-main-thm-inclusion}, and let \( g \colon A' \to A \) be a morphism in \( \cC \).
  Then the following diagrams of anafunctors commute up to isomorphism.
  \[
    \begin{tikzcd}
      \cC(A, \rS) \ar[r, "g^*"] \ar[d, "\El^{A}"'] & \cC(A', \rS) \ar[d, "\El^{A'}"] \\[5pt]
      \DOF(A) \ar[r, "g^*"] & \DOF(A')
    \end{tikzcd}
  \qquad \quad
  \begin{tikzcd}
    \cC(A, \rS) \ar[rd, "\El^A"] \ar[dd, "f_*"'] \\[-16pt]
    & \DOF(A) \\[-16pt]
    \cC(A, \rS') \ar[ru, "\El^A"']
  \end{tikzcd}
  \]
\end{propn}

\section{\alttext{\( \unein \)}{ε}-sets and \alttext{\( \unein \)}{ε}-objects} \label{sec:ein-objects}
We recall from the introduction that our approach to constructing a set theoretic universe in a 2-topos \( \cC \)---and thence also of implementing the Burali-Forti paradox---is based on the consideration of algebraic structures of the kind that are exhibited by transitive closures of hereditary (singleton) sets, endowed with the membership relation, and which we call hs-sets in the set-theoretic context, and hs-objects in the topos-theoretic context.
In this section, we formally introduce these objects and some related constructions, and prove various basic properties about them that we will use in what follows.

The definition of hs-object will be made in a piecemeal way: throughout the paper, we will need to consider structures which only have some of the three defining properties of hs-objects (namely extensionality, well-foundedness, and existence of a top element), so we name such structures as well.
What this amounts to is that we give a name to a structure not satisfying \emph{any} of these conditions---i.e., simply a set/object equipped with a binary relation.
There does not seem to be a well-established name for such structures (in contrast to the case of a set with a single binary \emph{operation}, which is called a \emph{magma}).
Or rather, there is one name---viz.\ \emph{directed graph}, if this is understood in the sense of allowing loops but not parallel edges---and indeed, in \cite{aczel-non-well-founded-sets}, this terminology is chosen for this purpose.
However, we find that the geometric/combinatorial associations of graphs make that name inappropriate for our context (though it is quite appropriate in \loccit, where general pointed graphs are taken as models for non-well-founded sets).
Instead, we use the name \emph{\( \unein \)-sets} or \emph{\( \unein \)-objects}, seeing as we usually use the symbol \( \unein \) to denote the given binary relation, and because our motivating example is that in which \( \unein \) is the actual membership relation \( \in \) on a hereditary set.

The most important related notion we introduce here is the analogue for hs-objects(/sets) of the membership relation on hereditary sets; this ends up being a morphism between the hs-objects satisfying certain conditions, and which we call an \emph{elemental morphism}.
Again, we will need to consider various weakenings of this notion, and so define it in a piecemeal manner (calling the intermediate notions \emph{consistent} and \emph{initial} morphisms).

The other main construction that we consider---which turns out to play a quite important role in what follows---is the \emph{downward closure} of an element (or subset) of an \( \unein \)-set \( (X, \unein) \) with respect to the relation \( \unein \).
In the topos-theoretic context, there are two versions of this construction, since we can speak either of subobjects of \( X \) or (generalized) elements of \( \pow X \); both are useful, and we consider them both.
We also study how \( \unein \)-objects in slice topoi \( \bC / U \) behave with respect to the natural functors \( \Sigma_f \colon \bC / U \toot \bC / V \colon f^* \) induced by a morphism \( f \colon U \to V \).

\subsection{\alttext{\( \unein \)}{ε}-sets} \label{subsec:ein-sets}
For the purpose of motivating the basic definitions, we begin in the set-theoretic context, i.e., by considering \( \unein \)-sets rather than \( \unein \)-objects.
We thus begin one of the two threads that will be developed in parallel through the paper, with, as we explained in \sref{subsec:sets-and-topoi}, the set-theoretic thread in each case providing motivation and clarification for the topos-theoretic one.
For reasons also explained there, all the definitions and propositions in this section will be encased in coloured boxes.
The claims made in this section, being entirely in the set-theoretic context, will have no bearing on the main results of the paper.
Moreover, several of the claims will simply be special cases of the analogous claims in the topos-theoretic context (namely by specializing to the topos \( \bC = \Set \)), and will therefore be proved later on.
However, some of the claims have no such analogue, being specific to the set-theoretic context and so, for completeness, we give their proofs (in \aref{sec:proofs}).

We use the terms ``class'' and ``large set'' interchangeably, and we may thus sometimes use ``set'' where we would elsewhere use ``small set''; we trust the reader to understand what is meant from context.
These can be understood (as is standard practice in category theory) to mean the subsets and elements, respectively, of some fixed Grothendieck universe.
(Of course, they can also simply be understood to mean ``class'' and ``set'', assuming we are working in, say, MK class-set theory; these two interpretations are closely related, see \cite[Theorem~13.1]{makkai-set-theory}.)
We recall that everything that we do in the set-theoretic context (i.e., everything in the coloured boxes) occurs in \emph{intuitionistic} set theory: all of the proofs must be intuitionistically valid.

\begin{settybox}
  \begin{defn}[cf.\ \jref{defn:topos-ein-hs}] \label{defn:hs-set}
    An \defword{\( \varepsilon \)-set} \( X = (X, \unein) \) is a set \( X \) together with a binary relation \( \unein \subset X \times X \).
    A subset \( S \subset X \) is \defword{(downward-)closed} if \( x \ein y \in S \) implies \( x \in S \), and is \defword{inductive} if the following holds: given \( y \in X \), if \( x \ein y \To x \in S \) for all \( x \in X \), then \( y \in S \).
    The \defword{downward closure} \( \dwnfn_X(S) \) or \( \dwnfn_X(x) \) (or just \( \dwnfn(S) \) and \( \dwnfn(x) \), respectively) of a subset \( S \subset X \) or element \( x \in X \) is the least downward-closed subset containing that subset or element.
    A \defword{top element} of \( X \) is an element \( \rt_X \in X \) with \( \dwnfn(\rt_X) = X \).

    We say that an \( \unein \)-set \( X \)
    \begin{enumerate}[(i)]
    \item is \defword{extensional} if: given any \( x,y \in X \), if \( w \ein x \ToT w \ein y \) for all \( w \in X \), then \( x = y \),
    \item is \defword{well-founded} if the relation \( \unein \) is well-founded, i.e., if any inductive subset \( S \subset X \) is equal to all of \( X \),
    \item is an \defword{hs-set} if it is extensional, well-founded, and has a top element.
    \end{enumerate}
  \end{defn}
\end{settybox}
As indicated in the definition, we will often conflate an \( \unein \)-set with its underlying set.
For an \( \unein \)-set \( X \), we will usually use \( \unein \), or for specificity \( \unein_X \), to denote its underlying relation.
The top element of an hs-set is uniquely determined (see \jref{rmk:on-topos-hs}), and we will always denote it by \( \rt_X \).
There is no reason in the definition of \( \unein \)-set to require the underlying set \( X \) to be \emph{small}, and in fact we will consider some large \( \unein \)-sets; of course, if \( X \) is large, the subsets ``\( S \)'' appearing in the definition must be understood to be possibly large as well.

As discussed in the introduction, the main examples of hs-sets come from the \( \unin \)-relation on the transitive closures of honest (singleton) hereditary sets.
Recall that the class \( \rV \) of \emph{hereditary sets} is the least class closed under set formation, i.e., such that \( x \subset \rV \) with \( x \) a set implies \( x \in \rV \).\footnote{
  It follows from the definition of the class \( \rV \) that it consists entirely of sets; indeed, taking \( P \subset \rV \) to be the class consisting of those \( x \in \rV \) which \emph{are} sets, we find that \( P \) is closed under set-formation and hence equal to \( \rV \).
  As mentioned in the introduction, it is common in axiomatic set theory (for example, in the axiomatization given in \sref{subsec:axioms-of-mk}---though notably not in \cite{zermelo-grundlagen-der-mengenlehre}) is to assume that \emph{everything} is a set (or rather, in class-set theory, that every element of a class is a set)---and moreover, in fact, via the axiom of foundation, that everything is a \emph{hereditary} set.
  Of course, under either of those assumptions, it is tautologous that \( \rV \) consists entirely of sets, and the above argument is empty.
  However, conceptually speaking, those assumptions are completely unnecessary, as one can simply restrict one's attention to \( \rV \), rather than assuming from the outset that everything belongs to it.
}

Recall also that a class \( A \) is \emph{transitive} if \( x \in y \in A \) implies \( x \in A \) (i.e., if \( y \in A \) implies \( y \subset A \)) for all sets \( y \).
For any class \( A \), we let \( \tc(A) \) denote its \emph{transitive closure}: the least transitive class \( B \) with \( A \subset B \).
For any class \( A \), we write \( \unin_A \) for the restriction \( \rstr{\unin}{A \times A} \) of the element relation to \( A \).

Here are some basic facts regarding \( \rV \) and the transitive closure operation:
\begin{settybox}
  \begin{lem} \label{lem:v-and-tc-facts}
    \enumbelow
    \begin{enumerate}[(i)]
    \item\label{item:v-and-tc-facts-h-if-h} A set is hereditary if and only if all of its elements are hereditary. \\
      (In particular, \( \rV \) is a transitive class.)
    \item\label{item:v-and-tc-facts-dwn-is-tc} A subclass \( A \subset \rV \) is downward-closed in \( (\rV, \unin_\rV) \) if and only if it is transitive. \\
      (More generally, if \( B \subset \rV \) is transitive, then \( A \subset B \) is downward-closed in \( (B, \unin_B) \) if and only if it is transitive.) \\
      (In particular \( \dwnfn_B(A) = \tc(A) \) for any \( A \subset B \).)
    \item\label{item:v-and-tc-facts-v-wf} \( (\rV, \unin_\rV) \) is a (large) extensional and well-founded \( \unein \)-set. \\
      (And so is \( (A, \unin_A) \) for any transitive \( A \subset \rV \).)
    \item\label{item:v-and-tc-facts-tc-is-set} The transitive closure of any hereditary set \( x \in \rV \) is again a (small) set. \\
      (Hence, it is the least transitive \emph{set} \( y \) with \( x \subset y \).)
    \item\label{item:v-and-tc-facts-tc-is-hs} \( \tc(\set{x}) \) is an hs-set for any \( x \in \rV \), with top element \( \rt_{\tc(\set{x})} = x \).
    \end{enumerate}
  \end{lem}

  \begin{proof}
    See \hyperref[proof:v-and-tc-facts]{\S\ref*{subsec:ein-sets-proofs}}.
  \end{proof}
\end{settybox}

Next, we have the converse to \jref{lem:v-and-tc-facts}~\ref{item:v-and-tc-facts-tc-is-hs}, discussed in the introduction.
For \( x \in \rV \), we write \( \unein_x \) as a shorthand for \( \unin_{\tc(\set{x})} \).
We note that the notion of \emph{isomorphism} (\( \cong \)) of \( \unein \)-sets is the obvious one (see also \jref{defn:ein-hs-mor}).
\begin{settybox}
  \begin{propn}[{Mostowski's collapsing lemma \cite[Theorem~3]{mostowski-undecidable-statement}}] \label{propn:mostowski-collapse}
    Given any hs-set \( X \), there is a unique hereditary set \( x \) such that \( \pbig{\tc(\set{x}), \ein_x} \cong (X, \ein_X) \).
  \end{propn}

  \begin{proof}
    Though this is well-known, it is normally considered in the context of classical set theory, so for completeness, we give a(n intuitionistically valid) proof in \hyperref[proof:mostowski-collapse]{\S\ref*{subsec:ein-sets-proofs}}.
  \end{proof}
\end{settybox}

\noindent
Moreover, the isomorphism \( X \toi \tc(\set{x}) \) is \emph{unique} because of:
\begin{settybox}
  \begin{propn}
    Any two isomorphisms \( X \toi Y \) between hs-sets are equal.
  \end{propn}

  \begin{proof}
    See \jref{propn:init-uniq}.
  \end{proof}
\end{settybox}
It follows that the (large) \defword{groupoid of hs-sets}, which we denote
\defword{\( \hsSet \)}, is in fact a \emph{setoid} (see \sref{subsec:cat-prelims}), and is equivalent to the class \( \rV \) of hereditary sets itself.
We define \( \rV_{\hs} \) to be the (super-large\footnote{%
  \label{footnote:super-large}
  It is super-large, being a collection of equivalence classes which are themselves large.
  If taking ``class'' and ``set'' to mean subset and element, respectively, of a Grothendieck universe, this is unproblematic.
  If working in class-set theory, however, this requires extending the theory to include ``super-classes''.
  On the other hand, the quotient setoid \( \rV_\hs \) is not strictly necessary; one can simply work with the setoid \( \hsSet \) instead.

  We also note that, though super-large, \( \rV_\hs \) is ``essentially (only) large'', i.e., isomorphic to a large set, namely \( \rV \).
})
quotient set of the setoid \( \hsSet \)---i.e., \( \rV_{\hs} \) is the collection of equivalence classes of hs-sets---so that \( \rV_{\hs} \cong \rV \).
We next introduce the element relation on \( \rV_\hs \):

\begin{settybox}
  \begin{defn}[cf.\ \jref{defn:ein-hs-mor}] \label{defn:ein-hs-set-mor}
    Given \( \unein \)-sets \( X \) and \( Y \), a function \( f \colon X \to Y \) is:
    \begin{enumerate}[(i)]
    \item \defword{consistent} if \( x \ein y \in X \) implies \( f(x) \ein f(y) \),
    \item \defword{conservative} if \( f(x) \ein f(y) \) implies \( x \ein y \) for any \( x,y \in X \),
    \item \defword{initial} if it is consistent and, given \( x_1 \in X \) and \( y_0 \in Y \) with \( f(x_1) \ein y_0 \), there exists \( x_0 \in X \) with \( x_0 \ein x_1 \) and \( f(x_0) = y_0 \) (thus, the image of \( f \) is an ``initial segment'' of \( Y \))\footnotemark
    \item \defword{elemental} if it is initial and \( f(\rt_X) \ein \rt_Y \) whenever \( \rt_X \in X \) and \( \rt_Y \in Y \) are top elements.
    \end{enumerate}

    For hs-sets \( X \) and \( Y \), we also write \( X \Ein Y \) to express that there \emph{exists} a (by the following proposition necessarily unique) elemental map \( X \to Y \).
  \end{defn}
\end{settybox}
\footnotetext{Initial morphisms are called \emph{simulations} in \cite{shulman-comparing}.
Also, note that the word ``initial'' here is not being used in the sense of ``initial object''.}

The main facts regarding elemental maps are:
\begin{settybox}
  \begin{propn} \label{propn:init-map-props}
    Suppose \( X \) and \( Y \) are \( \unein \)-sets.
    Then
    \begin{enumerate}[(i)]
    \item\label{item:init-map-props-inj} If \( X \) is extensional and well-founded, then any initial map \( X \to Y \) is injective.
    \item\label{item:init-map-props-unique} If \( X \) is well-founded and \( Y \) is extensional, then any two initial maps \( X \to Y \) are equal.
    \item\label{item:init-map-props-el-el} Given hereditary sets \( x \) \( y \), there exists a (by \ref{item:init-map-props-unique} unique) elemental \( \pbig{\tc(\set{x}), \ein_x} \to \pbig{\tc(\set{y}), \ein_y} \) iff \( x \in y \) (in which case this map is just the inclusion \( \tc(\set{x}) \hto \tc(\set{y}) \)).
    \end{enumerate}
  \end{propn}

  \begin{proof}
    For \ref{item:init-map-props-inj} and \ref{item:init-map-props-unique}, see \jrefs{propn:init-cons-mono}{propn:init-uniq}.
    For \ref{item:init-map-props-el-el}, see \hyperref[proof:init-map-props-el-el]{\S\ref*{subsec:ein-sets-proofs}}.
  \end{proof}
\end{settybox}

Let \( \unEin \) denote the \defword{groupoid of elemental maps}: it has objects triples \( (X,Y,f) \) with \( X,Y \in \hsSet \) and \( f \colon X \to Y \) an elemental map, and a morphism \( (X,Y,f) \to (X',Y',f') \) is a pair of isomorphisms \( (X \toi X', Y \toi Y') \) making a commutative square with \( f \) and \( f' \).
We have an obvious forgetful functor \( i_{\unEin} \colon \unEin \to \hsSet \times \hsSet \), and by definition, for \( X,Y \in \hsSet \), we have \( X \Ein Y \) if and only if \( (X,Y) \) is in the image of \( i_{\unEin} \).

It follows from \jref{propn:init-map-props} that \( \unEin \) is a setoid and that the functor \( i_{\unEin} \) is fully faithful (i.e., an ``injective map of setoids'').
Letting \( \unein_{\rV_\hs} \) be the (large) quotient set of \( \unEin \), we thus have an induced injection \( i_{\unein_{\rV_\hs}} \colon \unein_{\rV_\hs} \tto \rV_\hs \times \rV_\hs \), and the above proposition further shows that the map \( x \mapsto \tc(\set{x}) \) gives an isomorphism of \( \unein \)-sets \( (\rV, \unin_\rV) \toi (\rV_{\hs}, \ein_{\rV_{\hs}}) \).

We have thus established the reconstruction of the notion of hereditary set in terms of abstract structures (viz.\ hs-sets) described in \sref{subsec:abstract-concrete}.
The next things we will want to do are (in \sref{subsec:univ-ext-wf}) to show that \( \rV_{\hs} \) is extensional and well-founded, with a view to making it into an hs-set, and thus an ``element of itself'', leading to the Burali-Forti paradox; and then (in \sref{sec:v-is-a-model}) to show that \( \rV_{\hs} \) satisfies the remaining axioms of set theory.
But first, we leave the set-theoretic context aside and (in the rest of \sref{sec:ein-objects} and in \sref{sec:univs-in-2-topoi}) begin to develop the theory of \( \unein \)-objects in topoi and in 2-topoi.

\subsection{\alttext{\( \unein \)}{ε}-objects in a topos}
We now internalize the notions from the previous section in an arbitrary elementary topos.
The general procedure to do this is to express the definitions formally in higher-order logic, and then interpret them using the internal logic of the topos, which is done most conveniently using the KJ semantics.
Later (in \sref{subsec:hs-2-topos}), for the purpose of internalizing these notions in a \emph{2-topos}, we will need to make this formalization completely explicit.
But for now, we just explain the resulting topos-theoretic definitions as directly as possible, making use of logical formulas and the KJ semantics whenever convenient.
In that connection, we refer to \sref{subsec:logic-in-topoi} for our conventions regarding interpreting logic in topoi, and to \sref{subsec:cat-prelims} for our basic notation regarding power objects.

\begin{defn} \label{defn:topos-ein-hs}
  An \defword{\( \varepsilon \)-object} \( X = (X, \unein) \) in a topos \( \bC \) is an object \( X \in \bC \) together with a binary relation \( i_{\unein} \colon \unein \tto X \times X \).

  Given an \( \unein \)-object \( (X, \ein) \), we say that a morphism \( p \colon U \to \pow X \) is \defword{inductive} if
  \begin{equation} \label{eq:inductive-formula}
    U \Vdash_p \forall y \tolon X.\ (\forall x \tolon X.\ x \ein y \to x \in p) \to y \in p
  \end{equation}
  and is \defword{(downward-)closed} if
  \begin{equation} \label{eq:closed-formula}
    U \Vdash_p \forall x,y \tolon X.\ (x \ein y \wedge y \in p) \to x \in p.
  \end{equation}
  We also write
  \[
    \induc = \induc_X \tto \pow X
    \quad \AND \quad
    \clo = \clo_X \tto \pow X
  \]
  for the underobjects \( \set{ p \in \pow X \mid \phi } \tto \pow X \) defined by the formulas \( \phi \) appearing in \eqref{eq:inductive-formula} and \eqref{eq:closed-formula}, respectively, and we also say that an underobject \( P \tto X \) is inductive/closed if the morphism \( p = \top_P \exists_{i_P} \colon \tm \to \pow X \) classifying it is.

  A \defword{top element} of \( X \) is an element \( x \colon \tm \to X \) such that \( P \approx X \) for any closed \( P \tto X \) containing \( x \).
  (Equivalently, in the notation of \jref{defn:topos-dwn} below, \( x \) is a top element iff \( x \dwn = \top_X \colon \tm \to \pow X \).)

  We say that an \( \unein \)-object \( X \) is
  \begin{enumerate}[(i)]
  \item \defword{extensional} if
    \( \bC \vDash \forall x,y \tolon X.\ (\forall w \tolon X.\ w \ein x \tot w \ein y) \to x = y \)
  \item \defword{well-founded} if \( P \approx X \) for any inductive underobject \( P \tto X \)
  \item an \defword{hs-object} if it is extensional, well-founded, and has a top element \( \rt_X \colon \tm \to X \).
  \end{enumerate}
  As with \( \unein \)-sets, we often conflate an \( \unein \)-object with its underlying object, and if \( X \) is an \( \unein \)-object, we usually write \( i_{\unein} \colon \unein \tto X \times X \) or \( i_{\unein_X} \colon \unein_X \tto X \times X \) for its underlying relation.
\end{defn}

\begin{rmk} \label{rmk:on-topos-hs} \enumbelow
  \begin{enumerate}[1.]
  \item In \jref{cor:dwn-mono} below, we will show that an hs-object in fact has a \emph{unique} top element, which we will always denote by \( \rt_X \colon \tm \to X \).
  \item The condition that \( x \colon \tm \to X \) is a top element implies the following apparently stronger condition: any closed \( p \colon U \to \pow X \) containing \( x \) (that is, with \( x \in \unex_U p \)) is equal to \( \unex_U \top_X \):
    \[
      U \Vdash_x \forall p \tolon \pow X.\ (\clo(p) \wedge x \in p) \to p = \top_X.
    \]
    The definition of \( x \) being a top element says precisely that this condition holds for \emph{global} elements \( p \colon \tm \to \pow X \).
    The reason this implies the stronger condition is that it implies \( x \dwn = \top_X \colon \tm \to \pow X \) (see \jref{defn:topos-dwn}) and hence \( \unex_U x \dwn = \unex_U \top_X \colon U \to \pow X \) for any \( U \), which by definition means that any \( \unex_U \top_X \le p \) for any closed \( p \colon U \to \pow X \) with \( \unex_U x \in p \).
  \item Similarly, the definition of \( X \) being well-founded implies the apparently stronger
    \[
      \bC \vDash \forall p \tolon \pow X.\ \induc(p) \to p = \top_X,
    \]
    i.e., that any inductive \( p \colon U \to \pow X \) is equal to \( \top_X \), whereas again, this is by definition only required for global elements \( p \colon \tm \to \pow X \).
    The reason is essentially the same: the weaker condition implies that the underobject
    \[
      \set{ x \in X \mid \forall p \colon \pow X.\ \induc(p) \to x \in p } \tto X
    \]
    is equivalent to \( X \), and this in turn implies that \( x \in p \) for all \( x \colon U \to X \) and \( p \colon U \to \pow X \) inductive.
  \end{enumerate}
\end{rmk}

In \jref{defn:hs-set} (of ``hs-set''), we also introduced the operation of ``downward-closure'': the least downward-closed set containing a given element or set.
In a topos, we can formulate this as follows:
\begin{defn} \label{defn:topos-dwn}
  Given an \( \unein \)-object \( X \) in a topos \( \bC \) and a choice of power object \( \pow X \), we define the \defword{closure morphism} for \( X \), denoted \( \dwnfn = \dwnfn_X \colon X \to \pow X \), to be the morphism such that, for each \( x \colon U \to X \), \( x \dwn \) is the least closed \( p \colon U \to \pow X \) with \( x \in p \).

  The existence of the morphism \( \dwnfn \) is guaranteed by \jref{propn:comprehension-morphism}, since it is can be alternatively described as the unique morphism satisfying
  \[
    \bC \vDash
    \forall x, y \tolon X.\ x \in \dwnfn(y) \tot
    \pbig{\forall p \tolon \pow X.\ (\clo(p) \wedge y \in p) \to x \in p}.
  \]
  (We leave to the reader the easy verification that this is indeed equivalent.)

  Similarly, we define \( \dwnfn = \dwnfn_X \colon \pow X \to \pow X \) to be the morphism such that, for \( q \colon U \to \pow X \), \( q \dwn \) is the least closed \( p \colon U \to \pow X \) with \( q \le p \), or alternatively, the unique morphism satisfying
  \[
    \bC \vDash
    \forall x \tolon X.\ \forall q \tolon \pow X.\ x \in \dwnfn(q) \tot
    \pbig{\forall p \tolon \pow X.\ (\clo(p) \wedge q \le p) \to x \in p}.
  \]
  In particular, note that \( p \colon U \to \pow X \) is closed if and only if \( p \dwn = p \).

  We may write \( \dwnfn_X^X \colon X \to \pow X \) and \( \dwnfn_X^{\pow X} \colon \pow X \to \pow X \) if needed to distinguish between the two versions of the closure morphism.
  Note that \( \dwnfn_X^X = \sigma \dwnfn_X^{\pow X} \), with \( \sigma \colon X \to \pow X \) the singleton morphism.

  We also write \( {\preceq} \tto X \times X \) for the relation
  \[
    {\preceq}_X = \set{\br{x,y} \mid x \in \dwnfn(y)} \tto X \times X,
  \]
  which we call the \defword{ancestor relation}; thus, for \( x,y \colon U \to X \), we say that \( x \) is an \emph{ancestor} of \( y \), or \( y \) a \emph{descendent} of \( x \), if \( x \preceq y \), i.e., if \( x \in y \dwn \).
  Note that
  \[
    y \preceq y
  \]
  (by definition of \( y \dwn \)), that
  \[
    x \ein y \implies x \preceq y
  \]
  (by closedness of \( y \dwn \)), and that
  \[
    x \preceq y \preceq z \implies x \preceq z,
  \]
  i.e., that \( y \in z \dwn \) implies \( y \dwn \le z \dwn \) (since \( z \dwn \) is closed and \( y \dwn \) is least closed containing \( y \)).

  We will show in \jref{propn:no-loops} below that if \( X \) is well-founded, then \( \preceq \) is also antisymmetric, and thus a partial order.
\end{defn}

\begin{defn} \label{defn:ein-hs-mor}
  Given \( \unein \)-objects \( X \) and \( Y \) in a topos \( \bC \), a morphism \( f \colon X \to Y \) is:
  \begin{enumerate}[(i)]
  \item \defword{consistent} if
    \( \bC \vDash \forall x,y \colon X.\ x = y \to f(x) = f(y) \),
  \item \defword{conservative} if
    \( \bC \vDash \forall x,y \colon X.\ f(x) = f(y) \to x = y \),
  \item \defword{initial} if it is consistent 
    and
    \( \forall x_1 \tolon X.\ \forall y_0 \tolon Y.\ y_0 \ein f(x_1) \to
    (\exists x_0 \tolon X.\ f(x_0) = y_0 \wedge x_0 \ein x_1)  \)
  \item \defword{elemental} if it is initial and \( \rt_X f \ein \rt_Y \) for any\footnotemark{} top elements \( \rt_X \colon \tm \to X \) and \( \rt_Y \colon \tm \to Y \)
  \item an \defword{isomorphism of \( \unein \)-objects} if \( f \) is invertible in \( \bC \), and \( f \) and \( f\I \) are both consistent.
  \end{enumerate}

  As with hs-sets, we write \( X \Ein Y \) to express that there exists an elemental \( X \to Y \) between hs-objects \( X \) and \( Y \).
\end{defn}

\footnotetext{%
  We will only consider elemental morphisms between hs-objects, which have a unique top element, so we can say more simply: \( f \) is initial and \( \rt_X f = \rt_Y \).%
}
Note that consistent and conservative morphisms are closed under composition, and that the isomorphisms of \( \unein \)-objects are just the isomorphisms in the category of \( \unein \)-objects and consistent (or conservative) morphisms.
Note also that any isomorphism \( f \colon X \to Y \) of \( \unein \)-objects is initial.
(Indeed, given \( x_1 \colon U \to X \) and \( y_0 \colon U \to Y \) with \( y_0 \ein x_1 f \), then \( x_0 \defeq y_0 f\I \colon U \to X \) satisfies \( x_0 f = y_0 \) and \( x_0 = y_0 f\I \ein x_1 f f\I = x_1 \) by consistency of \( f\I \).)

\begin{defn}
  Given a topos \( \bC \), we write \( \hs(\bC) \) for the groupoid of hs-objects in \( \bC \) and isomorphisms between these, and \( \unEin_{\bC} \) for the groupoid of triples \( (X,Y,f) \) with \( X,Y \in \hs(\bC) \) and \( f \colon X \to Y \) elemental, and whose morphisms are pairs of isomorphisms between hs-objects forming a commutative square.
  There is an obvious forgetful functor \( i_{\unEin_\bC} \colon \unEin_\bC \to \hs(\bC) \times \hs(\bC) \).
\end{defn}
Note that, by \jref{propn:logical-fun-pres-satisf}, any logical functor \( F \colon \bC \to \bD \) between topoi induces a functor \( F_* \colon \hs(\bC) \to \hs(\bD) \), and similarly, any isomorphism 2-cell \( \beta \colon \bC \tocell \bD \) between logical functors induces a natural isomorphism \( \beta_* \colon \hs(\bC) \tocell \hs(\bD) \).

In \sref{subsec:univ-ext-wf}, we will be showing that the internalization of \( i_{\unEin_\bC} \colon \unEin_\bC \to \hs(\bC) \times \hs(\bC) \) in a 2-topos is an extensional, well-founded \( \unein \)-object; this reduces to the following ``external'' statement:
\begin{propn} \label{propn:v-ext-setoid-ff}
  Given any topos \( \bC \), the groupoids \( \hs(\bC) \) and \( \unEin_\bC \) are both setoids, and \( i_{\unEin_\bC} \) is fully faithful.
\end{propn}
\begin{proof}
  Unfolding the definitions, this immediately reduces to \jref{propn:init-uniq} below.
\end{proof}

\subsection{Basic properties of \alttext{\( \unein \)}{ε}-objects and morphisms} \label{subsec:ein-ob-props}
Throughout this section, let \( X \), \( Y \), and \( Z \) be \( \unein \)-objects in a topos \( \bC \).
The proofs of all statements in this section whose proof is not given here can be found in \sref{subsec:ein-ob-props-proofs}.

As per the discussion in \sref{subsec:sets-and-topoi}, several of the statements below have been proven elsewhere, e.g., \cite[\S8]{shulman-comparing} and \cite[\S9.2]{johnstone-topos-theory}, though mostly under somewhat different assumptions.

\prooflabel{proof:cons-wf-wf}
\begin{propn} \label{propn:cons-wf-wf}
  If \( f \colon X \to Y \) is a consistent morphism and \( Y \) is well-founded, then \( X \) is well-founded.
\end{propn}

\prooflabel{proof:init-uniq}
\begin{propn} \label{propn:init-uniq}
  If \( X \) is well-founded and \( Y \) is extensional, then there is at most one initial morphism \( X \to Y \).

  \noindent
  (Hence, since, as noted above, isomorphisms of \( \ein \)-objects are initial, there is also at most one isomorphism \( X \toi Y \) of \( \ein \)-objects if \( X \) is well-founded and \( Y \) is extensional.)
\end{propn}

\prooflabel{proof:init-cons-mono}
\begin{propn} \label{propn:init-cons-mono}
  If \( X \) is extensional and well-founded, then any initial morphism \( f \colon X \to Y \) is monic and conservative.

  \noindent
  (A more precise statement is: any initial mono \( f \colon X \to Y \) is conservative, without any assumptions on \( X \); and if \( X \) is extensional and well-founded, then any initial morphism \( f \colon X \to Y \) is monic.)
\end{propn}

\prooflabel{proof:dwn-dwn}
\begin{propn} \label{propn:dwn-dwn}
  If \( X \) is well-founded, then for any initial morphism \( f \colon X \to Y \), the following squares commute.
  \[
    \begin{tikzcd}
      \pow X \ar[r, "\exists_f"] \ar[d, "\dwn"'] & \pow Y \ar[d, "\dwn"] \\
      \pow X \ar[r, "\exists_f"] & \pow Y
    \end{tikzcd}
    \quad\AND\quad
    \begin{tikzcd}
      X \ar[r, "f"] \ar[d, "\dwn"'] & Y \ar[d, "\dwn"] \\
      \pow X \ar[r, "\exists_f"] & \pow Y
    \end{tikzcd}
  \]
\end{propn}

\prooflabel{proof:init-2-of-3}
\begin{propn} \label{propn:init-2-of-3}
  Given \( X \tox{f} Y \tox{g} Z \), if \( f \) and \( g \) are initial, so is \( f g \); and conversely, if \( Y \) is extensional and well-founded, and \( g \) and \( f g \) are initial, then so is \( f \).

  \noindent
  (More generally, without any assumptions on \( Y \): assuming \( g \) consistent and monic, if \( f g \) is initial, then so if \( f \); that this is indeed more general is because of \jref{propn:init-cons-mono}.)

  \noindent
  (Also: if \( g \) is conservative (resp.\ consistent), then if \( f g \) is consistent (resp.\ conservative), so is \( f \).)
\end{propn}

\subsubsection{The ancestor relation}
\( X \), \( Y \), and \( Z \) continue to denote \( \unein \)-objects in a topos \( \bC \).
\begin{defn} \label{defn:strict-ancestor}
  We write \( {\prec} \tto X \times X \) for the \defword{strict ancestor relation}, defined as
  \[
    {\prec} = \set{\br{x,z} \mid \exists y \tolon X\ x \preceq y \wedge y \ein z} \tto X \times X.
  \]
  Note that we have immediately, for \( x,y,z \colon U \to X \), that
  \[
    x \ein z \implies x \prec z \implies x \preceq z,
  \]
  since \( x \ein z \) implies \( x \preceq y \ein z \) with \( x = y \), and since, if \( x' \preceq y \ein z' \) for some \( U \oot U' \tox{y} X \), then \( x' \preceq z' \), hence \( x \preceq z \) by epicity of \( \delta_{U'} \).
  Note also that
  \[
    x \prec y \prec z \implies x \prec z,
  \]
  since if \( x' \preceq u \ein y \preceq v \ein z \) for some \( U \oot U' \tox{u,v} X \), then \( x' \preceq v \ein z \).
  We show in \jref{propn:no-loops} below that if \( X \) is well-founded, then \( \prec \) is irreflexive, and hence a strict partial order.
\end{defn}

\prooflabel{proof:dwn-opts}
\begin{propn} \label{propn:dwn-opts}
  \( \bC \vDash \forall x,y \tolon X.\ x \preceq y \tot (x = y \vee x \prec y) \)
\end{propn}

\prooflabel{proof:no-loops}
\begin{propn} \label{propn:no-loops}
  If \( X \) is well-founded, then
  \[
    \bC \vDash \forall x \tolon X.\ \neg (x \prec x)
  \]
  and hence
  \[
    \bC \vDash \forall x \tolon X.\ \neg (x \ein x)
    \qquad
    \AND
    \qquad
    \bC \vDash \forall x,y \tolon X.\ (x \preceq y \wedge y \preceq x) \to x = y.
  \]
\end{propn}

\begin{cor} \label{cor:dwn-mono}
  If \( X \) is well-founded, then \( \dwn \colon X \to \pow X \) is a monomorphism.

  \noindent
  (In particular, any \( x, y \colon U \to X \) with \( x \dwn = \unex \top_X = y \dwn \) are equal, so all top elements of \( X \) are equal.)
\end{cor}

\begin{proof}
  Given \( x,y \colon U \to X \) with \( x \dwn = y \dwn \colon U \to \pow X \), it follows that \( x \in x \dwn = y \dwn \) and \( y \in y \dwn = x \dwn \), hence \( x \preceq y \preceq x \), and hence, by \jref{propn:no-loops}, that \( x = y \).
\end{proof}

\subsubsection{Substuctures}
\( X \), \( Y \), and \( Z \) continue to denote \( \unein \)-objects in a topos \( \bC \).

\begin{defn}
  Any underobject \( i_A \colon A \tto X \) carries a (up to isomorphism) unique relation \( \unein_A \colon A \to A \times A \) for which \( i_A \) is consistent and conservative, namely \( \unein_A = (i_A \times i_A)^* \ein_X \).
  When regarding \( A \) as endowed with such a relation, we will use the term \defword{substructure}.
  We say that a substructure \( A \) of \( X \) is \defword{closed} if it is closed as an underobject (\jref{defn:topos-ein-hs})---or equivalently, by \jref{propn:closed-iff-init} below, if \( i_A \) is initial---and that it is \defword{elemental} if \( i_A \) is elemental.
  By an \defword{hs-substructure}, we mean a substructure which is also an hs-object.

  Note that if two substructure \( A,B \tto X \) are equivalent as underobjects, then the unique isomorphism \( A \to B \) commuting with \( i_A \) and \( i_B \) is automatically an isomorphism of \( \unein \)-objects; that is, an ``isomorphism of substructures'' is just the same as an isomorphism of underobjects.
\end{defn}

\prooflabel{proof:closed-iff-init}
\begin{propn} \label{propn:closed-iff-init}
  A substructure \( i_A \colon A \tto X \) is closed if and only if \( i_A \) is initial.

  \noindent
  (Hence, by \jref{propn:init-cons-mono}, ``closed substructure'' is essentially synonymous with ``initial mono''.)

  \noindent
  (Moreover, initiality holds here in the strong sense that given \( \bar x_1 \colon U \to A \) and \( x_0 \colon U \to X \) with \( x_0 \ein \bar x_1 i_A \), there is \( \bar x_0 \colon U \to A \) with \( \bar x_0 \ein \bar x_1 \) and \( \bar x_0 i_A = x_0 \).)
\end{propn}

\prooflabel{proof:sub-clo-induc}
\begin{propn} \label{propn:sub-clo-induc}
  Let \( i_A \colon A \tto X \) be a closed substructure.
  If \( P \tto X \) is closed or inductive, then so is \( i_A^* P \tto A \).
\end{propn}

\begin{propn} \label{propn:sub-ewf}
  Let \( A \tto X \) be a substructure.
  If \( X \) is well-founded, then so is \( A \); and assuming \( A \) is closed, if \( X \) is extensional, then so is \( A \).
\end{propn}

\begin{proof}
  The well-foundedness claim follows immediately from \jref{propn:cons-wf-wf}.
  For the extensionality claim, see \hyperref[proof:sub-ewf]{\S\ref*{subsec:ein-ob-props-proofs}}.
\end{proof}

\begin{defn} \label{defn:down-family}
  Given a global element \( x \colon \tm \to X \), we write
  \[
    i_{\Dwn_x} \colon \Dwn_x \tto X
  \]
for the underobject classified by \( x \dwn \colon \tm \to \pow X \), i.e., for the least closed underobject of \( X \) containing \( x \).
  More generally, for a generalized element \( x \colon U \to X \), we write \( i_{\Dwn_x} \colon \Dwn_x \to U \times X \) for the underobject classified by \( x \dwn \colon U \to \pow X \).
  We also set
  \[
    \delta_{\Dwn_x} \defeq i_{\Dwn_x} \pi_U \colon \Dwn_x \to U,
  \]
  and thus consider \( \Dwn_x \) as an object \( \und{\Dwn_x} = (\Dwn_x, \delta_{\Dwn_x}) \) in the slice \( \bC / U \).
\end{defn}

The construction \( \Dwn_x \) for generalized elements \( x \colon U \to X \) reduces to the one for global elements by passing to the slice \( \bC / U \), in the following sense:
\prooflabel{proof:closure-in-slice}
\begin{propn} \label{propn:closure-in-slice}
  Fix a morphism \( x \colon U \to X \), and write
  \[
    \und{\Dwn_{\br{\id_U, x}}} \tto \und{U \times X}
  \]
  for the underobject in \( \bC / U \) obtained by applying \jref{defn:down-family} to \( \br{\id_U, x} \colon \und U \to \und{U \times X} \) in \( \bC / U \) (where \( \und{U \times X} \) is given an \( \unein \)-object structure as in \jref{exm:pb-fun-logical}).

  Then the underobject
  \[
    i_{\Dwn_x} \colon \und{\Dwn_x} = (\Dwn_x, \delta_{\Dwn_x}) \tto \und{U \times X}
  \]
  in \( \bC / U \) is equivalent to \( \und{\Dwn_{\br{\id_U, x}}} \tto \und{U \times X} \).

  In other words, \( \und{\Dwn_x} \) is the least closed substructure of \( \und{U \times X} \) containing \( \br{\id_U, x} \).
\end{propn}
The proof of \jref{propn:closure-in-slice} depends on the following alternative characterization of \( \Dwn_x \):

\begin{lem} \label{lem:dwn-is-descendent}
  Writing \( ({\succeq}_X \tto X \times X) \) for the descendent relation
  \[
    ({\preceq}_X \tox{i_{{\preceq}_X}} X \times X \tox{\br{\pi_1, \pi_0}} X \times X),
  \]
  for any morphism \( x \colon U \to X \), there is a pullback square
  \[
    \begin{tikzcd}
      \Dwn_x \pb \ar[r, ""] \ar[d, "i_{\Dwn_x}"'] &[10pt] \succeq_X \ar[d, "i_{\succeq_X}"] \\
      U \times X \ar[r, "x \times \id_X"] & X \times X.
    \end{tikzcd}
  \]
  In particular, \( \Dwn_{\id_X} \approx {\succeq}_X \tto X \times X \).
\end{lem}

\begin{proof}
  Given \( y \colon U' \to X \), we must show that \( \br{\delta_U, y} \colon U' \to U \times X \) factors through \( \Dwn_x \) iff \( x' \succeq y \), i.e., iff \( y \in x' \dwn \); but this is just the statement that \( \Dwn_x \) is the underobject classified by \( x \dwn \colon U \to \pow X \).
\end{proof}

As an additional byproduct of \jref{lem:dwn-is-descendent}, we conclude:
\begin{cor} \label{cor:dwn-stable}
  Downward closure objects are stable under pullback in the following sense: for any morphisms \( U' \tox{\delta_{U'}} U \tox{x} X \), there is a pullback square
  \[
    \begin{tikzcd}
      \Dwn_{x'} \pb \ar[r, ""] \ar[d, >->] &[20pt] \Dwn_x \ar[d, >->] \\
      U' \times X \ar[r, "\delta_{U'} \times \id_X"] & U \times X
    \end{tikzcd}
  \]
\end{cor}
\begin{proof}
  Immediate from \jref{lem:dwn-is-descendent} and the 2-of-3 rule for pullback squares.
\end{proof}

The next propositions characterizes the downward closures \( \Dwn_x \tto X \) as precisely those substructures containing top elements, and relates the ``external'' element relation \( \unEin \) among these to the relation \( \unein_X \) on \( X \).\footnote{%
  In \cite[Definition~8.17]{shulman-comparing}, the characterization \jref{propn:dwn-hs}~\ref{item:dwn-hs-only-elems} of elemental substructures is taken to be the \emph{definition}.
}
\prooflabel{proof:dwn-hs}
\begin{propn} \label{propn:dwn-hs}\enumbelow
\begin{enumerate}[(i)]
\item\label{item:dwn-hs-top-elem} \label{item:dwn-hs-first} For any \( x \colon \tm \to X \), the unique element \( \bar x \colon \tm \to \Dwn_x \) with \( \bar x i_{\Dwn_x} = x \) is a top element of the closed substructure \( \Dwn_x \tto X \).

  \noindent
  (Hence by \jref{propn:sub-ewf}, if \( X \) is extensional and well-founded, \( \Dwn_x \) is an hs-object.)
\item \label{item:dwn-hs-closed-is-dwn}
  Conversely, if \( A \tto X \) is a closed substructure with a top element \( \bar x \colon \tm \to A \), then \( Z \approx \Dwn_{\bar x i_A} \).
\item \label{item:dwn-hs-only-elems}
  \label{item:dwn-hs-last}
  If \( X \) is an hs-object and \( x \colon \tm \to X \), then \( \Dwn_x \tto X \) is elemental if and only if \( x \ein \rt_X \).

  \noindent
  (Hence the elemental substructures of \( X \) are exactly those of the form \( \Dwn_x \) with \( x \ein \rt_X \).)
\end{enumerate}

We also have the following generalizations of these statements to generalized elements \( x \):
\begin{enumerate}[(i\('\))]
\item \label{item:dwn-hs-top-elem-gen} \label{item:dwn-hs-first-gen}
  For any \( x \colon U \to X \), \( \und{\Dwn_x} \tto \und{U \times X} \) has a top element \( \bar x \colon (U, \id) \to \und{\Dwn_x} \) with \( \bar x i_{\Dwn_x} = \br{\id, x} \).

    \noindent
    (Hence by \jref{propn:sub-ewf} and \jref{propn:pb-fun-is-logical}, if \( X \) is extensional and well-founded, \( \und{\Dwn_x} \) is an hs-object.)
  \item \label{item:dwn-hs-closed-is-dwn-gen}
    If \( \und A \tto \und{U \times X} \) is closed with a top element \( \bar x \colon \und U \to \und A \), then \( \und A \approx \und{\Dwn_{x}} \), where \( \br{\id, x} = \bar x i_A \).
  \item \label{item:dwn-hs-only-elems-gen} \label{item:dwn-hs-last-gen}
    If \( X \) is an hs-object, and \( x \colon U \to X \), then \( \und{\Dwn_x} \tto \und{U \times X} \) is elemental if and only if \( x \ein \unex \rt_X \).
  \end{enumerate}
\end{propn}

\prooflabel{proof:dwn-respects-ein}
\begin{propn} \label{propn:dwn-respects-ein}
  Given \( x,y \colon U \to X \), if \( x \ein y \), then \( \und{\Dwn_x} \Ein \und{\Dwn_y} \).
  Conversely, if \( X \) is extensional and well-founded, then \( \und{\Dwn_x} \Ein \und{\Dwn_y} \) implies \( x \ein y \).
\end{propn}

\subsection{Adding a top element} \label{subsec:top-element}
When we formed the (large) set \( \rV_{\hs} \) of hs-sets in \sref{subsec:ein-sets}, we obtained an \( \unein \)-set which is extensional and well-founded (as we will prove in \sref{subsec:univ-ext-wf}).
However, it fails to be an hs-set itself because it does not have a top element.
The same situation will arise when we construct the analogue of \( \rV_\hs \) in a 2-topos in \sref{subsec:hs-2-topos}.
As explained in \S\S\ref{subsec:unis-and-paradox}~and~\ref{subsec:ein-sets}, to carry out the Burali-Forti argument that, we want to a version of \( \rV_{\hs} \) which \emph{is} an hs-set.
We do this simply by adding a single element to \( \rV_{\hs} \) and declaring it maximal.
Here, we describe this procedure, in the general context of hs-objects in a topos.
In fact, we describe a somewhat more general procedure, which we will need in \sref{subsec:set-formation-closure}, in which the new element is declared not to be maximal, but only to be greater than a prescribed set \( S \) of elements.

First, we describe the set-theoretic construction:
\begin{settybox}
  \begin{defn}
    Given an \( \unein \)-set \( X \) and a subset \( S \subset X \), we define the \emph{\( S \)-top extension} of \( X \) to be the \( \unein \)-set \( \wh X = X \cup \set{t_{\wh X}} \) (where \( t_{\wh X} \) is some element not in \( X \)) with \( \unein_{\wh X} \) given by \( x \ein_{\wh X} y \) for \( x,y \in \wh X \) iff either \( x,y \in X \) and \( x \ein y \), or \( x \in S \) and \( y = t_{\wh X} \).
  \end{defn}
\end{settybox}

For the topos-theoretic version, we find it more convenient to first describe \( \wh X \) in terms of a universal property; in \jref{propn:top-ext-props}, we give an alternative characterization of \( \wh X \) analogous to the set-theoretic construction just described.

For the rest of this section, let \( X \) be an \( \unein \)-object in a topos \( \bC \).
\begin{defn} \label{defn:top-extension}
  Given a morphism \( i_S \colon S \to X \), we define a \defword{top extension along \( i_S \)} (or \defword{\( i_S \)-top extension}, or just \defword{\( S \)-top extension}) of \( X \) to be an \( \unein \)-object \( \wh X \) together with a consistent morphism \( i_X \colon X \to \wh X \) and a global element \( t_{\wh X} \colon \tm \to \wh X \) satisfying \( i_S i_X \ein \unex_S t_{\wh X} \), and which is universal in the following sense:

  For any \( \unein \)-object \( Y \) with a consistent morphism \( f \colon X \to Y \) and global element \( t_Y \colon \tm \to Y \) with \( i_S f \ein \unex_S t_Y \), there is a unique consistent morphism \( \bar f \colon \wh X \to Y \) with \( i_X \bar f = f \) and \( t_{\wh X} \bar f = t_Y \).
  \[
    \begin{tikzcd}
      S \ar[r, "i_S"] & X \ar[r, "i_X"] \ar[rd, "f"'] & \wh X \ar[d, "\bar f"', dashed] & \tm \ar[l, "t_{\wh X}"'] \ar[ld, "t_Y"] \\
      & & Y
    \end{tikzcd}
  \]
  We will see in \jref{propn:top-ext-props} that an \( S \)-top extension always exists, and it is obviously uniquely determined up to isomorphism.
  When \( i_S = \id_X \tolon X \to X \), then we also call an \( S \)-top extension simply a \defword{top extension}.
\end{defn}

\begin{rmk}
  Though it is convenient to state the definition in the above generality, we will in fact only employ it in the case where \( i_S \) is a monomorphism (and in fact a \emph{dense} monomorphism in the sense of \jref{defn:dense} below).

  In fact, taking the image \( \im(i_S) \tto X \), one immediately sees that an \( S \)-top extension is the same thing as an \( \im(i_S) \)-top extension (since the condition \( i_S f \ein \unex t_Y \) is equivalent to \( i_{\im(i_S)} f \ein \unex t_Y \)).
\end{rmk}

\begin{propn} \label{propn:top-ext-props}
  For any underobject \( S \tto X \), there exists an \( S \)-top extension \( \wh X \) of \( X \).

  Moreover:
  \begin{enumerate}[(i)]
  \item\label{item:top-ext-props-first} the morphism \( i_X \colon X \to \wh X \) is monic and initial (hence also conservative)
  \item\label{item:top-ext-props-coprod} \( X \tox{i_X} \wh X \xot{t_{\wh X}} \tm \) is a coproduct diagram
  \item\label{item:top-ext-props-el} \label{item:top-ext-props-last} given \( x \colon U \to \wh X \), we have \( x \ein \unex t_{\wh X} \) iff \( x \) factors through the underobject \( S \tox{i_S} X \tox{i_X} \wh X \).
  \end{enumerate}
\end{propn}
\begin{proof}
  See \hyperref[proof:top-ext-props]{\S\ref*{subsec:top-element-proofs}}.
\end{proof}

\begin{defn} \label{defn:dense}
  We say that \( S \tto X \) is \defword{dense} if its closure is \( X \), i.e., if \( \top_S \exists_{i_S} \dwn = \top_X \colon \tm \to \pow X \), or equivalently, if \( X \) is the least closed underobject of \( X \) containing \( S \).
\end{defn}

\begin{propn} \label{propn:top-ext-is-hs-struct}
  If \( X \) is extensional and well-founded and \( S \tto X \) is a dense underobject, then the \( S \)-top extension \( \wh X \) is an hs-object with top element \( t_{\wh X} \colon \tm \to \wh X \).
\end{propn}
\begin{proof}
  See \hyperref[proof:top-ext-is-hs-struct]{\S\ref*{subsec:top-element-proofs}}.
\end{proof}

\subsection{\alttext{\( \unein \)}{ε}-objects in slice topoi} \label{subsec:ein-in-slice}
Below, we will often need to consider \( \unein \)-objects in the slice categories \( \bC / U \) of a given topos \( \bC \).
Set-theoretically, this corresponds to a family of \( \unein \)-sets indexed by a given a set \( U \).
We will also need to relate \( \unein \)-objects in different slices.
Here, we collect some general results concerning these notions.

Recall our conventions regarding slice categories from \sref{subsubsec:under-slice-prelims}.
For the rest of this section, let \( f \colon A \to B \) be a morphism in a topos \( \bC \), and let \( \und X \in \bC / A \) and \( \und Y \in \bC / B \) be \( \ein \)-objects.
The proofs of all the statements in this section can be found in \sref{subsec:ein-in-slice-proofs}.

\begin{defn}
  We say that a morphism \( g \colon X \to Y \) over \( f \) in \( \bC \) is \defword{relatively consistent (over \( f \))}  if there exists a (necessarily unique) dotted morphism making the following square commute:
  \begin{equation} \label{eq:rel-consis-square}
    \begin{tikzcd}
      \ein_X \ar[r, dashed] \ar[d, "i_{\unein_X}"', >->] & \ein_Y \ar[d, "i_{\unein_Y}", >->] \\
      X \times_A X \ar[r, "g \times_f g"] & Y \times_B Y
    \end{tikzcd}
\end{equation}
(We recall from \sref{subsubsec:under-slice-prelims} that \( i_{\ein_Y} \), being a mono in \( \bC / B \), is also monic in \( \bC \).)
\end{defn}

\begin{rmk} \label{rmk:ein-a-notation}
  Relative consistency can be described in terms of generalized elements as follows.
  Given \( U \in \bC \) and \( x,y \colon U \to X \), let us write \( x \ein^A y \) to express that \( x \delta_X = y \delta_Y \colon U \to A \), and that the resulting morphism \( \br{x,y} \colon U \to X \times_A X \) factors through \( \unein_{X} \); this is the same as saying that the resulting morphisms \( x,y \colon (U, \delta_U) \to \und X \) in \( \bC / A \) satisfy \( x \ein_{X} y \), where \( \delta_U = x \delta_X = y \delta_Y \).

  Then the relative consistency of \( (f, u) \) means that for all \( x,y \colon U \to X \), if \( x \ein^A y \), then \( x u \ein^B y u \).
\end{rmk}

\begin{defn}
  We say that a morphism \( g \colon X \to Y \) over \( f \) in \( \bC \) is \defword{relatively conservative (over \( f \))} if for each \( U \in \bC \) and each \( x,y \colon U \to X \) with \( x \delta_X = y \delta_X \colon U \to A \), if \( x f \ein^B y f \), then \( x \ein^A y \).
  We say that \( g \) is \defword{strongly relatively conservative (over \( f \))} if for any \( x,y \colon U \to X \), \( x f \ein^B y f \) implies \( x \ein^A y \) (without assuming \( x \delta_X = y \delta_X \)).

  Note that if \( f \) is monic, then \( g \) is relatively conservative iff it is strongly relatively conservative.
\end{defn}

Given a topos \( \bC \), there is a category \( \bC^\to_{\ein} \) whose objects are \( \unein \)-objects in the slices \( \bC / A \) of \( \bC \) for varying \( A \), and whose morphisms are the relatively consistent morphisms.
This comes with an obvious forgetful functor \( \bC^\to_{\ein} \to \bC^\to \).
The upshot of the following two propositions is that \( \cod \colon \bC^\to_{\ein} \to \bC \) is a bifibration, and that \( \bC^\to_{\ein} \to \bC^\to \) is a morphism of bifibrations.

\prooflabel{proof:ein-pb-conds}
\begin{propn} \label{propn:ein-pb-conds}
  Given a morphism \( g \colon X \to Y \) over \( f \) in \( \bC \), the following are equivalent:
  \begin{enumerate}[(i)]
  \item\label{item:ein-pb-conds-pbs} \( g \) is relatively consistent over \( f \)---i.e., a dashed morphism as in \eqref{eq:rel-consis-square} exists---and the squares
    \begin{equation} \label{eq:ein-pb-square}
      \begin{tikzcd}
        X \ar[r, "g"] \ar[d, "\delta_X"'] & Y \ar[d, "\delta_Y"] \\
        A \ar[r, "f"] & B
      \end{tikzcd}
    \end{equation}
    and \eqref{eq:rel-consis-square} are both pullback squares.
  \item\label{item:ein-pb-conds-rel-consv} The square \eqref{eq:ein-pb-square} is a pullback square, and \( g \) is relatively consistent and relatively conservative over \( f \).
  \item\label{item:ein-pb-conds-univ} \( g \) is universal among relatively consistent morphisms to \( Y \) over \( f \), in the following sense:
    for any \( \unein \)-object \( \und W \) in \( \bC / A \) and any relatively consistent \( h \colon W \to Y \) over \( f \), there is a unique consistent morphism \( \bar h \colon \und W \to \und X \) in \( \bC / A \) with \( \bar h g = h \).
  \end{enumerate}
\end{propn}

\begin{defn} \label{defn:ein-cartesian}
  We say that a morphism \( g \colon X \to Y \) over \( f \) is a \defword{cartesian morphism of \( \unein \)-objects (over \( f \))} (or is \defword{\( \unein \)-cartesian}), and that the the square \eqref{eq:ein-pb-square} in \( \bC \) is a \defword{pullback square of \( \unein \)-objects} (or an \defword{\( \unein \)-pullback square}) if the equivalent conditions of \jref{propn:ein-pb-conds} are fulfilled.

  We note that by \jref{propn:ein-pb-conds}~\ref{item:ein-pb-conds-pbs}, given the \( \unein \)-object \( \und Y \in \bC / B \) and morphism \( f \colon A \to B \), the existence of an \( \unein \)-cartesian morphism \( g \colon X \to Y \) over \( f \) is guaranteed; indeed, the \( \unein \)-object \( \und X \in \bC / A \) is the just the image of the \( \unein \)-object \( \und Y \in \bC / B \) under the logical anafunctor \( f^* \colon \bC / B \to \bC / A \).
  In general, we will denote such a pullback \( (\und X, \ein_{\und X}) \) by \( (f^* \und Y, \ein_{f^* \und Y}) = \pbig{(f^* Y, \delta_{f^* Y}), \ein_{f^* \und Y}} \).
\end{defn}

\prooflabel{proof:ein-cocart-conds}
\begin{propn} \label{propn:ein-cocart-conds}
  Given a morphism \( g \colon X \to Y \) over \( f \) in \( \bC \), the following are equivalent:
  \begin{enumerate}[(i)]
  \item\label{item:ein-cocart-conds-iso} \( g \) is relatively consistent over \( f \)---i.e., a dashed morphism as in \eqref{eq:rel-consis-square} exists---and the morphism \( g \) as well as the dashed morphism in \eqref{eq:rel-consis-square} are both isomorphisms in \( \bC \).
  \item\label{item:ein-cocart-conds-strong} \( g \) is an isomorphism in \( \bC \), and is relatively consistent and strongly relatively conservative over \( f \).
  \item\label{item:ein-cocart-conds-univ} \( g \) is universal among relatively consistent morphisms to \( X \) over \( f \), in the following sense:
    for any \( \unein \)-object \( \und Z \) in \( \bC / B \) and any relatively consistent \( h \colon X \to Z \) over \( f \), there is a unique consistent morphism \( \und h \colon \und Y \to \und Z \) in \( \bC / B \) with \( g \und h = h \).
  \end{enumerate}
\end{propn}

\begin{defn} \label{defn:ein-cocart}
  We say that a morphism \( g \colon X \to Y \) in \( \bC \) is a \defword{cocartesian morphism of \( \unein \)-objects (over \( f \))} (or an \defword{\( \unein \)-cocartesian morphism}) if the equivalent conditions of \jref{propn:ein-cocart-conds} hold.

  We note that by \jref{propn:ein-cocart-conds}~\ref{item:ein-cocart-conds-iso}, given the \( \unein \)-object \( \und X \in \bC / A \) and morphism \( f \colon A \to B \), the existence of an \( \unein \)-cocartesian morphism \( g \colon X \to Y \) over \( f \) is guaranteed; indeed, we may take \( Y = X \) and \( g = \id_X \), and (fixing a product \( Y \times_B Y \)) let \( \unein_Y \) be the composite \( \unein_X \tox{i_{\ein_X}} X \times_A X \tox{\id \times_f \id} X \times_B X \), noting that \( \id \times_f \id \) is a monomorphism.
  In general, we will denote such a cocartesian lift by \( c \colon X \to X_U \), and the resulting \( \unein \)-object in \( \bC / A \) by \( (\und X_U, \ein_{\und X_U}) = \pbig{(X_U, \delta_{X_U}), \ein_{\und X_U}} \).
\end{defn}

The relatively consistent morphisms \( X \to Y \) of \( \unein \)-objects over \( f \colon A \to B \) are thus in bijection with the consistent morphisms \( \und X \to f^* \und Y \) in \( \bC / A \) and with the consistent morphisms \( \und X_B \to \und Y  \) in \( \bC / B \).
We next consider the preservation of various properties under these bijections.
\begin{defn} \label{defn:rel-init-elem}
  Fix a morphism \( g \colon X \to Y \) in \( \bC \) over \( f \).
  \begin{enumerate}[(i)]
  \item We say \( g \) is \defword{relatively initial (over \( f \))} if it is relatively consistent and for each \( U \in \bC \) and \( x_1 \colon U \to X \) and \( y_0 \colon U \to Y \) with \( y_0 \ein^B x_1 g \), there is \( U \oot U' \tox{x_0} X \) with \( x_0 \ein^A x_1 \) and \( x_0 g = y_0 \).
  \item We say \( g \) is \defword{relatively elemental (over \( f \))} if it is relatively initial over \( f \) and \( \rt_X g \ein^B f \rt_Y \) whenever \( \rt_X \colon (A, \id) \to \und X \) and \( \rt_Y \colon (B, \id) \to \und Y \) are top elements of \( \und X \) and \( \und Y \), respectively.
  \end{enumerate}
\end{defn}

\prooflabel{proof:rel-notions-compare}
\begin{propn} \label{propn:rel-notions-compare}
  Consider a commutative diagram as below, with \( c \) \( \unein \)-cocartesian over \( f \) and \( \pi_f \) \( \unein \)-cartesian over \( f \).
  \[
    \begin{tikzcd}
      X \ar[rrr, "c"] \ar[rdd, "\delta_X"', bend right] \ar[rrd, "g"] \ar[rd, "\bar g"']
      & & & X_B \ar[ddl, bend left] \ar[dl, "\und g"] \\
      & f^* Y \ar[d] \ar[r, "\pi_f"'] & Y \ar[d, "\delta_Y"] & \\
      & A \ar[r, "f"] & B &
    \end{tikzcd}
  \]
  Then:
  \begin{enumerate}[(i)]
  \item\label{item:rel-notions-compare-consv} \( \und g \) is conservative \( \iff \) \( g \) is strongly relatively conservative \( \To \) \( g \) is relatively conservative \( \iff \) \( \bar g \) is conservative
  \item\label{item:rel-notions-compare-init} \( \und g \) is initial \( \iff \) \( g \) is relatively initial \( \iff \) \( \bar g \) is initial
  \item\label{item:rel-notions-compare-elem}
    Assuming \( \und Y \) is well-founded:\footnote{%
      The well-foundedness is only needed to ensure that \( Y \), and hence \( f^* Y \), has at most one top element; and at any rate, we are only interested in the case where \( Y \) (as well as \( X \)) is in fact an hs-object.
    }
    \( g \) is relatively elemental \( \iff \) \( \bar g \) is elemental
  \end{enumerate}
\end{propn}

\begin{cor} \label{cor:rel-init-2-of-3}
  Relatively initial morphisms are closed under composition, and satisfy the 2-of-3 property: given \( g \colon X \to Y \) and \( k \colon Y \to Z \) over \( f \colon A \to B \) and \( h \colon B \to C \), respectively, if \( g k \) and \( k \) are relatively initial over \( f h \) and \( h \), respectively, and if \( \und Y \in \bC / B \) is extensional and well-founded, then \( g \) is relatively initial over \( f \).
\end{cor}

\begin{proof}
  This follows from contemplating the diagram
  \[
    \begin{tikzcd}
      X \ar[rd, dashed] \ar[rrddd, bend right=10pt] \ar[rrrrd, "g", bend left=10pt] \\
      & f^* Y \ar[rrd, dashed] \ar[rdd] \ar[rrr] &[-30pt] &[-20pt] &[-30pt] Y \ar[rrd, dashed] \ar[ddr] \ar[rrrd, "k", bend left=10pt] &[-30pt] &[-20pt] \\
      & & & f^* h^* Z \ar[rrr, crossing over] \ar[dl] \pb & & & h^* Z \ar[r] \ar[ld] \pb & Z \ar[d] \\
      & & A \ar[rrr, "f"] & & {} & B \ar[rr, "h"] & & C
    \end{tikzcd}
  \]
  and applying \jrefss{propn:rel-notions-compare}{propn:pb-fun-is-logical}{propn:init-2-of-3}.
\end{proof}

Finally, we compare the properties of \( \unein \)-objects in \( \bC / A \) and \( \bC / B \).
Since \( f^* \colon \bC / B \to \bC / A \) is logical, ``all'' properties of \( \und Y \in \bC / B \) (e.g., being extensional, well-founded, having a top element) are inherited by \( f^* \und Y \in \bC / A \).
In the other direction, we have the following.

\begin{defn} \label{defn:loc-ext}
  An \( \unein \)-object \( Z \) in a topos \( \bD \) is \defword{locally extensional} if the \( \unein \)-object \( \und{\Dwn_z} \) in \( \bD / U \) is extensional for each \( z \colon U \to Z \).
\end{defn}

\prooflabel{proof:loc-ext-char}
\begin{lem} \label{lem:loc-ext-char}
  An \( \unein \)-object \( Z \) in a topos \( \bD \) is locally extensional if and only if: for any \( x,y,z \colon U \to Z \), if \( x, y \preceq z \) and \( w \ein x' \ToT w \ein y' \) for all \( w \colon U' \to Z \), then \( x = y \).
\end{lem}

\prooflabel{proof:sum-props}
\begin{propn} \label{propn:sum-props}
  Suppose \( c \colon X \to X_B \) is an \( \unein \)-cocartesian morphism over \( f \).
  Then:
  \begin{enumerate}[(i)]
  \item If \( \und X \) is well-founded, so is \( \und X_B \)
  \item If \( \und X \) is extensional, then \( \und X_B \) is locally extensional
  \item If \( i_S \colon (S, \delta_S) \tto \und X \) is a dense underobject in \( \bC / A \) in the sense of \jref{defn:dense}, then \( i_S c \colon (S, \delta_S f) \tto \und X_B \) in \( \bC / B \) is again dense.
  \end{enumerate}
\end{propn}

\section{Universes of hereditary sets in 2-topoi} \label{sec:univs-in-2-topoi}
We now undertake the construction of the groupoid \( \hs(\rS) \) of \( \ein \)-objects in the DOF classifier \( \rS \in \cC \) in a 2-topos \( \cC \), which, as explained in the introduction, will serve as our model of MK, and also as our starting point for carrying out the Burali-Forti paradox.
As also explained there, we will more generally construct an object \( (E^\cT)^\iso \) of models for any finite higher-order theory \( \cT \); and we will do this by associating to any such \( \cT \) a \emph{topos sketch} (see below) having the same models as \( \cT \).

\subsection{Higher-order theory cotensors} \label{subsec:thy-cotensors}
We refer to \sref{subsec:logic-in-topoi} for the notion of a \emph{higher-order theory} \( \cT \) and of the \emph{groupoid of models} \( (\bC^\cT)^\iso \) of a higher-order theory \( \cT \) in a topos \( \bC \).

\begin{defn}
  Let \( \cC \) be a 2-category, let \( E \in \cC \) be an internal 1-topos, and let \( \cT \) be a higher-order theory.
  A \defword{strict cotensor} of \( E \) by \( \cT \) is a groupoid \( (E^\cT)^\iso \) together with a model \( \ev = \ev_{\cT,E} \colon \cT \to \cC\pbig{(E^{\cT})^\iso,E} \) of \( \cT \) in the topos \( \cC\pbig{(E^{\cT})^\iso,E} \) such that, for each groupoid \( A \in \cC \), the induced functor \( \cC\pbig{A,(E^{\cT})^\iso} \to \pbig{\cC(A,E)^{\cT}}^\iso \) taking \( f \) to \( \ev \circ f^* \) is an \emph{isomorphism}.
  A \defword{semi-strict cotensor} is defined in the same way, except that the functor \( \cC\pbig{A,(X^{\cT})^\iso} \to \cC\pbig{(A,X)^{\cT}}^\iso \) is only required to be a surjective-on-objects equivalence.
\end{defn}
Usually, we will omit ``semi-strict'', as this is the only kind of cotensor by higher-order theories we will be using (in particular, we will not be using the strict version for anything).

In \hyperref[proof:theory-cotensors]{\S\ref*{subsec:topos-sketches}}, we prove:
\begin{thm} \label{thm:theory-cotensors}
  If \( \cC \) is a corepita 2-category and \( \cT \) is a finite higher order theory, then any internal topos \( E \in \cC \) admits a semi-strict cotensor \( (E^{\cT})^\iso \).
\end{thm}
The proof consists of two steps.
One is a ``translation'' of models of ``syntactic'' theories into models of \emph{sketches} of an appropriate kind, in the manner described in \cite[\S9]{makkai-generalized-sketches-all}.
Once this is done, it then remains to show that cotensors by sketches of the appropriate kind exist, as was done in \cite[Appendix~A]{helfer-topoi-in-topoi} in the case of finite-limit sketches.
These two steps appear below as \jrefs{propn:thy-to-sketch}{propn:topos-sketch-cotensors}.

In \S\ref{subsec:hs-2-topos}, we will construct not only the groupoid \( \hs(E) \) of hs-sets in an internal topos \( E \), but also the groupoid \( \unein_E \to \hs(E) \times \hs(E) \) of elemental morphisms in \( E \).
For this purpose, we will also need to compare the cotensors of \( E \) by two different theories, for which we need the notion of an \emph{interpretation} of one theory into another.
We will only need a very simple special case of this notion, and we only define the version we will need.
\begin{defn} \label{defn:thy-embedding}
  Given higher-order theories \( \cT_0 \) and \( \cT_1 \) over languages \( \cL_0 \) and \( \cL_1 \), an \defword{embedding} \( i \colon \cT_0 \to \cT_1 \) consists of an injective map \( i_{\Ob} \colon \Ob \cL_0 \to \Ob \cL_1 \) and injective maps \( i_{\Fun} \colon \Fun \cL_0 \to \Fun \cL_1 \) and \( i_{\Rel} \colon \Rel \cL_0 \to \Rel \cL_1 \) respecting arities (relative to the map \( i \) on sorts) such that \( i(\phi) \in \cT_1 \) for each sentence \( \phi \in \cT_0 \) where \( i(\phi) \) is the sentence over \( \cL_1 \) resulting from applying \( i_{\Fun} \) and \( i_{\Rel} \) to the function and relation symbols appearing in \( \phi \).

An embedding \( i \colon \cT_0 \to \cT_1 \) can be composed with an interpretation \( M \colon \cT_1 \to \bC \) in a topos \( \bC \) to obtain an interpretation \( i M \colon \cT_0 \to \bC \), and this induces a functor \( i^* \colon \bC^{\cT_1} \to \bC^{\cT_0} \).

Given semi-strict cotensors \( E^{\cT_0},E^{\cT_1} \) for some internal topos \( E \) in a 2-category \( \cC \), the model \( \cT_0 \tox{i} \cT_1 \tox{\ev_{\cT_1,x}} \cC(E^{\cT_1},E) \) induces, by the universal property of \( E^{\cT_1} \) a morphism \( E^i \colon E^{\cT_1} \to E^{\cT_0} \), which is uniquely determined up to isomorphism.
\end{defn}

\subsection{Topos sketches} \label{subsec:topos-sketches}
In this section, we carry out the outline of the proof of \jref{thm:theory-cotensors} indicated above: we define the notion of \emph{topos sketch}, we define a translation of higher-order theories into topos sketches, and we prove the existence of cotensors by topos sketches.

We recall from \cite[Appendix~A]{helfer-topoi-in-topoi} that a \emph{pbt sketch} \( J \) consists of a category \( \abs{J} \) together with a set of commutative squares in \( \abs{J} \) and a set of objects in \( \abs{J} \), to be thought of as ``intended'' pullback squares and terminal objects, respectively.
A \emph{model} of \( J \) is then a functor from \( \abs{J} \) taking the specified squares and objects to actual pullback squares and terminal objects in the target category.
In \cite{makkai-generalized-sketches-all}, it is observed that for \emph{any} property of diagrams in a category, one can define a corresponding kind of sketch, and a general framework is established for introducing such definitions (as well as for making \emph{deductions} about the diagrams, which is the main point of that paper).
Here, we will only consider the particular kind of sketches we will be using, namely sketches for elementary topoi (which are also considered in \cite[p.~198]{makkai-generalized-sketches-all}).

To begin with, in addition to specifying finite limits, we should also specify \emph{power objects} in \( \abs{J} \).
Such a specification should consist of a diagram in \( \abs{J} \) of the shape
\begin{equation} \label{eq:power-diagram}
  \begin{tikzcd}[column sep=-5pt]
    & \in_X \ar[d, "i"] \\
    & X \times \pow X \ar[dl, "\pi_0"'] \ar[dr, "\pi_1"] \\
    X & & \pow X.
  \end{tikzcd}
\end{equation}
We can then define a \emph{model} of such a sketch to be a model of the underlying pbt sketch in which, additionally, each such diagram is taken to a genuine power object diagram, i.e., the span \( (\pi_0,\pi_1) \) is taken to a product diagram, and the morphism \( i \) to a universal relation (see \sref{subsubsec:topoi-prelims}).

In theory, this should in fact suffice as a notion of ``topos sketch'', since a topos is nothing but a category with finite limits and power objects---or, more to the point, since all other operations of interest in a topos, and in particular all the ``logical'' operations, can be derived from finite limits and power objects.
However, it would be highly impractical to proceed in this way, and it is better to simply endow a topos sketch with additional information corresponding to all these operations (cf.\ also the discussion of ``explicit definition'' in \cite[p.~193]{makkai-generalized-sketches-all}).
Namely, we want to include the specification of finite products, terminal and initial objects, and the logical operations (on underobjects) of conjunction, disjunction, implication, and of universal and existential quantification.

In order to organize this data, it is useful to make the following preliminary definition, in the spirit of the ``operations on diagrams'' of \cite[\S1.1]{makkai-ultraproducts} (an idea which was also invoked in \cite[Appendix~B]{helfer-topoi-in-topoi}).
We also note that we could use graphs as diagram-shapes (as in \cite{makkai-ultraproducts,makkai-generalized-sketches-all}) throughout, but we find it convenient to instead use categories (as was done in \cite{helfer-topoi-in-topoi}).

\begin{defn} \label{defn:topos-op}
  A \defword{topos operation} \( \omega \) is a triple \( (L_\omega, K_\omega, \sigma_\omega) \), where \( L_\omega \) is a category (always finite in the examples we will consider), \( i_{K_\omega} \colon K_\omega \hto L_\omega \) is a full subcategory, and \( \sigma_\omega = \set{\sigma_{\omega,\bC}}_\bC \) is a collection of diagrams \( \sigma_{\omega,_\bC} \subset \Ob \bC^{L_\omega} \) of shape \( L_\omega \) in each topos%
  \footnote{
    That is, each topos in some fixed universe.
  }
  \( \bC \), satisfying the following conditions:
  \begin{enumerate}[(i)]
  \item \label{item:topos-op-ext-invt}
    for each the topos \( \bC \), the set \( \sigma_{\omega,_\bC} \) is closed under isomorphism of diagrams,
  \item \label{item:topos-op-ext-ex}
    any diagram \( f \colon K_\omega \to \bC \) of shape \( K_\omega \) in a topos \( \bC \) extends to a diagram \( g \in \sigma_{\omega,\bC} \subset \Ob \bC^{L_\omega} \), and
  \item \label{item:topos-op-ext-uniq}
    given two diagrams \( f_0, f_1 \colon K_\omega \to \bC \) in a topos \( \bC \), and extensions \( g_0, g_1 \in \sigma_{\omega,\bC} \subset \Ob \bC^{L_\omega} \) of each of them, any isomorphism \( f_0 \toi f_1 \) extends to a unique isomorphism \( g_0 \toi g_1 \).
  \item \label{item:topos-op-ext-log-pres}
    given a logical functor \( F \colon \bC \to \bD \), if \( g \in \sigma_{\omega, \bC} \), then \( g F \in \sigma_{\omega, \bD} \).
  \end{enumerate}
  We define the following \defword{fundamental topos operations} \( \omega_0,\ldots,\omega_7 \) and \( \omega_{8+k} \) for \( k \ge 0 \):
  \begin{itemize}
  \item \( L_{\omega_0} \) is a poset of shape \eqref{eq:power-diagram}; \( K_{\omega_0} \) is the discrete subcategory on the objects labelled ``\( X \)'' and ``\( \pow X \)''; and \( \sigma_{\omega_0,\bC} \) is the set of power diagrams in \( \bC \).
  \item \( L_{\omega_1} \) is a commutative square; \( K_{\omega_1} \) is the underlying cospan; and \( \sigma_{\omega_1,\bC} \) is the set of pullback squares in \( \bC \).
  \item \( L_{\omega_2} \) is a single object; \( K_{\omega_1} \) is empty; and \( \sigma_{\omega_2,\bC} \) is the set of initial objects in \( \bC \).
  \item \( L_{\omega_i} \) for \( i = 3,4 \) are both the poset \( 0 \to 1 \to 2 \ot 3 \); the \( K_{\omega_i} \) are the sub-poset \( 0 \to 1 \to 2 \); and the \( \sigma_{\omega_i,\bC} \) are the sets of diagrams \( P \tox{i} X \tox{f} Y \xot{j} Q \), where \( i \) and \( j \) are monos, and \( Q \) is the universal quantification \( Q = \forall_f P \) or existential quantification \( Q = \exists_F P \) of \( P \) along \( f \), respectively.
  \item \( L_{\omega_i} \) for \( i = 5,6,7 \) are all the category consisting of three arrows \( f_0,f_1,f_2 \) with common codomain \( 3 \) and distinct domains \( 0,1,2 \); the \( K_{\omega_i} \) are all the full subcategory on \( 0,1,3 \); and the \( \sigma_{\omega_i,\bC} \) are the sets of triples \( (f_0,f_1,f_2) \) of monomorphisms \( f_i \colon P_i \tto X \) in \( \bC \), where \( P_2 \) is the conjunction \( P_2 = P_0 \wedge P_1 \), disjunction \( P_2 = P_0 \vee P_1 \), or implication \( P_2 = P_0 \To P_1 \) of \( P_0 \) and \( P_1 \), respectively.
  \item For each \( n \ge 0 \), \( L_{\omega_{8+n}} \) is the category consisting of \( n \) arrows \( f_1,\ldots,f_{n} \) with common codomain \( 0 \) and distinct domains \( 1,\ldots,n \); \( K_{\omega_{8+n}} \) is the discrete subcategory on \( \set{1,\ldots,n} \); and \( \sigma_{\omega_i,\bC} \) is the set of \( n \)-ary product diagrams in \( \bC \) (thus, of terminal objects for \( n = 0 \) and isomorphisms for \( n = 1 \)).
  \end{itemize}
  These are indeed all topos operations: in each case, \ref{item:topos-op-ext-invt} and \ref{item:topos-op-ext-uniq} are instances of the usual isomorphism-invariance of objects defined by universal properties, \ref{item:topos-op-ext-ex} is simply the statement that topoi \emph{do} have all of these operations, and \ref{item:topos-op-ext-log-pres} follows from \jrefs{propn:logical-fun-pres-colimits}{propn:logical-fun-pres-satisf}.

  We denote by \( \Omega_\fund = \set{\omega_i}_{i \ge 0} \) the set of all fundamental topos operations.
\end{defn}

\begin{defn}
  Let \( \Omega \) be a set of topos operations.
  An \defword{\( \Omega \)-topos sketch} is a category \( \abs{J} \) together with a set \( J_\omega \subset \Ob \abs{J}^{L_{\omega}} \) of diagrams in \( \abs{J} \) of shape \( L_{\omega} \) for each of the topos operations \( \omega \in \Omega \); we call the elements of \( J_\omega \) the set of \emph{specifications of type \( \omega \)} in \( J \).
  An \( \Omega \)-topos sketch is \defword{finite} if the category \( \abs{J} \) and each of the categories \( K_\omega \) and \( L_\omega \) is of finite presentation (see \sref{subsec:cat-prelims}), each of the sets \( J_\omega \) is finite, and only finitely many of the sets \( J_\omega \) are inhabited.
  By a \defword{topos sketch}, we mean an \( \Omega_\fund \)-topos sketch.

  A \defword{model} for an \( \Omega \)-topos sketch \( J \) in a topos \( \bC \) is a functor \( F \colon \abs{J} \to \bC \) such that \( f F \in \sigma_{\omega, \bC} \subset \Ob \bC^{L_{\omega}} \) for each \( f \in J_\omega \subset \abs{J}^{L_{\omega}} \).
  We write \( (\bC^J)^\iso \) for the full subcategory of \( (\bC^{\abs{J}})^\iso \) consisting of models.\footnote{%
    As with higher-order theory cotensors, the notation \( (\bC^J)^\iso \) is atomic, rather than being the core of some object \( E^J \); see \prerefpost{footnote~}{footnote:cotensor-atomic}{} on \pref{footnote:cotensor-atomic}.%
  }
  As with the categories of models for a theory (see \sref{subsec:logic-in-topoi}), logical functors \( F \colon \bC \to \bD \) and invertible 2-cells \( \beta \colon \bC \tocell \bD \) between such induce functors \( F_* \colon (\bC^J)^\iso \to (\bD^J)^\iso \) and natural isomorphisms \( \beta_* \colon (\bC^J)^\iso \tocell (\bD^J)^\iso \), respectively.

  If \( \cC \) is a 2-category and \( E \in \cC \) an internal topos, then a \defword{strict cotensor} of \( E \) by \( J \) is a groupoid \( (E^J)^\iso \in \cC \) together with a model \( \ev = \ev_{J,E} \colon J \to \cC\pbig{(E^J)^\iso,E} \) of \( J \) in the topos \( \cC(E^J,E) \) such that, for each groupoid \( A \in \cC \), the functor
  \[
    \cC\pbig{A,(E^J)^\iso} \to \pbig{\cC(A,E)^{\abs{J}}}^\iso
  \]
  defined on objects and morphisms by \( f \mapsto \ev \circ f^* \), is an isomorphism onto the subcategory \( \pbig{\cC(A,E)^J}^\iso \subset \pbig{\cC(A,E)^{\abs{J}}}^\iso \).
\end{defn}

Because any logical functor \( F \colon \bC \to \bC' \) between topoi preserves all of the fundamental topos operations, it follows that for any model \( M \colon \abs{J} \to \bC \) of a topos sketch \( J \) in \( \bC \), the composite \( M F \colon \abs{J} \to \bC' \) is a model of \( J \) in \( \bC' \).

Now the two aforementioned ingredients in the proof of \jref{thm:theory-cotensors} are:
\begin{propn} \label{propn:topos-sketch-cotensors}
  If \( \Omega \) is a set of topos operations, \( J \) is a finite \( \Omega \)-topos sketch, and \( \cC \) is a corepita 2-category, then any internal topos \( E \in \cC \) admits a strict cotensor \( (E^J)^\iso \) by \( J \).
\end{propn}
\begin{proof}
  See \sref{subsec:topos-sketch-cotensors-proof}.
\end{proof}

\begin{propn} \label{propn:thy-to-sketch}
  For any (finite) higher-order theory \( \cT \), there exists a (finite) topos sketch \( J \) and surjective equivalences \( F_{\bC} \colon (\bC^J)^\iso \to (\bC^\cT)^\iso \) for each topos \( \bC \), which are natural in the sense that, for any two topoi \( \bC \) and \( \bD \) and any invertible 2-cell \( \alpha \colon \bC \tocell \bD \) between logical functors, the following square of 2-cells strictly commutes.
  \begin{equation} \label{eq:cotensor-comparison-square}
    \begin{tikzcd}
      (\bC^J)^\iso \ar[d, "F_\bC"']
      \ar[r, bend left=15pt, ""' {name=fs}]
      \ar[r, bend right=15pt, "" {name=gs}]
      \ar[from=fs, to=gs, shorten=-2pt, Rightarrow, "\alpha_*"]
      & (\bD^J)^\iso \ar[d, "F_{\bD}"] \\
      (\bC^{\cT})^\iso
      \ar[r, bend left=15pt, ""' {name=fs}]
      \ar[r, bend right=15pt, "" {name=gs}]
      \ar[from=fs, to=gs, shorten=-2pt, Rightarrow, "\alpha_*"]
      & (\bD^{\cT})^\iso \\
    \end{tikzcd}
  \end{equation}
\end{propn}
\begin{proof}
  See \sref{subsec:thy-to-sketch-proof}.
\end{proof}

Assuming these two propositions, let us prove \jref{thm:theory-cotensors}:
\begin{proof}[Proof of \jref{thm:theory-cotensors}.] \label{proof:theory-cotensors}
  Within this proof, for readability, we drop the superscripts \( (-)^\iso \) from the notation \( (E^J)^\iso, (E^\cT)^\iso, (\bC^J)^\iso, (\bC^\cT)^\iso \) for cotensors of an internal topos \( E \) by---or the groupoid of models in a topos \( \bC \) of---a topos sketch \( J \) or higher-order theory \( \cT \).

  We are given a finite higher-order theory \( \cT \) and an internal topos \( E \) in a corepita 2-category \( \cC \), and we need to construct a cotensor \( E^{\cT} \).
  Let \( J \) and \( F_{\bC} \colon \bC^J \to \bC^{\cT} \) for each topos \( \bC \) be as in \jref{propn:thy-to-sketch}, and fix a cotensor \( \ev_{J,E} \colon J \to \cC(E^J,E) \) of \( E \) by \( J \), which exists by \jref{propn:topos-sketch-cotensors}.
  The model \( \ev_{J,E} \in \cC(E^J,E)^J \) of \( J \) in the topos \( \cC(E^J,E) \) induces a model \( \ev_{\cT,E} \defeq F_{\cC(E^J,E)}(\ev_{J,E}) \in \cC(E^J,E)^\cT \) of \( \cT \) in \( \cC(E^J,E) \), which we claim exhibits \( E^J \) as a semi-strict cotensor of \( E \) by \( \cT \).

  Briefly, the reason is that the 2-functors \( \cC(-,E)^\cT, \cC(-,E)^J \colon \cC^\op \to \Cat \) are by assumption naturally isomorphic, so that any object representing one also represents the other.
  In more detail, given a groupoid \( A \in \cC \), we have the composite
  \[
    \cC(A, E^J) \toi \cC(A, E)^J \tox{F_{\cC(A, E)}} \cC(A, E)^\cT,
  \]
  in which the first morphism is an isomorphism, and the second a surjective equivalence.
  It thus suffices to see that this composite is equal to the functor given by \( f \mapsto \ev_{\cT, E} \circ f^* \) (on objects and morphisms), which we wish to show is a surjective equivalence.
  It suffices to show this for morphisms in \( \cC(A, E^J) \), i.e., for 2-cells \( \alpha \colon A \tocell E^J \), and in this case, it follows from chasing \( \id_{E^J} \) around the following commutative diagram:
  \[
    \begin{tikzcd}[baseline=(bl.base)]
      \cC(A,E^J) \ar[r, "\sim"] & \cC(A,E)^J  \ar[r, "F_{\cC(A,E)}"] &[10pt] \cC(A,E)^\cT \\[10pt]
      \cC(E^J, E^J) \ar[r, "\sim"] \ar[u, bend left, ""' {name=f1}] \ar[u, bend right, "" {name=g1}] &
      \cC(E^J,E)^J \ar[r, "F_{\cC(E^J,E)}"] \ar[u, bend left, ""' {name=f2}] \ar[u, bend right, "" {name=g2}] &
      |[alias=bl]| \cC(E,E)^\cT \ar[u, bend left, ""' {name=f3}] \ar[u, bend right, "" {name=g3}].
      \ar[from=f1, to=g1, "\alpha^*", Rightarrow, shorten=-2pt]
      \ar[from=f2, to=g2, "\alpha^*", Rightarrow, shorten=-2pt]
      \ar[from=f3, to=g3, "\alpha^*", Rightarrow, shorten=-2pt]
    \end{tikzcd}
    \qedhere
  \]
\end{proof}

\subsection{\alttext{$\unein$}{ε}-objects in a 2-topos} \label{subsec:hs-2-topos}
We now apply \jref{thm:theory-cotensors}, which provides for cotensors of internal topoi by arbitrary higher-order theories, to the specific task of constructing an internal groupoid of hs-objects in an internal topos.
We begin by describing the universal property of the latter.

\begin{defn} \label{defn:int-topos-hs-classifier}
  Let \( \cC \) be a 2-category and \( E \in \cC \) an internal topos.
  A \defword{pre-hs-classifier} for \( E \) is a groupoidal object \( \hs(E) \in \cC \) equipped with a ``universal'' hs-object \( X_E \in \hs\pBig{\cC\pbig{\hs(E), E}} \) in the topos \( \cC\pbig{\hs(E), E} \), having the following universal property: for any groupoidal \( A \in \cC \), the functor
  \[
    \cC\pbig{A, \hs(E)} \to \hs\pbig{\cC(A, E)},
  \]
  taking \( f \colon A \to \hs(E) \) to the image \( f^* X_E \in \hs\pbig{\cC(A, E)} \) of the hs-object \( X_E \) under the logical functor \( f^* \colon \cC\pbig{\hs(E), E} \to \cC(A, E) \), and acting in the obvious way on morphisms, is an equivalence.

  An \defword{hs-classifier} for \( E \) is a pre-hs-classifier together with a groupoidal object \( \unEin_E \) and a morphism \( i_{\unEin_E} \colon \unEin_E \to \hs(E) \times \hs(E) \) such that, for each groupoidal \( A \in \cC \), the functor \( (i_{\Ein_E})_* \colon \cC(A, \unEin_E) \to \cC\pbig{A, \hs(E) \times \hs(E)} \) is fully faithful with essential image consisting of those pairs \( \br{f,g} \colon A \to \hs(E) \times \hs(E) \) such that \( f^* X_E \Ein g^* X_E \) in \( \cC(A, E) \).
\end{defn}

\begin{propn}
  If \( E \) is an internal topos in a groupoidant 2-category \( \cC \), then any pre-hs-classifier \( \hs(E) \) for \( E \) is setoidal, and for any hs-classifier for \( E \), the object \( \unEin_E \) is likewise setoidal, and the morphism \( i_{\unEin_E} \colon \unEin_E \to \hs(E) \times \hs(E) \) is fully faithful.
\end{propn}

\begin{proof}
  It follows from \jref{propn:v-ext-setoid-ff} that \( \cC\pbig{A, \hs(E)} \) is a setoid for each groupoid \( A \in \cC \), and we have by assumption that \( (i_{\unEin_E})_* \colon \cC(A, \unEin_E) \to \cC\pbig{A, \hs(E)} \times \cC\pbig{A, \hs(E)} \) is fully faithful, and hence that \( \cC(A, \unEin_E) \) is a setoid as well.
  Since \( \cC \) is groupoidant, it then follows from \cite[Theorem~6.2.1]{helfer-topoi-in-topoi} that \( \hs(E) \) and \( \unEin_E \) are setoidal and that \( i_{\unEin_E} \) is fully faithful.
\end{proof}

\begin{propn} \label{propn:hs-classifiers-exist}
  For any corepita 2-category \( \cC \) and internal topos \( E \in \cC \), there exists an hs-classifier \( \hs(E) \) for \( E \).
\end{propn}
\begin{proof}
  See \hyperref[proof:hs-classifiers-exist]{\S\ref*{subsec:topoculi-proofs}}.
\end{proof}

The proof makes use of the following lemma, which provides for the necessary higher-order theories.

\begin{lem} \label{lem:theories-exist}
  There exist finite higher-order theories \( \cT_\hs \), \( \cT_\hs \sqcup \cT_\hs \), and \( \cT_\unEin \), embeddings \( \cT_\hs \tox{i_0,i_1} \cT_{\hs} \sqcup \cT_\hs \tox{i} \cT_\unEin \) (in the sense of \jref{defn:thy-embedding}), and, for each topos \( \bC \), isomorphisms \( F_\bC \), \( G_\bC \), and \( H_\bC \) as follows, satisfying the following conditions \ref{item:theories-exist-conds-first}-\ref{item:theories-exist-conds-last}:
  \begin{enumerate}[(i)]
  \item \( F_\bC \colon (\bC^{\cT_\hs})^\iso \toi \hs(\bC) \)
  \item \( G_\bC \colon (\bC^{\cT_\hs \sqcup \cT_\hs})^\iso \toi \hs(\bC) \times \hs(\bC) \)
  \item \( H_\bC \colon (\bC^{\cT_\unEin})^\iso \toi \unEin_\bC \)
  \item \label{item:theories-exist-nat} \label{item:theories-exist-conds-first} \( f_\bC \) is natural in \( \bC \) in the same sense as in \jref{propn:thy-to-sketch}.
  \item \label{item:theories-exist-prod} The following square strictly commutes for \( j = 0,1 \):
    \[
      \begin{tikzcd}
         (\bC^{\cT_\hs \sqcup \cT_\hs})^\iso \ari[r, "G_\bC"] \ar[d, "i_j^*"'] & \hs(\bC) \times \hs(\bC) \ar[d, "\pi_j"] \\
        (\bC^{\cT_\hs})^\iso \ari[r, "F_\bC"] & \hs(\bC).
      \end{tikzcd}
    \]
  \item \label{item:theories-exist-inc} \label{item:theories-exist-conds-last} The following square strictly commutes:
    \[
      \begin{tikzcd}
         (\bC^{\cT_{\unEin}})^\iso \ari[r, "H_\bC"] \ar[d, "i^*"'] & \unEin_\bC \ar[d, "i_{\unEin_\bC}"] \\
        (\bC^{\cT_\hs \sqcup \cT_\hs})^\iso \ari[r, "G_\bC"] & \hs(\bC) \times \hs(\bC).
      \end{tikzcd}
    \]
  \end{enumerate}
\end{lem}
\begin{proof}
  See \hyperref[proof:theories-exist]{\S\ref*{subsec:topoculi-proofs}}.
\end{proof}

Our next, and main, task is to investigate the properties of the hs-classifier \( \hs(\rS) \) associated to a plentiful DOF classifier \( \rS \)---and first of all, to show that it is an extensional, well-founded \( \unein \)-object.
Now, the object \( \hs(\rS) \) resides in the 1-category \( \SOF(\tm) \) of setoidal objects in \( \cC \) (recall this notation from \sref{subsubsec:dof-sof-digression}).
However, we would like to consider them as lying in a \emph{topos}, for which purpose we now consider a second plentiful DOF classifier \( \rS' \in \cC \) with respect to which it is small, so that it lies in the topos \( \SOF_{\rS'}(\tm) \).
We note that, since \( \rS' \) is plentiful, the \( \rS' \)-smallness of \( \hs(\rS) \) automatically entails that of \( \hs(\rS) \times \hs(\rS) \) and of \( \Ein_{\rS} \).

\begin{thm} \label{thm:v-is-ewf}
  Let \( \rS \in \cC \) be a plentiful DOF classifier (hence an internal topos) in a groupoidant 2-category \( \cC \), and let \( \rS' \) be a plentiful DOF classifier for which \( p \colon \rS_* \to \rS \) and \( \hs(\rS) \) are \( \rS' \)-small.

  Then the \( \unein \)-object \( \hs(\rS) \) in the topos \( \SOF_{\rS'}(\tm) \) is extensional and well-founded (\jref{defn:topos-ein-hs}).
\end{thm}

\noindent
This will be proven in the following sections: via \jref{propn:2-topoi-give-topoculi}, it will be reduced to \jref{propn:v-is-ewf-topoculus}, which will be proven in \sref{subsec:univ-ext-wf}.

\begin{rmk}
  One could ask whether \jref{thm:v-is-ewf} is true if we assume only that \( \rS \in \cC \) is an internal topos, rather than a plentiful DOF classifier.
  Our proof relies crucially on this stronger assumption; we do not know if the more general statement is true.
\end{rmk}

The proof of \jref{thm:v-is-ewf} takes place almost entirely inside of the topos \( \DOF_{\rS'}(\tm) \), in which we have the additional data of the collection of \( \rS \)-small DOFs, and the object \( \hs(\rS) \) classifying \( \rS \)-small hs-objects.
As explained in \sref{subsec:intro-topoculi}, it is useful to axiomatize this situation, leading to the notion of \emph{topoculus}, which we turn to now.

\subsection{Topoculi} \label{subsec:topoculi}
As mentioned in \sref{subsec:intro-topoculi}, the notion of topoculus is a variant of that of ``category with small maps'' from \cite{joyal-moerdijk-ast} (of which there have since been other variants), which is an axiomatization of the collection of maps with small fibres in the category of classes.
As the motivating example of a topoculus is the set of \( \rS \)-small DOFs for a plentiful DOF classifier \( \rS \), the definition is closely related to that of plentifulness (plenitude?), see \sref{subsubsec:plenitude}.
Recall our conventions regarding slice categories from \sref{subsubsec:under-slice-prelims}.
\begin{defn} \label{defn:topoculus}
  Given a topos \( \bC \), a \defword{topoculus} \( \fS \) on \( \bC \) is a set of morphisms in \( \bC \), which we call the \defword{\( \fS \)-small} (or just \defword{small}) morphisms, satisfying the following:
  \begin{enumerate}[(i)]
  \item \label{item:topoculus-monos} all monomorphisms are small
  \item \label{item:topoculus-comp} any composition of small morphisms is small
  \item \label{item:topoculus-pb} small morphisms are stable under pullback (and in particular, under isomorphism)
  \item for each \( U \in \bC \), the full subcategory \( \fS_U \subset \bC / U \) consisting of the small morphisms is a sub-topos, i.e., \( \fS_U \) is itself a topos, and the inclusion \( \fS_U \hto \bC / U \) is logical
  \end{enumerate}
  We call an object \( (X, \delta_X) \in \bC / U \) in a slice of \( \bC \) (for example, in \( \bC \cong \bC / \tm \) itself) \defword{small} if the morphism \( \delta_X \colon X \to U \) is small, and we call an hs-object in \( \bC / U \) \defword{small} if its underlying object in \( \bC / U \) is small.

  Given a topoculus \( \fS \) on \( \bC \), a \defword{pre-hs-classifier} for \( \fS \) is an object \( \cV \) equipped with a small hs-object \( \und{X_\fS} \in \bC / \cV \), having the following universal property: for each \( U \in \bC \) and each small hs-object \( \und X \in \fS_U \subset \bC / U \), there is a unique morphism \( f \colon U \to \cV \) such that \( \und X \) is a pullback of \( \und{X_\fS} \) along \( f \) (in the sense of \jref{defn:ein-cartesian}).
  Note that there is at most one pre-hs-classifier up to isomorphism for a given topoculus.

  An \defword{hs-classifier} for \( \fS \) is an \( \unein \)-object \( (\cV, \ein_{\cV}) \) in \( \bC \) such that \( \cV \) is a pre-hs-classifier, and such that, given \( f,g \colon U \to \cV \), we have \( f \ein g \) iff \( f^* \und{X_\fS} \Ein g^* \und{X_\fS} \).

  A \defword{topoculized topos} \( (\bC, \fS) \) is a topos \( \bC \) together with a topoculus \( \fS \), and an \defword{hs-topoculized topos} \( (\bC, \fS, \cV) \) is a topoculized topos together with an hs-classifier \( \cV \) for \( \fS \).
\end{defn}

\begin{rmk} \label{rmk:topoculus-2-of-3}
  Given a topoculus \( \fS \) on a topos \( \bC \), it follows from axioms~\ref{item:wf-fp-props-first}-(iii) that given morphisms \( P \tox{f} Q \tox{g} A \) in \( \bC \), if \( f g \) is small, then so is \( f \).
  The argument is the same as that given in \cite[Proposition~2.3.4]{helfer-topoi-in-topoi}.\footnote{%
    \label{footnote:1-of-2-of-3}
    In \loccit, the statement that is being proven is the ``2-of-3'' property: that if \( g \) and \( f g \) are both small, then so is \( f \).
    However, the present stronger conclusion holds there too.
  }
  It follows that for each \( U \in \bC \), all the morphisms in \( \fS_U \subset \bC / U \) are themselves small, and hence that, for each \( \und{U'} \in \fS_U \) the isomorphism \( (\bC / U) / \und{U'} \toi \bC / U' \) restricts to an isomorphism
  \[
    \fS_U / \und{U'} \toi \fS_{U'}.
  \]

  We also note that, because the inclusion \( \fS_U \hto \bC / U \) is logical, any (topos-theoretic) property of some object/morphism/diagram in \( \fS_U \) will hold in \( \fS_U \) if and only if it does in \( \bC / U \).
\end{rmk}

We will also need the following further closure condition for topoculi:
\begin{propn} \label{propn:topoculus-epi-closed}
  If \( (\bC, \fS) \) is a topoculized topos, then given morphisms \( X \tox{f} Y \tox{g} Z \) in \( \bC \), if \( f \) is epi and \( f g \in \fS \), then \( g \in \fS \).
\end{propn}

\begin{proof}
  By assumption, we have that \( \und{X} = (X, f g) \) is in the subtopos \( \fS_Z \) of \( \bC / Z \), and we wish to show that \( \und Y = (Y, g) \) is as well.
  We have that \( f \colon \und X \to \und Y \) is an epi, since it is an epi in \( \bC \) (see \sref{subsubsec:under-slice-prelims}).
  Hence, forming the kernel pair \( \und X \times_{\und Y} \und X \), we have a coequalizer diagram \( \und X \times_{\und Y} \und X \toto \und X \to \und Y \) by \jref{propn:epis-eff}.
  Now \( \und X \times \und X \in \fS_Z \) since \( \fS_Z \), being a subtopos, is closed under finite limits (see \sref{subsec:logical-functors}), and hence \( \und X \times_{\und Y} \und X \in \fS_Z \) by \jref{defn:topoculus}~\ref{item:topoculus-monos}, being the domain of the monomorphism \( \und X \times_{\und Y} \und X \tto \und X \times \und X \).
  Hence, \( \und Y \in \fS_Z \) since \( \fS_Z \) is closed finite under colimits.
\end{proof}

\begin{exm}
  For any topos \( \bC \), the set \( \fS = \Ar \bC \) of \emph{all} arrows on \( \bC \) is a topoculus.
  By contrast, the set \( \fS \subset \Ar \bC \) of \emph{monomorphisms} in \( \bC \) is not, in general, a topoculus, since the category \( \fS_U \subset \bC / U \) will not have a subobject classifier in general.
\end{exm}

We now explain how a topoculus is obtained from a 2-topos.
Recall the notation \( \SOF(A) \) and the axiom (UA) from \sref{subsubsec:dof-sof-digression}.
\begin{propn} \label{propn:2-topoi-give-topoculi}
  Let \( \cC \) be a groupoidant 2-category.
  \begin{enumerate}[(i)]
  \item \label{item:2-topoi-give-topoculi-s-s} If \( \rS \) and \( \rS' \) are plentiful DOF classifiers in \( \cC \) such that \( p \colon \rS_* \to \rS \) is \( \rS' \)-small, then the set \( \fS \) of \( \rS \)-small SOFs is a topoculus on the topos \( \SOF_{\rS'}(\tm) \).
    Moreover, if the hs-classifier \( \hs(\rS) \) for \( \rS \) is \( \rS' \)-small, then it is an hs-classifier for \( \fS \).
  \item \label{item:2-topoi-give-topoculi-ua} If \( \cC \) satisfies Axiom~(UA), so that \( \SOF(\tm) \) is a 1-topos, then for any plentiful DOF classifier \( \rS \in \cC \), the set \( \fS \) of \( \rS \)-small SOFs is a topoculus on \( \SOF(\tm) \), and \( \hs(\rS) \) is an hs-classifier for \( \fS \).\footnote{%
      The same statements holds if the topoi \( \SOF_{\rS'(\tm)} \) and \( \SOF(\tm) \) in \ref{item:2-topoi-give-topoculi-s-s} and \ref{item:2-topoi-give-topoculi-ua} are replaced, respectively, by the equivalent topoi \( \DOF_{\rS'}(\tm) \) and \( \DOF(\tm) \), and \( \fS \) is taken to be the set of \( \rS \)-small \emph{DOFs}.
      However, in this case, \( \hs(\rS) \) does not itself lie in \( \SOF(\tm) \), and so it is rather an arbitrary \emph{DOF-collapse} of \( \hs(\rS) \)---which exists by the assumption that \( \hs(\rS) \) is \( \rS' \)-small for some \( \rS \)---which serves as an hs-classifier for \( \fS \). %
    }
  \end{enumerate}
\end{propn}

\begin{proof}
  See \hyperref[proof:2-topoi-give-topoculi]{\S\ref*{subsec:topoculi-proofs}}.
\end{proof}

By virtue of \jref{propn:2-topoi-give-topoculi}~\ref{item:2-topoi-give-topoculi-s-s}, \jref{thm:v-is-ewf} immediately reduces to the following proposition, which we prove in the next section.
\begin{propn} \label{propn:v-is-ewf-topoculus}
  Given a topoculized topos \( (\bC, \fS) \), any hs-classifier for \( \fS \) is extensional and well-founded.
\end{propn}

\subsection{Extensionality and well-foundedness of the universe} \label{subsec:univ-ext-wf}
Fix an hs-topoculized topos \( (\bC, \fS, \cV) \).
We now reformulate the extensionality and well-foundedness of the hs-classifier \( \cV \) asserted in \jref{propn:v-is-ewf-topoculus}, thus reducing the latter to \jrefs{propn:univ-ext}{propn:univ-wf} below.

The extensionality of \( \cV \) means that, for each \( U \in \bC \) and \( x,y \colon U \to \cV \), if \( w \ein x' \ToT w \ein y' \) for all \( w \colon U' \to \cV \), then \( x = y \).
By the universal property of \( \cV \), we are thus given hs-objects \( \und X, \und Y \in \fS_U \subset \bC / U \), and are supposing that for each \( \delta_{U'} \colon U' \to U \) in \( \bC \) each \( \und W \in \fS_{U'} \), we have \( \und W \Ein \delta_{U'}^* X \ToT \und W \Ein \delta_{U'}^* Y \); and we must show that \( \und X \cong \und Y \).

It turns out that in the proof, we only ever need to use the special case of the above assumption in which \( \delta_{U'} \) is itself small, and hence \( \und{U'} = (U', \delta_{U'}) \) is an object in \( \fS_{U} \).
We then have the equivalence \( \fS_{U'} \cong \fS_U / \und{U'} \) of \jref{rmk:topoculus-2-of-3}, under which the functor \( \delta_{U'}^* \colon \fS_U \to \fS_{U'} \) corresponds to the product (ana)-functor \( \und{U'} \times \colon \fS_U \to \fS_U / \und{U'} \).
The significance of this is that both the assumption and what we are trying to prove are stated entirely in terms of the topos \( \fS_U \) (and its slices \( \fS_U / \und{U'} \)).
We have thus reduced the question of extensionality in \jref{thm:v-is-ewf} to the following purely topos-theoretic (as opposed to topoculus-theoretic) proposition:
\begin{propn} \label{propn:univ-ext}
  Let \( \bC \) be a topos and let \( X, Y \in \hs(\bC) \).
  Suppose that for all \( U \in \bC \) and \( \und W \in \hs(\bC / U) \), we have that \( \und W \Ein (U \times X, \pi_U) \ToT \und W \Ein (U \times Y, \pi_U) \).
  Then \( X \cong Y \) (as hs-objects).
\end{propn}

\begin{proof}
  See \sref{subsec:univ-ext-proof}.
\end{proof}

Next, we turn to the well-foundedness of \( \cV \).
Again, we will translate this into a purely topos-theoretic question.
The well-foundedness of \( \cV \) means that for all \( P \tto \pow \cV \), if \( P \) is inductive, then \( P x \) for all \( U \in \bC \) and \( x \colon U \to \cV \).

Thus, fix \( P \tto \cV \) and \( U \in \bC \).
Let \( \fP \subset \Ob \hs(\fS_U / \und{U'}) \) be the set of hs-objects in \( \fS_{U} \) that are classified by a morphism \( x \colon U \to \cV \) with \( P x \).
We wish to show that if \( P \) is inductive, then \( \fP \) contains all hs-objects in \( \fS_{U} \).

For each \( \und{U'} \in \fS_{U} \), let \( \fP_{\und{U'}} \subset \Ob \hs(\fS_{U'}) \) be the set of hs-objects in \( \fS_{U} / \und{U'} \) whose image under the isomorphism \( \fS_{U} / \und{U'} \toi \fS_{U'} \) from \jref{rmk:topoculus-2-of-3} is classified by a morphism \( y \colon U' \to \cV \) with \( P y \).
(Thus, under the isomorphism \( \fS_{U} / \und U \toi \fS_U \), \( \fP_U \) is identified with \( \fP \).)
We now make some observations about the sets \( \fP_A \) for \( A \in \bD \defeq \fS_{U} \):
\begin{enumerate}[(i)]
\item\label{item:wf-fp-props-first} \label{item:wf-fp-props-stable} If an hs-object \( X \) in \( \bD / A \) is in \( \fP_A \), then \( f^*X \) is in \( \fP_{A'} \) for all \( f \colon A' \to A \) (and in particular, \( \fP_A \) is closed under isomorphisms of hs-objects).
\item\label{item:wf-fp-props-rep} For each \( A \in \bD \) and hs-object \( X \) in \( \bD / A \), there is an underobject \( A_X \tto A \) such that \( f \colon A' \to A \) factors through \( A_X \) iff \( f^*X \in \fP_{A'} \).\footnote{%
    \label{footnote:stack-semantics-repr}
    Cf.\ the notion of ``representable'' sentence in \cite{shulman-stack-semantics}.
  }
\end{enumerate}
The assumption that \( P \) is inductive gives, moreover:
\begin{enumerate}[(i), resume]
\item\label{item:wf-fp-props-last} \label{item:wf-fp-props-inductive} Given an hs-object \( X \) in \( \bD / A \), if \( W \in \fP_{A'} \) for each \( f \colon A' \to A \) and each \( W \) in \( \bD / A' \) with \( W \Ein f^* X \), then \( X \in \fP_A \).
\end{enumerate}
All of these properties are easily verified; we will comment only on \ref{item:wf-fp-props-rep}.
This property says that for each \( \fS \)-small \( \colon U' \to U \) in \( \bC \) and each hs-object \( X \) in \( \fS_{U'} \) classified by \( x \colon U' \to \cV \), there is an underobject \( U'_X \tto U' \) such that an \( \fS \)-small morphism \( f \colon U'' \to U' \) factors through \( U'_X \) if and only if \( P (f x) \); and indeed, this holds if we take \( U'_X = x^* P \tto U' \).

We have thus reduced the claim of well-foundedness in \jref{thm:v-is-ewf} the following:
\begin{propn} \label{propn:univ-wf}
  Let \( \bD \) be a topos, fix a subset \( \fP_A \subset \Ob \hs(\bD / A) \) for each \( A \in \bD \), and suppose that these satisfy the above properties \ref{item:wf-fp-props-first}-\ref{item:wf-fp-props-last}.
  Then \( \fP = \fP_\tm \) is equal to all of \( \Ob \hs(\bD) \) (and hence, by \ref{item:wf-fp-props-stable}, \( \fP_A = \Ob \hs(\bD / A) \) for all \( A \)).
\end{propn}
\begin{proof}
  See \sref{subsec:univ-wf-proof}.
\end{proof}

\subsection{The Burali-Forti paradox} \label{subsec:burali-forti}
Having established that the hs-classifier \( \cV \) for a topoculus is extension and well-founded, it follows from \jref{propn:top-ext-is-hs-struct} that its top-extension \( \wh \cV \) is an hs-object.
The main argument in the Burali-Forti paradox is now:

\begin{propn} \label{propn:too-universal}
  Let \( (\bC, \fS, \cV) \) be an hs-topoculized topos, and let \( \wh \cV \) be a top-extension of \( \cV \) as in \jref{defn:top-extension}.
  Then \( X \Ein \wh \cV \) for every small hs-object \( X \) in \( \bC \).
\end{propn}

\begin{proof}
 See \sref{subsec:too-universal-proofs}.
\end{proof}

\begin{thm}[Burali-Forti]\label{thm:world-explodes}
  Let \( (\bC, \fS, \cV) \) be an hs-topoculized topos.
  Suppose that every morphism in \( \bC \) is small---or even just that the morphism \( \unex_\cV \colon \cV \to \tm \) is small.
  Then the topos \( \bC \) is trivial (i.e., equivalent to the terminal category).
\end{thm}

\begin{proof}
  By \jref{propn:too-universal}, the hs-object \( \wh \cV \) in \( \bC \) satisfies \( X \Ein \wh \cV \) for every small hs-object \( X \) in \( \bC \).
  Now, note that the hs-object \( \wh \cV \) is itself small.
  Indeed, its underlying object is just a coproduct \( \cV + \tm \), and the sub-topos \( \fS_\tm \subset \bC \) has coproducts which are preserved by the logical inclusion functor \( \fS_\tm \hto \bC \); and it then follows that \( \wh \cV \times \wh \cV \) and \( \unein_\cV \tto \wh \cV \times \wh \cV \) are small as well, since small morphisms are closed under pullback and composition, and since all monos are small.

  Thus, if \( v \colon \tm \to \cV \) is the morphism classifying \( \wh \cV \), we conclude that \( x \ein v \) for all \( x \colon \tm \to \cV \), and in particular that \( v \ein v \).
  Thus, by \jref{propn:no-loops}, \( \tm \) is initial, which implies that \( \bC \) is trivial.
\end{proof}

We now return from the topoculus setting to the 2-topos setting and conclude:
\begin{cor} \label{cor:2-topos-world-explodes}
  Let \( \rS \) be a DOF classifier in a groupoidant 2-category \( \cC \) and suppose that every SOF in \( \cC \) is \( \rS \)-small
  Then \( \cC \) is trivial (i.e., equivalent to the terminal 2-category).
\end{cor}

\begin{proof}
  Assuming that every SOF in \( \cC \) is \( \rS \)-small, it follows that in the topos \( \SOF_\rS(\tm) = \SOF(\tm) \) with its topoculus and hs-classifier coming from \jref{propn:2-topoi-give-topoculi}, every morphism is small.
  Hence, by \jref{thm:world-explodes}, \( \SOF(\tm) \) is trivial.
  Since for each groupoidal \( A \in \cC \), the functor \( \SOF(\tm) \simeq \cC(\tm, \rS) \tox{\unex_A^*} \cC(A, \rS) \) is logical (by definition of \( \rS \) being an internal topos), it follows that \( \zset \cong \tm \) in \( \cC(A, \rS) \simeq \SOF(A) \) and hence that it is trivial as well.

  Next, we claim that every groupoid \( A \in \cC \) is equivalent to the terminal object.
  Note that \( \br{\partial_0,\partial_1} \colon A^\to \to A \times A \) is a DOF; this amounts to the claim that given a 2-cell \( \alpha \colon f_0 \To f_1 \) from \( A \) into \( X \) and 2-cells \( \beta_0 \colon f_0 \To g_0 \) and \( \beta_1 \colon f_1 \To g_1 \), there is a unique 2-cell \( \gamma \colon g_0 \To g_1 \) making
  \[
    \begin{tikzcd}
      f_0 \ar[r, "\alpha"] \ar[d, "\beta_0"'] & f_1 \ar[d, "\beta_1"] \\
      g_0 \ar[r, "\gamma", dashed] & g_1
    \end{tikzcd}
  \]
  commute, which is true since \( \beta_0 \) is invertible.
  It follows that \( \br{\partial_0, \partial_1} \) is an isomorphism, being a morphism in the trivial category \( \SOF(A \times A) \).
  Hence, for any two morphisms \( x,y \colon U \to A \), there is a unique 2-cell \( x \To y \).
  In particular, \( A \) is a setoid.
  Hence, \( A \to \tm \) is a morphism in the trivial 2-category \( \SOF(\tm) \), and hence must be an equivalence, as claimed.

  Next, we claim that \( \cC(A, X) \) is a groupoid for any \( X \) and any \emph{groupoid} \( A \).
  Indeed, consider a morphism \( \alpha \colon f \To g \) in \( \cC(A, X) \), i.e., a 2-cell into \( X \).
  This is classified by a morphism \( a \colon A \to X^\to \), which factors through some \( \bar a \colon A \to (X^\to)^\iso \), since \( A \) is a groupoid.
  Since the codomain of \( \bar a \) is equivalent to \( \tm \), as we have already shown, it follows that \( \bar a \cong \bar b \), where \( \bar b \colon A \to (X^\to)^\iso \) is the morphism with \( \bar b i_{(X^\to)^\iso} \partial = \id_f \) (and where \( \partial \colon X^\to \tocell X \) is the universal 2-cell).
  This gives an isomorphism \( \alpha = \bar a i \partial \cong \bar b i \partial = \id_f \) in \( \cC(A, X)^\to \), which amounts to a commutative square
  \[
    \begin{tikzcd}
      f \ar[r, "\alpha"] \ari[d, ""'] & g \ari[d, ""] \\
      f \ar[r, "\id_a"] & f,
    \end{tikzcd}
  \]
  with the vertical morphisms isos, as indicated.
  It follows that \( \alpha \) is an iso as desired.

  Next, we claim that the groupoid \( \cC(A, X) \) is trivial for any \( X \) and any groupoidal \( A \).
  First of all, this groupoid has an object, since \( X^\iso \simeq \tm \), and hence there is a morphism \( A \to X^\iso \), which we can compose with the inclusion \( i \colon X^\iso \to X \).
  Next, any two morphisms \( f,g \colon A \to X \) factor through \( X^\iso \), and hence are isomorphic, since \( X \simeq \tm \).
  And finally, any two isomorphisms \( \alpha,\beta \colon f \toi g \) both factor through \( X^\iso \), and hence are equal since \( X \simeq \tm \).

  Finally, we claim that for any object \( X \), the terminal morphism \( \unex_X \colon X \to \tm \) is an equivalence, which will complete the proof.
  We obtain a morphism \( \tm \to X \) in the other direction by using the triviality of \( \cC(\tm, X) \).
  Clearly the composite \( \tm \to X \to \tm \) is equal to the identity, so it remains to see that the other composite \( X \to \tm \to X \) is isomorphic to the identity.
  By \cite[Corollary~6.2.3]{helfer-topoi-in-topoi}, to show that two morphisms \( f,g \colon X \to X \) are isomorphic, it suffices to prove that they induce isomorphic natural transformations \( (\wh X)_\gpd \to (\wh X)_\gpd \), where \( (\wh X)_\gpd = \cC(-, X) \colon \cC_\gpd^\op \to \Cat \) is the 2-functor represented by \( X \), and where \( \cC_\gpd \) is the full sub-2-category on the groupoidal objects.
  Thus, we must produce an invertible modification \( \mu \colon f_* \toi g_* \); this consists of isomorphisms \( \mu \colon (f_*)_A \to (g_*)_A \) between \( (f_*)_A, (g_*)_A \colon \cC(A,X) \to \cC(A,X) \) for each \( A \in \cC_\gpd \), which are natural in \( A \).

  But since \( \cC(A,X) \) is trivial for all \( A \), there is a \emph{unique} isomorphism \( \mu_A \colon (f_*)_A \to (g_*)_A \) for each \( A \in \cC_\gpd \), and the naturality follows automatically from this uniqueness.
\end{proof}

\begin{rmk}
  There should be a simpler version of the above proof, avoiding the detours through groupoidal objects, by showing directly that for any object \( X \in \cC \), the internal hom fibration \( \mathrm{hom}_X \to X^\op \times X \) is a DOF, and hence an equivalence.
\end{rmk}

\section{The universe is  a model of set theory} \label{sec:v-is-a-model}
We established in \jref{thm:v-is-ewf} that the \( \unein \)-object \( \hs(\rS) \) of hs-objects in a plentiful DOF classifier \( \rS \)---and more generally, an hs-classifier \( \cV \) for a topoculus---is extensional and well-founded.
In other words, this says that \( \cV \) satisfies the axioms of extensionality and foundation from set theory.
We now begin the demonstration that \( \cV \) satisfies the remaining usual axioms of set theory---specifically, of MK (Morse-Kelley) set theory (see \sref{subsec:axioms-of-mk}).

Besides extensionality and foundation, these axioms fall into two groups.
The first group consists solely of the axiom schema of class-comprehension: for every first-order predicate, there is a class consisting of all \emph{sets} satisfying that predicate.
As we will be interpreting ``class'' using the power object \( \pow(\cV) \), this will turn out to be trivial to verify using the usual comprehension for power objects in a topos (see \sref{subsec:comprehension}).

The second group contains the ``set-existence axioms'': each of these say that some particular \emph{class} is in fact a \emph{set}.
For example, the axiom of pairing says that for any two sets, the class containing just those sets is a set; and the axiom of subsets says any subclass of a set is a set.
To explain how we will establish these axioms, let us first discuss how to do it for the universe \( \rV_\hs \) of hs-sets from \sref{subsec:ein-sets}.

Recall from there that that the universe \( \rV \) of hereditary sets is defined to be the least class closed under set-formation, and that we have an isomorphism of (large) \( \unein \)-set \( \rV_\hs \cong \rV \).
It follows immediately that \( \rV_\hs \) is also closed under set-formation:
\begin{settybox}
  \begin{propn} \label{propn:set-v-closed-under-sets}\enumbelow
    \begin{enumerate}[(i)]
    \item\label{item:set-v-closed-direct} For any \emph{small} subset \( A \subset \rV_\hs \), there is a unique \( x \in \rV_\hs \) such that \( w \ein x \) iff \( w \in A \) for all \( w \in \rV_\hs \).\footnotemark
    \item\label{item:set-v-closed-converse} Conversely, for any \( x \in \rV_\hs \), the class \( \set{w \in \rV_\hs \mid w \ein x} \) is small.
    \end{enumerate}
  \end{propn}
  \begin{proof}
    As we just wrote, this follows from the corresponding properties of \( \rV \) (which are immediate from its definition).
    However, as we are about to explain, we also want proofs of these which make no reference to \( \rV \).
    The (easy) proof of \ref{item:set-v-closed-converse} is given, together with the proof of its topos-theoretic analogue \jref{propn:topos-universe-small-closed-converse}, in \sref{subsec:proof-of-small-converse}.
    As for \ref{item:set-v-closed-direct}, we will state a sharper version in \jref{propn:set-v-closed-precise} below, which we will prove, along with it topos-theoretic analogue \jref{propn:topos-universe-small-closed-precise}, in \sref{subsec:universe-closure-proof}.
  \end{proof}
\end{settybox}
\footnotetext{%
  Since \( \rV_{\hs} \) is actually ``super-large'', as explained in \prerefpost{footnote~}{footnote:super-large}{} on \pref{footnote:super-large}, we should really say ``essentially small'' here---i.e., isomorphic to a small set---rather than small.
  We will suppress this distinction in what follows.%
}
From \jref{propn:set-v-closed-under-sets}, the satisfaction of the set-existence axioms by \( \rV_\hs \) follow easily, using the corresponding assumptions which we are making in our ambient (class-)set theory.
We note that, when interpreting the axioms of MK in the structure \( (\rV_{\hs}, \ein_{\rV_{\hs}}) \), we are interpreting ``class'' to mean an arbitrary subclass of \( \rV_{\hs} \), and ``set'' to mean a class of the form \( \set{ x \in \rV_{\hs} \mid x \ein y } \) for some \( y \in \rV_{\hs} \).
Thus, for example, to verify the axiom of subsets, given ``classes'' \( X,Y \subset \rV_\hs \) with \( X \subset Y \), and assuming \( Y \) is a ``set'', we must show that \( X \) is also a ``set''.
From \jref{propn:set-v-closed-under-sets}~\ref{item:set-v-closed-converse}, we conclude that the class \( Y \) is small, hence by the \emph{real} axiom of subsets, which we are assuming in our ambient class-set theory, we conclude that \( X \subset Y \) is small as well.\footnote{%
  When working in the context of a fixed Grothendieck universe \( \cU \) rather than in class-set theory (see \sref{subsec:ein-sets}), we rather use that \( \cU \) itself \emph{satisfies} the axiom of subsets; i.e., that \( X \subset Y \in \cU \) implies \( X \in \cU \).
}
Hence, by \jref{propn:set-v-closed-under-sets}~\ref{item:set-v-closed-direct}, we conclude that there is \( x \in \rV_\hs \) with \( X = \set{w \in \rV_\hs \mid w \ein x} \), as desired.
The proofs of the other set-existence axioms are similar.

Thus, we see that the key fact is \jref{propn:set-v-closed-under-sets}: the closure of \( \rV_\hs \) under set-formation.
Above, we deduced this from the isomorphism \( \rV_\hs \cong \rV \).
However, to better reflect the 2-topos-theoretic situation, we should not allow ourselves to make use of \( \rV \), and understand how to prove the closure of \( \rV_\hs \) under set-formation directly.

Given a small set \( A \subset \rV_\hs \), then, our goal is to construct the element \( x \in \rV_\hs \) with \( A = \set{w \in \rV_\hs \mid w \ein x} \)---that is, to construct an hs-set \( X \) such that any given hs-set \( W \) satisfies \( W \Ein X \) if and only if \( W \cong \floor{w} \) for some \( w \in A \).
Here, and below, for an isomorphism class \( w \in \rV_\hs \) of hs-sets, we write \( \floor{w} \) for the hs-set \( \pbig{\coprod_{Y \in w} Y} / \unsim \), where \( \unsim \) is the equivalence relation induced by the (by \jref{propn:init-uniq}) unique isomorphisms \( Y_0 \toi Y_1 \) for all \( Y_0, Y_1 \in w \), with the obvious induced \( \ein \)-set structure \( \ein_{\floor w} \).
We have a canonical isomorphism \( Y \toi \floor{w} \) for each \( Y \in w \), and hence \( w = \sqbig{\floor{w}} \).

As a first approximation to \( X \), we can form the disjoint union \( W = \coprod_{y \in A} \floor{y} \), with the \( \unein \)-set structure given by the disjoint union of the relations \( \unein_{\floor{y}} \) for each \( y \in W \).
This obviously has the property of admitting an initial maps \( \floor{y} \to W \) for each \( y \in A \).
If we furthermore adjoin a top element \( \rt_W \) as in \sref{subsec:top-element}, we can moreover make each of these initial maps elemental, and these will be the only elemental maps into \( W \).
The only problem is that \( W \) is not extensional in general, since there can be distinct elements \( y \in A \) such that the hs-set \( \floor{y} \) has a bottom element \( x \) (i.e., one with \( w \not\,\ein x \) for all \( w \in \floor{y} \)), and the resulting elements \( (y, x) \in W \) will violate extensionality.

However, we can correct this by passing to an appropriate quotient \( W / \unsim \) of \( W \), where \( x \sim y \) for \( x,y \in W \) if \( \dwnfn(x) \cong \dwnfn(y) \), with its induced \( \unein \)-set structure: for each \( U,V \ein W / \unsim \), we have \( U \ein V \) iff \( u \ein v \) for some \( u \in U \) and \( v \in V \).\footnote{%
  The same construction is given in \cite[Lemma~8.28]{shulman-comparing}.
}
Recalling that we still need to adjoin a top element, we thus define \( X \defeq \wh{W / \unsim} \) to be the \( S \)-top extension (\jref{defn:top-extension}), where \( S \subset W / \unsim \) is the image of \( \set{ \br{y, \rt_{\floor{y}}} \mid y \in A} \subset \coprod_{y \in A} \floor{y} = W \) under the quotient map \( W \to W / \unsim \).
We then have the following more precise version of \jref{propn:set-v-closed-under-sets}~\ref{item:set-v-closed-direct}:
\begin{settybox}
  \begin{propn} \label{propn:set-v-closed-precise}
    Given a small subset \( A \subset V_{\hs} \), and defining \( W = \coprod_{y \in A}\floor{y} \), \( \unsim \), \( S \), and \( X = \wh{(W / \unsim)} \) as above, then \( X \) is an hs-set, and its equivalence class \( x = [X] \in V_{\hs} \) satisfies \( A = \set{y \in \rV_\hs \mid y \ein x} \).
  \end{propn}
\end{settybox}

The proof of \jref{propn:set-v-closed-under-sets} and \jref{propn:set-v-closed-precise}, and of their topos-theoretic counterparts, thus involves two steps: in \sref{subsec:ewf-collapse}, we give the general construction which produces an extensional \( \unein \)-object from a non-extensional one by passing to a quotient; and then, in \sref{subsec:set-formation-closure}, we apply this general construction to \( W \) as above and show that (after adjoining a top element), the resulting hs-object has the desired property.
Finally, in \sref{subsec:axioms-of-mk}, we will then deduce that the axioms of MK hold in our universe.

\subsection{Extensional collapse of well-founded \alttext{\( \unein \)}{ε}-objects} \label{subsec:ewf-collapse}
The proofs all the statements in this section can be found in \sref{subsec:ewf-collapse-proof}.

We now explain how to take a quotient of a possibly non-extensional \( \unein \)-object to produce an extensional one.
In order to do this, we need the additional assumption---which clearly holds in the above example \( W \)---that the substructure \( \dwnfn(x) \) is extensional for each element \( x \), i.e., that our \( \unein \)-object is locally extensional in the sense of \jref{defn:loc-ext}.

Next, we characterize abstractly the above construction of the quotient \( W / \unsim \):
\prooflabel{proof:collapse-epi-iff-univ}
\begin{propn} \label{propn:collapse-epi-iff-univ}
  Let \( X \) and \( Y \) be \( \unein \)-objects in a topos \( \bC \) and \( q \colon X \to Y \) an initial morphism, and suppose that \( X \) is well-founded and \( Y \) is extensional.
  Then the following are equivalent:
  \begin{enumerate}[(i)]
  \item \( q \) is an epimorphism
  \item \label{item:collapse-epi-iff-univ-univ} \( q \) is the universal initial morphism from \( X \) to an extensional \( \unein \)-object, i.e., for any initial morphism \( f \colon X \to Z \) with \( Z \) extensional, there is a unique \( g \colon Y \to Z \) with \( q g = f \).
  \end{enumerate}
  Moreover, the morphism \( g \) in \ref{item:collapse-epi-iff-univ-univ} is again initial.
\end{propn}

\begin{defn}
  An \defword{extensional collapse} of a well-founded \( \unein \)-object \( X \) in a topos \( \bC \) is an extensional \( \unein \)-object \( Y \) equipped with an initial morphism \( X \to Y \) satisfying the equivalent conditions of \jref{propn:collapse-epi-iff-univ}.
\end{defn}

The main result of this section is:
\begin{propn} \label{propn:ewf-collapse}
  Any well-founded, locally extensional \( \unein \)-object \( X \) in a topos \( \bC \) admits an extensional collapse \( q \colon X \to Y \), and moreover \( Y \) is well-founded.
\end{propn}
\begin{proof}
  This is implied by the sharper statement \jref{propn:ewf-collapse-precise} below.
\end{proof}

As in the set-theoretic case, we will construct \( Y \) as a quotient of \( X \) by the equivalence relation which is described set-theoretically as \( x \sim y \ToT \dwnfn(x) \cong \dwnfn(y) \).
The counterpart of this for \( \unein \)-object \( X \) in a topos is the equivalence relation on generalized elements \( x,y \colon U \to X \) given by \( x \sim y \ToT \und{\Dwn_x} \cong \und{\Dwn_y} \).
In order to be able to take a quotient by this equivalence relation, we need to know that it is represented by an actual relation \( R \tto X \times X \):\footnote{%
  Cf.\ the notion of ``representable'' sentence in \cite{shulman-stack-semantics}; see also \prerefpost{footnote~}{footnote:stack-semantics-repr}{} on \pref{footnote:stack-semantics-repr}.
}
\prooflabel{proof:dwn-iso-rel-exists}
\begin{lem} \label{lem:dwn-iso-rel-exists}
  For any well-founded, locally extensional \( \unein \)-object \( X \) in a topos \( \bC \), there exists a (obviously unique up to equivalence) underobject \( R_X \tto X \times X \) such that \( x \Rrel_X y \) iff \( \und{\Dwn_x} \cong \und{\Dwn_y} \) for all \( x,y \colon U \to X \).

  \noindent
  (It follows that \( R_X \) is an equivalence relation, since \( {\cong} \) is.)
\end{lem}
For \( X \) well-founded and locally extensional, we will in general use \( R_X \tto X \times X \) to denote an underobject as in \jref{lem:dwn-iso-rel-exists}.
We can now state a more precise version of \jref{propn:ewf-collapse}:
\prooflabel{proof:ewf-collapse-precise}
\begin{propn} \label{propn:ewf-collapse-precise}
  Let \( X \) be a well-founded, locally extensional \( \unein \)-object in a topos, let \( q \colon X \to Y \) be a quotient by \( R_X \tto X \times X \)---i.e., a coequalizer of \( i_{R_X}\pi_0, i_{R_X} \pi_1 \colon R_X \toto X \)---and endow \( Y \) with the \( \unein \)-object structure \( \unein_Y \tto Y \times Y \) given by the image of \( \unein_X \tto X \times X \tox{q \times q} Y \times Y \).

  Then \( Y \) is well-founded and \( q \colon X \to Y \) is an extensional collapse.
\end{propn}

To prove \jref{propn:ewf-collapse-precise} and hence \jref{propn:ewf-collapse}, it remains to prove that the \( \unein \)-object \( Y \) constructed in the former is well-founded and extensional, and that \( q \) is initial.
All of these properties depend on the following ``initiality'' or ``bisimulation'' property of \( R_X \) (cf.\ \cite[\S2]{aczel-non-well-founded-sets}~and~\cite[Lemma~8.28]{shulman-comparing}):
\prooflabel{proof:r-is-bisimulation}
\begin{lem} \label{lem:r-is-bisimulation}
  For any well-founded and locally extensional \( \unein \)-object \( X \) in a topos \( \bC \), the relation \( R_X \tto X \) satisfies:
  \[
    \bC \vDash \forall x_0,x_1,y_1 \colon X.\
    x_0 \ein x_1 \wedge x_1 \Rrel_X y_1 \to \exists y_0.\ y_0 \ein y_1 \wedge x_0 \Rrel_X y_0.
  \]
\end{lem}

\noindent
The extensionality of \( Y \) depends on the following further lemma, showing that the relation \( R_X \) is trivial in case \( X \) is always extensional:
\prooflabel{proof:if-ext-r-is-equality}
\begin{lem} \label{lem:if-ext-r-is-equality}
  If \( X \) in any extensional, well-founded \( \unein \)-object in a topos \( \bC \), then for any \( x,y \colon U \to X \), if \( x \Rrel_X y \), then \( x = y \).

  \noindent
  (And more generally, if \( X \) is only locally extensional and well-founded, and given \( x,y,z \colon U \to X \) with \( x,y \preceq z \) and \( x \Rrel_X y \), then \( x = y \).)
\end{lem}

\subsection{Closure of the universe under set-formation} \label{subsec:set-formation-closure}
Continuing with the strategy for proving \jref{propn:set-v-closed-under-sets}~\ref{item:set-v-closed-direct} described above, having established the existence of the extensional collapse of well-founded, locally extensional \( \unein \)-objects, we go on to formulate the topos-theoretic analogue of the statement \jref{propn:set-v-closed-under-sets}~\ref{item:set-v-closed-direct} that the universe \( \rV_\hs \) is closed under set-formation.
We recall from there that the condition in the set-theoretic case, for each \emph{small} subset \( A \subset \rV_\hs \), there is \( x \in \rV_\hs \) with \( A = \set{y \in \rV_\hs \mid y \ein x} \).

To formulate the corresponding topos-theoretic condition, we put ourselves in the context of an hs-topoculized topos \( \bC \).
Now in place of a small subset of \( \rV_\hs \), we can consider an underobject \( A \tto \cV \)---or equivalently a global element \( a \colon \tm \to \pow \cV \)---for which the domain \( A \) is small (i.e., \( \fS \)-small).
More generally, we must consider \emph{generalized} elements \( a \colon U \to \pow \cV \).
These classify underobjects \( A \tto U \times \cV \), and the desired smallness condition is now that \( i_A \pi_U \colon A \to U \) is small (i.e., that the ``\( U \)-indexed family \( a \colon U \to \pow \cV \) of subsets of \( \cV \) is pointwise small'').

The conclusion in the set-theoretic setting, that \( A = \set{y \in \rV_\hs \mid y \ein x} \) for some \( x \in \rV_\hs \), is formulated in the topos-theoretic setting using the following definition:
\begin{defn} \label{defn:associated-class}
  Given an \( \unein \)-object \( V \) in a topos \( \bC \) and a power object \( \pow V \) for \( V \), we define the \defword{associated class morphism}
  \[
    \kappa = \kappa_V \colon V \to \pow V
  \]
  to be the unique morphism for which there is a pullback square
  \[
    \begin{tikzcd}
      \ein_V \ar[r, dashed] \ar[d, "i_{\ein_V}"'] \pb &[10pt] \in_V \ar[d, "i_{\in_V}"] \\
      V \times V \ar[r, "\id_V \times \kappa_V"] & V \times \pow V.
    \end{tikzcd}
  \]
  Thus, \( x \in y \kappa_V \) iff \( x \ein y \) for \( x,y \colon U \to V \).
  Note that \( \kappa_V \) is a monomorphism if and only if \( V \) is extensional.
\end{defn}

\begin{thm} \label{thm:topos-universe-small-closed}
  Let \( (\bC, \fS, \cV) \) be a hs-topoculized topos.
  Then any \( \alpha \colon U \to \pow \cV \) in \( \bC \) classifying an underobject \( i_A \colon A \tto U \times \cV \), for which \( i_A \pi_U \) is small, factors through the associated class morphism \( \kappa \colon \cV \to \pow \cV \).
\end{thm}
\begin{proof}
  This will follow from the sharper statement \jref{propn:topos-universe-small-closed-precise} below.
\end{proof}

To prove this theorem, given \( \alpha \colon U \to \pow \cV \), we must find \( x \colon U \to \cV \) with \( x \kappa = \alpha \).
Let us recall, from \jref{propn:set-v-closed-precise}, how we constructed \( x \) in the set-theoretic context.
There, given \( A \subset \rV_\hs \), we constructed the hs-set \( X \) which is to represent \( x \in \rV_\hs \) by adjoining a top element to the extensional collapse of the disjoint union \( W = \coprod_{y \in A} \floor{y} \).
We thus first need to identify the counterpart of this disjoint union in our present situation.

In the case of a global element \( \alpha \colon \tm \to \pow \cV \) and thus an underobject \( A \tto \cV \), the ``family of hs-objects indexed by \( A \)'' is precisely the hs-object \( \und W \) in \( \bC / A \) classified by \( i_A \).
The ``disjoint union'' should then simply be the domain \( W \in \bC \) of \( \delta_W \colon W \to A \), which is endowed with an \( \unein \)-object structure as in \jref{defn:ein-cocart}.
More generally, for a generalized element \( \alpha \colon U \to \pow \cV \), we have the underobject \( A \tto U \times \cV \) classified by \( \alpha \), and we consider the hs-object \( \und W \) in \( \bC / A \) classified by \( i_A \pi_{\cV} \colon A \to \cV \), and the resulting \( \unein \)-object \( \und W_U \) in \( \bC / U \) (again as in \jref{defn:ein-cocart}).

By \jref{propn:sum-props} \( \und W_U \) is well-founded and locally extensional, and hence, by \jref{propn:ewf-collapse}, we obtain an extensional collapse \( q \colon \und W_U \to \und Z \), giving an extensional, well-founded \( \unein \)-object \( \und Z \) in \( \bC / U \).

The next step is to adjoin a top element to \( \und Z \) at an appropriate dense underobject \( \und S \tto \und Z \), to obtain an hs-object.
Referring to the discussion above \jref{propn:ewf-collapse-precise}, where we had \( W = \coprod_{y \in A} \floor{y} \) and \( Z = W / \unsim \), we took \( S \) to be the image of \( \set{ (y, \rt_{\floor{y}}) \mid y \in A} \subset W \) under the quotient map.
In the present case, this subset of \( W \) corresponds to the top element \( \rt_{W} \colon (A, \id_A) \to \und W \) of \( \und W \in \bC / A \), viewed as a morphism \( \rt_{W} \colon (A, i_A \pi_U) \to \und W_U \) in \( \bC / U \).
We conclude that our underobject \( \und S \tto \und Z \) should be the image of the composite \( (A, \delta_A) \tox{\rt_{W}} \und W \tox{c} \und W_U \tox{q} \und Z \), where \( c \) is the \( \ein \)-cocartesian morphism underlying \( \und W_U \) as in \jref{defn:ein-cocart}.

We can now state a more precise version of \jref{thm:topos-universe-small-closed}, analogous to the sharpening \jref{propn:set-v-closed-precise} of \jref{propn:set-v-closed-under-sets}~\ref{item:set-v-closed-direct}:
\begin{propn} \label{propn:topos-universe-small-closed-precise}
  As in \jref{thm:topos-universe-small-closed}, let \( (\bC, \fS, \cV) \) be an hs-topoculized topos, let \( \alpha \colon U \to \pow \cV \) be a morphism classifying an underobject \( i_A \colon A \tto U \times \cV \), and suppose \( i_A \pi_A \) is small.
  Let \( \und W \in \bC / A \) be the small hs-object classified by \( i_A \pi_{\cV} \colon A \to \cV \), let \( c \colon W \to W_U \) be an \( \unein \)-cocartesian morphism over \( i_A \pi_U \colon A \to U \), let \( q \colon \und W_U \to \und Z \) be an extensional collapse, and let \( i_Z \colon \und Z \to \und{\wh Z} \) be a top extension along \( (A, i_A \pi_U) \tox{\rt_{W}} \und W \tox{c} \und W_U \tox{q} \und Z \) (\jref{defn:top-extension}) in \( \bC / U \).

  Then \( \wh Z \) is an hs-object in \( \bC / U \) for which \( \delta_{\wh Z} \colon \wh Z \to U \) is \( \fS \)-small, and the morphism \( x \colon U \to \cV \) classifying \( \wh Z \) satisfies \( x \kappa = \alpha \), with \( \kappa \) the associated class morphism (\jref{defn:associated-class}).
\end{propn}

\begin{proof}
  See \sref{subsec:universe-closure-proof}.
\end{proof}

Finally, we state the converse to \jref{thm:topos-universe-small-closed}, i.e., the topos-theoretic version of \jref{propn:set-v-closed-under-sets}~\ref{item:set-v-closed-converse}:
\begin{propn} \label{propn:topos-universe-small-closed-converse}
  Let \( (\bC, \fS, \cV) \) be an hs-topoculized topos.
  Then for any \( x \colon U \to \cV \), if \( i_A \colon A \tto U \times \cV \) is the underobject classified by \( x \kappa \colon U \to \pow \cV \) then \( i_A \pi_A \colon A \to U \) is small.
\end{propn}

\begin{proof}
  See \sref{subsec:proof-of-small-converse}.
\end{proof}

As we explain in the proof, this proposition immediately reduces to the ``universal'' case \( x = \id_\cV \colon \cV \to \cV \),  which is of independent interest.
In this case, the underobject of \( \cV \times \cV \) classified by \( \id_{\cV} \cdot \kappa \colon U \to \pow \cV \) is, by definition of \( \kappa \), simply the relation \( \unein_\cV \tto \cV \times \cV \).
The claim thus becomes:
\begin{cor} \label{cor:small-closed-converse-univ-case}
  Let \( (\bC, \fS, \cV) \) be an hs-topoculized topos.
  Then the morphism \( i_{\ein_\cV} \pi_1 \colon \unein_\cV \to \cV \) is small.
\end{cor}

\begin{rmk} \label{rmk:initial-algebra}
  Given an object \( X \) in a topoculized topos \( (\bC, \fS) \), we can define an \emph{\( \fS \)-power object} of \( X \) to be an object \( \pow_\fS X \) classifying relations \( R \tto X \times U \) for which \( i_A \pi_U \colon R \to U \) is \( \fS \)-small; this is like the ``power class'' in class-set theory (class of \emph{small} subsets of a class).
  Then \jrefs{thm:topos-universe-small-closed}{propn:topos-universe-small-closed-converse} together amount to the statement that if \( \cV \) is an hs-classifier for \( \fS \), then it \emph{has} an \( \fS \)-power object, given by \( \cV \) itself.

  In \cite[\S7]{van-den-berg-moerdisjk-aspects-1}, where a universe is constructed in a similar ``category of classes'' context, the existence of an isomorphism \( \pow_\fS \cV \toi \cV \) is deduced, via a classic observation of Lambek's, from the stronger statement that \( \cV \) is an \emph{initial algebra} for the functor \( \pow_\fS \); this parallels the definition of the class \( \rV \) of hereditary sets as the least class closed under set-formation.
  In the present context (specifically, in the context of \jref{propn:2-topoi-give-topoculi}), it remains true that \( \cV \) is an initial \( \pow_\fS \)-algebra.
  However, the proof is somewhat involved, and as we will not need this result, we do not pursue it further here.
\end{rmk}

\subsection{The axioms of MK} \label{subsec:axioms-of-mk}
We recall the axioms of MK set theory.
This can be and is often formulated in single-sorted first-order logic, with a single binary relation \( \unin \), and a single predicate \( M \) (``sethood'').
However, it will be convenient for us to formulate it in the two-sorted language, consisting of:
\begin{itemize}
\item two sorts, called \( \St \) and \( \Cls \)
\item a binary relation \( \in \) of arity \( \St \times \Cls \)
\item an operation \( \kappa \) of arity \( \St \to \Cls \)
\end{itemize}
We always use Latin letters for variables of sort \( \St \) and Greek letters \( \alpha,\beta,\theta \) for variables of sort \( \Cls \) (but not \( \phi \), which we use to denote a formula).
For \( x \) and \( y \) variables of sort \( \St \), we write \( x \ein y \) for \( x \in \kappa(y) \).
The axioms are the universal closures of the following formulas; some explanations follow.

\begin{tabular}{ll}
  (Sethood) & \( \kappa(x) = \kappa(y) \to x = y \) \\
  (Extensionality) & \( (\forall x.\ x \in \alpha \tot x \in \beta) \to \alpha = \beta \) \\
  (Comprehension\( _{\phi,w,\alpha} \)) & \( \exists \alpha\ \forall w\ (w \in \alpha \tot \phi) \) \\
  (Foundation) & \( (\forall x.\ \kappa(x) \subset \alpha \to x \in \alpha) \to (\forall x.\ x \in \alpha) \) \\
  \multicolumn{2}{l}{\textit{Set existence axioms:}} \\
  (Empty set) & \( \exists z\ \forall w\ \pBig{w \ein z \tot \bot} \) \\
  (Pairing) & \( \exists z\ \forall w\ \pBig{w \ein z \tot (w = x \vee w = y)} \) \\
  (Union) & \( \exists z\ \forall w\ \pBig{w \ein z \tot \pbig{\exists v.\ w \ein v \wedge v \ein x}}  \) \\
  (Power) & \( \exists z\ \forall w\ \pBig{w \ein z \tot \pbig{\kappa(w) \subset \kappa(x)}}  \) \\
  (Subsets) & \( \exists z\ \forall w\ \pBig{w \ein z \tot (w \ein x \wedge w \in \beta)} \) \\
  (Replacement) & \( \exists z\ \forall w\ \pBig{w \ein z \tot \pbig{\exists u.\ u \ein x \wedge u \beta w \wedge (\forall v.\ u \beta v \to v = w)}} \) \\
  (Infinity) & \( \exists z\ \forall w\ \pBig{w \ein z \tot \forall \beta.\ \pbig{\emptyset \in \beta \wedge \forall v.\ (v \in \beta \to v^+ \in \beta)} \to w \in \beta} \)
\end{tabular}

Some explanations:
\begin{itemize}
\item Comprehension is an axiom schema, with one instance for each formula \( \phi \) in the language of MK, each variable \( w \), and each variable \( \alpha \) not free in \( \phi \).
\item In Foundation and Power: we write \( \kappa(x) \subset \alpha \) as a shorthand for \( \forall v.\ v \in \kappa(x) \to v \in \alpha \) (and likewise with \( \kappa(w) \subset \kappa(x) \)).
\item In Infinity, we write \( v^+ \in \alpha \) as a shorthand for \( \exists u.\ u \in \alpha \wedge \forall t.\ \pbig{t \ein u \tot (t \ein v \vee t = v)} \); and we write \( \emptyset \in \alpha \) as a shorthand for \( \exists u.\ u \in \alpha \wedge \forall t.\ \neg t \in u \).
\item In Replacement: we write \( u \beta w \) as a shorthand for \( \exists p.\ p = (u,w) \wedge p \in \beta \) (and likewise with \( u \beta v \)).
  Here, in turn, we write \( p = (u,w) \) as a shorthand for
  \[
    \forall t.\ t \in p \tot \pbig{(\forall s.\ s \in t \tot s = u)
    \vee \pbig{\forall s.\ s \in t \tot (s = u \vee s = w)}}
  \]
  (and likewise with \( p = (u, v) \)).
\end{itemize}

A comment on the formulation of Replacement:
a common formulation is ``the image of any set \( x \) under any class function \( \beta \) is again a set''.
The above statement, which is equivalent in the presence of Comprehension and Subsets, essentially says ``given a set \( x \) and any class \( \beta \), if \( \beta' \subset \beta \) is the subclass of \( \beta \) consisting of ordered pairs, and \( x' \subset x \) is the subset of \( x \) consisting of elements in the domain of \( \beta' \) and having a unique image under \( \beta' \), then the image of \( x' \) under \( \beta' \) is a set''.

\subsubsection{Models of MK in topoi}
By definition, an interpretation of the language of MK in a topos \( \bC \) is a quadruple \( (V, P V, \in_V, \kappa) \) consisting of objects \( V \) and \( P V \) in \( \bC \), a binary relation \( \unein_V \tto V \times P V \), and a morphism \( \kappa \colon V \to P V \); let us call such a quadruple an \defword{MK-object}.
Thus, we may ask for any MK-object whether it is a \emph{model of MK} (or some fragment of MK), i.e., whether it satisfies the above axioms.
We will be interested in MK-objects of the following particular kind.
\begin{defn} \label{defn:associated-mk}
  Given an \( \unein \)-object \( V \) in a topos \( \bC \) and a power object \( \pow V \) for \( V \), we define the \defword{associated MK-object} of \( V \) to be the MK-object \( (V, \pow V, \in_{V}, \kappa_V) \), where \( \kappa_V \colon V \to \pow V \) is the associated class morphism of \jref{defn:associated-class}.
  When we say that an \( \unein \)-object \( V \) is a \textbf{model} of (some fragment of) MK, we mean that the associated MK-object is.
\end{defn}
\begin{rmk}
  When interpreting the language of MK in the MK-object associated to an \( \unein \)-object \( V \), the binary relations \( x \ein y \defiff \kappa x \in \kappa y \) and \( \alpha \subset \beta \defiff \forall x.\ x \in \alpha \to x \in \beta \) appearing in the axioms of MK are interpreted as the relation \( \unein_V \tto V \times V \) and the subset relation \( {\le}_V \tto \pow V \times \pow V \), respectively.
\end{rmk}

Let us write \( \MKo \) for MK minus the axiom of Infinity.
\begin{thm} \label{thm:model-of-mk}
  In any hs-topoculized topos \( (\bC, \fS, \cV) \), \( \cV \) is a model of \( \MKo \).

  In particular (by \jref{propn:2-topoi-give-topoculi}) if \( \cC \) is a groupoidant 2-category and \( \rS,\rS' \in \cC \) are plentiful DOF classifiers for which \( \hs(\rS) \) and \( p \colon \rS_* \to \rS \) are \( \rS' \)-small, then the \( \unein \)-object \( \hs(\rS) \) in the topos \( \SOF_{\rS'}(\tm) \) is a model for \( \MKo \); and if \( \cC \) satisfies (UA), then for any plentiful DOF classifier \( \rS \), the \( \unein \)-object \( \hs(\rS) \) in the topos \( \SOF(\tm) \) is a model of \( \MKo \).
\end{thm}

\begin{proof}
  The satisfaction of the first four axioms is immediate using the already established (in \sref{subsec:univ-ext-wf}) extensionality and well-foundedness of \( \cV \), and the usual extensionality and comprehension properties of power objects in topoi (\sref{subsec:comprehension}).
  The verification of the set existence axioms, which we give in \sref{subsec:proof-of-set-existence-axioms}, is more involved and makes essential use of the results of \sref{subsec:set-formation-closure}.
\end{proof}
\noindent
We will return to the Axiom of Infinity in \sref{subsec:infinity}.

\subsubsection{Relationship to ZF set theory} \label{subsubsec:relation-to-zf}
The theory ZF is a single sorted theory with a single binary relation \( \unin \); it is obtained from the above axioms for MK by making the following changes:
\begin{itemize}
  \item We treat all variables appearing in the formulas as variables of the same single sort, and replace \( \kappa(x) \) with \( x \) for any variable \( x \).
  \item Comprehension and sethood are removed.
  \item Foundation becomes: for all formulas \( \phi \) with \( y \) not free:
    \[
      [\forall x.\ (\forall v.\ v \in x \to \phi[v/x] ) \to \phi] \to \forall x.\ \phi,
    \]
where \( \phi[v/x] \) is the result of substituting \( v \) at all free occurrences of \( x \) in \( \phi \).
  \item Subsets becomes: for all formulas \( \phi \) with \( z \) not free: \( \exists z\ \forall w.\ w \in z \tot ( w \in x \wedge \phi ) \).
  \item Replacement becomes: for all formulas \( \phi \) with \( v \) and \( z \) not free:
    \[
      \exists z.\ \forall w.\ w \in z \tot \exists u\ \pbig{u \in x \wedge \phi \wedge (\forall v.\ \phi[v/w] \to v = w)}.
    \]
\end{itemize}
We write \( \ZFo \) for ZF minus infinity.

It is well-known, and straightforward (though somewhat involved) to verify, that given an MK-object \( (V, PV, \in_V, \kappa_V) \) which is a model for MK (or \( \MKo \)), then \( (V, \ein_V) \) is a model of ZF (or \( \ZFo \)), where \( \unein_V \tto V \times V \) is the pullback of \( \in_V \) along \( V \times V \tox{\id_V \times \kappa_V} V \times P V \).
(More fundamentally, ZF admits an \emph{interpretation} into MK, in which the single sort of ZF is interpreted as \( \St \), and the binary relation \( x \in y \) of ZF is interpreted \( x \in \kappa(y) \).)

In particular if an \( \unein \)-object \( V \) is a model of MK, then it is also a model of ZF.
From \jref{thm:model-of-mk}, we thus immediately conclude:
\begin{cor}
  Under the same assumptions as \jref{thm:model-of-mk}, \( \cV \) is a model for \( \ZFo \).
\end{cor}

\subsection{The axiom of infinity and natural numbers objects} \label{subsec:infinity}
In the previous section, we showed that, given a groupoidant 2-category \( \cC \) and a plentiful DOF classifier \( \rS \), its hs-classifier \( \hs(\rS) \) gives a model of \( \MKo \) (as soon as it is small with respect to some other plentiful DOF classifier \( \rS' \)).
We cannot expect in general that \( \hs(\rS) \) will also satisfy the axiom of infinity.
Indeed, in the case \( \cC = \Cat \), we might have \( \rS \) the category of finite sets, in which case \( \hs(\rS) \) will be the setoid of hereditarily finite \( \unein \)-sets, which does not satisfy the axiom of infinity.
However, note that the plentiful DOF classifier \( \rS' \), with respect to which \( \hs(\rS) \) is small, can no longer be the setoid of finite sets, since \( \hs(\rS) \) is infinite.
Therefore, we should expect \( \hs(\rS') \) to satisfy the axiom of infinity.

We can split this argument into two stages: first, the fact that \( \hs(\rS) \) is \( \rS' \)-small forces \( \rS' \) to have a natural numbers object (NNO); and second, if \emph{any} plentiful DOF classifier \( \rS \) has a natural numbers object, then \( \cV(\rS) \) satisfies the axiom of infinity.
\begin{thm} \label{thm:infinity}\enumbelow
  \begin{enumerate}[(i)]
  \item\label{item:infinity-nno-exists} Let \( (\bC, \fS) \) be a topoculized topos.
    If \( \fS \) has an hs-classifier, then \( \bC \) hs an NNO.
  \item\label{item:infinity-small-nno} Given an hs-topoculized topos \( (\bC, \fS, \cV) \), the NNO of \( \bC \) is small if and only if \( \cV \) satisfies the Axiom of Infinity (in the sense of \jref{defn:associated-mk}).

    \noindent
    (In particular, if the NNO of \( \bC \) is small, then \( \cV \) is a model for MK (in the sense of \jref{defn:associated-mk}), and hence also of ZF (as per \sref{subsubsec:relation-to-zf}).)
  \item\label{item:infinity-2cat-has-nno}
    Let \( \rS \), \( \rS' \), and \( \rS'' \) be plentiful DOF-classifiers in a groupoidant 2-category \( \cC \), and suppose that \( \hs(\rS) \) and \( p \colon \rS_* \to \rS \) are \( \rS' \)-small, and that \( \hs(\rS') \) and \( p' \colon \rS'_* \to \rS' \) are \( \rS'' \)-small.
    Then in the hs-topoculized topos \( \SOF_{\rS''}(\tm) \) of \jref{propn:2-topoi-give-topoculi} (or in \( \SOF(\tm) \), when \( \cC \) satisfies axioms (UA)), the NNO is small, and hence the hs-classifier \( \hs(\rS') \) is a model of MK.

    \noindent
    (In particular, if \( \cC \) is a groupoidant 2-category satisfying axiom (UA), then the topos \( \SOF(\tm) \) contains a model \( \hs(\rS') \) of MK.)
  \end{enumerate}
\end{thm}
\begin{proof}
  See \sref{subsec:infinity-proof}.
\end{proof}

\appendix

\section{Background on topoi and logic} \label{sec:topos-logic-background}
In this appendix, we continue what we started in \sref{subsec:cat-prelims}, namely recalling some basic notions and establishing notational conventions regarding topoi, but now with a particular view to the logical aspects of topoi.

\subsection{Interpreting logic in topoi} \label{subsec:logic-in-topoi}
\subsubsection{Logic}
We take the basic notions of formal logic (\emph{terms} and \emph{formulas}, \emph{free} and \emph{bound} variables, \emph{substitution}, etc.) for granted.
As a first remark concerning our conventions, we regard terms and formulas as abstract entities---elements of a certain set endowed with and freely generated by certain operations (namely the logical connectives)---rather than as sequences of elements (``symbols'') from a certain set (the ``alphabet'').
Thus, when we write down a formulas such as \( \forall x.\ \exists y.\ x = y \), this is simply a piece of notation to \emph{denote} a certain element of the set of formulas, rather than being the formula itself, \emph{qua} list of symbols.
Correspondingly, we can be flexible, and even inconsistent, with our notation for formulas, as long as it is always clear which formula is being denoted.

By a \defword{signature} \( \cL \), we always mean a multi-sorted first-order signature; thus, \( \cL \) is given by a set \( \Ob \cL \) of \emph{sorts}, together with sets \( \Fun \cL \) of function symbols and \( \Rel \cL \) of relation symbols of given sorted arities.
We write \( f \colon \vA \to B \) to indicate that \( f \) is a function symbol with domain \( \vA \) (a sequence of sorts) and codomain \( B \), and \( R \tolon \vA \) to indicate that \( R \) is a relation symbol of arity \( \vA \).

Given a set \( X \), write \( X^{< \omega} \) for the set of finite sequences in \( X \), and \( X^{< n} \) for the set of sequences of length less than \( n \in \N \).
Given a sequence \( \vx \in X^{<\omega} \), we write \( \len(\vx) \) for its length and \( x_i \) (\( 0 \le i < \len(\vx) \)) for its \( i \)-th element.
We write \( \vx \vy \) or \( \vx, \vy \) for the concatenation of sequences, and similarly \( \vx y \) or \( \vx, y \) for the extension of a sequence by a single element.

Variables are always sorted.
We write \( x \tolon A \) to express that \( x \) is a variable of sort \( A \), and \( \vx \tolon \vA \) to mean \( x_i \tolon A_i \) for all \( i \).
In formulas, when the sorts of variables aren't clear from context, we may include the sorts under the quantifiers as in ``\( \forall x \tolon A \)'' or ``\( \exists x \tolon A \)''.

We have the usual notion of \emph{first-order formula} over \( \cL \).
We define a \defword{higher-order formula} over \( \cL \) to be a first-order formula over the auxiliary signature \( \cL^\up \) extending \( \cL \), defined by adjoining sorts \( \cP^n A \) for each \( A \in \Ob \cL \) and \( n \ge 1 \), and relation symbols \( \in_{A,n} \colon \br{\cP^n A,\cP^{n+1} A} \) (which we usually just denote by \( \in \)) for each \( n \ge 0 \), where we write \( \cP^0 A \) for \( A \).
By a \defword{\( d \)-th order formula}, we mean a formula over the sub-signature \( \cL^{\up d} \) consisting only of the sorts \( \cP^n \) for \( n < d \).

\subsubsection{Interpretations} \label{subsubsec:interpretations}
Given a signature \( \cL \) and a topos (or more generally finite product category) \( \bC \), we have the usual notion of an \emph{interpretation} \( M \colon \cL \to \bC \) of \( \cL \) in \( \bC \), consisting of an object \( M A \) for each \( A \in \Ob \cL \), together with morphisms \( M f \colon M \vA \to M B \) for each function symbol \( f \colon \vA \to B \) and monomorphisms \( i_{M R} \colon M R \to M \vA \) for each relation symbol \( R \tolon \vA \).
Here, for a sequence \( \vA \in (\Ob \cL)^n \), we are writing \( M \vA \) for a fixed product \( M A_0 \times \cdots \times M A_{n-1} \) in \( \bC \), and we take as part of the data of \( M \) the choice of such a product for each sequence \( \vA \) appearing in the arity of a function or relation symbol in \( \cL \).
Any interpretation \( M \colon \cL \to \bC \) in a topos \( \bC \) can be extended to an interpretation \( M^\up \colon \cL^\up \to \bC \), in which \( M^\up (\unin_{A,n}) \tto M(\pow^n A) \times M(\pow^{n+1} A) \) exhibits \( M(\pow^{n+1} A) \) as a power object of \( M (\pow^n A) \) for all \( n \ge 0 \).
We call such an extension \( M^\up \) a \defword{power extension} of \( M \), and we note that it is unique up to a unique isomorphism (in the sense to be explained presently) extending \( \id_M \).\footnote{ %
  It is better to say that \( M^\up \) is a canonically defined ``ana-interpretation''.
  That is, for any particular sort \( \pow^n A \in \Ob \cL^\up \), we can say what we mean by \( M^\up(\pow^n A) \) without defining the full interpretation \( M^\up \); namely, it is a power object of \( M^\up(\pow^{n-1} A) \), where the latter object is, recursively, understood in the same way (the base case being \( n = 0 \), where \( M^\up(\pow^0 A) = M A \)).
  In accordance with our conventions in \sref{subsubsec:discipline}, we will always use \( M^\up(\pow^n A) \) in this way, rather than assuming we have chosen the whole interpretation \( M^\up \)---and moreover, we will just write \( M(\pow^n A) \) in place of \( M^\up (\pow^n A) \).
}

An \defword{isomorphism} of interpretations \( f \colon M_0 \toi M_1 \) consists of isomorphisms \( f \colon M_0 A \toi M_1 A \), commuting in the obvious sense with the morphisms \( M_0 f \) and \( M_1 f \) for each function symbol \( f \) and with \( i_{M_0(R)} \) and \( i_{M_1(R)} \) for each relation symbol \( R \).
There is an obvious category of interpretations and isomorphisms between them, which we denote by \( (\bC^\cL)^\iso \).
A logical functor (see \sref{subsec:logical-functors} below) \( F \colon \bC \to \bD \) can be composed with an interpretation \( M \colon \cL \to \bC \) in an evident way to produce an interpretation \( M F \colon \cL \to \bD \), and in this way we obtain a functor \( F_* \colon \bC^\cL \to \bD^\cL \).
Similarly an isomorphism \( \beta \colon F \Toi G \) between logical functors \( F,G \colon \bC \to \bD \) induces an isomorphism of models \( M \beta \colon M F \toi G F \), and a natural isomorphism \( F_* \Toi G_* \).

By a \defword{term-in-context} (resp.\ \defword{formula-in-context}) over \( \cL \), we mean a pair \( (\square, \vx) \) with \( \square \) a term (resp.\ higher-order formula) over \( \cL \) and \( \vx \) a sequence of distinct variables including all the free variables in \( \square \).
Given an interpretation \( M \colon \cL \to \bC \), a term or formula in context \( (t, \vx) \) or \( (\phi, \vx) \) with \( \vx \tolon \vA \in \Ob (\cL^\up)^n \), and given a product \( M(\vA) \) of \( M(A_0),\cdots,M(A_{n-1}) \) in \( \bC \), there is an associated \emph{interpretation} of \( (t, \vx) \) or \( (\phi, \vx) \), the former being a morphism \( M_{\vx}(t) \colon M(\vA) \to M(B) \) and the latter being a subobject (i.e., equivalence class of underobjects) of \( M(\vA) \); we typically write \( i_{M_{\vx}(\phi)} \colon M_{\vx}(\phi) \tto M(\vA) \) for a fixed underobject in this equivalence class.
For terms and first-order formulas, this is defined as in \cite{makkai-reyes-focl}, the definition being a straightforward, recursive one, and it is extended to higher-order formulas using the (essentially unique) extensions \( M^\up \colon \cL^\up \to \bC \) of \( M \).
In particular, if \( \phi \) is a \emph{sentence} (that is, a closed formula), then \( M_{\emptyset}(\phi) \) (where \( \emptyset \) is the empty sequence) is an underobject of the terminal object \( \tm \), and we say that \( M \) \emph{satisfies} \( \phi \), and write \( M \vDash \phi \)---or just \( \bC \vDash \phi \) when \( M \) is clear from context---if \( M(\phi, \emptyset) \approx \tm \).

The basic and central fact about the interpretation of logic in topoi is its \emph{soundness} with respect to intuitionistic logic: any sentence \( \phi \) which is provable in intuitionistic logic is satisfied by any interpretation \( M \) in any topos---however, we will not in fact be making any use of this in this paper (see the discussion in \S\ref{subsec:sets-and-topoi}).

A \defword{higher-order theory} (resp.\ \defword{\( d \)-th order theory}) \( \cT \) over \( \cL \) is a set of higher-order (resp.\ \( d \)-th order) sentences over \( \cL \).
When we say that a higher-order theory \( \cT \) is \emph{finite}, we mean not only that the set \( \cT \) is finite, but also the language \( \cL \) over which it is defined.
A \defword{model} for \( \cT \) in a topos \( \bC \) is an interpretation \( M \colon \cL \to \bC \) with \( M \vDash \phi \) for all \( \phi \in \cT \).
We denote by \( (\bC^\cT)^\iso \) the \defword{groupoid of models of \( \cT \) in \( \bC \)}, i.e., the full subcategory of \( (\bC^{\cL})^\iso \) consisting of models of \( \cT \).\footnote{%
  \label{footnote:cotensor-atomic}%
  Contrary to what is suggested by the notation, we do not introduce an object \( \bC^\cT \) of which \( (\bC^\cT)^\iso \) is the core.
}
The composite of a model of \( \cT \) with a logical functor \( F \colon \bC \to \bD \) is again a model, and there is thus induced functors \( F_* \colon (\bC^\cT)^\iso \to (\bD^\cT)^\iso \); and similarly, an invertible 2-cell \( \beta \colon \bC \tocell \bD \) induces a natural isomorphism \( \beta_* \colon (\bC^\cT)^\iso \tocell (\bD^\cT)^\iso \).

Associated to any topos \( \bC \), we have the \defword{internal language} \( \cL_{\bC} \) of \( \bC \): its sorts are the objects of \( \bC \); its function symbols of arity \( \vA \to B \) are pairs \( (P,f) \) with \( P \) a product \( (\pi_i \colon P \to A_i)_{i=0}^{\len \vA - 1} \) and \( f \colon P \to B \) a morphism; and its relations of arity \( \vA \) are pairs \( (P,R) \) with \( P \) a product as above, and \( R \tto P \) an underobject.
We always just write \( f \) for \( (P,f) \) and \( R \) for \( (P,R) \).

The internal language comes equipped with a tautological interpretation \( M_{\mathrm{taut}} \colon \cL_{\bC} \to \bC \) (and in fact, a general interpretation \( \cL \to \bC \) is nothing but an interpretation of languages \( \cL \to \cL_{\bC} \), suitably defined).
Whenever we speak of interpreting a certain sentence in a topos, and we have not explicitly specified the language, we always assume we are using the internal language.
Thus, for example, given an underobject \( R \tto X \times X \) of an object \( X \),
\[
  \bC \vDash \forall x,y,z \tolon X.\ (x R y \wedge y R z) \to x R z
\]
expresses that \( R \) is a transitive relation on \( X \).

Given a sequence of sorts \( \vA \in (\Ob C)^n \), a product \( \prod_{i=0}^{n-1} A_i \) in \( \bC \), and a formula-in-context \( (\phi, \vx) \) over \( \cL_{\bC} \) with \( \vx \tolon \vA \), we will write
\[
  \set{\vx \in \tprod_{i=0}^{n-1} A_i \mid \phi} \tto \prod_{i=0}^{n-1} A_i
  \quad
  \text{or just}
  \quad
  \set{\vx \mid \phi} \tto \prod_{i=0}^{n-1} A_i
\]
in place of \( (M_{\mathrm{taut}})_{\vx}(\phi) \tto \prod_{i=0}^{n-1} A_i \).

\begin{exm} \label{exm:im-as-setbuilder}
  The \emph{image} \( \im(f) \tto Y \) of a morphism \( f \colon X \to Y \) can be described as the underobject \( \set{ y \in Y \mid \exists x \tolon X.\ f(x) = y } \tto Y \); hence (by \jref{thm:kj-semantics} below), \( y \colon U \to Y \) factors through \( \im(f) \) if and only if there exists \( U \oot U' \tox{x} X \) with \( x f = y' \).
  In particular, \( f \colon X \to Y \) is epic iff \( \bC \vDash \forall y \tolon Y.\ \exists x \tolon X.\ f(x) = y \), and more generally, a family \( \set{f_j \colon X_j \to Y}_{j=0}^{n-1} \) is a cover iff \( \bC \vDash \forall y \tolon Y.\ \bigvee_{j=0}^{n-1} \exists x_j \tolon X_j.\ f_j(x) = y \).

  The \emph{subset} relation \( {\le_X} \tto \pow X \times \pow X \) can be described as the underobject \( \set{ \br{p, q} \in \pow X \times \pow X \mid \forall x.\ x \in p \to x \in q } \).
  Thus (again by \jref{thm:kj-semantics} below) \( p,q \colon U \to X \) satisfy \( p \le q \) iff \( x \in p' \To x \in q' \) for all \( x \colon U' \to X \); it follows immediately that \( {\le_X} \) is a \emph{preorder}; that it is a \emph{partial order} (i.e., anti-symmetric) follows from \emph{extensionality} (\jref{propn:topos-extensionality} below).
\end{exm}

\subsubsection{Kripke-Joyal semantics} \label{subsubsec:kj-semantics}
For the most part, whenever we interpret a sentence \( \phi \) in a topos \( \bC \) in this paper, we will understand the meaning of \( \bC \vDash \phi \) using the \emph{Kripke-Joyal} (KJ) semantics, which we now recall.

Fix a signature \( \cL \) and an interpretation \( M \colon \cL \to \bC \) in a topos \( \bC \).

Given a formula-in-context \( (\phi,\vx) \) over \( \cL \) with \( \vx \tolon \vA \in \Ob (\cL^\up)^n \), an object \( U \in \bC \), and morphisms \( a_i \colon U \to M(A_i) \) for \( 0 \le i < n - 1 \), we write \( U \Vdash^M_{\vx \mapsto \va} \phi \), or just \( U \Vdash_{\vx \mapsto \va} \phi \), and say that \defword{\( (U,\va) \) forces \( (\phi,\vx) \)}, if the morphism \( \br{a_i}_i \colon U \to M(\vA) \) factors through \( M_{\vx}(\phi) \tto M(\vA) \).
Note that, for a sentence \( \phi \), we have \( M \vDash \phi \) iff \( \tm \Vdash \phi \) iff \( U \Vdash \phi \) for all \( U \in \bC \).

It will often be convenient to use the same symbols to denote the variables \( x_i \tolon A_i \) and the morphisms \( x_i \colon U \to M(A_i) \), in which case we may just write \( U \vDash_{\vx} \phi \) in place of \( U \vDash_{\vx \mapsto \vx} \phi \).
We will use this abbreviated notation more generally if the names of the variables \( x_i \tolon A_i \) differ only from the names of the morphisms \( U \to M(A_i) \) by certain decorations on the latter, as per \jref{rmk:kj-bookkeeping} below; e.g., we will write \( U' \Vdash_{\vx'} \phi \) and \( U_0 \Vdash_{\vx_0} \phi \) in place of \( U' \Vdash_{\vx \mapsto \vx'} \phi \) and \( U_0 \Vdash_{\vx \mapsto \vx_0} \phi \), respectively.

The ``KJ semantics'' refers to the following theorem (for which see, e.g., \cite[\S\S VI.6-7]{mac-lane-moerdijk}); see \jref{rmk:kj-bookkeeping} for the notation \( \vx', \vx_0, \vx_1 \).
\begin{thm} \label{thm:kj-semantics}
  Given formulas-in-context \( (\phi, \vx) \) and \( (\psi,\vx) \) over a signature \( \cL \) with \( \vx \tolon \vA \), an interpretation \( M \colon \cL \to \bC \) in a topos \( \bC \), an object \( U \in \bC \) and morphisms \( x_i \colon U \to M(A_i) \), we have:

  \begin{tabular}{lll}
    \( U \Vdash_\vx \top \) & always & \\
    \( U \Vdash_\vx \bot \) & iff & \( U \) is initial \\
    \( U \Vdash_\vx \phi \wedge \psi \) & iff & \( U \Vdash_{\vx} \phi \) and \( U \Vdash_{\vx} \psi \) \\
    \( U \Vdash_\vx \phi \vee \psi \) & iff & \( U_0 \Vdash_{\vx_0} \phi \) and \( U_1 \Vdash_{\vx_1} \psi \) for some cover \( U_0,U_1 \to U \) \\
    \( U \Vdash_\vx \phi \to \psi \) & iff & \( U' \Vdash_{\vx'} \phi \) implies \( U' \Vdash_{\vx'} \psi \) for all \( U' \to U \) \\
    \( U \Vdash_\vx \neg \phi \) & iff & \( U' \Vdash_{\vx'} \phi \) implies that \( U' \) is initial for all \( U' \to U \) \\
    \( U \Vdash_\vx \forall y \tolon B.\ \phi \) & iff & \( U' \Vdash_{\vx' y} \phi \) for all \( U' \to U \) and all \( y \colon U' \to M(B) \) \\
    \( U \Vdash_\vx \exists y \tolon B.\ \phi \) & iff & \( U' \Vdash_{\vx' y} \phi \)  for some epi \( U' \too U \) and some \( y \colon U' \to M(B) \)
  \end{tabular}\\

  Moreover:

  If \( U \Vdash_\vx \phi \), then \( U' \Vdash_{\vx'} \phi \) for any \( U' \to U \);
  and conversely, if \( U_0, \ldots ,U_{n - 1} \to U \) is a cover (see \sref{subsec:cat-prelims}), and \( U_i \Vdash_{\vx_i} \phi \) for each \( i \), then \( U \Vdash_{\vx} \phi \).
\end{thm}

\begin{rmk} \label{rmk:kj-bookkeeping}
  By a \defword{generalized element} of an object \( X \) in a category \( \bC \), we simply mean a morphism \( x \colon U \to X \); we call the domain \( U \) the \defword{base} of the generalized element \( x \).
In case the base \( U \) is terminal, we call \( x \) a \defword{global element}.

  The chief utility of the KJ semantics is that it allows us to reason about the internal logic of topoi entirely in terms of generalized elements, and thus brings it very close to ordinary, set-theoretic arguments.
  The price one has to pay for this is a certain amount of ``bookkeeping'' related to changes of the base \( U \).
  To deal with this added complexity, it is helpful to introduce some notational conventions, which we now describe:
  \begin{itemize}
  \item Given a relation \( R \tto X_0 \times \cdots \times X_{n-1} \) and generalized elements \( x_i \colon U \to X_i \) for each \( i \), we write \( R(x_0,\ldots,x_{n-1}) \) for \( U \Vdash_{\vx} R(x_0,\ldots,x_{n - 1}) \)---in other words, to express that \( \br{x_0,\ldots,x_{n-1}} \to X_0 \times \cdots \times X_{n-1} \) factors through \( R \).
    When \( n = 2 \), we may instead use the infix notation \( x_0 R x_1 \).
    For instance, given a power object \( \pow X \) of an object \( X \) and generalized elements \( x \colon U \to X \) and \( p \colon U \to \pow X \), we would write \( x \in p \) to express that \( \br{x,p} \colon U \to X \times \pow X \) factors through \( \unin_X \tto X \times \pow X \).

    Incidentally, as mentioned in \sref{subsec:cat-prelims}, we use the ``geometric'' order \( x f \) for composing morphisms \( U \tox{x} X \tox{f} Y \) in a category.
    However, when writing first-order formulas, we use the conventional order \( f(x) \) for application; this leads to a certain reversal when interpreting formulas.
    For example, given morphisms \( x \) and \( f \) as above and an underobject \( P \tto Y \), the expression ``\( P(x f) \)'' expresses that \( U \Vdash_x P\pbig{f(x)} \).
  \item Given a generalized element \( x \colon U \to X \) and a new base \( \delta_{U'} \colon U' \to U \), we simply write \( x' \) as a shorthand for \( \delta_{U'} x \colon U' \to X \).
    For a sequence \( \vx \), we also write \( \vx' \) for \( (x_0',\ldots,x_{\len \vx - 1}') \).

    Similarly, given \( \delta_{U_0} \colon U_0 \to U \) or \( \delta_{U_1} \colon U_1 \to U \), we will write \( x_0 = \delta_{U_0} x \colon U_0 \to X \) and \( x_1 = \delta_{U_1} x \colon U_1 \to X \).

    If, as often happens, we have a morphism \( y \colon U' \to X \) and then we further introduce \( \delta_{U''} \colon U'' \to U' \), note that we write \( y'' \), and not \( y' \), for the morphism \( \delta_{U''} y \colon U'' \to X \).
  \item We will often introduce a new base \( U' \to U \) without explicitly naming the morphism.
    In such cases, we will denote this morphism, which we call the \defword{anchor} of the base \( U' \), by \( \delta_{U'} \) (this is consistent with our conventions for slice categories in \sref{subsec:cat-prelims}).
    Similarly, the morphism \( U'' \to U' \) will be called \( \delta_{U''} \), and hence the resulting morphism \( U'' \to U \) is \( \delta_{U''} \delta_{U'} \).
    Often, we will not even mention the anchor explicitly, but simply write something like ``let \( x \colon U \to X \) and \( y \colon U' \to X \)'', in which case, we are implicitly fixing an arbitrary object \( U' \) and morphism \( \delta_{U'} \colon U' \to U \).
    In some instances, it will be important that the anchor is an epi, in which case we may write something like ``let \( x \colon U \to X \) and \( U \oot U' \tox{y} X \)''.
  \end{itemize}
\end{rmk}

We emphasize again that, throughout this paper, we generally use the KJ semantics to understand the meaning of sentences interpreted in a topos, and we expect the reader to automatically carry out the KJ interpretations.
Thus, for example, given a morphism \( f \colon X \to Y \) in a topos \( \bC \), if we write \( \bC \vDash \forall y.\ \exists x.\ f(x) = y \), this should be read as: ``for every \( y \colon U \to Y \), there exists \( U \oot U' \tox{x} X \) such that \( x f = y' \)''.

\subsection{Unique existential quantification and effective epis} \label{subsec:uniq-ex-eff-epi}
Consider a formula-in-context \( (\phi, \vx y) \) over the internal language \( \bC_\cL \) of some topos \( \bC \) with \( \vx \colon \vX \) and \( y \colon Y \), and fix an object \( U \in \bC \) and morphisms \( x_i \colon U \to X_i \).
In general, from
\[
  U \Vdash_{\vx} \exists y.\ \phi
\]
we can only conclude the existence of \( U \oot U' \tox{y} X \) with
\[
  U' \Vdash_{\vx',y} \phi.
\]
However, we can sometimes draw the stronger conclusion that we can take \( U' = U \).
Namely, this happens in the context of a \emph{unique} existential quantifier.

The key point is the following general fact, known as ``effectiveness of epimorphisms'' (see, e.g., \cite[Corollary~C2.1.12~(d)]{johnstone-topos-theory}).
\begin{propn} \label{propn:epis-eff}
  Let \( \delta_{U'} \colon U' \to U \) be an epimorphism in a topos \( \bC \), and form the pullback \( U' \times_U U' \).
  Then
  \[
    \begin{tikzcd}
      U' \times_U U' \ar[r, shift left, "\pi_0"] \ar[r, shift right, "\pi_1"'] & U' \ar[r, "\delta_{U'}"] & U
    \end{tikzcd}
  \]
  is a coequalizer diagram.

  Thus, for any morphism \( x \colon U' \to X \), if \( \pi_0 x = \pi_1 x \colon U' \times_U U' \to X \), then there is a unique \( \bar x \colon U \to X \) with
  \[
    x = \delta_{U'} \bar x \colon U' \to X.
  \]

  (And more generally, given a cover \( U_0, \ldots, U_{n-1} \to U \) and pullbacks \( U_i \xot{\pi_i^{ij}} U_i \times_U U_j \tox{\pi_j^{ij}} \to U_j \), and given morphisms \( x_i \colon U_i \to X \), there is a morphism \( x \colon U \to X \) with \( \delta_{U_i} x = x_i \) for all \( i \) if and only if \( \pi^{ij}_i x_i = \pi^{ij}_j x_j \) for all \( i,j \).)
\end{propn}

We now explain the application of this proposition to unique existential quantification.
This will involve substitution (of a very simple kind).
For a formula \( \phi \) and a variable \( z \) which (for simplicity) does not appear in \( \phi \), we write \( \phi[z/y] \) for the resulting of substituting \( z \) for the free occurrences of \( y \) in \( \phi \).
As regards the KJ semantics, the key property regarding substitution is the following:
if \( (\phi, \vx y) \) is a formula-in-context over a signature \( \cL \) with \( \vx \tolon \vA \) and \( y \tolon B \), and \( M \colon \cL \to \bC \) is an interpretation in a topos \( \bC \), and \( z \tolon B \) is a variable distinct from the \( x_i \) and not appearing in \( \phi \), and given morphisms \( a_i \colon U \to M(A_i) \) and \( b \colon U \to M(B) \) in \( \bC \), we have:
\begin{equation}\label{eq:subst-and-kj}
  U \Vdash_{\vx,y \mapsto \va,b} \phi
  \qquad
  \IFF
  \qquad
  U \Vdash_{\vx,z \mapsto \va,b} \phi[z/y].
\end{equation}
Indeed, this follows from the fact that \( M_{\vx, y} \phi \approx M_{\vx, z} \phi[z / y] \tto M(\vA B) \); see \cite[Lemma~3.2.3]{makkai-reyes-focl}.

We thus have that
\[
  U \Vdash_{\vx \mapsto \va} \forall y,z.\ (\phi \wedge \phi[z / y]) \to y = z
\]
if and only if \( b = c \) for any \( b,c \colon U' \to M(B) \) with \( U' \Vdash_{\vx,y \mapsto \va',b} \phi \) and \( U' \Vdash_{\vx,y \mapsto \va',c} \phi \).

\begin{propn} \label{propn:ex-uniq}
  Let \( (\phi, \vx y) \), \( z \), \( \bC \), \( M \), \( \va \), and \( b \) be as above.
  Suppose that
  \[
    U \Vdash_{\vx \mapsto \va} (\exists y.\ \phi)
    \qquad\AND\qquad
    U \Vdash_{\vx \mapsto \va} \pbig{\forall y,z.\ (\phi \wedge \phi[z/y]) \to y = z}.
  \]
  Then there exists a (necessarily unique) \( b \colon U \to X \) with
  \[
    U \Vdash_{\vx,y \mapsto \va,b} \phi.
  \]
\end{propn}

\begin{proof}
  By the first assumption, we have \( U \oot U' \tox{b} X \) with \( U' \Vdash_{\vx,y \mapsto \va',b} \phi \).
  It follows that \( U' \times_U U' \Vdash_{\vx,y \mapsto \pi_j \va', \pi_j b} \phi \) for \( j = 0,1 \), and hence, by the second assumption, and since \( \pi_0 \va' = \pi_1 \va' \), that \( \pi_0 b = \pi_1 b \).
  Hence, by \jref{propn:epis-eff}, \( b \) factors through some \( \bar b \colon U \to X \), and we have \( U \Vdash_{\vx, y \mapsto \va, \bar b} \phi \) as desired, since \( \delta_U \colon U' \to U \) is epic.
\end{proof}

\begin{cor} \label{cor:ex-mono-uniq}
  If \( i \colon Y \to X \) is a monomorphism in a topos \( \bC \), and \( x \colon U \to X \) satisfies
  \[
    U \Vdash_x \exists y \colon A.\ y i = x,
  \]
  then there exists \( y \colon U \to Y \) with \( y i = x \).
\end{cor}

\begin{proof}
  Apply \jref{propn:ex-uniq} with \( \phi \) given by \( y i = x \).
\end{proof}

\subsection{Comprehension} \label{subsec:comprehension}
We discuss the two fundamental principles of higher-order logic, extensionality and comprehension, both of which are satisfied in any topos.
Let \( \cL \) be a signature.
For any sort \( X \in \Ob \cL \), extensionality for \( X \) is the sentence
\[
  \forall p, q \tolon \pow X.\ (\forall x \tolon X.\ x \in p \tot x \in q) \to p = q.
\]
The claim that this holds under any interpretation \( M \colon \cL \to \bC \) in a topos \( \bC \) is thus the following:
\begin{propn} \label{propn:topos-extensionality}
  Given any object \( X \in \bC \) in a topos \( \bC \), and given any \( p,q \colon U \to \pow X \), if \( U \Vdash_{p,q} \forall x \tolon X.\ x \in p \ToT x \in q \), then \( p = q \).
\end{propn}

\begin{proof}
  To see this, consider the underobjects \( P,Q \mapsto X \times U \) classified by \( p \) and \( q \) respectively.
  Then for \( x \colon U' \to X \), we have that \( x \in p' \) iff \( \br{x,\delta_{U'}} \colon U' \to U \times X \) factors through \( P \tto U \times X \), and likewise with \( q \) and \( Q \).
  The assumption thus implies that the underobjects \( P, Q \tto U \times X \) are equivalent and hence that \( p = q \).
\end{proof}

We next turn to comprehension, which is a sentence schema, with one sentence for each formula-in-context \( (\phi, \vx y) \) with \( \vx \tolon \vX \) and \( y \tolon Y \), and each variable \( p \tolon \pow X \) distinct from the \( \vx \) and \( y \) and not occurring in \( \phi \).
The statement of comprehension is then:
\[
  \forall \vx.\ \exists p \tolon \pow X.\ \forall y.\ y \in p \tot \phi.
\]
We note that, assuming extensionality, the \( p \) whose existence is being asserted here is unique.
The claim that this holds under in any topos is thus the following:
\begin{propn} \label{propn:topos-comprehension}
  Let \( \vX,Y \) be objects in a topos \( \bC \) with \( n = \len \vX \), and fix a formula-in-context \( (\phi, \vx y) \) over \( \cL_\bC \), where \( \vx \tolon \vX \) and \( y \tolon Y \).
  Then for any morphisms \( x_i \colon U \to X_i \) for \( i = 0,\ldots,n - 1 \), there exists a (necessarily unique) \( p \colon U \to \pow Y \) such that \( U \Vdash_{\vx p} \forall y.\ y \in p \tot \phi \).
\end{propn}

\begin{proof}
  Consider the underobject \( \set{\vx y \mid \phi} \tto \prod_{i=0}^{n-1} X_i \times Y \), and form the pullback
  \[
    \begin{tikzcd}
      R \pb \ar[r, ""] \ar[d, >>->, ""'] &[35pt] \set{\vx y \mid \phi} \ar[d, "", >>->] \\
      U \times Y \ar[r, "\br{x_i}_{i} \times \id_Y"] & \prod_{i=0}^{n-1} X_i \times Y.
    \end{tikzcd}
  \]
  The morphism \( p \colon U \to \pow Y \) classifying \( R \) then has the desired property.
\end{proof}

An equivalent statement is the following:
\begin{propn} \label{propn:comprehension-morphism}
  With \( \vX \), \( Y \), and \( (\phi, \vx y) \) as in \jref{propn:topos-comprehension}, there exists a unique morphism \( f \colon \prod_{i=0}^{n-1} X_i \to \pow Y \) with
  \[
    \bC \vDash \forall \vx.\ \forall y .\ y \in f(\vx) \tot \phi.
  \]
\end{propn}

\begin{proof}
  Take \( f \) to be the morphism classifying \( \set{\vx y \mid \phi} \tto \prod_{i=0}^{n-1} X_i \times Y \).
  Uniqueness follows from \jref{propn:topos-extensionality}, regarding \( f \) as a gentled element \( U = \prod_{i=0}^{n-1} X_i \to \pow Y \).
\end{proof}

\begin{exm} \label{exm:im-preim-by-comp}
  Given a morphism \( f \colon X \to Y \) and applying \jref{propn:comprehension-morphism} to the formulas-in-context
  \[
    \pbig{(\exists x \tolon X.\ x \in p \wedge f(x) = y), \br{p, y}}
    \quad\AND\quad
    \pbig{(f(x) \in q), \br{q, x}},
  \]
  respectively, we obtain the image and preimage morphisms \( \exists_f \colon \pow X \to \pow Y \) and \( f\I \colon \pow Y \to \pow X \).
  Thus, given \( p \colon U \to \pow X \) and \( q \colon U \to \pow Y \), as well as \( x \colon U \to X \) and \( y \colon U \to Y \), we have \( x \in q f\I \) if and only if \( x f \in q \), and \( y \in p \exists_f \) if and only there is \( U \oot U' \tox{x} X \) with \( x \in p' \) and \( x f = y' \).
\end{exm}

\subsection{Logical functors} \label{subsec:logical-functors}
We recall that a functor \( F \colon \bC \to \bD \) between topoi is \defword{logical} (or a \defword{morphism of elementary topoi}) if it preserves finite limits and power objects.
The significance of logical functors is that they preserve ``all'' of the structure of a topos.

For instance, we have the following:
\begin{propn}[{\cite[Corollary~1.37]{johnstone-topos-theory}}] \label{propn:logical-fun-pres-colimits}
  Any logical functor preserves finite colimits, and in particular epimorphisms.
\end{propn}

\begin{propn}[{\cite[p.~155]{johnstone-topos-theory}}]\label{propn:logical-fun-pres-satisf}
  Let \( \cL \) be a signature and \( M \colon \cL \to \bC \) an interpretation in a topos \( \bC \).
  Then for any sentence \( \phi \) over \( \cL \) and any logical functor \( F \colon \bC \to \bD \), if \( M \vDash \phi \), then \( M F \vDash \phi \).

  More generally, for any formula-in-context \( (\phi, \vx) \) over \( \cL \) with \( \vx \tolon \vA \),
  we have \( F\pbig{M_\vx(\phi)} \approx (M F)_\vx(\phi) \tto (M F)(\vA) \).
\end{propn}
\noindent
The proof of these propositions involves going over the proof that any elementary topos (i.e., category with finite limits and power objects) has finite colimits and supports the interpretation of full intuitionistic first-order logic (i.e., is a Heyting category); and checking that, not only does the \emph{existence} of finite limits and power objects imply the existence of finite colimits and the appropriate logical operations, but the \emph{preservation} of finite limits and power objects by a functor implies that these operations are preserved as well.

We will often be using \jref{propn:logical-fun-pres-satisf} to deduce that some property holds in \( \bD \) from the fact that a corresponding property holds in \( \bC \), leaving it to the reader to check that the property in question is expressible as the satisfaction of some formula.
Some sample applications include (with \( F \colon \bC \to \bD \) a logical functor):
\begin{itemize}
\item If \( X \) is an hs-object in \( \bC \) (in the sense of \jref{defn:topos-ein-hs}), then \( F X \) is an hs-object in \( \bD \), with underlying relation \( F i_{\unein_X} \colon F \unein_X \to F(X \times X) \), where \( F(X \times X) \) is given the product structure \( F X \xot{F \pi_0} F(X \times X) \tox{F \pi_1} F X \).
\item If \( X \) is an \( \unein \)-object in \( \bC \), then the closure morphism \( \dwn_{F X} \) for \( F X \) (\jref{defn:topos-dwn}) is given by \( F \dwn_X \colon F X \to F (\pow X) \)
\item If an underobject \( A \tto X \) is classified by a global element \( a \colon \tm \to \pow X \), then \( F A \tto F X \) is classified by \( F a \colon F \tm \to F X \).
\end{itemize}

In fact, almost all of our applications of \jref{propn:logical-fun-pres-satisf} will be to the following particular logical functors:
\begin{propn}[{\cite[Theorem~1.42]{johnstone-topos-theory}}]\label{propn:pb-fun-is-logical}
  If \( \bC \) is an elementary topos, then for any morphism \( f \colon U \to V \) in \( \bC \), the pullback anafunctor \( f^* \colon \bC / V \to \bC / U \) is logical.
  In particular, for any \( U \in \bC \), the product anafunctor \( (U \times -) \colon \bC \to \bC / U \) is logical.
\end{propn}
\noindent
We will spell out the relationship between the logic in \( \bC \) and \( \bC / U \) further in \jref{propn:slice-kj} below.

Sample applications of \jref{propn:pb-fun-is-logical} include:
\begin{exm} \label{exm:pb-fun-logical}
  \enumbelow
  \begin{itemize}
  \item If \( X \) is an hs-object in \( \bC \), then \( (U \times X, \pi_U) \) is an hs-object in \( \bC / U \) with underlying relation \( \id \times i_{\unein_X} \colon (U \times \unein_X, \pi_U) \to (U \times X \times X, \pi_U) \).
  \item If \( \pow X \) is a power object in \( \bC \) with universal relation \( \unin_X \tto X \times \pow X \), then \( (U \times \pow X, \pi_U) \) is a power object of \( (U \times X, \pi_U) \) in \( \bC / U \) with universal relation \( \id \times i_{\unin_X} \colon (U \times \unin_X, \pi_U) \to (U \times X \times \pow X, \pi_U) \).
    Here, the codomain is a product of \( (U \times X, \pi_U) \) and \( (U \times \pow X, \pi_U) \) via
    \[
      (U \times X, \pi_U) \xot{\id_U \times \pi_X} (U \times X \times \pow X, \pi_U) \tox{\id_U \times \pi_{\pow X}} (U \times \pow X, \pi_U).
    \]
  \item Given an \( \unein \)-object \( X \) in \( \bC \), the closure morphism \( \dwn_{(U \times X, \pi_U)} \) for \( (U \times X, \pi_U) \) is given by \( \id_U \times \dwn \colon (U \times X, \pi_U) \to (U \times \pow X, \pi_U) \).
  \end{itemize}
\end{exm}

As an application of \jrefs{propn:logical-fun-pres-satisf}{propn:pb-fun-is-logical}, we now describe a useful general theorem relating the KJ-semantics in a topos \( \bC \) and one of its slices \( \bC / U \).
(See also \cite{shulman-stack-semantics} for an extensive study of KJ semantics in slice topoi.)

Fix a language \( \cL \) and an interpretation \( M \colon \cL \to \bC \) in a topos \( \bC \), and fix an object \( U \in \bC \).
Then, given specifications for the logical anafunctor \( (U \times -) \colon \bC \to \bC / U \) at each of the objects \( M(A) \) and \( M(\vA) \) occurring in the interpretation \( M \), we have the composite interpretation \( \cL \to \bC \tox{(U \times -, \pi_U)} \bC / U \), which we denote by \( \wt M \).
Let \( (\phi, \vx) \) be a formula-in-context over \( \cL \) with \( \vx \tolon \vA \in \Ob (\cL^\up)^n \), and fix \( U' \to U \) and morphisms \( a_i \colon U' \to M A_i \) for \( 0 \le i < n \).
We then have induced morphisms \( b_i = \br{\delta_{U'}, a_i} \colon (U', \delta_{U'}) \to (U \times M A_i, \pi_U) \) in \( \bC / U \).
\begin{propn} \label{propn:slice-kj}
  With \( M \), \( \phi \), \( \va \), and \( \vb \) as above, we have
  \[
    U' \Vdash^M_{\vx \mapsto \va} \phi
    \quad
    \IFF
    \quad
    (U', \delta_{U'}) \Vdash^{\wt M}_{\vx \mapsto \vb} \phi.
  \]
\end{propn}

\begin{proof}
  By \jrefs{propn:logical-fun-pres-satisf}{propn:pb-fun-is-logical}, the interpretation \( \wt M_\vx (\phi) \tto \wt M(\vA) \) is given by
  \begin{equation} \label{eq:slice-kj-subob}
    \id \times i_{M_\vx(\phi)} \colon (U \times M_\vx(\phi), \pi_U) \to (U \times M(\vA), \pi_U).
\end{equation}
  Hence, by definition it follows that \( (U', \delta_{U'}) \Vdash^{\wt M}_{\vx \mapsto \vb} (\phi) \) if and only if
  \[
    \br{b_i}_i = \br{\delta_{U'}, \br{a_i}_i} \colon
    (U', \delta_{U'}) \to (U \times M(\vA), \pi_U)
  \]
  factors through \eqref{eq:slice-kj-subob}, which is equivalent to \( \br{a_i}_{i=0}^{n-1} \colon U' \to M(\vA) \) factoring through \( M_\vx(\phi) \tto M(\vA) \), and hence to \( U' \Vdash^M_{\vx \mapsto \va} \phi \), as desired.
\end{proof}

\jref{propn:slice-kj} allows us to relate properties of generalized elements in \( \bC \) and in the slice \( \bC / U \).
When applying it---as is the case with \jrefs{propn:logical-fun-pres-satisf}{propn:pb-fun-is-logical}---we will usually leave it to the reader to see that the property in question can indeed by described in terms of the satisfaction of an appropriate formula \( \phi \).
In most cases, we will apply this proposition in the rather trivial case where \( \phi \) is an atomic formula; thus, for example, given an \( \unein \)-object \( X \in \bC \) and morphisms \( x,y \colon U \to X \), we have the induced \( \unein \)-object \( \und{U \times X} \in \bC / U \) and global elements \( \br{\id_U, x}, \br{\id_U, y} \colon \und U \to \und{U \times X} \), and \jref{propn:slice-kj} tells us that \( x \ein y \) iff \( \br{\id_U, x} \ein \br{\id_U, y} \).

\clearpage

\etocdepthtag.toc{proofs}
\section{Proofs} \label{sec:proofs}
This appendix contains nearly all of the proofs in this paper.
The shorter proofs, for example of the various basic structural facts about \( \unein \)-objects, are collected together in subsections here according to the sections in which they appear in the paper.
The more substantial proofs, which may involve further subsidiary lemmas and definitions, are given a subsection of their own.

Throughout the proofs, we make fairly free use of our notation and conventions, as laid out in \sref{subsec:cat-prelims} and \aref{sec:topos-logic-background}, so the reader is advised to review those, especially when coming across unfamiliar notation.
In particular, the special notation related to the KJ semantics and generalized elements, introduced in \jref{rmk:kj-bookkeeping}, is used ubiquitously, and so familiarity with it is fairly indispensable to understanding the proofs.

On the other hand, as explained in \sref{subsec:sets-and-topoi}, most of the more complicated proofs are accompanied by the analogous set-theoretic arguments (set apart in coloured boxes), from which the main ideas of the proofs can be gleaned.
\vspace{1cm}

{
\renewcommand\hyperref[2][]{#2}
\etocsettagdepth{proofs}{subsection}
\etocsettocstyle{\textbf{Contents of Appendix~\thesection.}}{}
\localtableofcontents
}

\subsection{Proofs for \alttext{\hyperref[subsec:ein-sets]{Section~\ref*{subsec:ein-sets}}}{Section~\ref*{subsec:ein-sets}}: facts about hereditary sets} \label{subsec:ein-sets-proofs}
Here, we give the proofs of those claims in \sref{subsec:ein-sets} which will not be proven (and, in fact, generalized to an arbitrary topos \( \bC \) rather than just \( \bC = \Set \)) later on in the paper.
These all take place in the set-theoretic context and thus have no bearing on the main (i.e., topos-theoretic) results of the paper.
Accordingly, all of the proofs in this section are in coloured boxes (see \sref{subsec:sets-and-topoi}).

We note that in these proofs, we use some of the results from later in the paper.
There is no risk of circularity since, as we just said, none of those results depend on these.

\begin{settybox}
  \begin{proof}[Proof of \jref{lem:v-and-tc-facts}] \label{proof:v-and-tc-facts}
    We must prove the five claims from the proposition.

    For \ref{item:v-and-tc-facts-h-if-h}, we must show that a set is hereditary if and only if all of its elements are hereditary.
    In one direction, if \( x \subset \rV \) is a set of hereditary sets, then by definition of \( \rV \), we have \( x \in V \), i.e., \( x \) is hereditary.
    For the other direction, it suffices, by the definition of \( \rV \), to show that the class \( P = \set{ x \in \rV \mid x \subset \rV } \subset \rV \) is closed under set-formation.
    But if \( x \subset P \) is a set, then \( x \subset \rV \) (since \( P \subset \rV \)), hence \( x \in P \) by definition of \( P \).

    For \ref{item:v-and-tc-facts-dwn-is-tc}, we must show that if \( B \subset \rV \) is transitive, then \( A \subset B \) is downward-closed in \( B \) if and only if it is transitive.
    But if \( A \) is transitive, then \( x \in y \in A \) implies \( x \in A \) for all sets \( x \) and \( y \) whatsoever, hence in particular for all \( x,y \in B \).
    And conversely, if \( A \) is downward-closed, and given \( x \) and \( y \) with \( x \in y \in A \), it follows that \( y \in B \) (since \( A \subset B \)), hence that \( x \in B \) (since \( B \) is transitive), and hence that \( x \in A \) by closedness of \( A \).
    The parenthetical ``in particular'' statement follows since \( \tc(A) \), being the least transitive class whatsoever containing \( A \), is in particular the least transitive subclass of \( B \) containing \( A \), and hence, by what we just proved, the least downward-closed subclass of \( B \) containing \( A \).

    For \ref{item:v-and-tc-facts-v-wf}, we must show that the element relation on any transitive class \( A \subset \rV \) is extensional and well-founded; by \ref{item:v-and-tc-facts-dwn-is-tc} and \jref{propn:sub-ewf}, it suffices to prove this for \( \rV \).
    Regarding extensionality, if \( x,y \in \rV \) with \( w \in x \ToT w \in y \) for all \( w \in \rV \), then it follows from \ref{item:v-and-tc-facts-h-if-h} that \( w \in x \ToT w \in y \) for all \( w \) whatsoever, hence that \( x = y \) by the axiom of extensionality.
    As for well-foundedness, let \( P \subset \rV \) be inductive; we must show \( P = \rV \).
    By the definition of \( \rV \), it suffices to show that \( P \) is closed under set formation.
    Thus, let \( x \subset P \) be a set; we must show \( x \in P \).
    Since \( x \subset P \subset \rV \), we have \( x \in \rV \) by definition of \( \rV \).
    Moreover, for any \( w \in \rV \) with \( w \in x \), we have \( w \in P \) since \( x \subset P \), hence \( x \in P \) by inductivity of \( P \).

    For \ref{item:v-and-tc-facts-tc-is-set}, we must show that the transitive closure of a hereditary set is always a set.
    We show this by induction on \( x \in \rV \); i.e., we show that the subclass \( P \subset \rV \) consisting of those \( x \) with \( \tc(x) \) a set is inductive.
    Thus, fix \( x \in \rV \) and suppose that \( \tc(w) \) is a set for all \( w \in x \).
    Then, by the axioms of replacement and union, \( \bigcup_{w \in x} \tc(w) \) is again a set.
    But now observe that \( \tc(x) = x \cup \bigcup_{w \in x} \tc(w) \)---since the right-hand side is transitive and must be contained in any transitive class containing \( x \)---and is thus a set by pairing and union.

    Finally, for \ref{item:v-and-tc-facts-tc-is-hs}, we must show that \( \tc(\set{x}) \) is an hs-set for any \( x \in \rV \).
    Extensionality and well-foundedness follows from \ref{item:v-and-tc-facts-v-wf}, and the fact that \( x \) is a top element follows from \ref{item:v-and-tc-facts-dwn-is-tc}.
  \end{proof}
\end{settybox}

\begin{settybox}
  \begin{proof}[Proof of \jref{propn:mostowski-collapse}] \label{proof:mostowski-collapse}
    We must show that every hs-set \( X \) is isomorphic to \( \tc(\set{x}) \) for a unique \( x \).

    For uniqueness, we must prove that \( \tc(\set{y}) \cong \tc(\set{z}) \To y = z \) for all \( y, z \in \rV \), which we do by induction on \( y \) (using \jref{lem:v-and-tc-facts}~\ref{item:v-and-tc-facts-v-wf}).
    Thus, fix \( y \in \rV \), and suppose the inductive claim holds for all \( w \in y \), and that \( \tc(\set{y}) = \tc(\set{z}) \) for some \( z \in \rV \); we need to show \( y = z \).
    Given any \( w \in \tc(\set{y}) \), let \( x \in \tc(\set{z}) \) be the image of \( w \) under the given isomorphism \( f \colon \tc(\set{y}) \toi \tc(\set{z}) \).
    Since \( f \) is an isomorphism, it follows easily that \( f\pbig{\dwnfn(w)} = \dwnfn(x) \), hence \( f \) restricts to an isomorphism \( \tc(\set{w}) = \dwnfn(w) \toi \dwnfn(x) = \tc(\set{x}) \) (where the equalities are by \jref{lem:v-and-tc-facts}~\ref{item:v-and-tc-facts-dwn-is-tc}).
    Hence, by the inductive hypothesis, \( w = x \in X \).
    We thus conclude that \( \tc(\set{y}) \subset \tc(\set{z}) \), and the same argument gives \( \tc(\set{z}) \subset \tc(\set{y}) \), as desired.

    The proof of existence uses the principle of well-founded recursion, which is the following statement: fix a well-founded \( \unein \)-set \( (X, \ein_X) \) and a class \( A \).
    For each \( x \in X \), let \( X_x = \set{ w \in X \mid w \ein_X x } \) and consider the class \( \coprod_{x \in X} A^{X_x} \) of functions valued in \( A \) whose domain is the set \( X_x \) for some \( x \in X \).
    Then, for any function \( G \colon \coprod_{x \in X} A^{X_x} \to A \), there exists a unique function \( F \colon X \to A \) such that \( F(x) = G(x, \rstr{f}{X_x}) \) for each \( x \in X \).
    See \cite[Theorem~1.2.9]{nashaat-powell-valued-models} for a proof in intuitionistic set theory.

    Now, returning to our hs-set \( X \), we seek a hereditary set \( x \) and an isomorphism \( f \colon X \to \tc(\set{x}) \).
    By well-founded recursion, there exists a unique function \( F \colon X \to \rV \) valued in the class \( \rV \) of all hereditary sets satisfying \( F(x) = \set{ F(w) \mid w \ein_X x } \) for all \( x \in X \), as is seen by taking \( G \colon \coprod_{x \in X} \rV^{X_x} \to \rV \) to be given by \( G(x, h) = \set{ h(w) \mid w \ein_X x } \) (this is a hereditary set by \jref{lem:v-and-tc-facts}~\ref{item:v-and-tc-facts-h-if-h}).

    Let \( Y \subset \rV \) be the image of \( F \), endowed with the \( \unein \)-set structure \( \unin_Y \), and let \( f = F \colon X \to Y \).
    We claim that \( f \) is initial (see \jref{defn:ein-hs-set-mor}); indeed, given \( x_1 \in X \) and \( y_0 \in Y \) with \( y_0 \in f(x_1) = \set{ f(w) \mid w \ein_X y_0 } \), we have trivially that \( y_0 = f(x_0) \) for some \( x_0 \ein_X x_1 \), as required.
    It thus follows from \jref{propn:init-cons-mono} that \( f \) is injective and conservative, and it is surjective by its definition.
    Being a consistent and conservative bijection of \( \unein \)-sets, it is thus an isomorphism.

    It remains to see that \( Y = \tc\pbig{\set{f(\rt_X)}} \).
    First, obviously \( f(\rt_X) \in Y \), and \( Y \) is transitive, since if \( y \in Y \), then \( y = f\pbig{f\I(y)} = \set{ f(w) \mid w \ein_X f\I(y) } \subset Y \).
    Now let \( Z \in \rV \) be a transitive set with \( f(\rt_X) \in Z \); we need to show \( Y \subset Z \).
    We will show that \( W = \set{ x \in X \mid f(x) \in Z } \) is equal to all of \( X \), for which it suffices to show that \( W \) contains \( \rt_X \) and is downward closed.
    That \( \rt_X \in W \) is obvious, and if \( x \ein_X y \in W \), so that \( f(y) \in Z \), then \( f(x) \in f(y) \) by consistency of \( f \), hence \( f(x) \in Z \) by transitivity of \( Z \), so \( x \in W \), as desired.
  \end{proof}
\end{settybox}

\begin{settybox}
  \begin{proof}[Proof of \jref{propn:init-map-props}~\ref{item:init-map-props-el-el}] \label{proof:init-map-props-el-el}
    We must show, for \( x, y \in \rV \), that \( \tc(\set{x}) \Ein \tc(\set{y}) \) iff \( x \in y \), and that in this case the inclusion \( \tc(\set{x}) \hto \tc(\set{y}) \) is elemental.

    First suppose \( x \in y \).
    Then \( \tc(\set{x}) \), being transitive, is downward-closed in \( \tc(\set{y}) \) by \jref{lem:v-and-tc-facts}~\ref{item:v-and-tc-facts-dwn-is-tc}, and hence the inclusion \( i \colon \tc(\set{x}) \hto \tc(\set{y}) \) is initial by \jref{propn:closed-iff-init}.
    Moreover, \( i(\rt_{\tc(\set{x})}) = x \in y = \rt_{\tc(\set{y})} \) (where the equalities are by \jref{lem:v-and-tc-facts}~\ref{item:v-and-tc-facts-tc-is-hs}), hence \( i \) is an element morphism.

    Conversely, suppose that there exists an elemental map \( f \colon \tc(\set{x}) \to \tc(\set{y}) \).
    Arguing as above, since \( \tc(\set{x}) \) and \( \tc(\set{y}) \) are both transitive, they are both closed in \( \rV \), hence the inclusions \( i_x \colon \tc(\set{x}) \hto \rV \) and \( i_y \colon \tc(\set{y}) \hto \rV \) are initial.
    It then follows from \jrefs{propn:init-2-of-3}{propn:init-uniq} that \( f i_y = i_x \), and hence that \( \tc(\set{x}) \subset \tc(\set{y}) \), with \( f \) the inclusion.
    In particular, since \( f \) was assumed elemental, we have \( x = f(x) = f(\rt_{\tc(\set{x})}) \in \rt_{\tc(\set{y})} = y \), as desired.
  \end{proof}
\end{settybox}

\subsection{Proofs for \alttext{\hyperref[subsec:ein-ob-props]{Section~\ref*{subsec:ein-ob-props}}}{Section~\ref*{subsec:ein-ob-props}}: facts about \alttext{\( \unein \)}{ε}-objects and morphisms} \label{subsec:ein-ob-props-proofs}
As in \sref{subsec:ein-ob-props}, throughout this section, we let \( X \), \( Y \), and \( Z \) denote \( \unein \)-objects in a topos \( \bC \).

\begin{proof}[Proof of \jref{propn:cons-wf-wf}] \label{proof:cons-wf-wf}
  Assuming \( Y \) is well-founded, and given a consistent morphism \( f \colon X \to Y \), we must show that \( X \) is well-founded.

  Let \( P \tto X \) be inductive; we must show \( P \simeq X \).
  Consider the underobject
  \[
    Q = \set{ y \in Y \mid \forall x \tolon X.\ f(x) = y \to P x } \tto Y.
  \]
  It suffices to show \( Q \approx Y \)---i.e., by well-foundedness of \( Y \), that \( Q \) is inductive---since then, for any \( x \colon U \to X \), we have \( Q(x f) \), hence \( P x \) by definition of \( Q \).
  Thus, fix \( y_1 \colon U \to Y \) and suppose \( Q y_0 \) for all \( y_0 \colon U' \to Y \) with \( y_0 \ein y_1' \).
  Further, fix \( x_1 \colon U' \to X \) with \( x_1 f = y_1' \); we must show \( P x_1 \).

  Since \( P \) is inductive, it suffices to show \( P x_0 \) for all \( x_0 \colon U'' \to X \) with \( x_0 \ein x_1'' \).
  Fix such an \( x_0 \).
  We then have \( x_0 f \ein x_1'' f = y_1'' \) by consistency of \( f \).
  Hence, \( Q(x_0 f) \) by our hypothesis on \( y_1 \).
  Hence \( P(x_0) \), by definition of \( Q \), as desired.
\end{proof}

\begin{proof}[Proof of \jref{propn:init-uniq}] \label{proof:init-uniq}
  Assuming \( X \) is well-founded and \( Y \) extensional, we must show that any two initial \( f,g \colon X \to Y \) are equal.
  Let \( P \tto X \) be the equalizer of \( f \) and \( g \).
  By the well-foundedness of \( X \), it suffices to show that \( P \) is inductive.
  Fix \( x \colon U \to X \) and suppose that
  \begin{equation} \label{eq:init-uniq-pf-ih}
    w \ein x'
    \implies
    w f = w g \colon U' \to Y
    \FORALL
    w \colon U' \to X.
  \end{equation}
  We must show that \( x f = x g \).

  By extensionality of \( Y \), it suffices to show that
  \( z \ein x' f \ToT z \ein x' g \) for all \( z \colon U' \to Y \); we just show the \( (\To) \) direction, as the converse is the same argument.
  Given \( z \colon U' \to Y \) with \( z \ein x' f \), we have by initiality of \( f \) a morphism \( U' \oot U'' \tox{w} X \) with
  \( w f = z'' \quad \AND \quad w \ein x'' \).
  Hence, by \eqref{eq:init-uniq-pf-ih} and the consistency of \( g \), we have \( z'' = w f = w g \ein x'' g \), and hence \( z \ein x' g \), as desired, by epicity of \( \delta_{U''} \).
\end{proof}

\begin{proof}[Proof of \jref{propn:init-cons-mono}] \label{proof:init-cons-mono}
  We are given a morphism \( f \colon X \to Y \).
  We need to show that (i) if \( f \) is initial and monic, then it is conservative, and (ii) if \( f \) is initial and \( X \) is extensional and well-founded, then it \( f \) is monic.

  Regarding (i), assuming the hypotheses, and given \( x_1,x_2 \colon U \to X \) with \( x_1 f \ein x_2 f \), there exists by initiality \( U \oot U' \tox{x_0} X \) with \( x_0 \ein x_2' \) and \( x_0 f = x_1' f \).
  Thus, by monicity, \( x_0 = x_1' \), hence \( x_1' \ein x_2' \), hence \( x_1 \ein x_2 \), as desired, by epicity of \( \delta_{U'} \).

  For (ii), we first give the set-theoretic argument:
  \begin{settybox}
    Assuming the hypotheses, we prove by induction on \( x \in X \) that \( f(x) = f(y) \To x = y \) for all \( y \in X \).
    Thus, fix \( x \) and suppose that \( \forall y.\ f(w) = f(y) \To w = y \) for all \( w \ein x \); we want to show that the same holds for \( x \).
    Fixing \( y \) with \( f(x) = f(y) \), we must show \( x = y \).
    By extensionality, it suffices to prove \( w \ein x \ToT w \ein y \) for all \( w \in X \).

    Suppose \( w \ein x \).
    By consistency of \( f \), we have \( f(w) \ein f(x) = f(y) \).
    Hence, by initiality of \( f \), there is \( v \in X \) with \( v \ein y \) and \( f(v) = f(w) \).
    Hence, by the inductive hypothesis applied to \( w \), we conclude \( w = v \ein y \), as desired.

    Similarly, if \( w \ein y \), then \( f(w) \ein f(y) = f(x) \), hence we have \( v \in X \) with \( f(v) = f(w) \) and \( v \ein x \), and hence by the inductive hypothesis applied to \( v \), we have \( w = v \ein x \).
  \end{settybox}
  Consider the underobject
  \[
    P = \set{ w \in X \mid \forall y \colon X.\ f(w) = f(y) \to w = y} \tto X.
  \]
  Thus, for \( w \colon U \to X \), we have \( P w \) iff \( w' = y \) for all \( y \colon U' \to X \) with \( w' f = y f \).
  Hence, it suffices to show that \( P w \) for all \( w \colon U \to X \) and thus, by well-foundedness of \( X \), that \( P \) is inductive.

  Thus, fix \( x \colon U \to X \) and suppose that \( w \ein x' \To P w \) for all \( w \colon U' \to X \).
  We must show \( P x \).
  That is, given \( y \colon U' \to X \) with \( x' f = y f \), we must show \( x' = y \).
  By extensionality of \( X \), it suffices to show that
  \[
    w \ein x'' \iff w \ein y'' \FORALL w \colon U'' \to X.
  \]

  Assuming \( w \ein x'' \) (and hence \( P w \)), we have by consistency of \( f \) that \( w f \ein x'' f = y'' f \).
  By initiality of \( f \), we thus have \( U'' \oot U''' \tox{v} X \) with \( v f = w''' f \) and \( v \ein y''' \).
  Hence, since \( P w \), we have \( w''' = v \ein y''' \) and hence \( w \ein y'' \) by epicity of \( \delta_{U''} \), as desired.

  Similarly, if \( w \ein y'' \), we obtain \( w f \ein y'' f = x'' f \), hence by initiality obtain \( U'' \oot U''' \tox{v} X \) with \( v f = w''' f \) and \( v \ein x''' \); hence, \( P v \), and thus \( w''' = v \ein x''' \) and hence \( w'' \ein x'' \) by epicity of \( \delta_{U''} \).
\end{proof}

\begin{proof}[Proof of \jref{propn:dwn-dwn}] \label{proof:dwn-dwn}
  Recall from \sref{subsec:2-topoi-prelims} our notation related to the image morphism \( \exists_f \).

  Assuming \( X \) is well-founded, and given \( f \colon X \to Y \) initial, we must show that \( \exists_f \dwnfn_Y^{\pow Y} = \dwnfn_X^{\pow X} \exists_f \colon \pow X \to \pow Y \) and that \( f \dwnfn_Y^Y = \dwnfn_X^X \exists_f \colon X \to \pow Y \).
  Since \( \dwnfn_X^X = \sigma_X \dwnfn_X^{\pow X} \) and \( \sigma_X \exists_f = f \sigma_Y \) (as can be checked directly using \jref{propn:topos-extensionality}), the second of these follows from the first.

  For the first claim, here is the set-theoretic argument:
  \begin{settybox}
    Fix \( P \subset X \).
    We first show \( f\pbig{\dwnfn(P)} \subset \dwnfn\pbig{f(P)} \), i.e., that \( f(x) \in \dwnfn\pbig{f(P)} \) for all \( x \in \dwnfn(P) \).
    This is obvious if \( x \in P \), so it suffices to prove that the set of \( x \) for which it is true is downward closed.
    But by the consistency of \( f \), if \( f(x) \in \dwnfn\pbig{f(P)} \) and \( w \ein x \), then \( f(w) \ein f(x) \in \dwnfn\pbig{f(P)} \), hence \( f(w) \in \dwnfn\pbig{f(P)} \) since \( \dwnfn\pbig{f(P)} \) is closed.

    To show \( f\pbig{\dwnfn(P)} \supset \dwnfn\pbig{f(P)} \), it suffices to show that \( f\pbig{\dwnfn(P)} \) contains \( f(P) \) (which is obvious since \( \dwnfn(P) \supset P \)) and that it is closed.
    Thus, suppose that \( y \ein z \in f\pbig{\dwnfn(P)} \).
    Then \( z = f(x) \) for some \( x \in \dwnfn(P) \), hence by initiality of \( f \), \( y = f(w) \) for some \( w \in X \) with \( w \ein x \).
    It follows that \( w \in \dwnfn(P) \) since \( \dwnfn(P) \) is closed, and hence that \( y \in f\pbig{\dwnfn(P)} \), as desired.
  \end{settybox}

  Fix \( p \colon U \to \pow X \).
  We must show \( p \dwn_X \exists_f \le p \exists_f \dwn_Y \) and \( p \dwn_X \exists_f \ge p \exists_f \dwn_Y \).

  The first of these is equivalent to \( p \dwn_X \le p \exists_f \dwn_Y f\I \).
  By the definition of \( \dwn_X \), it suffices to show that \( p \le p \exists_f \dwn_Y f\I \) and that the right-hand side is closed.
  The last inequality is equivalent to \( p \exists_f \le p \exists_f \dwn_Y \), which follows since \( \id_{\pow Y} \le \dwnfn_Y \) by definition of \( \dwn_Y \).
  For the closedness claim, we claim more generally that \( \dwn_Y f\I \colon \pow Y \to \pow X \) is closed.
  Indeed, given \( \pow Y \xot{\delta_U} U \tox{w,x} X \) with \( w \ein x \in \delta_U \dwn_Y f\I \), we have by definition that \( x f \in \delta_U \dwn_Y \), hence by consistency of \( f \) and downward-closedness of \( \dwn_Y \) that \( w f \in \delta_U \dwn_Y \) and hence that \( w \in \delta_U \dwn_Y f\I \), as desired.

  As for \( p \dwn_X \exists_f \ge p \exists_f \dwn_Y \), it suffices by definition of \( \dwn_Y \) that \( p \exists_f \le p \dwn_X \exists_f \) (which is obvious since \( p \le p \dwn_X \) and \( \exists_f \) is monotone) and that the right-hand side is closed.
  Thus, fix \( y,z \colon U' \to Y \) with \( y \ein z \in p' \dwn_X \exists_f \).
  Thus by definition of \( \exists_f \), there is \( U' \oot U'' \tox{x} X \) with \( x f = z'' \) and \( x \in p'' \dwn_X \).
  By initiality of \( f \), we thus have \( U'' \oot U''' \tox{w} X \) with \( w f = y''' \) and \( w \ein x''' \).
  Since \( p''' \dwn_X \) is closed, it follows from \( w \ein x''' \in p''' \dwn_X \) that \( w \ein p''' \dwn_X \), and hence that \( y''' = w f \in p''' \dwn_X \exists_f \), and hence \( y \ein p' \) by the epicity of \( \delta_{U'''} \delta_{U''} \), as desired.
\end{proof}

\begin{proof}[Proof of \jref{propn:init-2-of-3}] \label{proof:init-2-of-3}
  We are given \( X \tox{f} Y \tox{g} Z \), and we must show four things: (i) if \( f \) and \( g \) are initial, then so is \( fg \); (ii) if \( g \) is consistent and monic and \( f g \) is initial, then so is \( f \); (iii) if \( g \) is conservative and \( f g \) is consistent, then so is \( f \); and (iv) if \( g \) is consistent and \( f g \) is conservative, then so is \( f \).

  Regarding (i): assume \( f \) and \( g \) are initial, and fix \( x_1 \colon U \to X \) and \( z_0 \colon U \to Z \) with \( z_0 \ein x_1 f g \).
  By initiality of \( g \) and \( f \), we have \( U \oot U' \tox{y_0} Y \) with \( y_0 \ein x_1' f \) and \( y_0 g = z_0' \), and \( U' \tox{x_0} X \) with \( x_0 \ein x_1' \) and \( x_0 f = y_0 \).
  Hence \( x_0 f g = y_0 g = z_0' \) and \( x_0 \ein x_1' \), as desired.

  Regarding (ii): suppose \( g \) is consistent and monic, and that \( f g \) is initial, and take \( x_1 \colon U \to X \) and \( y_0 \colon U \to Y \) with \( y_0 \ein x_1 f \).
  By consistency of \( g \), we have \( y_0 g \ein x_1 f g \).
  By initiality of \( f g \), we then have \( U \oot U' \tox{x_0} \) with \( x_0 \ein x_1' \) and \( x_0 f g = y_0' g \).
  By monicity of \( g \), this then gives \( x_0 f = y_0' \) and \( x_0 \ein x_1' \), as desired.

  Regarding (iii): suppose \( g \) is conservative and \( f g \) consistent, and take \( w,x \colon U \to X \) with \( w \ein x \).
  By consistency of \( f g \), we have \( w f g \ein x f g \).
  Hence, by conservativity of \( g \), \( w f \ein x f \), as desired.

  Finally, regarding (iv): suppose \( g \) is consistent and \( f g \) is conservative, and take \( w,x \colon U \to X \) with \( w f \ein x f \).
  By consistency of \( g \), we have \( w f g \ein x f g \), hence by conservativity of \( f g \), we have \( w \ein x \), as desired.
\end{proof}

\begin{proof}[Proof of \jref{propn:dwn-opts}] \label{proof:dwn-opts}
  Fix \( z \colon U \to X \), and let \( p \colon U \to \pow X \) be the morphism (given by \jref{propn:topos-comprehension}) such that
  \[
    U \Vdash_{p} \forall x \tolon X.\ x \in p \tot (x = z \vee x \prec z).
  \]
  We must show that \( p = z \dwn \), i.e., that \( p \) is least closed containing \( z. \)
  The set-theoretic argument:
  \begin{settybox}
    Fix \( z \in X \); we prove that \( P = \set{ x \in X \mid x = z \vee x \prec z } \) is least closed containing \( z \).
    First, it is clear that it contains \( z \).
    Second, it is closed: if \( v \ein x \in P \), then either \( x = z \), in which case \( v \prec x = z \) so \( v \in P \), or \( x \prec z \), hence \( v \prec x \prec z \) and hence \( v \prec z \), so \( v \in P \).

    Finally, \( P \) is least closed containing \( z \): suppose \( Q \subset X \) is closed and contains \( z \), and take \( x \in P \).
    Then either \( x = z \), hence \( x \in Q \) by assumption, or \( x \prec z \), i.e., \( x \preceq y \ein z \) for some \( y \in X \); thus \( y \in Q \) since \( Q \) is closed, and hence \( \dwnfn(y) \subset Q \), since \( \dwnfn(y) \) is least closed containing \( y \), and so \( x \in Q \), as desired, since \( x \in \dwnfn(y) \) by assumption.
  \end{settybox}
  First, it is clear that \( z \in p \).

  Second, \( p \) is closed: given \( v,x \colon U' \to X \) with \( v \ein x \in p' \), we have a cover \( U'_0,U'_1 \to U' \) with \( x'_0 = z'_0 \) and \( x'_1 \prec z_1' \).
  To show \( v \in p' \), it suffices to show that \( v_0' \in p_0' \) and \( v_1' \in p_1' \).
  We have \( v_0' \prec x_0' = z_0' \), hence \( v_0' \in p_0' \), and \( v_1' \prec x_1' \prec z_1' \), and hence \( v_1' \prec z_1' \), so \( v_1' \in p_1' \).

  Finally, given \( q \colon U' \to \pow X \) closed containing \( z \), we need \( p' \le q \).
  Fix \( x \colon U'' \to X \) with \( x \in p'' \).
  We thus have a cover \( U_0'', U_1'' \to U'' \) and \( y \colon U_1'' \to X \) with \( x_0'' = z_0'' \) and \( x_1'' \preceq y \ein z_1'' \).
  To show \( x \in q'' \), it suffices to show that \( x_0'' \in q_0'' \) and \( x_1'' \in q_1'' \).

  We have \( x_0'' \ein q_0'' \) since \( z_0'' \in q_0'' \) by assumption.
  And we have \( y \in q_1'' \) since \( y \ein z_1'' \in q_1'' \) and \( q \) is closed, and hence \( y \dwn \le q_1'' \) since \( y \dwn \) is least closed containing \( y \), and hence \( x_1'' \in q_1'' \) since \( x_1'' \in y \dwn \).
\end{proof}

\begin{proof}[Proof of \jref{propn:no-loops}] \label{proof:no-loops}
  Assuming \( X \) is well-founded, we need to show \( \bC \vDash \forall x \tolon X.\ \neg (x \prec x) \).
  The set-theoretic argument:
  \begin{settybox}
    We prove by induction on \( x \in X \) that \( x \not\prec x \).
    Thus, fix \( x \in X \) with \( w \not\preceq w \) for all \( w \ein x \), and suppose for the sake of contradiction that \( x \prec x \).
    Thus, by assumption, there is some \( w \) with \( x \preceq w \ein x \).
    By \jref{propn:dwn-opts}, either \( x = w \), in which case \( w = x \prec x = w \), or \( x \prec w \), in which case \( w \prec x \prec w \); hence in either case, \( w \prec w \), contradicting the induction hypothesis, since \( w \ein x \).
  \end{settybox}
  We consider the underobject \( P = \set{x \in X \mid \forall x.\ \neg (x \prec x)} \tto X \).
  We need to show that \( P \) is inductive.
  Fix \( x \colon U \to X \) and suppose
  \begin{equation} \label{eq:no-loops-pf-ih}
    w \prec w \implies U' \text{ is initial} \FORALL w \colon U' \to X \WITH w \ein x',
  \end{equation}
  and that \( x \prec x \).
  We need to show that \( U \) is initial.

  By definition, we have \( U \oot U' \tox{w} X \) with \( x' \preceq w \ein x' \).
  By \jref{propn:dwn-opts}, we then have a cover \( U_0',U_1' \to U' \) with \( x_0' = w_0' \) and \( x_1' \prec w_1' \).
  We will show that \( U_0',U_1' \) are initial, hence that \( U \) is initial.

  Indeed, we have \( w_0' = x_0' \prec x_0' = w_0' \), and hence that \( U_0' \) is initial by \eqref{eq:no-loops-pf-ih} since \( w_0' \ein x_0' \).
  And similarly \( w_1' \prec x_1' \prec w_1' \), hence that \( U_1' \) is initial by \eqref{eq:no-loops-pf-ih} since \( w_1' \ein x_1' \).

  Let us address the last two statements in \jref{propn:no-loops}.
  That \( U \vDash \neg (x \ein x) \) for all \( x \colon U \to X \) is obvious in light of what we just proved, since \( x \ein x \) implies \( x \prec x \).

  Finally, let us see that \( x \preceq y \preceq x \) implies \( x = y \) for \( x,y \colon U \to X \).
  There is a cover by four morphisms \( U' \to U \) such that for each one, we have either \( x' = y' \) or \( x' \prec y' \), and either \( y' = x' \) or \( y ' \prec x' \).
  To show that \( x = y \), it suffices to show that \( x' = y' \) for each \( U' \); i.e., we have reduced to checking the original claim \( x = y \) under each of these four assumptions.
  Of course there is nothing to do under the assumption \( x = y \), so we only need to consider the case \( x \prec y \prec x \), in which case \( U \) is initial by what we just showed, and hence \( x = y \) vacuously, as desired.
\end{proof}

\begin{proof}[Proof of \jref{propn:closed-iff-init}] \label{proof:closed-iff-init}
  We must show that a substructure \( A \tto X \), is closed iff \( i_A \) is initial.

  Suppose \( A \) is closed, and fix \( \bar x_1 \colon U \to X \) and \( x_0 \colon U \to X \) with \( x_0 \ein \bar x_1 i_A \).
  Then \( A x_0 \) by closedness of \( A \), hence \( x_0 \) factors through some \( \bar x_0 \colon U \to A \), and it follows that \( \bar x_0 \ein \bar x_1 \) by conservativity of \( i_A \); note that this proves initiality in the strong sense indicated in the proposition.

  Conversely, if \( i_A \) is initial, then given \( x,y \colon U \to X \) with \( x \ein y \) and \( A y \), we have \( \bar y \colon U \to A \) with \( \bar y i_A = y \), hence by initiality of \( i_A \) and \jref{cor:ex-mono-uniq}, there is \( \bar x \colon U \to A \) with \( x i_A = x \), i.e., \( A x \).
\end{proof}

\begin{proof}[Proof of \jref{propn:sub-clo-induc}] \label{proof:sub-clo-induc}
  Given \( A \tto X \) closed, we must show that if \( P \tto X \) is closed or inductive, then so is \( i_A^* P \tto A \).

  Suppose \( P \tto X \) is closed, and fix \( \bar x, \bar y \colon U \to A \) with \( \bar x \ein \bar y \) and \( (i_A^* P) \bar y \).
  Then \( \bar x i_A \ein \bar y i_A \) by consistency of \( i_A \) and \( P (\bar y i_A) \), hence \( P (\bar x i_A) \) by closedness of \( P \), hence \( (i_A^* P) \bar x \), as desired.

  Next, suppose \( P \tto X \) is inductive, and fix \( \bar y \colon U \to A \) with \( (i_A^* P) \bar x \) for all \( \bar x \colon U' \to A \) with \( \bar x \ein \bar y' \).
  Then any \( x \colon U' \to X \) with \( x \ein \bar y' i_A \) factors through \( \bar x \colon U' \to A \) with \( \bar x \ein \bar y' \) by closedness of \( A \), hence \( (i_A^* P) \bar x \) by assumption, and thus \( P x \).
  Hence \( P (\bar y i_A) \) by inductivity of \( P \), so \( (i_A^* P) \bar y \), as desired.
\end{proof}

\begin{proof}[Proof of \jref{propn:sub-ewf}] \label{proof:sub-ewf}
  Given a closed substructure \( A \tto X \) with \( X \) extensional, we must show that \( A \) is extensional.
  The set-theoretic argument:
  \begin{settybox}
    Suppose \( X \) is extensional, and fix \( x,y \in A \subset X \) with \( w \ein x \ToT w \ein y \) for all \( w \in A \); we need to show \( x = y \).
    (Note here that we do not need to distinguish between \( \unein_A \) and \( \unein_X \) since \( A \) is a substructure.)
    By extensionality of \( X \), it suffices to show that \( w \ein x \ToT w \ein y \) for all \( w \in X \).
    But if \( w \ein x \) with \( w \in X \), then \( w \in A \) since \( A \) is closed, hence \( w \ein y \) by assumption, and vice versa.
  \end{settybox}

  Suppose \( X \) is extensional, and fix \( \bar x, \bar y \colon U \to A \) with \( w \ein \bar x' \ToT w \ein \bar y' \) for all \( w \colon U \to A \).
  Set \( x = \bar x i_A \) and \( y = \bar y i_A \).
  To show \( \bar x = \bar y \), or equivalently \( x = y \), it suffices to show \( w \ein x' \ToT w \ein y' \) for any \( w \colon U' \to X \).
  But any \( w \colon U' \to X \) with \( w \ein x' \) factors through \( \bar w \colon U' \to A \) with \( \bar w \ein \bar x' \) by \jref{propn:closed-iff-init}, hence \( \bar w \ein \bar y' \) by assumption and \( w \ein y' \) by consistency of \( i_A \).
  And \( w \ein y' \To w \ein x' \) by the same argument.
\end{proof}

\begin{proof}[Proof of \jref{propn:closure-in-slice}] \label{proof:closure-in-slice}
  Given \( x \colon U \to X \), we need to show the agreement of the underobject \( \Dwn_x \tto U \times X \) and the underobject \( \und{\Dwn_{\br{\id, x}}} \tto \und{U \times X} \) in \( \bC / U \) coming from the global element \( \br{\id, x} \colon \und U \to \und{U \times X} \).
  (Note that these underobjects being equivalent in \( \bC / U \) is equivalent to their being equivalent in \( \bC \).)

  By \jref{propn:pb-fun-is-logical}, the descendent relation \( {\succeq}_{\und{U \times X}} \) is given by \( \id_U \times i_{\succeq_X} \colon (U \times {\succeq}_X, \pi_U) \to (U \times X \times X, \pi_U) \), and the product morphism \( \br{\id_U, x} \times \id_{U \times X} \colon \und{U} \times \und{U \times X} \to \und{U \times X} \times \und{U \times X} \) is given by \( \br{\pi_U, \pi_U x, \pi_X} \colon (U \times X, \pi_U) \to (U \times X \times X, \pi_U) \).
  Since pullbacks are preserved by the forgetful functor \( \bC / U \to \bC \), it thus follows from \jref{lem:dwn-is-descendent} that the following left-hand square is a pullback.
  \[
    \begin{tikzcd}
      \Dwn_{\br{\id_U, x}} \pb \ar[r] \ar[d, "i_{\Dwn_{\br{\id_U, x}}}"'] &[35pt] U \times {\succeq}_X \pb \ar[d, "\id_U \times i_{\succeq_X}"'] \ar[r, "\pi_{\succeq_X}"] &[20pt] {\succeq}_X \ar[d, "i_{\succeq_X}"] \\
      U \times X \ar[r, "\br{\pi_U, \pi_U x, \pi_X}"] & U \times X \times X \ar[r, "\br{\pi_X, \pi_{\pow X}}"] & X \times \pow_X.
    \end{tikzcd}
  \]
  The square on the right, and hence the outer square, is also a pullback, and hence, by \jref{lem:dwn-is-descendent} again, \( \Dwn_{\br{\id_U, x}} \approx \Dwn_x \tto U \times X \) as desired.
\end{proof}

\begin{proof}[Proof of \jref{propn:dwn-hs}] \label{proof:dwn-hs}
  There are three claims; it suffices to prove the more general ``primed'' versions \ref{item:dwn-hs-first-gen}-\ref{item:dwn-hs-last-gen}; though in each case, we will reduce it to the unprimed version.
  First, the set-theoretic argument:
  \begin{settybox}
    For \ref{item:dwn-hs-top-elem}, we must show that any \( x \in X \) is a top element in \( \dwnfn(x) \subset X \); but since \( \dwnfn(x) \) is closed, any closed \( P \subset \dwnfn(x) \) containing \( x \) is closed in \( X \) by \jrefs{propn:closed-iff-init}{propn:init-2-of-3}, hence is equal to all of \( \dwnfn(x) \) by definition.

    For \ref{item:dwn-hs-closed-is-dwn}, we must show that any closed \( A \tto X \) with a top element \( x \) is equal to \( \dwnfn(x) \); but for any closed \( P \subset X \) containing \( x \), \( P \cap A \subset A \) is also closed by \jref{propn:sub-clo-induc}, hence contains \( A \) since \( x \) is a top element.

    For \ref{item:dwn-hs-only-elems}, we must show that \( \dwnfn(x) \subset X \) is elemental iff \( x \ein \rt_X \).
    We already know that \( i \colon \dwnfn(x) \hto X \) is initial by \jref{propn:closed-iff-init}, hence the claim is just that \( i(\rt_{\dwnfn(x)}) \ein \rt_X \) iff \( x \ein \rt_X \), which is true, since we just proved that \( \rt_{\dwnfn(x)} = x \).
  \end{settybox}

  For \ref{item:dwn-hs-top-elem-gen}, given \( x \colon U \to X \), we must show that \( \bar x \colon (U, \id) \to \und{\Dwn_x} \) given by \( \bar x i_{\Dwn_x} = \br{\id, x} \) is a top element.
  Since \( \und{\Dwn_x} \) is, by \jref{propn:closure-in-slice}, the downward closure of \( \br{\id, x} \), this immediately reduces, by working in the slice \( \bC / U \), to \ref{item:dwn-hs-top-elem}, i.e., that \( \bar x \colon \tm \to \Dwn_x \) is a top element of \( \Dwn_x \) for any \( x \colon \tm \to X \), where \( \bar x i_{\Dwn_x} = x \).

  Thus, we must show \( P \approx \Dwn_x \) for any closed \( P \tto \Dwn_x \) containing \( \bar x \); this is an equivalence of underobjects of \( \Dwn_X \) or, equivalently, of underobjects of \( X \).
  But since \( P \tto \Dwn_X \) is closed, it follows from \jrefs{propn:closed-iff-init}{propn:init-2-of-3} that \( i_P i_{\Dwn_x} \colon P \tto X \) is as well; and since \( P \bar x \), it follows that \( P x \); hence \( \Dwn_x \le P \) by definition of \( \Dwn_x \).

  For \ref{item:dwn-hs-closed-is-dwn-gen}, we must show for any closed substructure \( \und A \tto \und{U \times X} \) with a top element \( \bar x \colon \und U \to \und A \) that \( \und A \approx \und{\Dwn_{x}} \) where \( \br{\id_U, x} = \bar x i_A \).
  Again, by passing to the slice \( \bC / U \), this immediately reduces to \ref{item:dwn-hs-closed-is-dwn}, i.e., that \( A \approx \Dwn_{\bar x i_A} \) for \( A \tto X \) with a top element \( \bar x \colon \tm \to A \).

  Thus, given \( P \tto X \) closed containing \( \bar x i_A \), we must show \( A \le P \).
  But \( i_A^* P \tto X \) is closed by \jref{propn:sub-clo-induc} and contains \( \bar x i_A \); hence \( i_A^* P \approx A \), and hence \( A \le P \), as desired.

  For \ref{item:dwn-hs-only-elems-gen}, if \( X \) is an hs-object and given \( x \colon U \to X \), we must show that \( \und{\Dwn_x} \tto \und{U \times X} \) is elemental iff \( x \ein \unex \rt_X \).
  Once again, now using that \( x \ein \unex \rt_X \) iff \( \br{\id_U, x} \ein \br{\id_U, \unex \rt_X} = \rt_{U \times X} \) by \jref{propn:slice-kj}, we reduce to \ref{item:dwn-hs-only-elems}, i.e., that \( \Dwn_x \tto X \) is elemental iff \( x \ein \rt_X \) for \( x \colon \tm \to X \).

  But we know that \( i_{\Dwn_x} \) is initial by \jref{propn:closed-iff-init}, so the claimed equivalence reduces to the claim that \( \rt_X i_{\Dwn_x} \ein \rt_X \) iff \( x \ein \rt_X \).
  This follows since \( x = \rt_X i_{\Dwn_X} \) by \jref{item:dwn-hs-top-elem}, which we already proved.
\end{proof}

\begin{proof}[Proof of \jref{propn:dwn-respects-ein}] \label{proof:dwn-respects-ein}
  Given \( x,y \colon U \to X \), we must show that \( x \ein y \) and implies \( \und{\Dwn_x} \Ein \und{\Dwn_y} \), and that the converse holds if \( X \) is extensional and well-founded.
  Since \( X \) being extensional and well-founded implies the same for \( \und{U \times X} \) by \jref{propn:pb-fun-is-logical}, since \( \und{\Dwn_x} \) and \( \und{\Dwn_y} \) are by \jref{propn:closure-in-slice} the downward closures of \( \br{\id_U,x}, \br{\id_U,y} \colon (U, \id) \to \und{U \times X} \), and since \( \br{\id_U, x} \ein \br{\id_U, y} \) iff \( x \ein y \) by \jref{propn:slice-kj}, we may reduce to the case in which \( U \) is terminal.

  Now, assuming \( x \ein y \), we have \( \Dwn_x \le \Dwn_y \) since \( \Dwn_y \) is closed and (hence) contains \( x \).
  It follows from \jrefs{propn:closed-iff-init}{propn:init-2-of-3} that the morphism \( f \colon \Dwn_x \to \Dwn_y \) with \( f i_{\Dwn_y} = i_{\Dwn_x} \) is initial, and since \( \rt_{\Dwn_x} f i_{\Dwn_y} = \rt_{\Dwn_x} i_{\Dwn_y} = x \ein y = \rt_{\Dwn_y} i_{\Dwn_y} \) by \jref{propn:dwn-hs}, it follows that \( \rt_{\Dwn_x} f \ein \rt_{\Dwn_y} \) by conservativity of \( i_{\Dwn_y} \).

  Conversely, if \( X \) is extensional and well-founded, then so are \( \Dwn_x \) and \( \Dwn_y \) by \jref{propn:sub-ewf}, hence if \( f \colon \Dwn_x \to \Dwn_y \) is elemental, then \( f i_{\Dwn_y} = i_{\Dwn_x} \) by \jref{propn:init-uniq}.
  Hence, by \jref{propn:dwn-hs} and by the consistency of \( i_{\Dwn_y} \), we have \( x = \rt_X i_{\Dwn_x} = \rt_X f i_{\Dwn_y} \ein \rt_Y i_{\Dwn_y} = y \).
\end{proof}

\subsection{Proofs for \alttext{\hyperref[subsec:top-element]{Section~\ref*{subsec:top-element}}}{Section~\ref*{subsec:top-element}}: existence of top extensions} \label{subsec:top-element-proofs}
As in \sref{subsec:top-element}, throughout this section, we let \( X \) denote an \( \unein \)-object in a topos \( \bC \).

\begin{proof}[Proof of \jref{propn:top-ext-props}] \label{proof:top-ext-props}
  Given an underobject \( S \tto X \), we must produce an \( S \)-top extension \( X_S \) and show that it satisfies \ref{item:top-ext-props-first}-\ref{item:top-ext-props-last}.
  In fact, these properties essentially dictate the construction:

  We take \( \wh X \) to be a coproduct \( \wh X = X + \tm \), and write \( t_{\wh X} \colon X \to \wh X \) and \( i_X \colon X \to \wh X \) for the coprojections.
  We then define \( \unein_{\wh X} \tto \wh X \times \wh X \) as the underobject
  \[
    \unein_{\wh X} =
    \set{ \br{x, y} \in \wh X \times \wh X \mid
      x \mathrel{\wh{\ein}} y \vee (y = t_{\wh X} \wedge \wh S x)
    }
    \tto \wh X \times \wh X
  \]
  where we write \( \wh{\unein} \tto \wh X \times \wh X \) for the underobject \( \unein_X \tto X \times X \tox{i_X \times i_X} \wh X \times \wh X \) and \( \wh S \) for the underobject \( S \tox{i_S} X \tox{i_X} \wh X \).

  It is obvious that \( \wh X \), thus defined, satisfies \ref{item:top-ext-props-coprod} and \ref{item:top-ext-props-el}, and that \( i_X \) is monic and consistent; it remains to see that \( i_X \colon X \to \wh X \).
  Thus, fix \( \bar y \colon U \to X \) and \( x \colon U \to \wh X \) with \( x \ein y \), where \( y = \bar y i_X \); we seek \( U \oot U' \tox{\bar x} X \) with \( \bar x \ein \bar y' \) and \( \bar x i_X = x' \).
  We have a cover \( U_0,U_1 \to U \) with \( x_0 \mathrel{\wh\ein} y_0 \) and \( y_1 = \unex t_{\wh X} \) and \( \wh S x_1 \).

  Now \( \unex t_{\wh X} = y_1 = \bar y_1 i_X \colon U_1 \to X \) factors both through \( i_X \) and through \( t_{\wh X} \), hence \( U_1 \) is initial by disjointness of coproducts in topoi.
  Hence \( \delta_{U_0} \colon U_0 \to U \) is an epi, and since \( x_0 \mathrel{\wh\ein} y_0 \), we have by definition a morphism \( \bar x_0 \colon U_0 \to X \) with \( \bar x_0 \ein \bar y_0 \) and \( \bar x_0 i_X = x_0 \), as desired.

  It remains to see that \( \wh X \) satisfies the universal property of the \( S \)-top extension.
  Fix an \( \unein \)-object \( Y \), a global element \( t_Y \colon \tm \to Y \), and a consistent morphism \( f \colon X \to Y \) with \( i_S f \ein \unex_S t_Y \).
  We must show there is a unique consistent morphism \( \bar f \colon \wh X \to Y \) with \( i_X \bar f = f \) and \( t_{\wh X} \bar f = t_Y \).
  In fact, there is a unique \( \bar f \) with the last two properties by the universal property of the coproduct \( \wh X \), so it remains to see that \( \bar f \) is consistent.

  Thus, fix \( x,y \colon U \to X \) with \( x \ein y \).
  First, if \( y \) factors through \( i_X \), then so does \( x \) since \( i_X \colon X \tto \wh X \) is a closed substructure, and hence \( x f = x \bar f \ein y \bar f = y f \) by consistency of \( f \colon X \to Y \), where \( \bar x i_X = x \).
  Second, if \( y = \unex t_X \), then \( x \) must factor through \( S \tox{i_S} X \tox{i_X} \wh X \), and \( y \bar f = \unex t_Y \) by assumption, and hence \( x \bar f \ein y \bar f \) since \( i_S f \ein \unex t_Y \).
  Finally, in the general case, we have a cover of \( U \) by two morphisms \( U' \to U \) such that one of the two above assumptions holds on each \( U' \), and hence \( x' \bar f \ein y' \bar f \) in each case, and hence \( x \bar f \ein y \bar f \).
\end{proof}

\begin{proof}[Proof of \jref{propn:top-ext-is-hs-struct}] \label{proof:top-ext-is-hs-struct}
  We are assuming \( X \) is well-founded and extensional, and are given a dense underobject \( S \tto X \), and must show that the \( S \)-top extension \( \wh X \) is an hs-object, and that \( t_{\wh X} \) is its top element.
  The set-theoretic argument:
  \begin{settybox}
    We first check that \( t_{\wh X} \in \wh X \) is a top element.
    Let \( P \subset \wh X \) be closed containing \( t_{\wh X} \).
    Then by definition of \( \unein_{\wh X} \) and closedness of \( P \), it follows that \( S \subset P \).
    Now by \jref{propn:sub-clo-induc}, \( X \cap P \subset X \) is closed in \( X \), hence is equal to \( X \) since it contains the dense subset \( S \).
    Thus, \( X \subset P \) and hence, since \( t_{\wh X} \in P \) by assumption, we conclude that \( P = \wh X \) as desired.

    Next, we check well-foundedness.
    Let \( P \subset \wh X \) be inductive.
    Then by \jref{propn:sub-clo-induc} again, \( P \cap X \) is still inductive, hence equal to \( X \) by well-foundedness of \( X \).
    Hence \( X \subset P \).
    In particular, \( P \) contains all \( x \in \wh X \) with \( x \ein t_{\wh X} \) (by definition of \( \unein_{\wh X} \)) and hence also \( t_{\wh X} \) by inductivity of \( P \), and hence \( P = \wh X \), as desired.

    Finally, we check extensionality.
    Fix \( x,y \in \wh X \) and suppose \( w \ein x \) iff \( w \ein y \) for all \( w \in \wh X \); we want \( x = y \).
    If \( x = t_{\wh X} = y \), there is nothing to do, so we may assume at least one of \( x \) and \( y \) is in \( X \).
    If \( x \in X \) and \( y = t_{\wh X} \), then \( x \preceq y \) since \( t_{\wh X} \) is a top element, as we just showed.
    Hence, by \jref{propn:dwn-opts}, either \( x = y \), in which case we are done, or \( x \prec y \), i.e., \( x \preceq w \ein y \) for some \( w \in \wh X \).
    It follows by assumption that \( w \ein x \), hence that \( x \prec x \), contradicting \jref{propn:no-loops}.

    If \( y \in X \) and \( x = t_{\wh X} \), the argument is the same.
    Finally, if \( x,y \in X \), then the claim follows from the extensionality of \( X \) and its closedness in \( \wh X \) (\jref{propn:sub-ewf}).
  \end{settybox}

  First, let us see that \( t_{\wh X} \) is indeed a top element.
  Fix \( P \tto \wh X \) closed containing \( t_{\wh X} \); we need to show \( P \approx \wh X \).
  Since by \jref{propn:top-ext-props}, the underobjects \( \tm \tox{t_{\wh X}} \wh X \xot{i_X} X \) is a cover (see \sref{subsubsec:topoi-prelims}), and since \( (\tm, t_{\wh X}) \le P \) by assumption, it suffices to prove \( X \le P \)---or equivalently, \( i_X^* P \simeq X \tto X \).

  Now, any \( x \colon U \to \wh X \) with \( S x \) satisfies \( x \ein \unex t_{\wh X} \) by the assumption \( i_S i_X \ein \unex_S t_{\wh X} \), and hence \( P x \) since \( P \) is closed.
  Thus \( S \le P \).
  Since \( P \tto \wh X \) is closed and \( X \tto \wh X \) is a closed substructure by \jrefs{propn:closed-iff-init}{propn:top-ext-props}, the intersection \( i_X^* P \tto X \) is closed by \jref{propn:sub-clo-induc}.
  Since \( S \le P \), it follows that \( S \approx i_X^* S \le i_X^* P \tto X \), hence that \( i_X^* P \approx X \), as desired, since \( i_X^* P \) is closed and \( S \tto X \) is dense.

  Next, well-foundedness: fix \( P \tto \wh X \) inductive; we want \( P \approx \wh X \).
  Again, it suffices to prove that \( X \le P \), since we then obtain \( P t_{\wh X} \) by inductivity of \( P \) and by the fact that \( x \ein t_{\wh X} \) implies \( X x \), for \( x \colon U \to \wh X \).
  Thus, it suffices to show that \( i_X^* P \approx X \tto X \); but this follows by well-foundedness of \( X \), since \( i_X^* P \) is inductive by \jref{propn:sub-clo-induc}.

  Finally, extensionality: fix \( x,y \colon U \to \wh X \) with \( w \ein x' \) iff \( w \ein y' \) for all \( w \colon U' \to \wh X \); we want \( x = y \).
  We have a cover of \( U \) by four morphisms \( U' \to U \) such that for each one, we have that either \( x' = \unex t_{\wh X} \) or \( X x' \), and that either \( y' = \unex t_{\wh X} \) or \( X y' \).
  To show that \( x = y \), it suffices to show that \( x' = y' \) for each \( U' \); i.e., we have reduced to checking the original claim \( x = y \) under each of these four assumptions.

  When \( x = \unex t_{\wh X} = y \), there is nothing to do.
  Suppose \( X x \) and \( y = \unex t_{\wh X} \).
  Since by assumption we have \( w \ein x' \) iff \( w \ein y' \) for all \( w \colon U' \to \wh X \), it follows (from the definition of \( \prec \), \ref{defn:strict-ancestor}) that
  \begin{equation} \label{eq:ext-precs-agree}
    w \prec x'
    \IFF
    w \prec y'
    \FORALL
    w \colon U' \to \wh X.
  \end{equation}
  Since \( t_{\wh X} \) is a top element, as we have just proven, we have that \( x \preceq y \) and hence, by \jref{propn:dwn-opts} that there is a cover of \( U \) by two morphisms \( U' \to U \) such that, for each one, either \( x' \prec y' \) or \( x' = y' \).
  In the second case, there is nothing to do.
  In the first case, since \( x' \prec y' \), it follows from \eqref{eq:ext-precs-agree} that \( x' \prec x' \), hence by \jref{propn:no-loops} that \( U' \) is initial, and hence that \( x' = y' \), as desired.

  In the case where \( x = \unex t_{\wh X} \) and \( X y \), the argument is the same.
  Finally, in the case where \( X x \) and \( X y \), we have \( \bar x, \bar y \colon U \to X \) with \( \bar x i_X = x \) and \( \bar y i_X = y \).
  Now, for any \( \bar w \colon U \to X \), by consistency and conservativity of \( i_X \), we have \( \bar w \ein \bar x \) iff \( \bar w \ein \bar y \).
  Hence \( \bar x = \bar y \) by extensionality of \( X \), and hence \( x = y \), as desired.
\end{proof}

\subsection{Proofs for \alttext{\hyperref[subsec:ein-in-slice]{Section~\ref*{subsec:ein-in-slice}}}{Section~\ref*{subsec:ein-in-slice}}: facts about \alttext{\( \unein \)}{ε}-objects in slices} \label{subsec:ein-in-slice-proofs}
As in \sref{subsec:ein-in-slice}, throughout this section, let \( f \colon A \to B \) be a morphism in a topos \( \bC \), and let \( \und X \in \bC / A \) and \( \und Y \in \bC / B \) be \( \ein \)-objects.

\begin{proof}[Proof of \jref{propn:ein-pb-conds}] \label{proof:ein-pb-conds}
  We must show that the three conditions in the proposition are equivalent.

  Conditions \ref{item:ein-pb-conds-pbs} and \ref{item:ein-pb-conds-rel-consv} both include that \( g \) is relatively consistent and that \( g \) exhibits \( X \) as the pullback of \( Y \) along \( f \).
  Thus, it only remains to see that, under these circumstances, \( g \) is relatively conservative if and only if \( \unein_X \) is the pullback of \( \unein_Y \) along \( g \times_f g \colon X \times_A X \to Y \times_B Y \), and this follows immediately from the definition of relative conservativity.

  Next, assuming \ref{item:ein-pb-conds-rel-consv}, to prove \ref{item:ein-pb-conds-univ}, we need to show that given a relatively consistent morphism \( h \colon W \to Y \) from an \( \unein \)-object \( \und W \in \bC / A \), the unique induced morphism \( \bar h \colon \und W \to \und X \) in \( \bC / A \) with \( \bar h g = h \) is consistent.
  We can check this directly: given \( y,z \colon U \to W \) with \( y \ein^A z \), we have that \( y \bar h g = y h \ein^B z h = z \bar h g \) by relative consistency of \( h \), hence \( y \bar h \ein^A z \bar h \) by relative consistency of \( g \) (and since \( y \bar h \delta_X = y \delta_W = z \delta_W = z \bar h \delta_X \) by the assumption \( y \ein^A z \)), as desired.

  Finally, the implication \ref{item:ein-pb-conds-univ}\( \To \)\ref{item:ein-pb-conds-pbs} now follows formally, since we already have the reverse implication: assuming \( \und X \) satisfies \ref{item:ein-pb-conds-pbs}, let \( \und{\wt X} \in \bC / A \) be the image of \( \und Y \in \bC / B \) under the pullback anafunctor \( g^* \colon \bC / B \to \bC / A \), so that \( \und{\wt X} \), with its projection morphism \( \tilde g \colon \wt X \to Y \) satisfies \ref{item:ein-pb-conds-pbs} by definition, and hence also satisfies \ref{item:ein-pb-conds-univ}.
  Since \ref{item:ein-pb-conds-univ} determines the \( \unein \)-object \( \und X \) together with the morphism \( g \colon X \to Y \) uniquely up to isomorphism, it follows that there is an isomorphism \( \und X \to \und{\wt X} \) in \( \hs(\bC / A) \) commuting with \( g \) and \( \tilde g \), and hence that \( (\und X, g) \) satisfies \ref{item:ein-pb-conds-pbs} since \( (\und{\wt X}, \tilde g) \) does.
\end{proof}

\begin{proof}[Proof of \jref{propn:ein-cocart-conds}] \label{proof:ein-cocart-conds}
  As in the previous proof, we must show that the three conditions in the proposition are equivalent.

  Assuming \ref{item:ein-cocart-conds-iso}, to prove \ref{item:ein-cocart-conds-strong}, we just need to show that \( g \) is strongly relatively conservative.
  Thus, given \( x,y \colon U \to X \) with \( x g \ein^B y g \), we need to show \( x \ein^A y \); in other words, we are assuming that \( \br{x g, y g} \colon U \to Y \times Y \) factors through \( \unein_Y \tto Y \times_B Y \tto Y \times Y \), and we want \( \br{x,y} \colon U \to X \times X \) to factor through \( \unein_X \tto X \times_A X \tto X \times X \).
  But this follows immediately from the assumption that \( g \) and the dashed morphism in \eqref{eq:rel-consis-square} on \pref{eq:rel-consis-square} are both isomorphisms.

  Next, assuming \ref{item:ein-cocart-conds-strong}, to prove \ref{item:ein-cocart-conds-univ}, we need to show that, given an \( \ein \)-object \( \und Z \in \bC / B \) and a relatively consistent morphism \( h \colon X \to Z \), the unique morphism \( \und h \colon \und Y \to \und Z \) in \( \bC / B \) with \( g \bar h = h \) (which exists since \( g \) is an isomorphism) is consistent.
  Given \( y,z \colon U \to Y \) with \( y \ein^B z \), we have that the morphisms \( \tilde y, \tilde z \colon U \to X \) with \( \tilde y g = y \) and \( \tilde z g = z \), satisfy \( \tilde y \ein^A \tilde z \) by the strong relative conservativity of \( g \).
  Hence, \( y \und h = \tilde y h \ein^B \tilde z h = z \und h \) by consistency of \( h \), as desired.

  Finally, the implication \ref{item:ein-cocart-conds-univ}\( \To \)\ref{item:ein-cocart-conds-iso} follows from the reverse implication in the same way as in the proof of \jref{propn:ein-pb-conds} above.
\end{proof}

\begin{proof}[Proof of \jref{propn:rel-notions-compare}] \label{proof:rel-notions-compare}
  We are given the various data shown in the diagram in the statement of the proposition, and we must prove the various implications relating the properties of \( g \), \( \bar g \), and \( \und g \).

  We begin with \ref{item:rel-notions-compare-consv}.

  First, if \( \und g \) is conservative (and hence strongly relatively conservative over \( \id_B \)), then since \( c \) is strongly relatively conservative, the strong relative conservativity of \( g \) follows from the quite obvious fact that this property is preserved under composition.
  Conversely, suppose that \( g \) is strongly relatively conservative, and fix \( x,y \colon U \to X_B \) with \( x \und g \ein^B y \und g \).
  Then, considering the morphisms \( \tilde x = x c \I, \tilde y = y c \I \colon U \to X \), since \( \tilde x g = x \und g \ein^B y \und g = \tilde y g \) it follows from strong relative conservativity of \( g \) that \( \tilde x \ein^A \tilde y \), and hence from consistency of \( c \) that \( x \ein^B y \), as desired.

  That strong relative conservativity implies conservativity is obvious

  Next, suppose that \( g \) is relatively conservative; we want to show that \( \bar g \) is conservative.
  Hence, fix \( x,y \colon U \to X \) with \( x \delta_X = y \delta_X \) and \( x \bar g \ein^A y \bar g \).
  Then \( x g = x \bar g \pi_f \ein^B y \bar g \pi_f = y g \) by relative consistency of \( \pi_f \), and hence \( x \ein^A y \) by relative conservativity of \( g \), as desired.
  Conversely, suppose \( \bar g \) is conservative, and fix \( x,y \colon U \to X \) with \( x \delta_X = y \delta_X \) and \( x g \ein^B y g \).
  Then \( x \bar g \pi_f \ein^B y \bar g \pi_f \) and hence \( x \bar g \ein^A y \bar g \) and \( x \ein^A y \) by relative conservativity of \( \pi_f \) and conservativity of \( \bar g \).

  We next check the equivalence of the conditions from \ref{item:rel-notions-compare-init}, namely that \( \bar g \) being initial, \( g \) being relatively initial, and \( \und g \) being initial are all equivalent.
  Suppose \( g \) is relatively initial and fix \( x_1 \colon U \to X \) and \( y_0 \colon U \to f^* Y \) with \( x_1 \delta_X = y_0 \delta_Y \) and \( y_0 \ein^A x_1 \bar g \).
  Then \( y_0 \pi_f \ein^B x_1 \bar g \pi_f = x_1 g \) by relative consistency of \( \pi_f \).
  Hence, by relative initiality of \( g \) there is \( U \oot U' \tox{x_0} X \) with \( x_0 \ein^A x_1' \) and \( x_0 \bar g \pi_f = x_0 g = y_0' \pi_f \), and hence \( x_0 \bar g = y_0' \) by relative consistency of \( \pi_f \).
  Conversely, supposing \( \bar g \) is initial, and fixing \( x_1 \colon U \to X \) as above and \( y_0 \colon U \to Y \) with \( y_0 \ein^B x_1 g \), we obtain by \( \ein \)-cartesianness of \( \pi_f \) a morphism \( \tilde y_0 \colon U \to f^* Y \) with \( \tilde y_0 \pi_f = y_0 \) and \( \tilde y_0 \ein^A x_1 \bar g \).
  Then, by initiality of \( \bar g \), we obtain \( U \oot U' \tox{x_0} X \) with \( x_0 \ein x_1' \) and \( x_0 \bar g = \tilde y_0' \) and hence \( x_0 g = y_0' \), as desired.
  The proof that \( g \) being relatively and \( \und g \) being initial are equivalent is very similar, and we leave it to the reader.

  Finally, regarding \ref{item:rel-notions-compare-elem}, suppose \( \rt_X \colon \und A \to \und X \) and \( \rt_Y \colon \und B \to \und Y \) are top elements.
  Applying the logical anafunctor \( f^* \colon \bC / B \to \bC / A \), we obtain a morphism \( \tilde \rt_Y \colon \und A \to f^* \und Y \) which is again a top element by \jref{propn:pb-fun-is-logical} (hence the unique top element of \( f^* \und Y \) by \jref{cor:dwn-mono}).
  Since we have already proven that \( g \) is relatively initial if and only if \( \bar g \) is initial, and since \( \rt_X \bar g \ein_{f^* Y} \tilde \rt_Y \ToT \rt_X g = \rt_X \bar g \pi_f \ein_{Y} \tilde \rt_Y \pi_f = \rt_Y \) by \( \ein \)-cartesianness of \( \pi_f \), it follows that \( \bar g \) is elemental iff \( g \) is relatively elemental, as desired.
\end{proof}

\begin{proof}[Proof of \jref{lem:loc-ext-char}] \label{proof:loc-ext-char}
  We are given an \( \unein \)-object \( Z \) in a topos \( \bD \) and must show that \( Z \) is locally extensional if and only if for any \( x,y,z \colon U \to Z \) satisfying
  \begin{equation*} \label{eq:loc-ext-char-pf-cond} \tag{\( \star \)}
    x,y \preceq z
    \AND
    w \ein x' \ToT w \ein y'
    \FORALL
    w \colon U' \to Z,
  \end{equation*}
  we have \( x = y \).

  Suppose \( Z \) is locally extensional and fix \( x,y,z \colon U \to Z \) satisfying \eqref{eq:loc-ext-char-pf-cond}; we need to show \( x = y \).
  We thus have global elements \( \tilde x = \br{\id, x}, \tilde y = \br{\id, y}, \tilde z = \br{\id,z} \colon \und U \to \und{U \times Z} \) in \( \bD / U \), and \jref{propn:slice-kj} implies that \( \tilde x, \tilde y \ein \tilde z \), and that \( \tilde w \ein \tilde x' \ToT \tilde w \ein \tilde y' \) for all \( \tilde w \colon \und U' \to \und{U \times Z} \).
  Hence, by definition of \( \und{\Dwn_z} \), \( \tilde x \) and \( \tilde y \) factor through morphisms \( \bar x, \bar y \colon \und U \to \und{\Dwn_z} \), and since \( \und{\Dwn_z} \) is a substructure, we have \( \bar w \ein \bar x' \ToT \bar w \ein \bar y' \) for all \( \bar w \colon \und U' \to \und{\Dwn_z} \).
  By extensionality of \( \und{\Dwn_z} \), it follows that \( \bar x = \bar y \), hence that \( \tilde x = \tilde y \) and \( x = y \), as desired.

  Conversely, assume that \eqref{eq:loc-ext-char-pf-cond} implies \( x = y \) for all \( x,y,z \colon U \to Z \), and fix \( z \colon U \to Z \) and \( \bar x, \bar y \colon \und U' \to \und{\Dwn_z} \) in \( \bD / U \) with \( \bar w \ein \bar x'' \ToT \bar w \ein \bar y'' \) for all \( \bar w \colon \und{U''} \to \und{\Dwn_z} \).
  We need to prove that \( \bar x = \bar y \).

  Set \( \br{\id, x} = \tilde x = \bar x i_{\Dwn_z} \colon \und{U'} \to \und{U \times Z} \) and \( \br{\id, y} = \tilde y = \bar y i_{\Dwn_z} \) and \( \tilde z = \br{\id, z} \).
  We then have \( \tilde x, \tilde y \preceq \tilde z' \) by definition of \( \Dwn_z \), and hence that \( x,y \preceq z' \) by \jref{propn:slice-kj}.
  Next, given any \( w \colon U'' \to Z \), it follows from \jref{propn:slice-kj}, the closedness of the substructure \( \und{\Dwn_z} \to \und{U \times Z} \), and our assumption on \( \bar x \) and \( \bar y \), that \( w \ein x'' \ToT w \ein y'' \) for all \( w \colon U'' \to Z \).
  Hence \( x = y \) by our assumption on \( Z \), and hence \( \tilde x = \tilde y \) and \( \bar x = \bar y \), as desired.
\end{proof}

\begin{proof}[Proof of \jref{propn:sum-props}] \label{proof:sum-props}
  We are given an \( \unein \)-cocartesian morphism \( c \colon X \to X_B \) over \( f \).
  We must show that the properties of \( \und X \) indicated in the proposition are inherited by \( \und X_B \).

  We first make a general observation, that we will use several times: given an underobject \( i_P \colon \und P \tto \und X \), a morphism \( x \colon \und U \to \und X \) in \( \bC / A \) factors through \( \und P \) if and only if \( x \), regarded as a morphism in \( \bC \), factors through \( i_P \colon P \tto X \); thus, we may express this unambiguously by writing \( P x \).
  Note also that, since \( c \) is an iso, we have \( P x \) iff \( P (x c) \), where the ``\( P \)'' in the second expression refers to the underobject \( i_P c \colon P \tto X_B \).

  Now, suppose \( \und X \) is well-founded, and fix an inductive underobject \( i_P \colon \und P \tto \und X_B \) in \( \bC / B \); we need to show that \( \und P \approx \und X_B \), i.e., that \( i_P \) is an iso.
  By the well-foundedness of \( \und X \), it suffices to prove that \( i_P \colon (P, i_P c\I \delta_X) \tto \und X \) is still inductive.

  Thus, fix \( x \colon U \to X \) and suppose that \( P w \) for all \( w \colon U' \to X \) with \( w \ein^A x' \).
  Then for any \( w \colon U' \to Y \) with \( w \ein^B x' c \), it follows from strong relative conservativity of \( c \) that \( w c \I \ein^A x' \), hence that \( P (w c \I) \), hence that \( P w \).
  Since \( \und P \tto \und X_B \) is inductive, it follows that \( P (x c) \), hence that \( P x \), as desired.

  Next, supposing \( \und X \) is extensional, we need to show that \( \und X_B \) is locally extensional.
  We verify the condition from \jref{lem:loc-ext-char}.
  Thus, fix \( x,y,z \colon \und U \to \und X_B \) in \( \bC / B \) with \( x, y \preceq z \), and suppose that \( w \ein x' \) iff \( w \ein y' \) for all \( w \colon \und U' \to \und X_B \); we need to show \( x = y \).

  We first claim that \( w c\I \delta_X = z' c\I \delta_X \colon U' \to A \) for all \( w \colon \und{U'} \to \und X_B \) with \( w \preceq z \).
  In other words, letting \( p \colon \und U \to \pow \und X_B \) be such that \( w \in p' \) iff \( w c\I \delta_X = z' c\I \delta_X \colon \und{U'} \to (A, f) \) for \( w \colon \und{U'} \to \und X_B \) (which we have by \jref{propn:comprehension-morphism}), we claim that \( z \dwn \le p \).
  By definition of \( z \dwn \), it suffices to show that \( z \in p \) (which is clear) and that \( p \) is closed.
  But given \( v,w \colon \und{U'} \to \und X_B \) with \( v \ein w \in p' \), we then have \( v \ein^B w \) (considering \( v \) and \( w \) as morphisms \( U' \to X_B \) in \( \bC \)), hence \( v c\I \ein^A w c\I \) by strong relative conservativity of \( c \), hence that \( v c\I \delta_X = w c \I \delta_X = z c\I \delta_X \), as desired.

  We thus have that \( x c\I \delta_X = y c\I \delta_X = z c\I \delta_X \), so that we have morphisms \( x c\I,y c\I \colon (U, z c\I \delta_X) \to \und X \) in \( \bC / A \).
  Hence, by extensionality of \( \und X \), it suffices to show that \( w \ein^A x' c\I \ToT w \ein^A y' c\I \) for all \( w \colon U' \to X \).
  But this follows immediately from relative consistency and conservativity of \( c \) and our assumption on \( x \) and \( y \).

  Finally, supposing that \( i_S \colon (S, \delta_S) \tto \und X \) is a dense underobject in \( \bC / A \), we must show that \( i_S c \colon (S, \delta_s f) \tto \und X_B \) is a dense underobject in \( \bC / B \).
  Thus, we must show that any closed underobject \( i_P \colon \und P \tto \und X_B \) containing \( (S, \delta_S f) \) is equivalent to \( \und X_B \), i.e., that, \( i_P \) is an iso.
  But for any such underobject, \( i_P c\I \colon (P, i_P c\I \delta_X) \tto \und X \) is an underobject in \( \bC / A \) which contains \( (S, \delta_S) \), and which is closed: given \( x,y \colon U \to X \) with \( x \ein^A y \) and \( P y \), we have \( x c \ein^B y c \) (by relative consistency of \( c \)) and \( P (y c) \), and hence \( P (x c) \) and \( P c \).
  Hence, by density of \( (S, \delta_S) \), the morphism \( i_P c\I \), and hence also \( i_P \), is an iso, as desired.
\end{proof}

\subsection{Proof of \alttext{\hyperref[propn:topos-sketch-cotensors]{Proposition~\ref*{propn:topos-sketch-cotensors}}}{Proposition~\ref*{propn:topos-sketch-cotensors}}: existence of sketch cotensors} \label{subsec:topos-sketch-cotensors-proof}
We need to prove that any internal topos \( E \) in a corepita 2-category \( \cC \) admits a cotensor by any finite \( \Omega \)-topos sketch.
The proof closely parallels that of the corresponding statement \cite[Proposition~A.0.1]{helfer-topoi-in-topoi} concerning finite-limit sketches (and, indeed, generalizes it, except for the fact that here, we consider only \emph{internal topos} objects \( E \), rather than general lex objects, and only consider the \emph{groupoidal} object of models of \( J \) in \( E \)).
As was the case there, we first prove the existence of cotensors \( (E^J)^\iso \) when \( J \) is a ``simple'' sketch containing only one specification.

The following lemma is proven in exactly the same way as \cite[Lemma~A.0.2]{helfer-topoi-in-topoi}:
\begin{lem} \label{lem:one-spec-sketch-cotensor}
  Let \( \Omega \) be a set of topos operations, and suppose \( J \) is a finite \( \Omega \)-topos sketch such that \( \abs{J} = L_{\omega} \) for some \( \omega \in \Omega \), and \( J_\omega = \set{\id_{\abs{J}}} \) and \( J_{\omega'} = \emptyset \) for all \( \omega' \ne \omega \).

  Then there exists a cotensor \( (E^J)^\iso \) for any internal topos \( E \in \cC \) in any corepita 2-category \( \cC \).
\end{lem}

\begin{proof}
  We first consider the following general construction.
  Let \( L \) be a category and \( i_K \colon K \hto L \) a full subcategory.
  We form a new category \( L' \) by adjoining to \( L \) an isomorphic copy of each object in \( \Ob L - \Ob K \).
  Thus, \( \Ob L' = \Ob L \sqcup (\Ob L - \Ob K) \), and the morphisms are given by \( L'\pbig{(i, x), (j, y)} = L(x, y) \) (where \( i,j \in \set{0,1}) \), with the obvious composition.
  We have two evident embeddings \( i_0,i_1 \colon L \to L' \) (with \( i_0(x) = i_1(x) = (0, x) \) for \( x \in K \) and \( i_j(x) = (j, x) \) for \( x \in \Ob L - \Ob K \)) and a natural isomorphism \( \alpha \colon i_0 \toi i_1 \) with \( \alpha_x = \id_x \in L'\pbig{i_0(x), i_1(x)} \) for all \( x \in L \).
  This gives a functor \( \bar \alpha \colon L \times ( \toi ) \to L' \), and we have a strict pushout square:
  \[
    \begin{tikzcd}
      K \times (\toi) \ar[r, "i_K \times \id"] \ar[d, "\pi_K"'] &[15pt] L \times ( \toi ) \ar[d, "\bar \alpha"]
      \\
      K \ar[r, "i_K i_0 = i_K i_1"] & L'. \po
    \end{tikzcd}
  \]
  Now, another way to state the condition \jref{defn:topos-op}~\ref{item:topos-op-ext-uniq} in the definition of a topos operation \( (L, K, \sigma) \) is that any two extensions \( f_0,f_1 \in \sigma_{\bC} \subset \Ob \bC^L \) of a diagram \( f \colon K \to \bC \) can be further uniquely extended to a diagram \( f \colon L' \to \bC \):
  \begin{equation} \label{eq:doubled-diag-uniq-ext}
    \begin{tikzcd}
      L \ar[rd, "f_0"'] \ar[r, "i_0"] & L' \ar[d, "f", dashed]  & L \ar[ld, "f_1"] \ar[l, "i_1"'] \\
      & \bC. &
    \end{tikzcd}
  \end{equation}

  We now proceed exactly as in \loccitnodot.
  Apply the above construction to our given topos operation \( \omega = (L_\omega, K_\omega, \sigma_\omega) \qefed (L, K, \sigma) \) to obtain a finite category \( L' \).
  Given an object \( E \in \cC \) in a corepita 2-category \( \cC \), form the cotensor-core\footnote{%
    This just means the core of the cotensor, though it also has a simple universal property that can be described directly, see \cite[Remark~3.2.2]{helfer-topoi-in-topoi}.%
}
  \( (E^{K})^\iso \), which exists since \( K \) is by assumption of finite presentation (see \sref{subsubsec:corepita-prelims}), so that we have a universal diagram \( \wt \ev_{K, E} \colon K \to \cC\pbig{(E^{K})^\iso, E} \) of shape \( K \) in the topos \( \cC\pbig{(E^{K})^\iso, E} \).
  By \jref{defn:topos-op}~\ref{item:topos-op-ext-ex}, there exists an extension \( g \colon L \to \cC\pbig{(E^K)^\iso, E} \) of \( \wt \ev_{K, E} \) with \( g \in \sigma_{\omega,\cC((E^{K})^\iso, E)} \).
  The diagram \( g \) is represented by some morphism \( p \colon (E^{K})^\iso \to (E^{L})^\iso \).
  Now form \( (E^J)^\iso \) as the pullback of the isofibration \( (E^{i_1})^\iso \colon (E^{L'})^\iso \to (E^{L})^\iso \) along \( p \):
  \[
    \begin{tikzcd}
      (E^J)^\iso \pb \ar[r, "\pi_1"] \ar[d, "\pi_0"'] & (E^{L'})^\iso \ar[d, "E^{i_1}"] \ar[r, "E^{i_0}"] & (E^L)^\iso \\
      (E^K)^\iso \ar[r, "p"] & (E^L)^\iso
    \end{tikzcd}
  \]
  Note that \( (E^J)^\iso \), being a pullback of groupoids, is itself a groupoid.
  The same argument as in \loccitspace (now using \jref{defn:topos-op}~\ref{item:topos-op-ext-invt}) shows that the diagram \( \pi_1 E^{i_0} \) of shape \( L \) in the topos \( \cC\pbig{(E^J)^\iso, E} \) is in \( \sigma_{\omega,\cC((E^J)^\iso, E)} \), and thus produces a model for \( J \) in \( \cC\pbig{(E^J)^\iso, E} \).
  The verification that this satisfies the requisite universal property, making \( (E^J)^\iso \) a cotensor \( E \) by \( J \), now proceeds exactly as in \cite[Lemma~A.0.2]{helfer-topoi-in-topoi} (now using \eqref{eq:doubled-diag-uniq-ext}, and using \jref{defn:topos-op}~\ref{item:topos-op-ext-uniq}~and~\ref{item:topos-op-ext-log-pres} to obtain the full-faithfulness of \( p \)).
\end{proof}

Still following the proof of \cite[Proposition~A.0.1]{helfer-topoi-in-topoi}, we can now finish the proof of \jref{propn:topos-sketch-cotensors} using \jref{lem:one-spec-sketch-cotensor}.
Fix a finite \( \Omega \)-topos sketch \( J \).
  For each \( \omega \in \Omega \), let \( J^\omega \) be the \( \Omega \)-topos sketch with \( \abs{J^\omega} = L_\omega \) as in \jref{lem:one-spec-sketch-cotensor}.

  Now fix an internal topos \( E \) in a groupoidant 2-category \( \cC \).
  By \jref{lem:one-spec-sketch-cotensor}, for each \( \omega \in \Omega \), we have a cotensor \( (X^{J^\omega})^\iso \), with a universal model \( \abs{J^\omega} = L_\omega \to \cC\pbig{(X^{J^\omega})^\iso, X} \), classified by a morphism \( p_\omega \colon (X^{J^\omega})^\iso \to (X^{L_\omega})^\iso \).
  We also have a cotensor-core \( (X^{\abs{J}})^\iso \) and, for each \( \omega \in \Omega \) and \( s \in J_\omega \subset \Ob \abs{J}^{L_\omega} \), a morphism \( (X^s)^\iso \colon (X^{\abs{J}})^\iso \to (X^{L_\omega})^\iso \).
  Now form the strict pullback square
  \[
    \begin{tikzcd}
      (X^J)^\iso \pb \ar[r, "\pi_1"] \ar[d, "\pi_0"'] &[35pt] \prod_{\omega \in \Omega} \prod_{s \in J_\omega} (X^{J^\omega})^\iso \ar[d, "\prod_\omega \prod_s p_\omega"] \\[10pt]
      (X^{\abs J})^\iso \ar[r, "\br{\br{(X^s)^\iso}_s}_\omega"] & \prod_{\omega \in \Omega} \prod_{s \in J_\omega} (X^{L_\omega})^\iso,
    \end{tikzcd}
  \]
  noting that the displayed products are all finite, and that the right-hand vertical morphism is an isofibration (the latter by \jref{defn:topos-op}~\ref{item:topos-op-ext-invt} and because products of isofibrations are again isofibrations).
  The left-hand vertical morphism \( \pi_0 \) represents a diagram of shape \( \abs{J} \) in \( \cC\pbig{(X^J)^\iso, X} \), which we claim is a model for \( J \), and exhibits \( (X^J)^\iso \) as a cotensor of \( X \) by \( J \).

  That \( \pi_0 \) is a model for \( J \) means that for each \( \omega \in \Omega \) and \( s \in J_\omega \), the diagram \( L_\omega \to \cC\pbig{(X^J)^\iso, X} \) represented by \( (X^J)^\iso \tox{\pi_0} (X^{\abs{J}})^\iso \tox{(X^s)^\iso} (X^{L_\omega})^\iso \) is in \( \sigma_{\omega, \cC((X^J)^\iso, X)} \).
  This follows from the commutativity of the above square, the fact that \( p_\omega \) represents a diagram in \( \sigma_{\omega, \cC((X^{J^\omega})^\iso)} \), and \jref{defn:topos-op}~\ref{item:topos-op-ext-log-pres}.
  That this exhibits \( (X^J)^\iso \) as a cotensor of \( X \) by \( J \) amounts to the claim that for any groupoid \( A \in \cC \), the functor \( (\pi_0)_* \colon \cC\pbig{A, (X^J)^\iso} \to \cC\pbig{A, (X^{\abs J})^\iso} \) is an isomorphism onto the full subcategory consisting of morphisms \( A \to (X^{\abs J})^\iso \) representing functors \( \abs{J} \to \cC(A, X) \) which are models of \( J \).
  The full-faithfulness and injectivity-on-objects of \( \pi_0 \) follows from the right-hand vertical morphism in the above pullback square having these properties.
  And the fact that \( (\pi_0)_* \) has the correct image follows by inspecting the above pullback square, and noting that a morphism \( f \colon A \to (X^{\abs J})^\iso \) represents a model of \( J \) if and only if \( f \cdot (X^s)^\iso \colon A \to (X^{L_\omega})^\iso \) factors through \( p_\omega \colon (X^{J^\omega})^\iso \to (X^{L_\omega})^\iso \) for each \( \omega \in \Omega \) and \( s \in J_\omega \).

\subsection{Proof of \alttext{\hyperref[propn:thy-to-sketch]{Proposition~\ref*{propn:thy-to-sketch}}}{Proposition~\ref*{propn:thy-to-sketch}}: every theory gives a topos sketch} \label{subsec:thy-to-sketch-proof}
Given a higher order theory \( \cT \), we need to produce a corresponding topos sketch \( J \), which is finite if \( \cT \) is.
The construction is the same as that in \cite[\S9]{makkai-generalized-sketches-all}.

We refer to \aref{sec:topos-logic-background} for our conventions concerning logical formulas.
In everything that follows, we may take \( d = \omega \), in which case by a \( d \)-th order theory, we simply mean a higher-order theory.

\begin{defn}
  Let \( \cL \) be a signature.
  A multi-sorted \defword{\( d \)-th order \( \cL \)-fragment} is a pair \( \cF = (T_\cF, F_\cF) \) satisfying the conditions \ref{item:frag-conds-first}~\ref{item:frag-conds-last} below, where
  \begin{enumerate}[(i)]
  \item \( T_\cF \) is a set of terms-in-context \( (t, \vx) \) over \( \cL^{\up d} \),
  \item \( F_\cF \) is a set of formulas-in-context \( (\phi, \vx) \) over \( \cL^{\up d} \).
  \end{enumerate}
  These must satisfy:
  \begin{enumerate}[(i), resume]
  \item \label{item:frag-conds-first} If \( \pbig{f(t_0,\ldots,t_{n-1}), \vx} \in T_\cF \), then \( (t_i, \vx) \in T_\cF \) for all \( i \).
  \item If \( (t_0 = t_1, \vx) \in F_\cF \) or \( (R t_0,\ldots,t_{n-1}, \vx) \in F_\cF \) for some relation symbol \( R \) in \( \cL \), then \( (t_i, \vx) \in T_\cF \) for all \( i \)
  \item If \( (\phi_0 \square \phi_1, \vx) \in F_\cF \), where \( \square \) is one of \( \wedge \), \( \vee \), or \( \To \), then \( (\phi_0, \vx), (\phi_1, \vx) \in F_\cF \).
  \item \label{item:frag-conds-last} If \( (\square y.\ \phi, \vx) \in F_\cF \), where \( \square \) is one of \( \forall \) or \( \exists \) then \( (\phi, \vx') \in F_\cF \), where \( \vx' = \vx \) if \( y \in \vx \) and \( \vx = \vx y \) otherwise.
  \end{enumerate}

  We will sometimes identify a fragment \( (T_\cF, F_\cF) \) with the disjoint union \( T_\cF \sqcup F_\cF \), so that we can talk about, e.g., the intersection or union of fragments.

  Note that the intersection of any two \( d \)-th order \( \cL \)-fragments is again one.
  Thus, we may define the \( d \)-th order \( \cL \)-fragment \defword{generated} by any set \( \cT \) of \( d \)-th order sentences over \( \cL \) to be the least \( d \)-th order \( \cL \)-fragment with \( (\phi, \emptyset) \in F_\cF \) for all \( \phi \in \cT \).

  The following technical definition is not strictly needed, but makes things simpler in what follows: we say a formula \( \phi \) is \defword{clean} if the free and bound variables in \( \phi \) and its subformulas are distinct; thus, if we call a formula \( \phi \) \emph{pre-clean} if no variable appear both free and bound in \( \phi \), then \( \phi \) is (recursively) defined to be clean if it is pre-clean and each of its immediate subformulas are clean.
  It is clear, by renaming bound variables, that any theory \( \cT \) is equivalent (i.e., has the same models as) a clean theory \( \cT' \) (i.e., a set of clean sentences), thus giving isomorphisms \( (\bC)^{\cT'} \toi (\bC)^{\cT} \) for each topos \( \bC \), which are natural in \( \bC \); thus, to prove \jref{propn:topos-sketch-cotensors}, it suffices to consider clean theories \( \cT \).

  A \( d \)-th order \( \cL \)-fragment \( \cF \) is \defword{clean} if for each \( (\phi, \vx) \in F_\cF \), the set of bound variables of \( \phi \) is disjoint from \( \vx \).
  It is easy to see that the \( d \)-th order \( \cL \)-fragment generated by a set of clean sentences is clean.
\end{defn}

\begin{defn} \label{defn:fragment-sketch}
  \def\arraystretch{1.3}

  Let \( \cL \) be a signature.
  To a clean \( d \)-th order \( \cL \)-fragment \( \cF \), we will associate a topos sketch \( J_\cF \).\footnote{%
    Though the definition makes sense for a general fragment, it will only produce the ``right'' sketch for clean fragments. %
  }
  We will first describe the underlying category \( \abs{J_\cF} \) as a quotient of a free category on a certain graph; that is, we will give the objects, the generating arrows, and then the generating relations.
  We will then describe the specifications of \( J_{\cF} \).

  Each of these entities (objects, arrows, relations, and specifications) will come in different families, each family corresponding to one row of one of the following tables.
  The family depends on certain parameters, which are specified in the last column.
  For the objects and generating arrows, the first column introduces a notation for the entities being introduced.

  In what follows, let \( N \) be the least integer \( N \ge 2 \) such that all of the arities of function and relation symbols in \( \cL \), and all the sequences \( \vx \) of variables appearing in all \( (t, \vx) \in T_\cF \) and all \( (\phi, \vx) \in F_\cF \), have length \( < N \) (we may have \( N = \omega \)).
  Let us also write \( T_\cF^\star \) for the set of all pairs \( (\vt, \vx) \), where \( (f \vt, \vx) \in T_\cF \) for some \( f \in \Fun \cL^{\up d} \), or \( (R \vt, \vx) \in F_\cF \) for some \( R \in \Rel \cL^{\up d} \), or \( \pbig{(t_0 = t_1), \vx} \in F_\cF \), where \( \vt = (t_0,t_1) \).
  \\

  \noindent
  \textbf{Objects of \( \abs{J_{\cF}} \)}:
  \begin{center}
    \begin{tabular}{l|l}
      Notation & Parameters \\
      \hline
      \( A \) & \( A \in \Ob \cL^{\up d} \) \\
      \( \Pi \vA \) & \( \vA \in (\Ob \cL^{\up d})^{< N} \) \\
      \( R \) & \( R \in \Rel \cL \) \\
      \( [\vx \tolon \phi] \) & \( (\phi, \vx) \in F_\cF \)
    \end{tabular}
  \end{center}

  \noindent
  \textbf{Generating arrows of \( \abs{J_{\cF}} \)}:
  \begin{center}
    \begin{tabular}{l|l|l}
      Notation & Domain and codomain & Parameters \\
      \hline
      \( \pi_i^{\vA} \) & \( \Pi \vA \to A_i \) & \( \vA = \in (\Ob \cL^{\up d})^{< N} \) and \( 0 \le i < n \) \\
      \( \pi_{\vA}^{\vA B} \) & \( \Pi (\vA B) \to \Pi \vA \) & \( \vA \in (\Ob \cL^{\up d})^{< N - 1} \) and \( B \in \Ob \cL^{\up d} \) \\
      \( i_R \) & \( R \to \Pi \vA \) & \( R \in \Rel \cL \) with \( R \tolon \vA \) \\
      \( i_{[ \vx \tolon \phi ]} \) & \([\vx \tolon \phi] \to \vA \) & \( (\phi, \vx) \in F_\cF \) with \( \vx \tolon \vA \) \\
      \( f \) & \( \Pi \vA \to B \) & \( f \in \Fun \cL \) with \( f \colon \vA \to B \) \\
      \( \lambda \vx.\ t \) & \( \Pi \vA \to B \) & \( (t, \vx) \in F_\cF \) with \( \vx \tolon \vA \) and \( t \tolon B \) \\
      \( \Delta_B \) & \( B \to \Pi \br{B, B} \) & \( B \in \Ob \cL \) \\
      \( \lambda \vx .\ \vt \) & \( \Pi \vA \to \Pi \vB \) & \( (\vt, \vx) \in T_\cF^\star \) with \( \vx \tolon \vA \) and \( \vt \tolon \vB \) \\
      \( \pi_{R, \vt, \vx} \) & \( [\vx \tolon R \vt] \to R \) & \( (R \vt, \vx) \in F_\cF \) where \( R \in \Rel \cL^{\up d} \) \\
      \( \pi_{t_0, t_1, \vx} \) & \( [ \vx \tolon t_0 = t_1] \to B \) & \( (t_0 = t_1, \vx) \in F_\cF \) where \( t_0, t_1 \tolon B \)
    \end{tabular}
  \end{center}

  \noindent
  \textbf{Generating relations of \( \abs{J_{\cF}} \)}:
  \begin{center}
    \begin{tabular}{l|l}
      Relation & Parameters \\
      \hline
      \( \lambda \vx.\ x_i = \pi_i^{\vA} \) & \( (x_i, \vx) \in T_\cF \) with \( \vx \colon \vA \) and \( 0 \le i < \len \vx \) \\
      \( \pi_{\vA}^{\vA B} \pi_i^{\vA} = \pi_i^{\vA B} \) & \( \vA \in (\Ob \cL^{\up d})^{< N - 1} \) and \( B \in \Ob \cL^{\up d} \) and \( 0 \le i < \len \vA \) \\
      \( (\lambda \vx.\ \vt) \pi_i^{\vB} = \lambda \vx.\ t_i \) & \( (\vt, \vx) \in T_\cF^\star \) and \( 0 \le i < \len \vt \), where \( \vt \tolon \vB \) \\
      \( \Delta_B \pi_0^{\br{B,B}} = \Delta_B \pi_1^{\br{B,B}} = \id_B \) & \( B \in \Ob \cL \) \\
      \( (\lambda \vx.\ f \vt) = (\lambda \vx.\ \vt) f \colon \vA \to C \) & \( (f \vt, \vx) \in T_\cF \) where \( \vx \tolon \vA \) and \( f \in \Fun \cL \) with \( f \colon \vB \to C \) \\
      \prerefpost{(}{eq:frag-sketch-rel-square}{a)} below commutes & \( (R \vt, \vx) \in F_\cF \) where \( R \in \Rel \cL^{\up d} \) with \( R \tolon \vB \) and \( \vx \tolon \vA \) \\
      \prerefpost{(}{eq:frag-sketch-rel-square}{b)} below commutes & \( (t_0 = t_1, \vx) \in F_\cF \) where \( t_0,t_1 \tolon B \) and \( \vx \tolon \vA \)
    \end{tabular}
  \end{center}

  \begin{equation} \label{eq:frag-sketch-rel-square}
    (a)\quad
    \begin{tikzcd}
      {[\vx \tolon R \vt]} \ar[r, "\pi_{R, \vt, \vx}"] \ar[d, "i_{[\vx \tolon R \vt]}"'] &[15pt] R \ar[d, "i_R"] \\
      \Pi \vA \ar[r, "\lambda \vx.\ \vt"] & \Pi \vB
    \end{tikzcd}
    \qquad \qquad
    (b)\quad
    \begin{tikzcd}
      {[\vx \tolon t_0 = t_1]} \ar[r, "\pi_{t_0, t_1, \vx}"] \ar[d, "i_{[\vx \tolon t_0 = t_1]}"'] & B \ar[d, "\Delta_B"] \\
      \Pi \vA \ar[r, "\lambda \vx.\ \br{t_0, t_1}"] & \Pi \br{B, B}
    \end{tikzcd}
  \end{equation}

  \noindent
  \textbf{Specifications of \( J_{\cF} \)}:

  For each family of specifications, we indicate which type of specification it is, i.e., which fundamental topos operation from \jref{defn:topos-op} it corresponds to; in the third line, we recall that a pullback square can be used to specify that a given morphism \( f \colon X \to Y \) is monic, using the square with two instances of \( f \) and two instances of \( \id_X \).
 \begin{center}
   \begin{tabular}{l|l|l}
     Fundamental topos operation & Diagram & Parameters \\
     \hline
     \( \omega_{8 + n} \) (product diagram) & \( (\pi_0^{\vA}, \ldots, \pi_{\len \vA-1}^{\vA}) \) & \( \vA \in (\Ob \cL^{\up})^{< N} \) \\
     \( \omega_0 \) (power diagram) & \prerefpost{(}{eq:frag-sketch-spec-diags}{a)} below & \( A \in \Ob \cL \) and \( 1 \le i < d \) \\
     \( \omega_1 \) (monomorphism) & \( i_R \) & \( R \in \Rel \cL \) \\
     \( \omega_1 \) (pullback square) & \prerefpost{(}{eq:frag-sketch-rel-square}{a)} above &
     same as for \( \pi_{R, \vt, \vx} \) above \\
     \( \omega_1 \) (pullback square) & \prerefpost{(}{eq:frag-sketch-rel-square}{b)} above &
     same as for \( \pi_{t_0, t_1, \vx} \) above \\
     \( \omega_{8 + 1} \) (isomorphism) & \( i_{[\vx \tolon \top]} \colon [\vx \tolon \top] \to \Pi \vA \) & \( (\top, \vx) \in F_\cF \), where \( \vx \tolon \vA \) \\
     \( \omega_2 \) (initial object) & \( [\vx \tolon \bot] \) & \( (\bot, \vx) \in F_\cF \) \\
     \( \omega_5,\omega_6,\omega_7 \) (\( \vee \), \( \wedge \), or \( \to \)) & \( (i_{[\vx \tolon \phi_0]}, i_{[\vx \tolon \phi_1]}, i_{[\vx \tolon \phi]}) \) & \( (\phi, \vx) \in F_\cF \) where \( \phi = \phi_0 \square \phi_1 \), and \( \square \) is \( \wedge \), \( \vee \), or \( \to \) \\
     \( \omega_3,\omega_4 \) (\( \forall \) or \( \exists \)) &
     \prerefpost{(}{eq:frag-sketch-spec-diags}{b)} below &
     \( (\square y.\ \phi, \vx) \in F_\cF \) where \( \square \) is \( \forall \) or \( \exists \)
   \end{tabular}
 \end{center}

  \begin{equation} \label{eq:frag-sketch-spec-diags}
    (a)
    \begin{tikzcd}
      &[-40pt] \in_{\pow^{i-1} A} \ar[d, "i_{\in_{\pow^{i-1} A}}"] &[-40pt] \\
      & \Pi \br{\pow^{i-1} A, \pow^{i} A} \ar[ld, "\pi_0^{\br{\pow^{i-1} A, \pow^{i} A}}"'] \ar[rd, "\pi_1^{\br{\pow^{i-1} A, \pow^{i} A}}"] & \\[10pt]
      \pow^{i-1} A & & \pow^{i} A
    \end{tikzcd}
    \qquad\qquad
    (b) \quad
    \begin{tikzcd}
      {[\vx y \tolon \phi]} \ar[d, "i_{[\vx y \tolon \phi]}"'] &
      {[\vx \tolon \square y.\ \phi]} \ar[d, "i_{[\vx \tolon \square y. \phi]}"] \\
      \Pi (\vA B) \ar[r, "\pi_{\vA}^{\vA B}"] & \Pi \vA
    \end{tikzcd}
  \end{equation}

  \noindent
  (End of \jref{defn:fragment-sketch}.)
\end{defn}

\begin{defn} \label{defn:thy-sketch}
  Let \( \cT \) be an \( d \)-th order theory over \( \cL \), and let \( \cF \) be the \( d \)-th order fragment generated by \( \cT \).
  We define a sketch \( J_\cT \), an extension of \( J_\cF \) with the same underlying category \( \abs{J_\cT} = \abs{J_\cF} \) by adding, for each \( \phi \) in \( \cT \) an isomorphism specification on the morphism \( i_{[\emptyset \colon \phi]} \colon [\emptyset \colon \phi] \to \Pi \emptyset \).
\end{defn}

\begin{lem} \label{lem:finite-thy-finite-sketch}
  The higher-order \( \cL \)-fragment generated by any finite list of sentences is finite; and if \( \cF \) is a finite \( d \)-th order \( \cL \)-fragment or \( \cT \) a finite set of \( d \)-th order sentences over \( \cL \), then \( J_{\cF} \) and \( J_{\cT} \) are finite topos sketches.
\end{lem}
\begin{proof}
  The first claim is proven by an induction showing more generally that the fragment generated by any finite set of formulas-in-context is finite, where the induction step uses that the fragment generated by a set of formulas-in-context is the union of the fragments generated by each one, and that the union of finitely many finite fragments is finite.
  The second claim follows simply by observing that, if \( \cF \) is finite, then in the definition of \( J_{\cF} \), there are only finitely many objects, generating morphisms, generating relations, and specifications.
\end{proof}

Note that if \( \cF_0 \subset \cF_1 \) are two \( d \)-th order \( \cL \)-fragments, then \( J_{\cF_0} \subset J_{\cF_1} \), in the sense that \( \abs{J_{\cF_0}} \subset \abs{J_{\cF_1}} \) and \( (J_{\cF_0})_\omega \subset (J_{\cF_1})_\omega \) for all \( \omega \).
In particular, if \( \cF \) is the fragment generated by any \( d \)-th order theory \( \cT \) over \( \cL \), then \( J_{\cF} \) contains \( J_{\cF_0} \), where \( \cF_0 \) is the \emph{minimal} \( \cL \)-fragment, i.e., the one generated by the empty theory.
We denote this fragment by \( J_\cL \)

\begin{defn}
  Let \( \cL \) be a signature and \( J_\cL \) the associated \( d \)-th order topos sketch for some \( d \).
  To each model \( M \colon J_\cL \to \bC \) of \( J_\cL \) in a topos \( \bC \) we can associate in an obvious a way an interpretation \( \rstr{M}{\cL} \colon \cL \to \bC \), since the data involved in specifying a model \( J_\cL \to \bC \) is a strict superset of that involved in specifying an interpretation \( \cL \to \bC \), and thus from a model of \( J_\cL \), we obtain an interpretation of \( \cL \) simply by restriction.

  Indeed, inspecting \jref{defn:fragment-sketch} in the case \( T_\cF = F_\cF = \emptyset \), we see that a model \( M \colon J_\cL \to \bC \) is given by objects \( M(A) \in \bC \) for each \( A \in \Ob \cL \); products \( M(\Pi \vA) = \prod_i M(A_i) \) in \( \bC \) for each \( \vA \in (\Ob \cL)^{< N} \); power objects \( M(\pow^{i} A) \) of \( M(\pow^{i-1} A) \) for each \( A \in \Ob \cL \) and \( 1 \le i < d \); morphisms \( M(f) \colon M(\Pi \vA) \to M(B) \) for each \( f \colon \vA \to B \) in \( \cL \); and underobjects \( M(i_R) \colon M(R) \to M(\vA) \) for each \( R \tolon \vA \) in \( \cL \).
  (There is furthermore the data of morphisms \( M(\Delta_B) \colon M(B) \to M(\Pi\br{B,B}) \) and \( M(\pi_{\vA}^{\vA B}) \colon M\pbig{\Pi(\vA B)} \to \Pi \vA \), but the values of these are forced by the relations imposed on \( \abs{J_\cL} \).)

  An interpretation \( \cL \to \bC \) is given by the same data, except that we do not need the power objects \( M(\pow^{i} A) \), and we only need those products \( M(\Pi \vA) \) for which \( \vA \) appears as the arity of some function or relation symbol in \( \cL \).
  In fact, we see that, by retaining the power objects \( M(\pow^{i} A) \) as well, \( M \) also gives rise to an interpretation \( \rstr{M}{\cL^{\up d}} \colon \cL^{\up d} \to \bC \) which is a \emph{power extension} of \( \rstr{M}{\cL} \) in the sense explained in \sref{subsec:logic-in-topoi}.

  Similarly, to each isomorphism \( M_0 \toi M_1 \) of models \( J_\cL \to \bC \), we obtain an isomorphism of interpretations \( \rstr{M_0}{\cL} \toi \rstr{M_1}{\cL} \) and \( \rstr{M_0}{\cL^{\up d}} \toi \rstr{M_1}{\cL^{\up d}} \), again by restriction.
  We thus obtain restriction functors \( \rstr{(-)}{\cL} \colon (\bC^{J_\cL})^\iso \to (\bC^\cL)^\iso \) and \( \rstr{(-)}{\cL^{\up d}} \colon (\bC^{J_\cL})^\iso \to (\bC^{\cL^{\up d}})^\iso \).

  Since for any \( d \)-th order \( \cL \)-fragment \( \cF \) or theory \( \cT \), the sketch \( J_\cF \) or \( J_\cT \) is an extension of \( J_\cL \), we similarly obtain restriction functors \( \rstr{(-)}{\cL} \colon \bC^{J_\cF} \to (\bC^\cL)^\iso \) and \( \rstr{(-)}{\cL} \colon \bC^{J_\cT} \to (\bC^\cL)^\iso \) (and likewise \( \rstr{(-)}{\cL^{\up d}} \)).
\end{defn}

We next want to show that the functor \( \rstr{(-)}{\cL} \colon (\bC^{J_{\cT}})^\iso \to (\bC^{\cL})^\iso \) factors through \( (\bC^\cT)^\iso \hto (\bC^\cL)^\iso \), and that the resulting functor \( (\bC^{J_{\cT}})^\iso \to (\bC^\cT)^\iso \) is a surjective equivalence, thus proving the first half of \jref{propn:thy-to-sketch}.

\begin{lem} \label{lem:sketch-model-no-surprises}
  Let \( \cF \) be a \emph{clean} \( d \)-th order \( \cL \)-fragment, let \( J_{\cF} \) be the associated topos sketch, and fix a model \( M \colon J_\cF \to \bC \) of \( J_\cF \) in a topos \( \bC \).
  Then:
  \begin{enumerate}[(i)]
  \item For any \( (t, \vx) \in T_\cF \), the morphisms \( M_\vx(t), (\rstr{M}{\cL^{\up d}})_{\vx}(t) \colon M(\Pi \vA) \to M(B) \) agree, where \( \vx \tolon \vA \) and \( t \colon B \).
  \item For any \( (\phi, \vx) \in F_\cF \), the underobjects \( M(i_{[\vx \colon \phi]}) \colon M([\vx \colon \phi]) \tto M(\Pi \vA) \) and \( (\rstr{M}{\cL^{\up d}})_{\vx}(\phi) \tto M(\Pi \vA) \) are equivalent (noting that the latter is in fact only defined up to equivalence), where \( \vx \tolon \vA \).
  \end{enumerate}
\end{lem}

\begin{proof}
  The proofs are by induction on \( t \) and \( \phi \), referring to the various elements of \jref{defn:fragment-sketch} for the corresponding induction steps.

  When \( t \) is a variable \( x_i \), this follows from the relation \( \lambda \vx.\ x_i = \pi_i^{\vA} \) in \jref{defn:fragment-sketch}.
  When \( t = f \vt \) for some \( f \colon \vB \to C \) in \( \cL \), this follows from the relations \( \lambda \vx.\ f \vt = (\lambda \vx.\ \vt) f \) and \( (\lambda \vx. \vt) \pi_i^{\vB} = \lambda \vx.\ t_i \) in \jref{defn:fragment-sketch} and the inductive hypothesis applied to the \( t_i \).

  For atomic \( \phi \), the value of \( M([\vx \colon \phi]) \) is forced to be the correct one by the two pullback specifications in \jref{defn:fragment-sketch} displayed in \eqref{eq:frag-sketch-rel-square}, by the relations \( \Delta_B \pi_0^{\br{B,B}} = \Delta_B \pi_1^{\br{B,B}} = \id_B \), by part (i) which we already proved, and by the \( \omega_2 \) and \( \omega_{8+1} \) specifications.

  For non-atomic \( \phi \), the value \( M([\vx \colon \phi]) \) is forced to be the correct one by the \( \omega_i \) specifications for \( i = 3,4,5,6,7 \) in \jref{defn:fragment-sketch}, by the relations \( \pi_{\vA}^{\vA B} \pi_i^{\vA} = \pi_i^{\vA B} \), and by the inductive hypothesis applied to the components of \( \phi \); we note that this is where we use that the fragment \( \cF \) is clean.
\end{proof}

\begin{propn} \label{propn:thy-sketch-surj-eq}
  For any clean \( d \)-th order theory \( \cT \) over \( \cL \) generating the \( d \)-th order \( \cL \)-fragment \( \cF \), and any topos \( \bC \):
  \begin{enumerate}[(i)]
  \item The functor \( \rstr{(-)}{\cL} \colon (\bC^{J_\cF})^\iso \to (\bC^{\cL})^\iso \) is a surjective equivalence.
  \item The functor \( \rstr{(-)}{\cL} \colon (\bC^{J_\cT})^\iso \to (\bC^{\cL})^\iso \) factors through \( (\bC^\cT)^\iso \hto (\bC^\cL)^\iso \), and the resulting functor \( (\bC^{J_{\cT}})^\iso \to (\bC^\cT)^\iso \) is a surjective equivalence.
  \end{enumerate}
\end{propn}

\begin{proof}
  We first verify the full-faithfulness in (i).
  Fix models \( M_0, M_1 \colon J_\cF \to \bC \) and an isomorphism \( \alpha \colon \rstr{M_0}{\cL} \toi \rstr{M_1}{\cL} \); we want to show that \( \alpha \) extends to a unique isomorphism \( \alpha \colon M_0 \to M_1 \).
  Referring to the objects of \( \abs{J_\cF} \) as specified \jref{defn:fragment-sketch}, since \( \alpha \) is already defined on the objects \( A \), \( \Pi \vA \), and \( R \), it remains to define it on the objects \( [\vx \colon \phi] \).
  By \jref{lem:sketch-model-no-surprises}, for each such object, we have that \( M_j(i_{[\vx \colon \phi]}) \colon M_j([\vx \colon \phi]) \tto M_j(\Pi \vA) \), where \( \vx \tolon \vA \), is equivalent to the underobject \( (\rstr{M_j}{\cL})_{\vx}(\phi) \tto M_j(\Pi \vA) \).
  Now it is a general fact (proven by induction on \( \phi \)) that for isomorphic interpretations \( N_0,N_1 \colon \cL \to \bC \), the underobjects \( (N_j)_{\vx}(\phi) \tto N_j(\vA) \) for \( j = 0,1 \) agree under the given isomorphism \( N_0(\vA) \toi N_1(\vA) \).
  Thus, we may, and are forced to, define \( \alpha_{[\vx \colon \phi]} \) to be the unique isomorphism making
  \[
    \begin{tikzcd}
      M_0([\vx \colon \phi]) \ar[r, dashed, "\alpha_{[\vx \colon \phi]}"] \ar[d, >->, "i_{[\vx \colon \phi]}"] &[10pt]
      M_1([\vx \colon \phi]) \ar[d, >->, "i_{[\vx \colon \phi]}"]
      \\
      M_0(\Pi \vA) \ar[r, "\alpha_{\vA}"', "\sim'"] &
      M_1(\Pi \vA)
    \end{tikzcd}
  \]
  commute, where \( \alpha_{\vA} \) is the unique morphism with \( \alpha_{\vA} \pi_i = \pi_i \alpha_{A_i} \) for all \( i \).
  It remains to see that \( \alpha \), thus defined, is natural, which involves a condition for each of the generating arrows of \( \abs{J_\cF} \) given in \jref{defn:fragment-sketch}.
  For the arrows \( i_R \) and \( f \), this follows from \( \alpha \) being an isomorphism of interpretations; for \( i_{[\vx \colon \phi]} \), it follows from the definition of \( \alpha_{[\vx \colon \phi]} \); for \( \pi_{R, \vt, \vx} \) and \( \pi_{t_0, t_1, \vx} \), it follows as a bi-product of the base-cases in the proof that \( (\rstr{M_0}{\cL})_{\vx}(\phi) \cong (\rstr{M_1}{\cL})_{\vx}(\phi) \); for \( \pi_i^{\vA} \), it follows from the definition of \( \alpha_{\vA} \); for \( \pi_{R, \vt, \vx} \) and \( \pi_{t_0, t_1, \vx} \), it follows as a bi-product of the definition of \( \alpha_{[\vx \colon \phi]} \); and in the remaining cases \( \pi_{\vA}^{\vA B} \), \( \Delta_B \), \( \lambda \vx.\ t \), and \( \lambda \vx.\ \vt \), it follows from the fact that the set of morphisms with respect to which \( \alpha \) is natural is closed under composition and forming the induced morphism into a product, and using induction on \( t \).
  We leave the remaining details to the reader.

  Next, we address the essential surjectivity in (i).
  We are given an interpretation \( \rstr{M}{\cL} \colon \cL \to \bC \), and want to extend it to a model \( M \colon J_{\cF} \to \bC \).
  On objects, we only need to define \( M \) on the objects \( [\vx \colon \phi] \), and of course we set \( M([\vx \colon \phi]) = (\rstr{M}{\cL})_{\vx}(\phi) \); and this also gives the definitions of \( M(i_{[\vx \colon \phi]}) \defeq i_{(\rstr{M}{\cL})_{\vx}(\phi)} \), and of \( M(\pi_{R, \vt, \vx}) \) and \( M(\pi_{t_0, t_1, \vx}) \) as the projections in the pullback squares defining \( (\rstr{M}{\cL})_{\vx}(R \vt) \) and \( (\rstr{M}{\cL})_{\vx}(t_0 = t_1) \).
  It remains to define \( M \) on the morphisms \( \pi_{\vA}^{\vA B} \), \( \lambda \vx.\ t \), \( \Delta_B \), \( \lambda \vx.\ \vt \); for all of these, the definition (using recursion on \( t \)) is obvious, and is in fact uniquely determined by the imposed relations in \( \abs{J_\cF} \).

  The fact that all the relations and sketch specifications imposed in \( J_\cF \) are satisfied by \( M \) is immediate from the definitions of the action of \( M \) on the relevant objects and morphisms.

  Finally, we turn to (ii).
  The fact that \( \rstr{(-)}{\cL} \colon (\bC^{J_\cT})^\iso \to (\bC^{\cL})^\iso \) factors through \( (\bC^{\cT})^\iso \) follows immediately from \jref{lem:sketch-model-no-surprises} and the definition of \( J_\cT \) as an extension of \( J_{\cF} \).
  We now have a commutative square of restriction functors
  \[
    \begin{tikzcd}
      (\bC^{J_\cT})^\iso \ar[r, ""] \ar[d, ""'] & (\bC^{J_\cF})^\iso \ar[d, ""] \\
      (\bC^{\cT})^\iso \ar[r, ""] & (\bC^{\cL})^\iso
    \end{tikzcd}
  \]
  where the right vertical functor is a surjective equivalence, as we just showed.
  It is easy to see that the top horizontal functor is fully faithful, since we only introduce new specifications (and no new objects or morphisms) in extending \( J_\cF \) to form \( J_\cT \), and that the bottom horizontal functor is fully faithful for the same reason.
  It follows that the left-hand vertical functor is fully faithful.
  That this functor is also surjective on objects follows from the fact that the right-hand vertical functor is surjective on objects, and the obvious fact that if \( \rstr{M}{\cL} \colon \cL \to \bC \) is a model of \( \cT \) for some \( M \colon J_\cF \to \bC \), then \( M \) can be extended to a model of \( J_\cT \).
\end{proof}

We are now nearly done with the proof of \jref{propn:thy-to-sketch}: given a \( d \)-th order theory \( \cT \) over \( \cL \) (with \( 0 \le d \le \omega \)), we take \( J \) to be the sketch \( J_\cT \) of \jref{defn:thy-sketch}---which is finite if \( \cT \) is by \jref{lem:finite-thy-finite-sketch}---and for each topos \( \bC \), we take \( \alpha \colon (\bC^J)^\iso \to (\bC^\cT)^\iso \) to be the restriction functor \( \rstr{(-)}{\cL} \colon (\bC^{J_\cT})^\iso \to (\bC^\cT)^\iso \), which is a surjective equivalence by \jref{propn:thy-sketch-surj-eq}.
It remains only to see that for any invertible 2-cell \( \beta \colon \bC \tocell \bD \) between logical functors, the square \eqref{eq:cotensor-comparison-square} on \pref{eq:cotensor-comparison-square} strictly commutes.
However, this is immediate from the definition of the restriction functor \( \rstr{(-)}{\cL} \).

\subsection{Proofs for Sections~\ref{subsec:hs-2-topos}~and~\ref{subsec:topoculi}: construction of topoculi and hs-classifiers} \label{subsec:topoculi-proofs}
In this section, we give the proofs of \jrefs{propn:hs-classifiers-exist}{lem:theories-exist} from \sref{subsec:hs-2-topos} and \jref{propn:2-topoi-give-topoculi} from \sref{subsec:topoculi}.

\begin{proof}[Proof of \jref{propn:hs-classifiers-exist}] \label{proof:hs-classifiers-exist}
  As in the \hyperref[proof:theory-cotensors]{proof of Theorem~\ref*{thm:theory-cotensors}}, we leave out the ``\( \iso \)'' superscripts for readability.

  Fix higher-order theories \( \cT_{\hs} \), \( \cT_{\hs} \sqcup \cT_{\hs} \) and \( \cT_{\unEin} \), embeddings \( i_0 \), \( i_1 \), and \( i \), and isomorphisms \( F_\bC \), \( G_\bC \), and \( H_\bC \) as in \jref{lem:theories-exist}.
  By \jref{thm:theory-cotensors}, there exist cotensors \( E^{\cT_{\hs}} \), \( E^{\cT_{\hs} \sqcup \cT_{\hs}} \), and \( E^{\cT_{\unEin}} \) of \( E \) by each of these theories.

  We claim, first of all, that \( E^{\cT_\hs} \) is a pre-hs-classifier for \( E \).
  As universal hs-object, we take the image \( X_E \) of the universal model \( \cT_\hs \to \cC(E^{\cT_\hs}, E) \) under the isomorphism \( F_{\cC(E^{\cT_\hs}, E)} \colon \cC(E^{\cT_\hs}, E)^{\cT_\hs} \toi \hs\pbig{\cC(E^{\cT_\hs}, E)} \).
  To see that this is indeed a universal hs-object, we note that for each groupoid \( A \in \cC \), we have an isomorphism \( \cC(A, E^{\cT_\hs}) \toi \cC(A, E)^{\cT_\hs} \toix{F_{\cC(A,E)}} \hs\pbig{\cC(A, E)} \), and we need only see that this is the intended functor \( f \mapsto f^* X_E \); this is done by the same argument as in \hyperref[proof:theory-cotensors]{proof of Theorem~\ref*{thm:theory-cotensors}} (now using \jref{lem:theories-exist}~\ref{item:theories-exist-nat}).

  Next, we endow \( E^{\cT_\hs} \) with an hs-classifier structure via the morphism \( \br{E^i E^{i_0},E^i E^{i_1}} \colon E^{\cT_\unEin} \to E^{\cT_\hs} \times E^{\cT_\hs} \).
  That this indeed makes \( E^{\cT_\hs} \) into an hs-classifier follows from the commutativity of the following diagram, in which the lower squares commute by \jref{lem:theories-exist}~\ref{item:theories-exist-prod}~and~\ref{item:theories-exist-inc}, and the upper squares commute by the definition of \( E^i \), \( E^{i_0} \), and \( E^{i_1} \).
  \[
    \begin{tikzcd}[baseline=(bl.base)]
      \cC(A, E^{\cT_{\Ein}}) \ar[r, "(E^i)_*"] \ar[d, "\sim" sloped, "f \mapsto \ev \circ f^*"'] & \cC(A, E^{\cT_\hs \sqcup \cT_\hs}) \ar[r, "(E^{i_0})_*", shift left] \ar[r, "(E^{i_1})_*"', shift right] \ar[d, "\sim" sloped, "f \mapsto \ev \circ f^*"'] & \cC(A, E^{\cT_\hs}) \ar[d, "\sim"' sloped, "f \mapsto \ev \circ f^*"] \\
      \cC(A, E)^{\cT_{\Ein}} \ar[r, "i^*"] \ar[d, "H_{\cC(A,E)}"', "\sim" sloped] & \cC(A, E)^{\cT_\hs \sqcup \cT_\hs} \ar[r, "i_0^*", shift left] \ar[r, "i_1^*"', shift right] \ar[d, "G_{\cC(A,E)}", "\sim"' sloped] & \cC(A, E)^{\cT_\hs} \ar[d, "F_{\cC(A,E)}", "\sim"' sloped] \\
      |[alias=bl]| \unEin_{\cC(A, E)} \ar[r, "i_{\unEin_{\cC(A,E)}}"] & \hs\pbig{\cC(A, E)} \times \hs\pbig{\cC(A, E)} \ar[r, "\pi_0", shift left] \ar[r, "\pi_1"', shift right] & \hs\pbig{\cC(A, E)}
    \end{tikzcd}
    \qedhere
  \]
\end{proof}

\begin{proof}[Proof of \jref{lem:theories-exist}] \label{proof:theories-exist}
  We must produce the finite higher-order theories \( \cT_\hs \), \( \cT_{\hs} \sqcup \cT_{\hs} \), and \( \cT_{\unEin} \), the embeddings \( i_0 \), \( i_1 \), and \( i \), and the isomorphisms \( F_\bC \), \( G_\bC \), and \( H_\bC \), and verify that they satisfy the properties \ref{item:theories-exist-nat},~\ref{item:theories-exist-prod},~and~\ref{item:theories-exist-inc} listed in the statement of the lemma.

  The existence of the theory \( \cT_{\hs} \) and the isomorphisms \( F_\bC \) is almost immediate from the way that hs-objects are defined: an \( \unein \)-object in \( \bC \) is nothing but an interpretation of the signature \( \cL_\unein \) consisting of a single sort \( X \) and a single binary relation symbol \( \unein \), and an hs-object is a model for the theory \( \cT_\hs \) over \( \cL_\unein \) consisting of the sentences expressing the extensionality, well-foundedness, and the existence of a top element:
  \begin{itemize}
  \item \( \forall x,y \tolon X.\ (\forall z \tolon X.\ z \ein x \tot z \ein y) \to x = y \)
  \item \( \forall p \tolon \pow X.\ \pbig{\forall x \tolon X.\ (\forall y \tolon X.\ y \ein x \to y \in p) \to x \in p} \to \forall x.\ x \in p \)
  \item \(
    \exists \rt_X \tolon X.\
    \forall p \tolon \pow X.\
    \pBig{\rt_X \in p \wedge \forall x,y \tolon X.\ \pbig{(x \ein y \wedge y \in p) \to x \in p}}
    \to \forall x \tolon X.\ x \in p.
    \)
  \end{itemize}
  Moreover, an isomorphism of models of \( \cT_\hs \) is the same thing as an isomorphism of hs-objects.
  That the naturality property \ref{item:theories-exist-nat} holds for the isomorphisms \( F_\bC \) is obvious.

  (There is a slight subtlety here, in that the last sentence displayed above a priori only guarantees the existence of a \emph{generalized} element \( \tm \oot U \tox{\rt_X} X \), rather than a global element \( \rt_X \colon \tm \to X \) as required by \jref{defn:topos-ein-hs}; however, the latter is guaranteed by \jrefs{cor:dwn-mono}{propn:epis-eff}.)

  The existence of the theory \( \cT_{\hs} \sqcup \cT_{\hs} \) and isomorphisms \( G_\bC \) satisfying \ref{item:theories-exist-prod} then follows readily: as the notation suggests, we simply take \( \cT_{\hs} \sqcup \cT_{\hs} \) to be ``two copies'' of \( \cT_\hs \): the underlying signature \( \cL_\unein \sqcup \cL_\unein \) has two sorts \( X \) and \( Y \) and two relations \( \unein_X \tolon \br{X,X} \) and \( \unein_Y \tolon \br{Y,Y} \), and two copies of the axioms of \( \cT_\hs \), with one copy referring to \( \unein_X \) and the other \( \unein_Y \).

  Finally, we define the signature \( \cL_\unEin \) as an extension of \( \cL_\unein \sqcup \cL_\unein \) by adding a single function symbol \( f \colon X \to Y \), and \( \cT_\hs \) by adding to the axioms of \( \cT_{\hs} \sqcup \cT_{\hs} \) the following sentences expressing that \( f \) is elemental, whereupon the existence of \( H_\bC \) satisfying \ref{item:theories-exist-inc} is immediate.
  \begin{itemize}
  \item \( \forall x_0,x_1 \tolon X.\ x_0 \ein_X x_1 \to f(x_0) \ein_Y f(x_1) \)
  \item \( \forall x_1 \tolon X.\ \forall y_0 \tolon Y.\ y_0 \ein_Y f(x_1) \to \exists x_0 \tolon X.\ f(x_0) = y_0 \wedge x_0 \ein_X x_1 \)
  \item \( \forall \rt_X \tolon X.\ \forall \rt_Y \tolon Y.\ \pbig{\topp_X(\rt_X) \wedge \topp_Y(\rt_Y)} \to f(x) \ein_Y y \).
  \end{itemize}
  Here \( \topp_X(\rt_X) \) is the formula appearing above under the quantifier ``\( \exists \rt_X \tolon X \)'' in the sentence expressing existence of a top element, and \( \topp_Y(\rt_Y) \) is the same formula, but with \( Y \) substituted for \( X \).
\end{proof}

\begin{proof}[Proof of \jref{propn:2-topoi-give-topoculi}] \label{proof:2-topoi-give-topoculi}
  We verify the claims in \jref{propn:2-topoi-give-topoculi}~\ref{item:2-topoi-give-topoculi-s-s}; the arguments for \ref{item:2-topoi-give-topoculi-ua} are similar.

  We must verify that the given \( \fS \) satisfies the defining conditions of a topoculus.
  That the morphisms in \( \fS \) include all monomorphisms and are closed under composition and pullback follows immediately from \( \rS' \) being a (pre-)plentiful DOF classifier (see \sref{subsubsec:plenitude}); here, we note that any monomorphism in \( \DOF_{\rS'}(\tm) \) is still monic in \( \cC \), since the inclusion of \( \DOF_{\rS'}(\tm) \) into the underlying 1-category of \( \cC \) preserves pullbacks.

  That \( \fS_U \) is a topos follows from the isomorphism \( \fS_U \cong \SOF_{\rS}(U) \) and \jref{thm:topoi-in-topoi-main-thm}.
  Similarly, that \( \fS_{U} \hto \SOF_{\rS'}(\tm) / U \) is logical for each \( U \in \SOF_{\rS'}(\tm) \) also follows from that theorem, since this inclusion is isomorphic to the inclusion \( \SOF_{\rS}(U) \hto \SOF_{\rS'}(U) \).

  Finally, supposing that \( \hs(\rS) \in \SOF_{\rS'}(\tm) \), we must verify that it is an hs-classifier for \( \fS \)---or more precisely, that it, together with the relation \( \unEin_{\rS} \tto \hs(\rS) \times \hs(\rS) \) coming from \( \hs(\rS) \) being an hs-classifier for \( \rS \), are an hs-classifier for \( \fS \).
  (Note that, by definition, the morphisms in \( \SOF_{\rS'}(\tm) \) are really \emph{isomorphism classes} of morphisms in \( \cC \) (see \sref{subsubsec:dof-sof-digression}), so here we are dealing with the isomorphism class of \( i_{\unEin_{\rS}} \); however, we will not distinguish notationally between morphisms and their isomorphism classes.)

  As universal hs-object \( \und{X_\fS} \in \fS_{\hs(\rS)} \), we of course take the image of \( X_\rS \) (in the notation of \jref{defn:int-topos-hs-classifier}) under the equivalence \( \El^{\hs(\rS)} \colon \cC(\hs(\rS), \rS) \toi \SOF_{\rS}(\rS) \cong \fS_{\hs(\rS)} \).
  We then need to check, for each \( A \in \SOF_{\rS'}(\tm) \) and \( X \in \hs(\fS_A) \), that there is a unique \( f \colon A \to \hs(\rS) \) with \( A \cong f^* X_\rS \).
  Since \( \El^A \colon \cC(A, \rS) \to \fS_A \) is an equivalence, we have an hs-object \( Y \in \cC(A, \rS) \) mapping to \( X \) under \( \El^A \), and hence, by the universal property of \( \hs(\rS) \), there is a unique up to isomorphism \( f \colon A \to \hs(\rS) \) with \( f^* X_\rS \cong Y \).
  By the commutativity of the square of anafunctors
  \[
    \begin{tikzcd}
      \cC(\hs(\rS), \rS) \ar[r, "f^*"] \ar[d, "\El^{\hs(\rS)}"', "\sim" sloped] & \cC(A, \rS) \ar[d, "\El^{A}", "\sim"' sloped] \\
      \fS_{\hs(\rS)} \ar[r, "f^*"] & \fS_A
    \end{tikzcd}
  \]
  (\jref{propn:el-pb-compat}), this is equivalent to the desired condition \( f^* \und{X_{\fS}} \cong \und X \) (and the uniqueness up to isomorphism of \( f \) in \( \cC \) gives the required uniqueness of \( f \) in \( \SOF_{\rS'}(\tm) \)).
\end{proof}

\subsection{Proof of \alttext{\hyperref[propn:univ-ext]{Proposition~\ref*{propn:univ-ext}}}{Proposition~\ref*{propn:univ-ext}}: extensionality of the universe} \label{subsec:univ-ext-proof}
We are given hs-objects \( X \) and \( Y \) in a topos \( \bC \).
We are assuming that \( W \Ein \und{U \times X} \) iff \( W \Ein \und{U \times Y} \) for all \( U \in \bC \) and all \( W \in \hs(\bC / U) \), and we need to show that \( X \cong Y \).

Thus, our first task is to construct a morphism \( X \to Y \) (and, symmetrically, \( Y \to X \)).

The set-theoretic argument:
\begin{settybox}
  Fix hs-sets \( X \) and \( Y \), and assume \( Z \Ein X \) iff \( Z \Ein Y \) for all hs-sets \( Z \).
  We need to show \( X \cong Y \).

  Given any \( b \in X \) with \( b \ein \rt_X \), we have \( \dwnfn(b) \Ein X \) by \jref{propn:dwn-hs}, and thus by assumption \( \dwnfn(b) \Ein Y \), so we have a (by \jref{propn:init-uniq} unique) elemental map \( f_b \colon \dwnfn(b) \to Y \).
  Then \( f_b(x) = f_{b'}(x) \) for \( x \in \dwnfn(b) \cap \dwnfn(b') \) by \jref{propn:init-uniq} since the restrictions of \( f_b \) and \( f_{b'} \) to \( \dwnfn(b) \cap \dwnfn(b') \) are both initial by \jrefs{propn:closed-iff-init}{propn:init-2-of-3}.
  Thus we can define \( f \colon X \to Y \) by putting \( f(x) = f_b(x) \) if \( x \in \dwnfn(b) \), and \( f(\rt_X) = \rt_Y \).
  This is well-defined by \jrefs{propn:dwn-opts}{propn:no-loops}.

  It then remains to show that \( f \) is initial, and hence that \( g \colon Y \to X \), constructed symetrically, is as well.
  Once this is done, it follows from \jref{propn:init-2-of-3} that \( f g \) and \( g f \) are initial, hence equal to the identity by \jref{propn:init-uniq}, and we're done.

  We first show \( f \) is consistent.
  Fix \( x,y \in X \) with \( x \ein y \); we need to show \( f(x) \ein f(y) \).
  If \( y = \rt_X \), then \( f(y) = \rt_Y \) and \( f(x) = f_x(x) = f_x(\rt_{\dwnfn(x)}) \ein \rt_Y \) since \( f_x \colon \dwnfn(x) \to Y \) is elemental.
  If \( y \in \dwnfn(b) \) for some \( b \ein \rt_X \), then \( x \in \dwnfn(b) \) as well, and \( f(x) = f_b(x) \ein f_b(y) = f(y) \), since \( f_b \) is consistent.

  Finally, we show \( f \) is initial.
  Fix \( x_1 \in X \) and \( y_0 \in Y \) with \( y_0 \ein f(x_1) \); we need to find \( x_0 \in X \) with \( x_0 \ein x_1 \) and \( f(x_0) = y_0 \).
  If \( x_1 \in \dwnfn(b) \) for some \( b \ein \rt_X \), then \( f(x_1) = f_b(x_1) \), and the claim follows from the initiality of \( f_b \).

  If \( x_1 = \rt_X \), then \( y_0 \ein \rt_Y = f(x_1) \), hence \( \dwnfn(y_0) \Ein Y \), and hence by assumption we have an elemental morphism \( h \colon \dwnfn(y_0) \to X \).
  Letting \( b = h(y_0) \), we have \( b \ein \rt_X = x_1 \) since \( h \) is elemental, and \( h\pbig{\dwnfn(y_0)} = \dwnfn(b) \) by \jref{propn:dwn-dwn}.
  The composite \( \dwnfn(y_0) \tox{h} \dwnfn(b) \tox{f_b} Y \) must be equal to the inclusion \( \dwnfn(y_0) \hto Y \) by \jrefs{propn:init-2-of-3}{propn:init-uniq}.
  Hence \( f(b) = f_b(b) = f_b\pbig{h(y_0)} y_0 \) and \( b \ein x_1 \), so \( x_0 = b \) is as desired.
\end{settybox}

We now turn to the topos-theoretic proof, starting with the construction of \( f \colon X \to Y \).

Consider the underobjects \( \set{\rt_X},\oring X \tto X \) given by
\[
  \set{\rt_X} = \set{x \in X \mid x = \rt_X} \tto X
  \quad\AND\quad
  \oring X = \set{x \in  X \mid x \prec \rt_X } \tto X.
\]
It follows from \jref{propn:dwn-opts} that these form a cover of \( X \).
Moreover, \( \set{\rt_X} \wedge \oring X \approx \zset \) by \jref{propn:no-loops}.
Hence, by \jref{propn:epis-eff}, to define a morphism \( f \colon X \to Y \), it suffices to define morphisms \( f_\rt \colon \set{\rt_X} \to Y \) and \( f_\circ \colon \oring X \to Y \).
Obviously, we take \( f_\rt = \unex \rt_Y \).

In analogy to the set-theoretic construction, we next want to define \( f_\circ \) by ``defining it separately on each subset \( \dwnfn(b) \) for \( b \ein \rt_X \)''.
However, this element-wise construction is not available to us, so instead we pass to the ``generalized element consisting of all \( b \ein \rt_X \)'':

Let \( \abs{X} \tto X \) be the underobject given by
\begin{equation*}
  \abs X = \set{ b \in X \mid b \ein \rt_X } \tto X,\footnote{%
    This notation for the ``set of elements'' of an hs-object comes from \cite{shulman-comparing}; note that the corresponding definition for hs-sets satisfies \( \abs{\tc(\set{x})} = x \) for a hereditary set \( x \) (with \( \tc \) as in \sref{subsec:ein-sets}).
  }
\end{equation*}
and consider the underobject \( \Dwn_{i_{\abs X}} \tto \abs X \times X \) of \jref{defn:down-family}.
By \jref{lem:dwn-is-descendent}, we have a pullback square
\begin{equation} \label{eq:univ-ext-pf-ances-sq}
  \begin{tikzcd}
    \Dwn_{i_{\abs X}} \pb \ar[r, ""] \ar[d, "i_{\Dwn_x}"'] &[20pt] \succeq_X \ar[d, "i_{\succeq_X}"] \\
    \abs{X} \times X \ar[r, "i_{\abs X} \times \id_X"] & X \times X.
  \end{tikzcd}
\end{equation}

Since clearly \( i_{\abs X} \ein \unex_{\abs X} \rt_X \), it follows from \jref{propn:dwn-hs} that
\( \und{\Dwn_{i_{\abs X}}} \tto \und{\abs X \times X} \)
is an elemental substructure in \( \bC / \abs{X} \).
Hence, by assumption \( \und{\Dwn_{i_{\abs X}}} \Ein \und{\abs X \times Y} \) so that we have a (unique) elemental morphism
\[
  \tilde f_\circ \colon \und{\Dwn_{i_{\abs X}}} \to \und{\abs X \times Y}.
\]

\begin{claim} \label{claim:rho-exists}
  \( \oring X \tto X \) is the image of the \( \Dwn_{i_{\abs X}} \tox{i_{\Dwn_{\abs{X}}}} \abs X \times X \) under \( \pi_X \colon \abs{X} \times X \to X \), i.e., there exists a (necessarily unique) \( \rho \) making the following square commute, and \( \rho \) is epic.
  \begin{equation}\label{eq:rho-defn}
    \begin{tikzcd}
      \Dwn_{i_{\abs X}} \ar[r, "i_{\Dwn_{\abs{X}}}", >->] \ar[d, "\rho"', ->>, dashed] & \abs X \times X \ar[d, "\pi_X"] \\
      \oring X \ar[r, "i_{\oring X}", >->] & X
    \end{tikzcd}
\end{equation}
\end{claim}

\begin{proof}
  The image in question is equivalently the image of
  \begin{equation} \label{eq:univ-ext-pf-prec-rel}
    \Dwn_{i_{\abs X}} \tox{i_{\Dwn_{\abs{X}}}} \abs X \times X \tox{i_{\abs X} \times X} X \times X
  \end{equation}
  under \( \pi_1 \colon X \times X \to X \).
  By \eqref{eq:univ-ext-pf-ances-sq} and the definition of \( \abs{X} \), the underobject \eqref{eq:univ-ext-pf-prec-rel} is equivalently
  \[
    \set{ \br{b, x} \mid x \preceq b \wedge b \ein \rt_X } \tto X \times X.
  \]
  The claim thus follows from the definitions of \( \oring X \) and \( \prec_X \).
\end{proof}

\begin{claim} \label{claim:rho-init}
  The morphism \( \rho \) is relatively initial over \( \unex \colon \abs X \to \tm \) (in the \jref{defn:rel-init-elem}).
\end{claim}

\begin{proof}
  First note that the remaining three morphisms in \eqref{eq:rho-defn} are relatively initial: \( i_{\oring X} \) and \( i_{\Dwn_{\abs X}} \) are initial morphisms in \( \bC \) and \( \bC / \abs{X} \) by \jref{propn:closed-iff-init}, and \( \pi_X \) is relatively initial over \( \unex_{\abs X} \) by \jref{propn:rel-notions-compare} (since \( \id_{\abs{X} \times X} \) is initial in \( \bC / \abs{X} \)).
  The claim thus follows from \jref{cor:rel-init-2-of-3}.
\end{proof}

Now, the morphism \( \Dwn_{i_{\abs X}} \tox{\tilde f_\circ} \abs X \times Y \tox{\pi_Y} Y \) is a single entity corresponding to all of the morphisms ``\( f_b \colon \dwnfn(b) \to Y \)'' from the above set-theoretic proof.
There, we next showed that all of the ``\( f_b \)'' agreed on their common domain.
Here, this corresponds to the claim that \( \tilde f_\circ \pi_Y \colon \Dwn_{i_{\abs X}} \to Y \) descends along the epimorphism \( \rho \) of \jref{claim:rho-exists} to a morphism \( f_\circ \colon \oring X \to Y \).
By \jref{propn:epis-eff}, this amounts to:
\begin{claim} \label{claim:f-well-def}
  With \( \pi_0 \) and \( \pi_1 \) the projections \( \pi_0,\pi_1 \colon \Dwn_{i_{\abs X}} \times_{\oring X} \Dwn_{i_{\abs X}} \to \Dwn_{i_{\abs X}} \), we have \( \pi_0 \tilde f_\circ \pi_Y = \pi_1 \tilde f_\circ \pi_Y \colon \Dwn_{i_{\abs X}} \times_{\oring X} \Dwn_{i_{\abs X}} \to Y \).
\end{claim}

\begin{proof}[Proof of \jref{claim:f-well-def}]
  Given \( w_0, w_1 \colon U \to \Dwn_{i_{\abs X}} \) with \( w_0 \rho = w_1 \rho \colon U \to \oring X \), we need to show \( w_0 \tilde f_\circ \pi_Y = w_1 \tilde f_\circ \pi_Y \colon U \to Y \).
  For \( j = 0,1 \), we have \( w_j i_{\Dwn_{i_{\abs X}}} = \br{b_j, x} \colon U \to \abs{X} \times X \) for some \( b_j \colon U \to \abs{X} \) and some \( x \colon U \to X \).

  We now consider the substructures \( \und{\Dwn_{b_j i_{\abs X}}} \to \und{U \times X} \) induced by the morphisms \( b_j i_{\abs X} \colon U \to X \).
  By \jref{cor:dwn-stable}, we have pullback squares
  \begin{equation} \label{eq:wj-induces-vj}
    \begin{tikzcd}
      U \ar[rrd, "w_j", bend left=10pt] \ar[rd, "v_j", dashed] \ar[ddr, "\br{\id, x}"', bend right] & &[20pt] \\
      & \Dwn_{b_j i_{\abs X}} \pb \ar[r, "\lambda_j"] \ar[d, >->] & \Dwn_{i_{\abs X}} \ar[d, >->] \\
      & U \times X \ar[r, "b_j \times \id_X"] & \abs{X} \times X.
    \end{tikzcd}
  \end{equation}
  and hence induced morphisms \( v_j \) as shown.

  Since \( b_j i_{\abs X} \ein \unex \rt_X \), it follows from \jref{propn:dwn-hs} that each \( \und{\Dwn_{b_j i_{\abs X}}} \to \und{U \times X} \) is an elemental substructure.
  Hence, by assumption, we have elemental morphisms
  \( \tilde f_{j} \colon \und{\Dwn_{b_j i_{\abs X}}} \to \und{U \times Y} \).
  It follows from \jref{propn:pb-fun-is-logical} and \jref{propn:init-uniq} that
  \begin{equation} \label{eq:f-induces-fj}
    \begin{tikzcd}
      \Dwn_{b_j i_{\abs X}} \pb \ar[r, "\lambda_j"] \ar[d, "\tilde f_j"'] &[30pt] \Dwn_{i_{\abs X}} \ar[d, "\tilde f_\circ"] \\
      U \times Y \ar[r, "b_j \times \id_Y"] & \abs{X} \times Y.
    \end{tikzcd}
  \end{equation}
  is a pullback square.
  Indeed, the logical anafunctor \( b_j^* \colon \bC / \abs{X} \to \bC / U \) preserves initial morphisms, and \( \tilde f_j \) is the unique initial morphism \( \und{\Dwn_{b_j i_{\abs X}}} \to \und{U \times Y} \).

  By \eqref{eq:wj-induces-vj} and \eqref{eq:f-induces-fj}, we have
  \[
    w_j \tilde f_\circ \pi_Y = v_j \lambda_j \tilde f_\circ \pi_Y = v_j \tilde f_j \pi_Y.
  \]
  Hence, to prove \( w_0 \tilde f_\circ \pi_Y = w_1 \tilde f_\circ \pi_Y \), it suffices to prove
  \[
    v_0 \tilde f_0 \pi_Y = v_1 \tilde f_1 \pi_Y \colon U \to Y.
  \]

  Consider the intersection \( \und{\Dwn_{b_0 i_{\abs X}}} \wedge \und{\Dwn_{b_1 i_{\abs X}}} \tto \und{U \times X} \),
  and write \( u_j \colon \Dwn_{b_0 i_{\abs X}} \wedge \Dwn_{b_1 i_{\abs X}} \tto \Dwn_{b_j i_{\abs X}} \) for the inclusions, which are initial by \jrefs{propn:closed-iff-init}{propn:init-2-of-3} and since the intersection of closed substructures is closed.
  We thus have an equality of initial morphisms
  \begin{equation} \label{eq:uf-init-eq}
    u_0 \tilde f_0 = u_1 \tilde f_1 \colon \und{\Dwn_{b_0 i_{\abs X}}} \wedge \und{\Dwn_{b_1 i_{\abs X}}} \to \und{U \times Y}
  \end{equation}
  by \jrefs{propn:init-2-of-3}{propn:init-uniq}.

  Since \( v_0 i_{\Dwn_{b_0 i_{\abs X}}} = \br{\id_U, x} = v_1 i_{\Dwn_{b_0 i_{\abs X}}} \), there exists \( v \colon U \to \Dwn_{b_0 i_{\abs X}} \wedge \Dwn_{b_1 i_{\abs X}} \) with \( v u_0 = v_0  \) and \( v u_1 = v_1 \).
  \[
    \begin{tikzcd}
      &[30pt] \Dwn_{b_0 i_{\abs X}} \ar[rd, >->, ""] \\
      U \ar[ru, "v_0"] \ar[r, "v", near end] \ar[rd, "v_1"'] &
      \Dwn_{b_0 i_{\abs X}} \wedge \Dwn_{b_1 i_{\abs X}} \ar[u, "u_0"'] \ar[d, "u_1"] \ar[r, >->, ""] &
      U \times X \\
      & \Dwn_{b_1 i_{\abs X}} \ar[ru, >->, ""]
    \end{tikzcd}
  \]
  We thus conclude from \eqref{eq:uf-init-eq} that
  \( v_0 \tilde f_0 \pi_Y
  = v u_0 \tilde f_0 \pi_Y
  = v u_1 \tilde f_1 \pi_Y
  = v_1 \tilde f_1 \pi_Y \)
  as desired.
\end{proof}

By virtue of \jref{claim:f-well-def} and \jref{propn:epis-eff},
there exists a unique morphism \( f_\circ \colon \oring X \to Y \) with \( \rho f_\circ = \tilde f_\circ \pi_Y \colon \Dwn \to Y \), and thence a unique morphism \( f \colon X \to Y \) with \( i_{\oring X} f = f_\circ \) and \( \rt_X f = \rt_Y \).
We now want to show \( f \) is in initial.
We summarize the situation with the following commutative diagram.
\begin{equation} \label{eq:univ-ext-pf-mega-diag}
  \begin{tikzcd}
    & &[-15pt] &[-15pt] \abs X \times Y \ar[dddd, "\pi_Y"] \\
    \Dwn_{i_{\abs X}} \ar[r, "i_{\Dwn_{\abs{X}}}"', >->] \ar[dd, "\rho"', ->>] \ar[rrru, "\tilde f_\circ", bend left=10pt] & \abs X \times X \ar[dd, "\pi_X"'] \ar[rru, "\id \times f"'] \\[-12pt]
    & & \tm \ar[dl, "\rt_X"] \ar[ddr, "\rt_Y"] \\[-12pt]
    \oring X \ar[r, "i_{\oring X}", >->] \ar[rrrd, "f_\circ"', bend right=10pt] & X \ar[rrd, "f"] \\
    & & & Y
  \end{tikzcd}
\end{equation}

\begin{claim} \label{claim:f-consis}
  \( f \) is consistent.
\end{claim}

We first prove a special case:
\begin{claim} \label{claim:f-consis-global-top}
  The morphism \( i_{\abs X} \colon \abs{X} \to X \) satisfies \( i_{\abs X} f \ein \unex_{\abs X} \rt_Y \).
\end{claim}

\begin{proof}
  Since \( \abs X \le \oring X \), we have a morphism \( \ibar_{\abs X} \colon \abs X \tto \oring X \) with \( \ibar_{\abs X} i_{\oring X} = i_{\abs X} \colon \abs X \to X \).

  By \jref{propn:dwn-hs}, the top element \( \rt_{\Dwn_{i_{\abs X}}} \colon \abs X \to \Dwn_{i_{\abs X}} \) of the hs-object \( \und{\Dwn_{i_{\abs X}}} \in \bC / \abs X \) satisfies \( \rt_{\Dwn_{i_{\abs X}}} i_{\Dwn_{i_{\abs X}}} \pi_X = i_{\abs X} \colon \abs X \to X \) and hence
  \[
    \rt_{\Dwn_{i_{\abs X}}}\ \rho = \ibar_{\abs X} \colon \abs X \to \oring X.
  \]
  By \jref{propn:pb-fun-is-logical}, the top element of \( \und{\abs{X} \times Y} \) in \( \bC/ \abs X  \) is \( \br{\id_{\abs X}, \unex \rt_Y} \).
  Since \( \tilde f_\circ \) is elemental, it follows that \( \rt_{\Dwn_{i_{\abs X}}} \tilde f_\circ \ein \br{\id_{\abs X}, \unex \rt_Y} \) and hence from
  relative consistency of \( \pi_Y \colon \abs{X} \times Y \to Y \) (\jref{propn:ein-pb-conds}), and referring to \eqref{eq:univ-ext-pf-mega-diag}, we obtain, as desired:
  \[
    i_{\abs X} f =
    \ibar_{\abs X} f_\circ =
    \rt_{\Dwn_{i_{\abs X}}} \rho f_\circ =
    \rt_{\Dwn_{i_{\abs X}}} \tilde f_\circ \pi_Y \ein
    \br{\id_{\abs X}, \unex \rt_Y} \pi_Y =
    \unex \rt_Y.
    \qedhere
  \]
\end{proof}

\begin{proof}[Proof of \jref{claim:f-consis}] \label{proof:f-consis}
  Fix \( x,y \colon U \to X \) with \( x \ein y \).
  We need to show \( x f \ein y f \).

  There is a cover of \( U \) by two morphisms \( U' \to U \) such that for each one, we have either \( y' = \unex_{U'} \rt_X \) or \( y' \prec \unex_{U'} \rt_X \).
  To show that \( x f \ein y f \), it suffices to show that \( x' f \ein y' f \) in each of these cases.
  In other words, we have reduced to checking the original claim \( x f \ein y f \) under each of these assumptions.

  When \( y = \unex_U \rt_X \), we have \( x \ein \unex \rt_X \), hence \( x = \bar x i_{\abs X} \) for some \( \bar x \colon U \to \abs X \).
  But we have \( i_{\abs X} f \ein \unex_{\abs X} \rt_Y = \unex_{\abs X} \rt_X f \) by \jref{claim:f-consis-global-top} and hence \( x f = \bar x i_{\abs X} f \ein \bar x \unex_{\abs X} \rt_X f = \unex_U \rt_X f \), as desired.

  When \( y \prec \unex_U \rt_X \), it follows that \( x \prec \unex_U \rt_X \), so that we have \( \bar x,\bar y \colon U \to \oring X \) with \( \bar x i_{\oring X} = x \) and \( \bar y i_{\oring X} = y \).
  By \jref{claim:rho-init}, there exists \( U \oot U' \tox{\tilde x, \tilde y} \Dwn_{i_{\abs X}} \) with \( \tilde x \rho = x \) and \( \tilde y \rho = y \) and \( \tilde x \ein^{\abs X} \tilde y \) (in the notation of \jref{rmk:ein-a-notation}).
  Since \( \tilde f_\circ \) is elemental and in particular consistent, using the relative conservativity of \( \pi_Y \colon \abs{X} \times Y \to Y \), (\jref{propn:ein-pb-conds}) and referring to \eqref{eq:univ-ext-pf-mega-diag}, we thus have, as desired:
  \(
  x f
  = \bar x f_\circ
  = \tilde x \tilde f_\circ \pi_y
  \ein \tilde y \tilde f_\circ \pi_y
  = \bar y f_\circ
  = y f.
  \)
\end{proof}

\begin{claim} \label{claim:f-init}
  \( f \) is initial.
\end{claim}
\begin{proof}[Proof of \jref{claim:f-init}]
  Fix \( x_1 \colon U \to X \) and \( y_0 \colon U \to Y \) with \( y_0 \ein x_1 f \).
  As in the \hyperref[proof:f-consis]{Proof of Claim~\ref*{claim:f-consis}}, we may assume that either \( x_1 = \unex_U \rt_X \) or \( x_1 \prec \unex_U \rt_X \).

  Suppose first that \( x_1 \prec \unex_U \rt_X \), so that we have \( \bar x_1 \colon U \to \oring X \) with \( \bar x_1 i_{\oring X} = x_1 \).
  We can then find \( U \oot U' \tox{w_1} \Dwn_{i_{\abs X}} \) with \( w_1 \rho = \bar x_1 \), and we set \( \br{b, \tilde x_1} = w_1 i_{\Dwn_{i_{\abs X}}} \colon U' \to \abs{X} \times X \).

  Now, we have \( \br{b, y_0} \pi_Y = y_0 \ein x_1 f = \bar x_1 f_\circ = w_1 \tilde f_\circ \pi_Y \) and hence, by the relative conservativity of \( \pi_Y \), that \( \br{b, y_0} \ein^{\abs X} w_1 \tilde f_\circ \).
  Hence, since \( \tilde f_\circ \) is elemental and hence initial, there exists \( U \oot U' \tox{w_0} \Dwn_{i_{\abs X}} \) with \( w_0 \ein^{\abs X} w_1' \) and \( w_0 \tilde f_\circ = \br{b', y_0'} \).
  Setting \( x_0 = w_0 \rho \colon U' \to \oring X \), it follows from the relative consistency of \( \rho \) (~\jref{claim:rho-init}) that \( x_0 \ein x_1' \), and we have, as desired:
  \[
    x_0 f = w_0 \tilde f_\circ \pi_Y = \br{b', y_0'} \pi_Y = y_0'.
  \]

  Finally, consider the case in which \( x_1 = \unex_{U} \rt_X \), so that \( x_1 f = \unex_{U} \rt_Y \) by definition.
  We then have \( y_0 \ein \unex_{U} \rt_Y \) and hence by \jref{propn:dwn-hs}, the top element of the elemental structure \( \und{\Dwn_{y_0}} \tto \und{U \times Y} \) in \( \bC / U \) satisfies
  \[
    \rt_{\Dwn_{y_0}} i_{\Dwn_{y_0}} = \br{\id_U, y_0} \colon U \to U \times Y.
  \]
  Since \( \und{\Dwn_{y_0}} \Ein \und{U \times Y} \), it follows from our assumption (made in the statement of \jref{propn:univ-ext}) that we have an elemental morphism
  \[
    h \colon \und{\Dwn_{y_0}} \to \und{U \times X}.
  \]
  Set \( x_0 = \rt_{\Dwn_{y_0}} h \pi_X \colon U \to X \) so that
  \[
    \br{\id_U, x_0} = \rt_{\Dwn_{y_0}} h \colon U \to U \times X.
  \]
  Since \( h \) is elemental, we then have
  \[
    \br{\id_U, x_0} = \rt_{\Dwn_{y_0}} h \ein \unex \rt_{U \times X} = \br{\id_{U}, \unex \rt_X} = \br{\id_{U}, x_1}
  \]
  and hence \( x_0 \ein x_1 \) by \jref{propn:slice-kj}.
  It remains to see that \( x_0 f = y_0 \).
  By \jref{propn:dwn-hs}, the underobject \( h \colon \und{\Dwn_{y_0}} \tto \und{U \times X} \) is equivalent to \( \und{\Dwn_{x_0}} \tto \und{U \times X} \).
  By \jref{cor:dwn-stable}, we thus have a pullback diagram
  \[
    \begin{tikzcd}
      \Dwn_{y_0} \pb \ar[r, "h", >->] \ar[d] &[10pt] \pb U \times X \ar[r, "\id \times f"] \ar[d, "\bar x_0 \times \id"] &[10pt] U \times Y \ar[d, "\bar x_0 \times \id"] \\
      \Dwn_{i_{\abs X}} \ar[r, "i_{\Dwn_{i_{\abs X}}}", >->] & \abs X \times X \ar[r, "\id \times f"] & \abs X \times Y,
    \end{tikzcd}
  \]
  where \( \bar x_0 \colon U \to \abs{X} \) satisfies \( \bar x_0 i_{\abs X} = x_0 \) (this exists since \( x_0 \ein x_1 = \unex \rt_X \)).
  Since the bottom composite is the elemental morphism \( \tilde f_\circ \colon \und{\Dwn_{i_{\abs X}}} \to \und{\abs X \times Y} \) in \( \bC / \abs X \), it follows from \jref{propn:pb-fun-is-logical} that the top composite is likewise initial in \( \bC / U \), and hence equal to \( i_{\Dwn_{y_0}} \) by \jref{propn:init-uniq}.
  We thus have, as desired:
  \[
    \br{\id_U, x_0 f}
    = \br{\id_U, x_0} (\id_X \times f)
    = \rt_{\Dwn_{y_0}} h (\id_X \times f)
    = \rt_{\Dwn_{y_0}} i_{\Dwn_{y_0}}
    = \br{\id_U, y_0}.
    \qedhere
  \]
\end{proof}

We have thus produced an initial morphism \( f \colon X \to Y \).
By the symmetry of the situation, we likewise obtain an initial \( g \colon Y \to X \).
By \jref{propn:init-2-of-3}, \( f g \) and \( g f \) are initial, hence equal to the identity by \jref{propn:init-uniq}, and thus \( f \) is an isomorphism of hs-objects.
This completes the proof of \jref{propn:univ-ext}.

\subsection{Proof of \alttext{\hyperref[propn:univ-wf]{Proposition~\ref*{propn:univ-wf}}}{Proposition~\ref*{propn:univ-wf}}: well-foundedness of the universe} \label{subsec:univ-wf-proof}
We recall the setup:
we are given a topos \( \bC \), a set \( \fP_A \) of hs-objects in \( \bC / A \) for each \( A \in \bC \), and are supposing that they satisfy the properties \ref{item:wf-fp-props-first}-\ref{item:wf-fp-props-last} listed above \jref{propn:univ-wf}.
We must show that \( \fP \defeq \fP_\tm \) consists of all hs-objects in \( \bC \).

The set-theoretic argument:
\begin{settybox}
  Here the claim is: if \( P \) is a class of \( \unein \)-sets which is closed under isomorphism and inductive with respect to \( \unEin \) (in the sense that \( X \in P \) whenever \( W \in P \) for all \( W \Ein X \)), then \( P \) is the class of \emph{all} \( \unein \)-sets.

  To prove this, we take an arbitrary \( \unein \)-set \( X \) and prove by induction on \( x \in X \) that \( \dwnfn(x) \in P \) (and hence that \( X = \dwnfn(\rt_X) \in P \), as desired).
  Thus, fix \( x \in X \) and suppose \( \dwnfn(w) \in P \) for all \( w \ein x \).
  To show that \( \dwnfn(x) \in P \), it suffices by assumption to show that \( W \in P \) for all \( W \Ein \dwnfn(x) \).
  Given such a \( W \), we have by \jref{propn:dwn-hs} that \( W \cong \dwnfn(w) \) for some \( w \ein x \), and hence that \( W \in P \) by the inductive hypothesis and isomorphism invariance of \( P \).
\end{settybox}

Now, let us take an hs-object \( X \) in \( \bC \).
We need to show it is in \( \fP \).
We need to imitate the induction in the above set-theoretic proof, which involves the consideration of the family of sets \( \dwnfn(x) \) for each \( x \in X \), which we do by the same device used in the beginning of \sref{subsec:univ-ext-proof}.

We consider the hs-substructure \( \und{\Dwn_{\id_X}} \tto \und{X \times X} \) in \( \bC / X \).
(By \jref{lem:dwn-is-descendent}, this is the same as the descendent relation \( {\succeq} \tto X \times X \), though we will not need this.).

Now, by our assumption \ref{item:wf-fp-props-rep} on \pref{item:wf-fp-props-rep}, it follows that there is an underobject \( X_{\Dwn} \tto X \) with the property that, for \( f \colon U \to X \), we have that \( f^* \und{\Dwn_{\id_X}} \in \fP_U \) if and only if \( f \) factors through \( X_{\Dwn} \).
We also note that by \jref{cor:dwn-stable}, we have \( f^* \und{\Dwn_{\id_X}} = \und{\Dwn_f} \in \bC / U \), so we can equivalently express this as \( \und{\Dwn_{f}} \in \fP_U \) if and only if \( f \) factors through \( X_{\Dwn} \).

The induction from the above set-theoretic proof will now take the form as an induction over \( X_{\Dwn} \) showing that \( X_\Dwn \approx X \).
This will finish the proof, since from \( X_\Dwn \approx X \), it follows by definition that \( \und{\Dwn_{f}} \in \fP_U \) for all \( f \colon U \to X \) and, taking \( U = \tm \) and \( f = \rt_X \colon \tm \to X \), we conclude \( X = \Dwn_{\rt_X} \in \fP_\tm \), as desired.

It remains to show that the underobject \( X_\Dwn \tto X \) is inductive.
Let \( x \colon U \to X \) be such that \( X_\Dwn w \) for all \( w \colon U' \to X \) with \( w \ein x' \).
We want to show \( X_\Dwn x \).
Thus, we are assuming that \( \und{\Dwn_{w}} \in \fP_{U'} \) whenever \( w \ein x' \), and want to show \( \und{\Dwn_{x}} \in \fP_U \).

By our assumption \ref{item:wf-fp-props-inductive}, to show that \( \und{\Dwn_{x}} \in \fP_U \), it suffices to show \( \und W \in \fP_{U'} \) for all \( \delta_{U'} \colon U' \to U \) and all \( \und W \in \hs(\bC / {U'}) \) with \( \und W \Ein \delta_{U'}^* \und{\Dwn_{x}} = \und{\Dwn_{x'}} \).
Fix such a \( U' \) and \( \und W \), and write \( i_W \colon \und W \to \und{\Dwn_{x'}} \) for the unique elemental morphism.

The composite
\[
  \und W \tox{i_W} \und{\Dwn_{x'}} \tox{i_{\Dwn_{x'}}} \und{U' \times X}
\]
is initial by \jref{propn:init-2-of-3}, and hence by \jref{propn:dwn-hs}, we have an equivalence \( \und W \approx \und{\Dwn_w} \) for the global element \( \br{\id_{U'}, w} = {\rt_W i_W i_{\Dwn_{x'}}} \colon \und{U'} \to \und{U' \times X} \).
In order to invoke our induction hypothesis and draw the desired conclusion \( W \in \fP_{U'} \), it thus suffices to prove \( w \ein x' \).
But since \( \und W \Ein \und{\Dwn_{x'}} \), this follows from \jref{propn:dwn-respects-ein}.

\subsection{Proof of \alttext{\hyperref[propn:too-universal]{Proposition~\ref*{propn:too-universal}}}{Proposition~\ref*{propn:too-universal}}: the hs-set containing all hs-sets} \label{subsec:too-universal-proofs}
We are given an hs-topoculized topos \( (\bC, \fS, \cV) \), and must show that the top extension \( \wh \cV \) satisfies \( X \Ein \wh \cV \) for every small hs-object \( X \) in \( \bC \).
We begin with the set-theoretic argument:
\begin{settybox}
  Here, the statement we want to prove is the following: let \( (\rV_{\hs}, \unein_{\rV_{\hs}}) \) be the (large) \( {\unein} \)-set of all isomorphism-classes of (small) hs-sets as in \sref{subsec:ein-sets}; this is extensional and well-founded by the set-theoretic versions of the results in \sref{subsec:univ-ext-wf}.
  Let \( \wh \rV_\hs \) be the top-extension of \( \rV_\hs \), so that \( \wh \rV_\hs = \rV_\hs \cup \set{t_{\wh \rV_\hs}} \) which is thus a hs-set by \jref{propn:top-ext-is-hs-struct}.
  Then \( X \Ein \rV_{\hs} \) for every (small) \( \unein \)-set \( X \).

  Fix a small hs-set \( X \).
  To find an elemental map \( X \to \wh \rV_\hs \), it suffices to find an initial morphism \( f \colon X \to \rV_{\hs} \subset \wh \rV_\hs \), since the condition \( f(\rt_X) \in \rt_{\rV} \) is then automatic.
  We claim that \( f \) given by \( f(x) = [\dwnfn(x)] \) is initial, where we write \( [X] \) for the equivalence class of the hs-set \( X \).
  It is consistent since, if \( x \ein y \), then \( \dwnfn(x) \Ein \dwnfn(y) \) by \jref{propn:dwn-respects-ein}.
  And it is initial: if \( [W] \ein f(x) \), so that \( W \Ein \dwnfn(x) \), then \( W \cong \dwnfn(w) \) for some \( w \ein x \) by \jref{propn:dwn-hs}, and we thus have \( f(w) = [W] \), as desired.
\end{settybox}

Fix a small hs-object \( X \) in \( \bC \).
By \jref{propn:top-ext-props}, we have that \( i_{\cV} \colon \cV \to \wh \cV \) is initial and satisfies \( i_{\cV} \ein \unex_{\cV} \rt_{\wh \cV} \).
Since initial morphisms are closed under composition (\jref{propn:init-2-of-3}), it thus suffices to find an initial morphism \( f \colon X \to \cV \).

Consider the hs-substructure \( \und{\Dwn_{\id_X}} \tto \und{X \times X} \) in \( \bC / X \) from \jref{defn:down-family}; it is again small by closure of small morphisms under pullback and composition, and since monomorphisms are small.
By the universal property of \( \cV \), there is thus a (unique) morphism \( f \colon X \to \cV \) in \( \bC \) with \( f^* \und{X_\fS} \cong \und{\Dwn_{\id_X}} \), where \( \und{X_\fS} \in \fS_\cV \) is the universal small hs-object associated to the hs-classifier \( \cV \).
We claim that \( f \) is initial, which will prove the claim.

Note that for \( x \colon U \to X \) in \( \bC \), the morphism \( x f \colon U \to \cV \) is the one classifying \( x^* f^* \und{X_\fS} \cong x^* \und{\Dwn_{\id_X}} \cong \und{\Dwn_x} \in \hs(\fS_U) \) by \jref{cor:dwn-stable}.

We now show consistency of \( f \).
Fix \( x,y \colon U \to X \) in \( \bC \) with \( x \ein y \).
We must prove \( x f \ein y f \), i.e., that \( \und{\Dwn_x} \Ein \und{\Dwn_y} \) in \( \fS_U \); but this follows directly from \jref{propn:dwn-respects-ein}.

Finally, we prove that \( f \) is initial.
Fix \( x_1 \colon U \to X \) in \( \bC \) and \( y_0 \colon U \to \cV \) with \( y_0 \ein x_1 f \).
We need to find \( U \oot U' \tox{x_0} X \) with \( x_0 f = y_0 \) and \( x_0 \ein x_1 \); in fact, we will be able to take \( U' = U \).
Thus, we have an hs-object \( \und{Y_0} = y_0^* \und{X_{\fS}} \in \fS_U \) with \( \und{Y_0} \Ein \und{\Dwn_{x_1}} \), so that we have an initial morphism \( h \colon \und{Y_0} \to \und{\Dwn_{x_1}} \).
Defining \( x_0 \colon U \to X \) by \( \br{\id_U, x_0} = \rt_{Y_0} h i_{\Dwn_{x_1}} \colon (U, \id) \to (U \times X, \pi_U) \), it follows from \jref{propn:dwn-hs} that \( \und{Y_0} \approx \und{\Dwn_{x_0}} \), and hence that \( x_0 f = y_0 \), and by \jref{propn:dwn-respects-ein}, we have \( x_0 \ein x_1 \), as desired.

\subsection{Proofs for \alttext{\hyperref[subsec:ewf-collapse]{Section~\ref*{subsec:ewf-collapse}}}{Section~\ref*{subsec:ewf-collapse}}: extensional collapse of well-founded \alttext{\( \unein \)}{ε}-objects} \label{subsec:ewf-collapse-proof}
In this section, we give the proofs of \jrefs{propn:collapse-epi-iff-univ}{propn:ewf-collapse-precise} and \jrefss{lem:dwn-iso-rel-exists}{lem:r-is-bisimulation}{lem:if-ext-r-is-equality} from \sref{subsec:ewf-collapse}.

\begin{proof}[Proof of \jref{propn:collapse-epi-iff-univ}] \label{proof:collapse-epi-iff-univ}
  We are given an initial morphism \( q \colon X \to Y \) from a well-founded \( \ein \)-object \( X \) to an extensional \( \ein \)-object \( Y \), and must show that \( q \) is epic if and only if it is universal among initial morphisms from \( X \) to extensional \( \ein \)-objects.

  Suppose first that \( q \) has the stated universality property, and take the image \( X \toox{\bar q} \im(q) \tto Y \) of \( q \); we must show \( \im(q) \approx Y \).
  By \jref{propn:dwn-dwn}, \( \im(q) \tto Y \) is a closed substructure of \( Y \), and is hence extensional by \jref{propn:sub-ewf}.
  It follows from \jrefs{propn:closed-iff-init}{propn:init-2-of-3} that \( \bar q \) is initial, and hence, by assumption, must factor through \( q \) via a morphism \( r \colon Y \to \im(q) \).
  Since \( q \) is epic and \( q r i_{\im(q)} = \bar q i_{\im(q)} = q \), it follows that \( r i_{\im(q)} = \id_Y \), and hence that \( \im(q) \approx Y \), as desired.

  Next, we prove the second implication: suppose that \( q \) is epic, and fix \( f \colon X \to Z \) initial with \( Z \) extensional.
  We wish to find a morphism \( g \colon Y \to Z \) with \( q g = f \); the uniqueness of \( g \) is clear since \( q \) is epi, and that \( g \) is initial then follows from \jref{propn:init-2-of-3}.
  By \jref{propn:epis-eff}, it suffices to show that, if we form the pullback \( X \times_Y X \), then \( \pi_0 f = \pi_1 f \colon X \to Z \), i.e., that for \( x_0, x_1 \colon U \to X \), if \( x_0 q = x_1 q \), then \( x_0 f = x_1 f \).

  The proof is by induction: let \( P \tto X \) be given by
  \[
    P = \set{x \in X \mid \forall y \tolon X.\ q(x) = q(y) \to f(x) = f(y)} \tto X.
  \]
  We must show that \( P \) is inductive.
  Thus, fix \( x \colon U \to X \) and suppose that \( P w \) for all \( w \colon U' \to X \) with \( w \ein x \), and fix \( y \colon U' \to X \) with \( x' q = y q \); we need to show \( x' f = y f \).

  By extensionality of \( Z \), it suffices to show that \( w \ein x'' f \ToT w \ein y'' f \) for all \( w \colon U'' \to Z \).
  Assuming \( w \ein x'' f \), then since \( f \) is initial, we have \( U'' \oot U''' \tox{v} X \) with \( v f = w''' \) and \( v \ein x''' \).
  Thus, \( v q \ein x''' q = y''' q \) by consistency of \( q \) and assumption.
  Hence, by initiality of \( q \), there is \( U''' \oot U'''' \tox{u} X \) with \( u \ein y'''' \) and \( u q = v'''' q \).
  By the induction hypothesis applied to \( v \) and consistency of \( f \), it follows that \( w'''' = v'''' f = u f \ein y'''' f \), and hence that \( w \ein y'' f \), as desired, by epicity of \( \delta_{U''''} \delta_{U'''} \).
  The converse \( w \ein y'' f \To w \ein x'' f \) is proven by a similar argument.

  Next, we show that the resulting morphism \( g \colon Y \to Z \) with \( g q = f \) is consistent.
  Let \( y_0, y_1 \colon U \to Y \) with \( y_0 \ein y_1 \).
  By epicity and initiality of \( q \), we then have \( U \oot U' \tox{x_0,x_1} X \) with \( x_0 \ein x_1 \) and \( x_0 q = y_0'  \) and \( x_1 q = y_1' \).
  Hence \( y_0' g = x_0 q g = x_0 f \ein x_1 f = x_1 q g = y_1' g \) by initiality of \( f \), and hence \( y_0 g \ein y_1 g \), as desired, by epicity of \( \delta_{U'} \).

  Finally, we show that \( g \) is initial.
  Let \( y_1 \colon U \to Y \) and \( z_0 \colon U \to Z \) with \( z_0 \ein y_1 g \).
  By epicity of \( q \) and initiality of \( f \), we then have \( U \oot U' \tox{x_0,x_1} X \) with \( x_0 \ein x_1 \) and \( x_0 f = z_0' \) and \( x_1 q = y_1' \).
  Hence \( y_0 \defeq x_0 q \) satisfies \( y_0 = x_0 q \ein x_1 q = y_1' \) and \( y_0 g = x_0 q g = x_0 f = z_0' \), as desired.
\end{proof}

\begin{proof}[Proof of \jref{lem:dwn-iso-rel-exists}] \label{proof:dwn-iso-rel-exists}
  We are given a well-founded, locally extensional \( \ein \)-object \( X \) in a topos \( \bC \), and must produce a relation \( R_X \tto X \times X \) with \( x \Rrel_X y \) iff \( \und{\Dwn_x} \cong \und{\Dwn_y} \) for \( x,y \colon U \to X \).

  To begin with, fix more generally any \( \ein \)-objects \( A \) and \( B \) in \( \bC \).
  Consider the following formulas\footnote{%
      We have not introduced the pairing notation \( (a,b) \) in our logical syntax in \sref{subsec:logic-in-topoi}; rather, \( (a, b) \) is simply the application of the function symbol in the internal language of \( \bC \) given by the identity morphism \( A \times B \to A \times B \). %

      The unique existential quantifier \( \exists \unex \) is of course an abbreviation for the conjunction of the existence and uniqueness formulas as usual.
    } with free variables \( f \tolon \pow(A \times B) \), \( p \tolon \pow  A \), and \( q \tolon \pow B \):
\begin{center}
  \begin{tabular}{ll}
    ``\( f \) has domain in \( p \)''
    & \( \forall a \tolon A.\ \forall b \tolon B.\ (a,b) \in f \to a \in p \)
    \\
    ``\( f \) has codomain in \( q \)''
    & \( \forall a \tolon A.\ \forall b \tolon B.\ (a,b) \in f \to b \in q \)
    \\
    ``\( f \) is functional on \( p \)''
    & \( \forall a \tolon A.\ a \in p \to \exists\unex b \tolon B.\ (a,b) \in f \)
    \\
    ``\( f \) is cofunctional on \( q \)''
    & \( \forall b \tolon B.\ b \in q \to \exists\unex a \tolon A.\ (a,b) \in f \)
    \\
    ``\( f \) is a partial bijection from \( p \) to \( q \)''
    & The conjunction of the previous four.
    \\
    ``\( f \) respects \( \unein \)''
    & \makecell[lt]{\( \forall a_0,a_1 \tolon A.\ \forall b_0,b_1 \tolon B.\) \\
      \( \pbig{(a_0,b_0) \in f \wedge (a_1,b_1) \in f} \to (a_0 \ein a_1 \tot b_0 \ein b_1) \)}
    \\
    ``\( f \) is a partial isomorphism from \( p \) to \( q \)''
    & The conjunction of the previous two.
  \end{tabular}
\end{center}
Now we set
\[
  R_X = \set{\br{x,y} \mid \exists f \tolon \pow(X \times X).\ f\text{ is a partial isomorphism from } \dwnfn(x) \text{ to } \dwnfn(y)} \tto X \times X.
\]
Now, we first establish our claim in the case of global elements \( x, y \colon \tm \to X \).
Returning to the case of general \( \ein \)-objects \( A,B \in \bC \) as above, this comes down the following claim: given global elements \( p \colon \tm \to \pow A \) and \( q \colon \tm \to \pow B \) classifying underobjects \( P \tto A \) and \( Q \tto B \), and given \( f \colon \tm \to \pow(X \times X) \) classifying \( F \tto A \times B \), then \( f \) is a partial isomorphism from \( p \) to \( q \) if and only if \( F \) factors through \( P \times Q \tto A \times B \), and \( F \tto P \times Q \) is the graph of an isomorphism of \( \unein \)-objects \( P \to Q \).

That this is so is just a matter of inspecting the defintions: \( F \) factoring through \( P \times Q \) is equivalent to \( f \) having domain in \( p \) and codomain in \( q \), and \( F \) being the graph of an isomorphism (of objects) \( P \to Q \)---which just means that the composites of \( F \tto P \times Q \) with the two projections are isomorphisms---is then equivalent to \( f \) being functional on \( p \) and cofunctional on \( q \) by \jref{propn:ex-uniq}; and finally, \( f \) respecting \( \ein \) is then equivalent to this isomorphism \( P \to Q \) being an isomorphism of \( \ein \)-objects.

Next, to pass from global elements to general elements, we need the following claim:
 for any \( x,y \colon U \to X \) and \( f,g \colon U \to \pow(X \times X) \), if \( f,g \) are partial isos from \( x \dwn \) to \( y \dwn \), then \( f = g \).

 Indeed, given such \( x \), \( y \), \( f \), and \( g \), we have by \jref{propn:slice-kj} that \( \br{\id,f},\br{\id,g} \colon \und{U} \to \und{U \times \pow(X \times X)} \) are partial isos between \( \br{\id, x} \dwn,\br{\id, y} \dwn \colon \und{U} \to \und{U \times \pow X} \).
  Moreover, \( \und{U \times X} \) is still locally extensional by \jref{propn:pb-fun-is-logical} and the characterization of local extensionality given in \jref{lem:loc-ext-char}.
  Since \( \br{\id,f} = \br{\id,g} \) implies \( f = g \), we have thus reduced the claim to the case where \( U \) is terminal, so that \( x,y,f,g \) are global elements.

  In this case, we have underobjects \( \Dwn_x,\Dwn_y \tto X \) and, by what we just proved above in the global case, \( f \) and \( g \) each classifies the graph of an isomorphism \( \Dwn_x \to \Dwn_y \).
  By local extensionality and well-foundedness of \( X \) and \jref{propn:sub-ewf}, \( \Dwn_x \) and \( \Dwn_y \) are extensional and well-founded, hence by \jref{propn:init-uniq}, there is at most one isomorphism \( \Dwn_x \to \Dwn_y \), hence \( f = g \), as required.

  Finally, we prove the proposition in the general case.
  Fix \( x,y \colon U \to X \).
  By definition, we have \( x \Rrel_X y \) if and only if there is \( U \oot U' \tox{f} \pow(X \times X) \) which is a partial isomorphism from \( x' \dwn \) to \( y' \dwn \); and by \jref{propn:ex-uniq} and the uniqueness of \( f \) just proven, we can take \( U' = U \).

  Now, by \jref{propn:slice-kj}, \( f \colon U \to \pow(X \times X) \) being a partial isomorphism from \( x \dwn \) to \( y \dwn \) is equivalent to \( \br{\id_{U},f} \colon \und{U} \to \und{U \times X \times X} \) being a partial isomorphism from
  \[
    \br{\id_{U}, x \dwn_X}
    =
    \br{\id_{U}, x} \dwn_{\und{U \times X}}
    \quad\text{ to }\quad
    \br{\id_{U}, y \dwn_X}
    =
    \br{\id_{U}, y} \dwn_{\und{U \times X}}
    \colon \und{U} \to \und{U \times \pow X}.
  \]
  And by \jref{propn:closure-in-slice}, we have that \( \und{\Dwn_x} \cong \und{\Dwn_y} \) iff \( \und{\Dwn_{\br{\id,x}}} \cong \und{\Dwn_{\br{\id,y}}} \).
  Thus, we have reduced the claim to the case in which \( U \) is terminal---i.e., in which \( x \) and \( y \) are global elements---which we have already treated above.
\end{proof}

\begin{proof}[Proof of \jref{lem:r-is-bisimulation}] \label{proof:r-is-bisimulation}
  We are given a well-founded, locally extensional object \( X \) in a topos \( \bC \), and \( x_0,x_1,y_1 \colon U \to X \) with \( x_0 \ein x_1 \Rrel_X y_1 \), and we must find \( U \oot U' \tox{y_0} X \) with \( x_0 \Rrel y_0 \ein y_1 \); in fact, we will be able to take \( U' = U \).
  Our assumption \( x_1 \Rrel y_1 \) says that \( \und{\Dwn_{x_1}} \cong \und{\Dwn_{y_1}} \).
  By \jref{propn:dwn-respects-ein}, the assumption \( x_0 \ein x_1 \) gives \( \und{\Dwn_{x_0}} \Ein \und{\Dwn_{x_1}} \), and hence \( \und{\Dwn_{x_0}} \Ein \und{\Dwn_{y_1}} \).
  Let \( f \colon \und{\Dwn_{x_0}} \to \und{\Dwn_{y_1}} \) be the unique elemental morphism, and set \( \br{\id_U, y_0} = \rt_{\Dwn_{x_0}} f i_{\Dwn_{y_1}} \colon \und U \to \und{U \times X} \).
  Then since \( f \) is elemental and \( i_{\Dwn_{y_1}} \) is consistent, and using \jref{propn:dwn-hs}, we have
  \[
    \br{\id, y_0} = \rt_{\Dwn_{x_0}} f i_{\Dwn_{y_1}} \ein \rt_{\Dwn_{y_1}} i_{\Dwn_{y_1}} = \br{\id, y_1},
  \]
  and hence \( y_0 \ein y_1 \) by \jref{propn:slice-kj}.
  Moreover, since by \jrefs{propn:init-cons-mono}{propn:init-2-of-3}, \( f i_{\Dwn_{y_1}} \colon \und{\Dwn_{x_0}} \to \und{U \times X} \) is initial and monic, \jrefs{propn:closed-iff-init}{propn:dwn-hs} give \( \und{\Dwn_{x_0}} \cong \und{\Dwn_{y_0}} \), hence \( x_0 \Rrel y_0 \), as desired.
\end{proof}

\begin{proof}[Proof of \jref{lem:if-ext-r-is-equality}] \label{proof:if-ext-r-is-equality}
  We are given an extensional, well-founded object \( X \) in a topos \( \bC \), and must show that \( R_X \tto X \times X \) is equivalent the equality relation \( \Delta \colon X \to X \times X \).
  Let \( P \tto X \) be the underobject given by
  \[
    P = \set{ x \in X \mid \forall y \colon X.\ x \Rrel_X y \to x = y }.
  \]
  We need to show that \( P \) is inductive.
  Thus, fix \( x \colon U \to X \), suppose that \( P w \) for all \( w \colon U' \to X \) with \( w \ein x' \), and suppose \( x' \Rrel_X y \) for some \( y \colon U' \to X \); we need to show \( x' = y \).
  By extensionality, it suffices to show \( w \ein x'' \) iff \( w \ein y'' \) for all \( w \colon U'' \to X \).
  By \jref{lem:r-is-bisimulation}, for each \( w \colon U'' \to X \) with \( w \ein x'' \), there is \( U'' \oot U''' \tox{v} X \) with \( w''' \Rrel_Z v \ein y''' \), and hence \( w''' = v \) by the induction hypothesis applied to \( w \).
  Hence \( w''' \ein y''' \) and hence \( w'' \ein y'' \), as desired, by epicity of \( \delta_{U'''} \).
  The converse \( w \ein y'' \To w \ein z'' \) is proven by a similar argument.

  Now suppose \( X \) is only locally extensional and well-founded, and fix \( x,y,z \colon U \to X \) with \( x,y \preceq z \) with \( x \Rrel_X y \); again, we must show \( x = y \).
  We have the induced global elements \( \br{\id, x}, \br{\id, y}, \br{\id, z} \colon \und U \to \und{U \times X} \), and we have that \( \br{\id, x} \Rrel_{U \times X} \br{\id, y} \) (since \( \und{\Dwn_{\br{\id, x}}} = \und{\Dwn_{x}} \) by \jref{propn:closure-in-slice} and likewise with \( y \)), and, by \jref{propn:slice-kj}, we have \( \br{\id, x}, \br{\id, y} \preceq \br{\id, z} \).
  Moreover, if \( \br{\id, x} = \br{\id, y} \), then \( x = y \).
  Hence, we have reduced to proving the claim for global elements, i.e., when \( U \) is terminal.

  In this case, \( x \) and \( y \) factor as morphisms \( \tm \tox{\bar x, \bar y} \Dwn_z \tto X \).
  Moreover, considering \( \Dwn_{\bar x} \tto \Dwn_z \) as an underobject of \( X \), it is equivalent to \( \Dwn_x \), the latter being least closed containing \( x \), and likewise with \( y \).
  Hence, our assumption \( x \Rrel_X y \) gives \( \Dwn_{\bar x} \cong \Dwn_{\bar y} \) and hence, since \( \Dwn_z \) is extensional by assumption and well-founded by \jref{propn:sub-ewf}, it follows that \( \bar x = \bar y \) from the first part of the proposition proved above, and hence that \( x = y \), as desired.
\end{proof}

\begin{proof}[Proof of \jref{propn:ewf-collapse-precise}] \label{proof:ewf-collapse-precise}
We are given a well-founded, locally extensional \( \unein \)-object \( X \) in a topos \( \bC \), we form the quotient \( q \colon X \to Y \) by the equivalence relation \( R_X \tto X \times X \), and we define \( \unein_Y \tto Y \times Y \) as the image of \( \unein_X \tto X \times X \tox{q \times q} Y \times Y \).
Thus, for \( y_0,y_1 \colon U \to Y \), we have \( y_0 \ein y_1 \) iff there is \( U \oot U' \tox{x_0, x_1} X \) with \( x_0 \ein x_1 \) and \( y_i' = x_i q \) for \( i = 0,1 \).
Also, by effectiveness of equivalence relations in topoi (see \cite[Proposition~1.23]{johnstone-topos-theory}), we have that \( R_X \) is equivalent to the kernel pair of \( q \), and hence that for \( x_0, x_1 \colon U \to X \), we have \( x_0 \Rrel_X x_1 \) iff \( x_0 q = x_1 q \).

We must prove that \( Y \) is an extensional collapse---i.e., that \( Y \) is extensional and \( q \) is initial---and that \( Y \) is well-founded.
Note that the consistency of \( q \) is immediate from the definition of \( \unein_Y \).

\textit{Initiality of \( q \)}:
given \( x_1 \colon U \to X \) and \( y_0 \colon U \to Y \) with \( y_0 \ein x_1 q \), we wish to find \( U \oot U' \tox{x_0} X \) with \( x_0 \ein x_1' \) and \( x_0 q = y_0' \).
By definition of \( \unein_Y \), we have \( U \oot U' \tox{\tilde x_0, \tilde x_1} X \) with \( \tilde x_0 q = y_0' \) and \( \tilde x_1 q = x_1' q \) (hence \( \tilde x_1 \Rrel_X x_1' \)) and \( \tilde x_0 \ein \tilde x_1 \).
  Hence, by \jref{lem:r-is-bisimulation}, there exists \( U' \oot U'' \tox{x_0} X \) with \( x_0 \Rrel_X \tilde x_0'' \) (hence \( x_0 q = y_0'' \)) and \( x_0 \ein x_1'' \), as desired.

\textit{Well-foundedness of \( Y \)}:
Let \( P \tto Y \) be inductive.
  Since \( q^* P \tto X \tox{q} Y \) factors through \( P \) and \( q \) is epic, to prove \( P \approx Y \), it suffices to show \( q^* P \approx X \)---i.e., since \( X \) is well-founded, that \( q^* P \) is again inductive.
  Thus, let \( x \colon U \to X \) and suppose \( (q^* P) w \) for all \( w \colon U' \to X \) with \( w \ein x' \); we must show \( (q^* P) x \), i.e., that \( P (x q) \).

  By the assumption that \( P \) is inductive, it suffices to show that \( P y \) for any \( y \colon U' \to Y \) with \( y \in x' q \).
  Fixing such a \( y \), we have by definition of \( \unein_Y \) some \( U' \oot U'' \tox{\tilde w, \tilde x} X \) with \( \tilde w \ein \tilde x \) and \( \tilde w q = y'' \) and \( \tilde x q = x'' q \), hence \( \tilde x \Rrel_X x'' \).
  Hence, by \jref{lem:r-is-bisimulation}, we have \( U'' \oot U''' \tox{w} X \) with \( \tilde w''' \Rrel_X w \ein x''' \), and hence \( (q^* P) w \) by our assumption on \( x \).
  Hence \( P (w q) \), and hence \( P y \), as desired, since \( w q = \tilde w''' q = y''' \) and \( \delta_{U'''} \delta_{U''} \delta_{U'} \) is epic.

\textit{Extensionality of \( Y \):}
fix \( y,z \colon U \to Y \) with \( x \ein y' \ToT x \ein z' \) for all \( x \colon U' \to Y \); we must show \( y = z \).
Thus, fixing \( U \oot U' \tox{\tilde y, \tilde z} X \) with \( \tilde y q = y' \) and \( \tilde z q = z' \), we must prove \( \tilde y \Rrel_X \tilde z \).
Our assumption on \( x \) and \( y \) together with the already established initiality of \( q \) gives that for any \( v \colon U'' \to X \) with \( v \ein \tilde y \), there is \( U'' \oot U''' \tox{w} X \) with \( v''' \Rrel_X w \ein z''' \) (and vice versa).
Moreover, this \( w \) is unique by \jref{lem:if-ext-r-is-equality}, and hence we can take \( U''' = U'' \) by \jref{propn:ex-uniq}.

Now, to show that \( \tilde y \Rrel_X \tilde z \), i.e., that \( \und{\Dwn_{\tilde y}} \cong \und{\Dwn_{\tilde z}} \), it suffices by \jrefs{cor:dwn-stable}{propn:univ-ext} to show that for each \( U'' \to U' \) and \( \und W \in \hs(\bC / U'') \), we have \( \und W \Ein \und{\Dwn_{\tilde y''}} \) iff \( \und W \Ein \und{\Dwn_{\tilde z''}} \).

Fixing such a \( U'' \) and \( W \) with \( \und W \Ein \und{\Dwn_{\tilde y''}} \), and letting \( f \colon \und W \to \und{\Dwn_{\tilde y''}} \) be the unique elemental morphism, we have by \jrefs{propn:init-2-of-3}{propn:dwn-hs} that \( \und{W} = \und{\Dwn_{v}} \), where \( \br{\id_{U''}, v} = \rt_W f i_{\Dwn_{\tilde y''}} \colon \und{U''} \to \und{U'' \times X} \); and from \jref{propn:dwn-respects-ein}, it follows that \( v \ein \tilde y'' \).
Hence, by assumption, there is \( w \colon U'' \to X \) with \( v \Rrel_X w \ein \tilde z'' \), and hence, by \jref{propn:dwn-respects-ein} again, \( \und W = \und{\Dwn_v} \cong \und{\Dwn_w} \Ein \und{\Dwn_{z''}} \), as desired.
And the converse \( \und W \Ein \und{\Dwn_{z''}} \To \und W \Ein \und{\Dwn_{y''}} \) is of course proven the same way.
\end{proof}

\subsection{Proof of \alttext{\hyperref[propn:topos-universe-small-closed-precise]{Proposition~\ref*{propn:topos-universe-small-closed-precise}}}{Proposition~\ref*{propn:topos-universe-small-closed-precise}}: closure of the universe under set-formation} \label{subsec:universe-closure-proof}
We recall the setup:
\begin{itemize}
\item \( (\bC, \fS, \cV) \) is an hs-topoculized topos
\item \( \alpha \colon U \to \pow \cV \) is a morphism in \( \bC \), classifying an underobject \( i_A \colon A \tto U \times \cV \)
\item \( \und W \in \bC / A \) is the small hs-object classified by \( i_A \pi_{\cV} \colon A \to \cV \)
\item \( \rt_W \colon (A, \id) \to \und W \) is the top element of \( \und W \)
\item \( c \colon W \to W_U \) is an \( \unein \)-cocartesian morphism of \( \unein \)-objects over \( i_A \pi_U \)
\item \( q \colon \und W_U \to \und Z \) is an extensional collapse
\item \( i_Z \colon \und Z \to \und{\wh Z} \) is a top extension along \( A \tox{\rt_{W}} W \tox{c} W_U \tox{q} Z \)
\end{itemize}
Thus, \( \und W_U \), \( \und Z \), and \( \und{\wh Z} \) are all \( \unein \)-objects in \( \bC / U \).
Moreover, we are assuming:
\begin{itemize}
\item \( i_A \pi_A \colon A \to U \) is small
\end{itemize}
The data is summarized in the following diagram
\[
    \begin{tikzcd}
      W \ar[rr, "c"] \ar[d, "\delta_W"] & &
      W_U \ar[r, "q"] \ar[d, "\delta_{W_U}"] &
      Z \ar[r, "i_Z"] &
      \wh Z \ar[d, "\delta_{\wh Z}"] \\
      A \ar[r, "i_A"] \ar[u, bend left=40pt, "\rt_W" pos=0.45]
      & U \times \cV \ar[r, "\pi_U"]
      & U \ar[rr, "\id"]
      && U,
    \end{tikzcd}
  \]
What we want to show is:
\begin{enumerate}[(i)]
\item \( \und{\wh Z} \) is an hs-object
\item the morphism \( \delta_{\wh Z} \colon \wh Z \to U \) is small
\item the morphism \( x \colon U \to \cV \) classifying \( \und{\wh Z} \) satisfies \( y \ein x' \) iff \( y \in \alpha' \) for all \( y \colon U' \to \cV \).
\end{enumerate}

It will also be good to keep in mind the analogous set-theoretic situation (\jref{propn:set-v-closed-precise}):
\begin{settybox}
  We have a small subset \( A \subset \rV_{\hs} \) of the class of hs-objects.
  We form \( W = \coprod_{y \in A} \floor{y} \), which we endow with the disjoint union \( \unein \)-set structure.
  We then form an extensional collapse \( q \colon W \to Z \), and a top extension \( \wh Z \) of \( Z \) along \( q\pbig{ \set{(y, \rt_{\floor y}) \mid y \in A} } \subset Z \).

  We then want to show that \( \wh Z \) is an hs-object, and that its equivalence class \( x = [\wh Z] \in \rV_\hs \) satisfies \( A = \set{y \in \rV_\hs \mid y \ein x} \).
\end{settybox}

\begin{claim}
  \( \und W_U \) is locally extensional and well-founded, and \( \und{Z} \) is extensional and well-founded.
\end{claim}
\noindent
This follows from \jref{propn:sum-props} and \jref{propn:ewf-collapse}.

\begin{claim}
  \( \rt_{W} c q \colon (A, i_A \pi_U) \to \und Z \) is a monomorphism.
\end{claim}
Here is the set-theoretic analogue:
\begin{settybox}
  Claim: The restriction of \( q \colon W \to Z \) to \( \set{(y, \rt_{\floor y}) \mid y \in A} \) is injective.

  Proof:
  If \( q(y, \rt_{\floor y}) = q(z, \rt_{\floor z}) \), then by definition of the extensional collapse \( q \colon W \to Z = W / \unsim \) (\jref{propn:ewf-collapse-precise}), this means that \( \dwnfn(y, \rt_{\floor y}) \cong \dwnfn(z, \rt_{\floor z}) \).
  We have \( \dwnfn(y, \rt_{\floor y}) = \set{y} \times \floor{y} \subset W \), as follows from \jref{propn:dwn-dwn}, since the inclusion \( \floor{y} \hto W \) is initial (by the definition of \( \unein_W \)).
  Similarly, \( \dwnfn(x, \rt_{\floor z}) = \set{z} \times \floor{z} \subset W \).
  We thus conclude that \( \floor{y} \cong \floor{z} \), and hence that \( y = z \) and \( (y, \rt_{\floor y}) = (z, \rt_{\floor z}) \), as desired.
\end{settybox}

\begin{proof}
  As most of the objects involved in this proof, namely \( (A, i_A \pi_U), W_U, Z \in \bC / U \) and \( W \in \bC / A \cong (\bC / U) / (A, i_A \pi_U) \) lie in the topos \( \bC / U \), it will be convenient to temporarily change notation: we set \( \bD \defeq \bC / U \), \( B \defeq (A, i_A \pi_U) \in \bD \), \( \und X \defeq \und W \in \bD / B \), \( X_\tm \defeq \und W_U \in \bD \), and \( Y \defeq Z \in \bD \).

  Thus, we have the hs-object \( \und X \) in \( \bD / B \), the \( \ein \)-cocartesian morphism \( c \colon X \to X_\tm \) over \( !_B \), and the extensional collapse \( q \tolon X_\tm \to Y \) in \( \bD \), and we need to show that \( \rt_X c q \colon B \to Y \) is a monomorphism.
  Fix \( u,v \colon U' \to B \) with \( u \rt_X c q = v \rt_X c q \colon U \to Y \).
  By \jref{propn:ewf-collapse-precise}, this means that \( \und{\Dwn_{u \rt_{X} c}} \cong \und{\Dwn_{v \rt_{X} c}} \in \bD / U' \).

  We have the morphism \( \br{\delta_X, c} \colon \und X \to \und{B \times X_\tm} \) in \( \bC / B \), and and applying \jref{propn:rel-notions-compare} twice, we find that \( c \) is relatively initial over \( !_B \), and thence that \( \br{\delta_X,c} \) is initial in \( \bD / B \); and it is also monic since \( c \) is an iso.
  It follows from \jref{propn:dwn-hs} that \( \und X = \und{\Dwn_{\rt_X c}} \).
  \[
    \begin{tikzcd}
      X \ar[dr, "\br{\delta_X, c}", near end] \ar[drr, "c", bend left=15pt] \ar[ddr, "\delta_X"', bend right=10pt] \\
      & B \times X_\tm \pb \ar[d, "\pi_B"] \ar[r, "\pi_{X_\tm}"'] & X_\tm \ar[d] \\
      & B \ar[r, "!_B"] & \tm
    \end{tikzcd}
    \qquad \qquad
    \begin{tikzcd}
      W \ar[dr] \ar[drr, "c", bend left=15pt] \ar[ddr, "\delta_W"', bend right=10pt] \\
      & A \times_U W_U \pb \ar[d] \ar[r] & W_U \ar[d] \\
      & A \ar[r, "i_A \pi_U"] & U
    \end{tikzcd}
  \]
Hence, by \jref{cor:dwn-stable} and our assumption on \( u \) and \( v \), we have
\[
  u^* \und{X} =
  u^* \und{\Dwn_{\rt_X c}} =
  \und{\Dwn_{u \rt_X c}} \cong
  \und{\Dwn_{v \rt_X c}} =
  v^* \und{\Dwn_{\rt_X c}} =
  v^* \und{X}.
\]

We now return to the original notation: thus, we have \( \und{U'} \in \bC / U \) and \( u,v \colon \und{U'} \to \und{A} \), and are trying to show that \( u = v \), and so far, we have shown \( u^* \und{W} \cong v^* \und{W} \in \hs(\bC / U') \).
Recalling that \( \und{W} \in \hs(\bC / A) \) is a small hs-object, classified by \( i_A \pi_\cV \colon A \to \cV \), we conclude that the morphisms \( u i_A \pi_{\cV}, v i_A \pi_{\cV} \colon U' \to \cV \), classifying the small hs-objects \( u^* \und{W} \) and \( v^* \und{W} \) in \( \bC / U' \), are equal.
  Since also \( u i_A \pi_U = \delta_{U'} = v i_A \pi_U \colon U' \to U \), we conclude that \( u i_A = v i_A \colon U' \to U \times \cV \), and hence that \( u = v \), as desired, since \( i_A \) is monic.
\end{proof}

\begin{claim} \label{claim:dense-mono-success}
  The monomorphism \( \rt_{W} c q \colon (A, i_A \pi_U) \to \und Z \) is dense.
\end{claim}

\begin{proof}
  We have that \( \rt_W c \colon (A, i_A \pi_U) \to \und W_U \) is dense by \jref{propn:sum-props}.
  Since \( q \colon \und W_U \to \und Z \) is initial, it follows from \jref{propn:dwn-dwn} that \( \rt_W c q \) is also dense.
\end{proof}

We now come to the main claims (those listed as ``(i)-(iii)'' above):
\begin{claim}
  \( \und{\wh Z} \) is an hs-object.
\end{claim}

\begin{proof}
  This follows from \jref{claim:dense-mono-success} and \jref{propn:top-ext-is-hs-struct}.
\end{proof}

\begin{claim}
  \( \delta_{\wh Z} \colon\wh Z \to U \) is small.
\end{claim}

\begin{proof}
  We assumed that \( i_A \pi_A \colon A \to U \) and \( \delta_W \colon W \to A \) are small.
  It follows from \jref{defn:topoculus}~\ref{item:topoculus-comp} that \( \delta_W i_A \pi_A \colon W \to U \) is small; i.e., that \( (W, \delta_W i_A \pi_A) \in \bC / U \) lies in the sub-topos \( \fS_U \subset \bC / U \).
  We have \( (W, \delta_W i_A \pi_A) \cong \und W_U \) by \jref{propn:ein-cocart-conds}, and hence that \( \und W_U \in \fS_U \) as well by \jref{defn:topoculus}~\ref{item:topoculus-pb}.
  Since \( q \colon \und W_U \to \und Z \) is an epi, it follows from \jref{propn:topoculus-epi-closed} that \( \und Z \in \fS_U \).
  Finally, since \( \und {\wh Z} \cong \und Z + \tm \), it follows that \( \wh Z \in \fS_U \) by the closure of the subtopos \( \fS_U \subset \bC / U \) under finite colimits.
\end{proof}

\begin{claim} \label{claim:rel-elem-crit}
  Let \( x \colon A \to \cV \) be a morphism classify an hs-object \( X \in \bC / A \) and let \( y \colon A' \to \cV \) be a morphism classifying \( Y \in \cV / A' \).
  Then \( y \ein x' \) iff there exists a relatively elemental morphism \( Y \to X \) over \( \delta_{A'} \colon A' \to A \) (see \jref{defn:rel-init-elem}).
\end{claim}

\begin{proof}
  By definition, \( x' \) classifies the pullback hs-object \( \delta_{A'}^* X \).
  Hence, the claim follows from \jref{propn:rel-notions-compare}.
\end{proof}

\begin{claim}
  The morphism \( x \colon U \to \cV \) classifying \( \und{\wh Z} \) satisfies \( y \ein x' \) iff \( y \in \alpha' \) for all \( y \colon U' \to \cV \).
\end{claim}

We begin with the set-theoretic proof:
\begin{settybox}
  First, note that \( y \ein x = [\wh Z] \) is equivalent to \( \floor{y} \Ein \wh Z \).

  Thus, we need to show that \( \floor{y} \Ein \wh Z \) iff \( y \in A \) for all \( y \in \rV_\hs \).
  Now, if \( y \in A \), then the composite
  \[
    \floor{y}
    \toi \set{y} \times \floor{y}
    \hto W
    \tox{q} Z
    \tox{i_Z} \wh Z,
  \]
  is initial by \jrefs{propn:init-2-of-3}{propn:top-ext-props} and by the definition of \( \unein_W \).
  Since \( q(y, \rt_{\floor{y}}) \) is contained in the set along which we formed the top-extension \( \wh Z \), it follows that \( i_Z\pbig{q(y, \rt_{\floor y})} \ein \rt_{\wh Z} \) and hence that \( \floor{y} \Ein \wh Z \), as desired.

  Conversely, suppose \( \floor{y} \Ein \wh Z \), so that \( \floor{y} \cong \dwnfn(w) \) for some \( w \in \wh Z \) with \( w \ein \rt_{\wh Z} \) by \jref{propn:dwn-hs}.
  By the definition of \( \unein_{\wh Z} \), it follows that \( w = i_Z\pbig{q(a, \rt_{\floor a})} \) for some \( a \in A \), and thus by \jref{propn:dwn-hs} and the argument we just gave, \( \dwnfn(w) \cong \floor{a} \).
  Hence \( \floor{y} \cong \floor{a} \), so \( y = a \in A \), as desired.
\end{settybox}

\begin{proof}
  Fix \( y \colon U' \to \cV \), and suppose \( y \in \alpha' \); we want to prove that \( y \ein x' \).
  Since \( \alpha \colon U \to \pow \cV \) classifies the underobject \( i_A \colon A \tto U \times \cV \), it follows by definition that \( \br{\delta_{U'}, y} \colon U' \to U \times \cV \) factors through \( A \tto U \times \cV \) via some morphism \( \tilde y \colon U' \to A \).
  Writing \( \und Y \in \fS_{U'} \subset \bC / U' \) for the hs-object classified by \( y = \tilde y i_A \pi_{\cV} \), it follows that we have a \( \unein \)-pullback square
  \begin{equation} \label{eq:closure-proof-pb-sq}
    \begin{tikzcd}
      Y \pb \ar[r, "i_Y"] \ar[d] & W \ar[d] \\
      U' \ar[r, "\tilde y"] & A
    \end{tikzcd}
\end{equation}
  for some morphism \( i_Y \colon Y \to W \).
  By \jref{claim:rel-elem-crit}, what we want to prove is that we have a relatively elemental morphism
  \begin{equation} \label{eq:closure-proof-rel-elem-sq}
    \begin{tikzcd}
      Y \ar[r, "f"] \ar[d] & \wh Z \ar[d] \\
      U' \ar[r, "\delta_{U'}"] & U.
    \end{tikzcd}
  \end{equation}

  By \jref{cor:rel-init-2-of-3}, the composite
  \[
    \begin{tikzcd}
      Y \ar[r, "i_Y"] \ar[d, "\delta_Y"'] &
      W \ar[r, "c"] \ar[d, "\delta_W"] &
      W_U \ar[r, "q"] \ar[d, "\delta_{W_U}"] &
      Z \ar[r, "i_Z"] &
      \wh Z \ar[d, "\delta_{\wh Z}"] \\
      U' \ar[r, "\tilde y"] &
      A \ar[r, "i_A \pi_U"]
      & U \ar[rr, "\id"]
      && U,
    \end{tikzcd}
  \]
  is relatively initial.
  It remains to see that \( \rt_Y i_Y c q i_Z \ein^U \tilde y i_A \pi_U \rt_{\wh Z} \).
  Now, since \( i_Y \) is \( \ein \)-cartesian, we have by \jref{propn:pb-fun-is-logical} that \( \rt_Y i_Y = \tilde y \rt_W \), and since \( \wh Z \) is the top-extension of \( Z \) along \( \rt_W c q \), we have \( \rt_W c q i_Z \ein^U i_A \pi_U \rt_{\wh Z} \), and hence \( \rt_Y i_Y c q i_Z \ein \tilde y i_A \pi_U \rt_{\wh Z} \), as desired.

  Now suppose, conversely, that \( y \ein x' \), i.e.\ (by \jref{claim:rel-elem-crit}), that we have a relatively elemental morphism \( f \colon Y \to \wh Z \) over \( \delta_{U'} \) as in \eqref{eq:closure-proof-rel-elem-sq}; we wish to show that \( y \in \alpha' \).
  We thus have that \( \rt_Y f \ein^U \delta_{U'} \rt_{\wh Z} \), and hence, by \jref{propn:top-ext-props} and the definition of \( \wh Z \), that \( \rt_Y f \colon  U' \to \wh Z \) factors through \( \rt_W c q i_Z \colon A \tto \wh Z \), via some morphism \( \tilde y \colon U' \to A \), so that
  \begin{equation} \label{eq:image-arg-eq-1}
    \rt_Y f = \tilde y \rt_W c q i_Z.
  \end{equation}

  Now, fix an \( \unein \)-cartesian morphism \( \pi_{i_A \pi_U} \colon (i_A \pi_U)^* \wh Z \to \wh Z \), so that we have an induced morphism \( \ol{c q i_Z} \colon W \to (i_A \pi_U)^* \wh Z \).
  \[
    \begin{tikzcd}
      W \ar[rrd, bend left=15pt, "c q i_Z"] \ar[rdd, bend right] \ar[rd, "\ol{c q i_Z}"', dashed, near end] \\
      & (i_A \pi_U)^* \wh Z \ar[r, "\pi_{i_A \pi_U}"'] \ar[d, ""'] &[10pt] \wh Z \ar[d, ""] \\
      & A \ar[r, "i_A \pi_U"] & U
    \end{tikzcd}
  \]
  Since \( c q i_Z \) is relatively initial, \( \ol{c q i_Z} \) is initial by  \jref{propn:rel-notions-compare}, and hence a mono by \jref{propn:init-cons-mono} (in \( \bC / A \), hence also in \( \bC \)).

  We also have an induced morphism \( \bar f \colon Y \to (i_A \pi_U)^* Z \), which is relatively initial over \( \bar y \) by \jref{cor:rel-init-2-of-3}.
  \[
    \begin{tikzcd}
      Y \ar[rrd, bend left=15pt, "f"] \ar[dd, "\delta_Y"'] \ar[rd, "\bar f"', dashed] \\
      & (i_A \pi_U)^* \wh Z \ar[r, "\pi_{i_A \pi_U}"'] \ar[d, ""'] &[10pt] \wh Z \ar[d, ""] \\
      U' \ar[r, "\bar y"] & A \ar[r, "i_A \pi_U"] & U
    \end{tikzcd}
  \]
  We next claim that \( \bar f \colon Y \to (i_A \pi_U)^* \wh Z \) factors through the monomorphism \( \ol{c q i_Z} \colon W \to (i_A \pi_U)^* \wh Z \) (say, via a morphism \( i_Y \colon Y \to W \)).
  \[
    \begin{tikzcd}
      Y \ar[r, "i_Y", dashed] \ar[dd, "\delta_Y"'] \ar[rd, "\bar f"'] & W \ar[d, "\ol{c q i_Z}"] \\
      & (i_A \pi_U)^* \wh Z \ar[d, ""'] \\
      U' \ar[r, "\bar y"] & A
    \end{tikzcd}
  \]
  This amounts to showing that the image of the morphism \( \bar f \colon \und Y = (Y, \delta_Y \bar y) \to (i_A \pi_U)^* \und {\wh Z} \) in \( \bC / A \) lies inside that of \( \ol{c q i_Z} \colon \und W \to (i_A \pi_U)^* \und{\wh Z} \), i.e., that
  \begin{equation} \label{eq:image-arg-eq-3}
    \top_{\und Y} \exists_{\bar f} \le \top_{\und W} \exists_{\ol{c q i_Z}} \colon (A, \id) \to \pow \pbig{(i_A \pi_U)^* \und{\wh Z}}.
  \end{equation}
  And indeed, since \( \rt_Y \colon \und U' = (U', \bar y) \to \und Y \) is dense by \jref{propn:sum-props}, we have, using \eqref{eq:image-arg-eq-1}~and~\eqref{eq:image-arg-eq-3} as well as \jref{propn:dwn-dwn} (twice):
  \[
    \begin{split}
      \top_{\und Y} \exists_{\bar f}
      & = \top_{\und U'} \exists_{\rt_Y} \dwn \exists_{\bar f}
      = \top_{\und U'} \exists_{\rt_Y} \exists_{\bar f} \dwn
      = \top_{\und U'} \exists_{\rt_Y \bar f} \dwn \\
      & = \top_{\und U'} \exists_{\tilde y \rt_W \ol{c q i_Z}} \dwn
      = \top_{\und U'} \exists_{\tilde y \rt_W} \exists_{\ol{c q i_Z}} \dwn
      = \top_{\und U'} \exists_{\tilde y \rt_W} \dwn \exists_{\ol{c q i_Z}}
      \le \top_{\und W_U} \exists_{\ol{c q i_Z}}.
    \end{split}
  \]

  We next claim that the resulting square \eqref{eq:closure-proof-pb-sq} is an \( \unein \)-pullback square.
  We have that \( i_Y \) is relatively initial by \jref{cor:rel-init-2-of-3} and hence relatively consistent and conservative by \jrefs{propn:rel-notions-compare}{propn:init-cons-mono}.
  Hence, it remains only to prove that \eqref{eq:closure-proof-pb-sq} is a pullback square, or in other words, that the induced morphism \( \ol{i_Y} \colon \und Y \to \tilde y^* \und W \) into the \( \unein \)-pullback of \( W \) along \( \tilde y \) is an isomorphism.
  We know by \jref{propn:init-cons-mono} that it is a mono, since it is initial by \jref{propn:rel-notions-compare}, so we need only see it is an epimorphism.

  This will follow at once from \jref{propn:dwn-dwn} if we can show that \( \rt_Y \ol{i_Y} = \rt_{\tilde y^* W} \), or what amounts to the same (by \jref{propn:pb-fun-is-logical}), that \( \rt_Y i_Y = \tilde y \rt_W \).
  But this follows since \( \rt_Y i_Y \ol{c q i_Z} = \rt_Y \ol f = \tilde y \rt_W \ol{c q i_Z} \) by \eqref{eq:image-arg-eq-1}, and \( \ol{c q i_Z} \) is monic.

  Finally, since \eqref{eq:closure-proof-pb-sq} is an \( \unein \)-pullback square, it follows that \( \tilde y i_A \pi_{\cV} = y \colon U' \to \cV \).
  And we also have \( \tilde y i_A \pi_U = \delta_{U'} \) by post-composing \eqref{eq:image-arg-eq-1} on both sides with \( \delta_{\wh Z} \).
  Hence \( \br{\delta_{U'}, y} \colon U' \to U \times \cV \) factors through \( i_A \colon A \tto U \times \cV \), which means that \( y \in \alpha' \), as desired.
\end{proof}

\subsection{Proof of \alttext{\hyperref[propn:topos-universe-small-closed-converse]{Proposition~\ref*{propn:topos-universe-small-closed-converse}}}{Proposition~\ref*{propn:topos-universe-small-closed-converse}}: converse of closure of the universe under set formation} \label{subsec:proof-of-small-converse}
We are given an hs-topoculized topos \( (\bC, \fS, \cV) \) and a morphism \( x \colon U \to \cV \), and must show that the underobject \( i_A \colon A \tto U \times \cV \) classified by \( x \kappa \colon U \to \pow \cV  \) is such that \( i_A \pi_U \) is small.

We begin with the set-theoretic argument, i.e., the proof of \jref{propn:set-v-closed-under-sets}~\ref{item:set-v-closed-converse}:
\begin{settybox}
  We need to show that for any \( x \in \rV_{\hs} \), the class \( A = \set{y \in \rV_{\hs} \mid y \ein x} \) is a set.

  Now by definition, for each such \( y \), we have an elemental morphism \( i_y \colon \floor{y} \to \floor{x} \).
  By \jref{propn:dwn-hs}, the map \( y \mapsto i_y(\rt_{\floor{y}}) \) is a bijection \( A \to \set{w \in \floor{x} \mid w \ein \rt_{\floor{x}}} \), and hence \( A \) is small since the codomain is a small.
\end{settybox}

We first treat the ``universal'' case \( x = \id_{\cV} \) (i.e., we first prove \jref{cor:small-closed-converse-univ-case}).
That is, we show that the morphism \( {\ein_{\cV}} \tox{i_{\ein_\cV}} \cV \times \cV \tox{\pi_1} \cV \) is small.

Writing \( i_j = i_{\unein_{\cV}} \pi_j \colon \unein_{\cV} \to \cV \) for \( j = 0,1 \), we have the pullback \( \unein \)-objects \( i_0^* \und X_\fS, i_1^* \und X_\fS \in \bC / \unein_{\cV} \) (in the notation of \jref{defn:topoculus}), and since \( i_0 \ein_\cV i_1 \), it follows from \( (\cV, \unein_\cV) \) being an hs-classifier that we have an elemental morphism \( f \colon i_0^* X \to i_1^* X \).
\[
  \begin{tikzcd}
    i_0^* X_\fS \ar[r, "f"] \ar[dr] & i_1^* X_\fS \ar[rr, "\pi"] \ar[d] \pb & & X_\fS \ar[d, "\delta_{X_\fS}"] \\
    {} & \ein_{\cV} \ar[r, "i_{\ein_{\cV}}"] & \cV \times \cV \ar[r, "\pi_1"] & \cV
  \end{tikzcd}
\]
We thus have a morphism \( g = \rt_{i_0^* X_\fS} f \pi \colon \unein_{\cV} \to X_{\fS} \).

\begin{claim}
  The morphism \( g \) is a monomorphism.
\end{claim}

\begin{proof}
  We must show that any two morphisms \( U \to {\ein_{\cV}} \) with the same composite with \( g \) are equal.
  A morphism \( u \colon U \to {\ein_{\cV}} \) classifies a pair of small hs-objects \( \und X,\und Y \in \bC / U \) for which there exists an elemental morphism \( h \colon \und X \to \und Y \).
  A morphism \( v \colon U \to X_{\fS} \) classifies a small hs-object \( \und Y \in \bC / U \) together with a section \( U \to Y \).
  Given \( u \colon U \to {\ein_{\cV}} \) classifying \( h \colon \und X \to \und Y \), the morphism \( u g \) classifies \( \und Y \) with the section \( \rt_X h \colon U \to Y \).

  Hence, to give two morphisms \( u_0,u_1 \colon U \to {\ein_{\cV}} \) with \( u_0 g = u_1 g \) is to give small hs-objects \( \und{X_0},\und{X_1},\und Y \in \bC / U \) with elemental morphisms \( h_j \colon \und{X_j} \to \und Y \) (for \( j = 0,1 \)) with \( \rt_{X_0} h_0 = \rt_{X_1} h_1 \).
  It thus follows from \jref{propn:dwn-hs} that \( \und{X_0} \cong \Dwn_{\rt_{X_0} h_0} \cong \Dwn_{\rt_{X_0} h_1} \cong \und{X_1} \), and hence that \( u_0 = u_1 \), as desired.
\end{proof}

Since \( g \) is a monomorphism, it is also small by \jref{defn:topoculus}~\ref{item:topoculus-monos}.
Hence, since \( \delta_{X_\fS} \) is small by the definition hs-classifier, it follows from \jref{defn:topoculus}~\ref{item:topoculus-comp} that \( i_{\ein_\cV} \pi_1 = g \delta_{X_{\fS}} \) is small as well.
This completes the proof of \jref{cor:small-closed-converse-univ-case}.

We now prove the general case.
We have pullback squares
\[
  \begin{tikzcd}[column sep=30pt]
    A \ar[r, dashed] \ar[d, "i_A"', >->] \pb & \ein_{\cV} \pb \ar[d, "i_{\ein_{\cV}}", >->] \ar[r] & \in_\cV \ar[d, >->] \\
    \cV \times U \ar[r, "\id_{\cV} \times x"] \ar[d, "\pi_1"] \pb & \cV \times \cV \ar[d, "\pi_1"] \ar[r, "\id \times \kappa"] & \cV \times \pow \cV \\
    U \ar[r, "x"] & \cV.
  \end{tikzcd}
\]
Thus, we are done by the smallness of \( i_{\ein_{\cV}} \pi_1 \) and the stability of small morphisms under pullbacks (\jref{defn:topoculus}~\ref{item:topoculus-pb}).

\subsection{Proof of \alttext{\hyperref[thm:model-of-mk]{Theorem~\ref*{thm:model-of-mk}}}{Theorem~\ref*{thm:model-of-mk}}: satisfaction of the set-existence axioms} \label{subsec:proof-of-set-existence-axioms}
We are given an hs-topoculized topos \( (\bC, \fS, \cV) \), and need to show that the MK-object associated to \( \cV \) as in \jref{defn:associated-mk} satisfies all the set existence axioms of MK, as given in \sref{subsec:axioms-of-mk}.

Each axiom has the form:
\[
  \forall x_1,\ldots,x_m,\beta_1,\ldots,\beta_n\
  \exists z\ \forall w\ \pbig{w \in \kappa(z) \tot \phi}
\]
for some formula \( \phi \), specifically with \( (m,n) \) equal to \( (0,0) \) (for Empty set), \( (1, 0) \) (for Union and Power), \( (2, 0) \) (for Pairing), or \( (1, 1) \) (for Subsets and Replacement).

Thus, in each case, we are given some morphisms \( x_i \colon U \to \cV \) and \( \beta_j \colon U \to \pow \cV \), and we need to find a morphism \( z \colon U \to \cV \) so that \( \alpha \defeq z \kappa \colon U \to \pow \cV \) satisfies
\begin{equation} \label{eq:mk-proof-alpha-cond}
  U \Vdash_{\vx, \vek \beta, \alpha} \forall w\ (w \in \alpha \tot \phi).
\end{equation}
Now, by Comprehension, we know that there is a unique \( \alpha \colon U \to \pow \cV \) satisfying \eqref{eq:mk-proof-alpha-cond}, so it only remains to show that \( \alpha \) factors through \( \kappa \colon \cV \to \pow \cV \).
In each case, we will show this using \jref{thm:topos-universe-small-closed}: we consider the underobject \( i_A \colon A \tto U \times \cV \) classified by \( \alpha \), and verify that \( i_A \pi_U \colon A \to U \) is small.

In fact, it suffices to prove this in the ``universal'' case, in which \( U = \cV^m \times (\pow \cV)^n \), and \( x_i = \pi_{\cV^m} \pi_i \) and \( \beta_i = \pi_{(\pow \cV)^m} \pi_i \).
In this case, \( A \) is simply the underobject
\begin{equation} \label{eq:univ-underob}
  A = \set{ \vx \vek \beta w \mid \phi } \tto \cV^m \times (\pow \cV)^n \times \cV.
\end{equation}
To see that this universal case indeed suffices, suppose that we have already shown that the morphism \( \alpha_{\univ} \colon \cV^m \times (\pow \cV)^n \to \pow \cV \) arising from the universal case factors through \( \kappa \colon \cV \to \pow \cV \).
Then, for general morphisms \( x_i \colon U \to \cV \) and \( \beta_i \colon U \to \pow \cV \), we claim that the morphism \( \alpha \colon U \to \pow \cV \) of interest (i.e., the one satisfying \eqref{eq:mk-proof-alpha-cond}) is given by
\[
  \alpha = \br{\vx, \vek \beta} \alpha_{\univ},
\]
and hence that \( \alpha \) factors through \( \kappa \) since \( \alpha_{\univ} \) does.
Indeed, writing \( \pi_\cV \vek \pi \) and \( \pi_{\pow \cV} \vek \pi \) for the sequences consisting, respectively of the morphisms \( \pi_{\cV} \pi_i \) and \( \pi_{\pow \cV} \pi_i \): from
\[
  \cV^m \times \cV^n \Vdash_{\vx, \vek \beta, \alpha \mapsto \pi_\cV \vek \pi,\pi_{\pow \cV}\vek \pi, \alpha_{\univ}} \forall w\ (w \in \alpha \tot \phi),
\]
it follows that
\[
  U \Vdash_{\vx, \vek \beta, \alpha \mapsto \vx,\vek \beta,\br{\vx, \vek \beta} \alpha_\univ} \forall w\ (w \in \alpha \tot \phi),
\]
since \( \br{\vx, \vek \beta} \pi_{\cV} \pi_i = x_i \), \( \br{\vx, \vek \beta} \pi_{\pow \cV} \pi_i = \beta_i \), and \( \br{\vx, \vek \beta, \alpha_{\univ}} \pi_{\pow \cV} \pi_i = \beta_i \).

Thus, to summarize, for each of the axioms, we must specify the formula \( \phi \), and then prove that the underobject \eqref{eq:univ-underob} is such that \( i_A \br{\pi_{\cV^m}, \pi_{(\pow \cV)^n}} \colon A \to \cV^m \times (\pow \cV)^n \) is small.
We now proceed to do this for each of the axioms in turn.

\subsubsection{Empty set}
In this case, \( \phi \) is \( \bot \).
Thus, the underobject \( A = \set{w \in \cV \mid \phi} \tto \cV \) is the initial underobject \( \zset \to \cV \), and hence \( i_A \unex \colon A \to \tm \) is just \( \zset \to \tm \), which is small, since it is a monomorphism.

\subsubsection{Pairing}
In this case, \( \phi \) is \( w = x \vee w = y \).
Hence, \( A = \set{\br{x,y,w} \in \cV^3 \mid \phi} \tto \cV \) is the union of the underobjects \( \pi_{010},\pi_{011} \colon \cV^2 \tto \cV^3 \).
This can equally be described as the image of the morphism
\[
  [\pi_{010},\pi_{011}] \colon \cV^2 + \cV^2 \to \cV^3.
\]
We thus have a commutative diagram
\[
  \begin{tikzcd}
    \cV^2 + \cV^2 \ar[r, "q", ->>] \ar[rr, bend left, "{[\pi_{010},\pi_{011}]}"] \ar[rrd, "{[\id, \id]}"', bend right] &
    A \ar[r, "i_A", >->] \ar[rd, "i_A \pi_{01}"'] & \cV^3 \ar[d, "\pi_{01}"] \\
    & & \cV^2
  \end{tikzcd}
\]
with \( q \) epic, as indicated.
Since \( [\pi_{010},\pi_{011}] \colon \cV^2 + \cV^2 \to \cV^2 \) is small (the subtopos \( \fS_{\cV^2} \subset \bC / \cV^2 \) being closed under terminal objects and coproducts), it follows from \jref{propn:topoculus-epi-closed} that \( i_A \pi_{01} \) is small, as desired.

\subsubsection{Unions}
In this case, \( \phi \) is \( \exists v.\ w \ein v \wedge v \ein x \).
The resulting underobject \( A \tto \cV \times \cV \) is the image under \( \pi_{20} \) of the intersection of \( \pi_{01}^* \ein_{\cV} \) and \( \pi_{12}^* \ein_{\cV} \).
It follows that we have a commutative diagram
\[
  \begin{tikzcd}
    & &[-15pt] \pi_{01}^* \unein \wedge \pi_{12}^* \unein \ar[dl, >->] \ar[dr, >->, "\rho"] \ar[rrr, "g", ->>] \pbr &[-15pt] & &[-25pt] A \ar[dddd, "i_A \pi_0"] \\
    & \pi_{01}^* \ein_\cV \ar[dl] \ar[dr, >->] \pbr & & \pi_{12}^* \ein_\cV \ar[dl, >->] \ar[dr, "\sigma"] \pbr \\
    \ein_{\cV} \ar[dr, >->, "i_{\ein_{\cV}}"'] & & \cV \times \cV \times \cV \ar[dl, "\pi_{01}"'] \ar[dr, "\pi_{12}"] \pbr & & \ein_{\cV} \ar[dl, >->, "i_{\ein_{\cV}}"] \\
    & \cV \times \cV \ar[dr, "\pi_1"'] & & \cV \times \cV \ar[dl, "\pi_0"] \ar[drr, "\pi_1"] \\
    & & \cV & & & \cV,
  \end{tikzcd}
\]
with pullback squares as indicated, and with \( g \) an epi.

Since \( i_{\ein_{\cV}} \pi_1 \colon \unein_{\cV} \to \cV \) is small by \jref{cor:small-closed-converse-univ-case}, it follows that the morphism \( \rho \sigma \colon \pi_{01}^* \unein \wedge \pi_{12}^* \unein \to \unein_{\cV} \) is small.
Using the smallness of \( i_{\ein_{\cV}} \pi_1 \) again and closure of small morphisms under composition, it follows that \( \rho \sigma i_{\ein_{\cV}} \pi_1 \colon \pi_{01}^* \unein \wedge \pi_{12}^* \unein \to \cV \) is small.
Since \( g \) is epi, it thus follows from \jref{propn:topoculus-epi-closed} that \( i_A \pi_0 \) is small, as desired.

\subsubsection{Power set}
In this case, \( \phi \) is \( \kappa(w) \subset \kappa(x) \).
The resulting underobject \( A \tto \cV \times \cV \) is thus simply the pullback
\[
  \begin{tikzcd}
    A \pb \ar[d, >->] \ar[r] & \le_{\cV} \ar[d, "i_{\le_{\cV}}", >->] \\
    \cV \times \cV \ar[r, "\kappa \times \kappa"] & \pow \cV \times \pow \cV.
  \end{tikzcd}
\]
Our task is to prove that \( i_A \pi_1 \colon A \to \cV \) is small, or in other words that \( (A, i_A \pi_1) \in \bC / \cV \) is contained in the full subcategory \( \fS_{\cV} \subset \bC / \cV \).
We will show that \( (A, i_A \pi_1) \) admits a monomorphism into the power object \( \pow (\ein_{\cV}, i_{\ein_{\cV}} \pi_1) \) in \( \bC / \cV \), which implies that \( (A, i_A \pi_1) \in \fS_\cV \), since \( (\ein_{\cV}, i_{\ein_{\cV}} \pi_1) \in \fS_\cV \) by \jref{cor:small-closed-converse-univ-case} and \( \fS_\cV \subset \bC / \cV \) is a sub-topos.

Letting \( \wt A = (\id_{\cV} \times \kappa)^* {\le}_{\cV} \tto \pow_\cV \times \cV \), we have the diagrams
\[
  \begin{tikzcd}
    A \pb \ar[r, "j", >->] \ar[d, "i_A"'] &
    \wt A \pb \ar[d, >->, "i_{\wt A}"] \ar[r] & \le_{\cV} \ar[d, "i_{\le_{\cV}}", >->] \\
    \cV \times \cV \ar[r, "\kappa \times \id"] \ar[rd, "\pi_1"'] & \pow \cV \times \cV \pb \ar[r, "\id \times \kappa"] \ar[d, "\pi_1"'] & \pow \cV \times \pow \cV \ar[d, "\pi_1"] \\
    & \cV \ar[r, "\kappa"] & \pow \cV
  \end{tikzcd}
  \qquad \AND \qquad
  \begin{tikzcd}
    \unein_\cV \pb \ar[r, ""] \ar[d, "i_{\ein_\cV}"', >->] &[20pt] \unin_\cV \ar[d, "i_{\unin_\cV}", >->] \\
    \cV \times \cV \pb \ar[r, "\id \times \kappa"] \ar[d, "\pi_1"'] & \cV \times \pow \cV \ar[d, "\pi_1"] \\
    \cV \ar[r, "\kappa"] &\pow \cV,
  \end{tikzcd}
\]
with pullback squares as indicated, and where \( j \) is monic since \( \kappa \) is.
Thus, to prove our assertion, it suffices to show that \( (\wt A, i_{\wt A} \pi_1) \) is a power object of \( (\unein_\cV, i_{\ein_\cV} \pi_1) \) in \( \bC / \cV \); and hence, since the logical anafunctor \( \kappa^* \colon \bC / \pow \cV \to \bC / \cV \) preserves power objects, it suffices to show that \( (\le_{\cV}, i_{\le_{\cV}} \pi_1) \) is a power object of \( (\unin_\cV, i_{\unin_\cV} \pi_1) \) in \( \bC / \pow \cV \).
We have thus reduced our problem to the following general topos-theoretic one:

\begin{claim}
  Let \( X \in \bC \) be any object in a topos \( \bC \), let \( \pow X \) be a power object for \( X \), and let \( {\le_X} \tto \pow X \times \pow X \) be the subset relation.
  Then \( (\le_X, i_{\le_X} \pi_1) \) is a power object of \( (\in_{X}, i_{\in_{X}} \pi_1) \) in \( \bC / \pow X \).
\end{claim}

This, in turn, follows immediately from the following two general claims.
Given an underobject \( A \tto X \) of an object \( X \) in a topos \( \bC \), classified by the global element \( a \colon \tm \to \pow X \), by the \defword{\boldmath \( A \)-part of \( \pow X \)}, we mean the underobject
  \[
    B = \set{p \in \pow X \mid p \le a} \tto \pow X.
  \]
  Thus, a morphism \( p \colon U \to \pow X \) factors through \( B \) if and only if there exists a dashed morphism making the following square commute, where \( P \tto X \times U \) is the underobject classified by \( p \).
  \[
    \begin{tikzcd}
      P \ar[r, dashed] \ar[d, >->] & A \ar[d, >->] \\
      X \times U \ar[r, "\pi_X"] & X
    \end{tikzcd}
  \]
\begin{claim} \label{claim:part-is-pow}
  Let \( A \tto X \) be an underobject of an object \( X \) in a topos \( \bC \), and let \( B \tto \pow X \) be the \( A \)-part of \( \pow X \).
  Then \( B \) is a power object of \( A \).
\end{claim}

\begin{claim} \label{claim:sub-is-el-part}
  Let \( X \) be an object in a topos \( \bC \), let \( \pow X \) be a power object of \( X \), and let \( {\le_X} \tto \pow X \times \pow X \) be the subset relation.
  Then, regarding \( (\pow X \times \pow X, \pi_1) \) as a power object of \( (X \times \pow X, \pi_1) \) in \( \bC / \pow X \) as in \jref{exm:pb-fun-logical}, we have that \( (\le_X, i_{\le_X} \pi_1) \) is the \( (\in_{X}, i_{\in_{X}} \pi_1) \)-part of \( (\pow X \times \pow X, \pi_1) \).
\end{claim}

\begin{proof}[Proof of \jref{claim:sub-is-el-part}]
  We must prove that, given \( \und U = (U, \delta_U) \in \bC / \pow X \), a morphism \( \br{p, \delta_U} \colon \und U \to (\pow X \times \pow X, \pi_1) \) factors through \( (\le_X, i_{\le_X} \pi_1) \tto (\pow X \times \pow X, \pi_1) \) iff there exists a dashed morphism making the following square commute, where \( P \tto (X \times \pow X, \pi_1) \times \und U \) is the underobject classified by \( \br{p, \delta_U} \).
  \[
    \begin{tikzcd}
      P \ar[r, dashed] \ar[d, ""', >->] &[40pt] (\in_X, i_{\in_X} \pi_1) \ar[d, "i_{\in_X}"] \\
      (X \times \pow X, \pi_1) \times \und U \ar[r, "\pi_{(X \times \pow X, \pi_{\pow X})}"] & (X \times \pow X, \pi_1).
    \end{tikzcd}
  \]
  Now, the product \( (X \times \pow X, \pi_1) \times \und U = \pbig{(X \times \pow X) \times_{\pow X} (U, \pi_1 \delta_U)} \) is given by \( (X \times U, \pi_1 \delta_U) \) with projections
  \[
    (X \times \pow X, \pi_1) \xot{\id_X \times \delta_U}
    (X \times U, \pi_1 \delta_U) \tox{\pi_1}
    \und U.
  \]
  Moreover, under this identification, as well as the similar identification
  \[
    (X \times \pow X, \pi_1) \times (\pow X \times \pow X, \pi_{1}) = (X \times \pow X \times \pow X, \pi_2),
  \]
  the morphism \( \id_{(X \times \pow X, \pi_1)} \times \br{p, \delta_U} \colon (X \times \pow X, \pi_1) \times \und U \to (X \times \pow X, \pi_1) \times (\pow X \times \pow X, \pi_{1}) \) is given by
  \[
    \br{\pi_0, \pi_1 p, \pi_1 \delta_U} \colon (X \times U, \pi_1 \delta_U) \to (X \times \pow X \times \pow X, \pi_2).
  \]
  We thus have pullback squares
  \[
    \begin{tikzcd}
      P \pb \ar[d, >->] \ar[r] &[30pt] \in_X \times \pow_X \pb \ar[r] \ar[d, "i_{\in_X} \times \id"] &[30pt] \in_X \ar[d, >->] \\
      X \times U \ar[r, "\br{\pi_0, \pi_1 p, \pi_1 \delta_U}"] & X \times \pow_X \times \pow_X \ar[r, "\br{\pi_0, \pi_1}"] & X \times \pow_X,
    \end{tikzcd}
  \]
  showing that \( P \tto X \times U \) is the underobject classified by \( p \) in \( \bC \).
  Thus, the claim is reduced to showing that \( \br{p, \delta_U} \colon U \to \pow X \times \pow X \) factors through \( {\le}_X \tto \pow X \times \pow X \) if and only there exists a dashed morphism making
  \[
    \begin{tikzcd}
      P \ar[r, dashed] \ar[d, >->] &[10pt] \in_X \ar[d, >->] \\
      X \times U \ar[r, "\id \times \delta_U"] & X \times \pow_X
    \end{tikzcd}
  \]
  commute, where now \( P \) is the underobject classified by \( p \).
  But this is the same as saying that \( P \le Q \), where \( Q \tto X \times U \) is the underobject classified by \( \delta_U \).
  And indeed, this is the case if and only if \( p \le \delta_U \), as required.
\end{proof}

\begin{proof}[Proof of \jref{claim:part-is-pow}]
  A reformulation of the defining condition for \( B \) to be the \( A \)-part of \( \pow X \) is that a morphism \( p \colon U \to \pow X \) classifying \( i_P \colon P \tto X \times U \) factors through \( B \tto \pow X \) if and only if \( i_P \) factors through \( i_A \times \id \colon A \times U \to X \times U \).
  It follows that the natural bijection \( \bC(U, \pow X) \toi \Sub(X \times U) \) restricts to a (likewise natural) bijection
  \[
    \begin{tikzcd}
      \bC(U, B) \ar[r, "\sim", dashed] \ar[d, "(i_B)_*"', >->] & \Sub(A \times U) \ar[d, "(i_A \times \id)_*", >->] \\
      \bC(U, \pow X) \ar[r, "\sim"] & \Sub(X \times U),
    \end{tikzcd}
  \]
  showing that \( B \) is a power object of \( A \), as desired.
\end{proof}

\subsubsection{Subsets}
In this case, \( \phi \) is \( w \ein x \wedge w \in \beta \).
Hence the underobject \( A = \set{ \br{x, \beta, w} \mid \phi } \tto \cV \times \pow \cV \times \cV \) is the intersection of \( \pi_{20}^* \unein_\cV \) and \( \pi_{21}^* \unin_{\cV} \).
We thus have the following diagram, with pullback squares as indicated and \( f \) a mono.
\[
  \begin{tikzcd}
    A \ar[r, >->, "f"] \ar[rd, "i_A"'] & \pi_{20}^* \unein_\cV \pb \ar[r] \ar[d, "i_{\pi_{20}^* \unein_\cV}"] & \unein_\cV \ar[d, >->, "i_{\unein_\cV}"] \\
    & \cV \times \pow \cV \times \cV \pb \ar[r, "\pi_{20}"] \ar[d, "\pi_{01}"] & \cV \times \cV \ar[d, "\pi_1"] \\
    & \cV \times \pow \cV \ar[r, "\pi_0"] & \cV
  \end{tikzcd}
\]
Since \( i_{\unein_{\cV}} \pi_1 \colon \unein_{\cV} \to \cV \) is small by \jref{cor:small-closed-converse-univ-case}, its pullback \( i_{\pi_{20}^* \unein_\cV} \pi_{01} \) along \( \pi_0 \) is small as well, and hence \( i_A \) is the composite of a small morphism with the monomorphism \( f \), and thus small.

\subsubsection{Replacement}
In this case, \( \phi \) is
\[
  \exists u.\ u \ein x \wedge u \beta w \wedge (\forall v.\ u \beta v \to v = w).
\]
We recall here from \sref{subsec:axioms-of-mk} that \( u \beta w \) and \( u \beta v \) are shorthands for certain complicated formulas.
However, we will in fact not need to know anything about these formulas, except for the obvious fact that they are related to each other by a substitution \( (u \beta v) \equiv (u \beta w)[w / v] \), and hence that, for \( u,w \colon U \to \cV \) and \( \beta \colon U \to \pow \cV \), we have
\[
  U \Vdash_{\beta, w, u} u \beta w
  \quad \iff \quad
  U \Vdash_{\br{\beta, v, u} \mapsto \br{\beta, w, u}} u \beta w,
\]
(cf.\ \eqref{eq:subst-and-kj} in \sref{subsec:uniq-ex-eff-epi}), and hence, finally, that, writing \( \psi \) for the formula \( u \beta w \wedge (\forall v.\ u \beta v \to v = w) \), we have, for any \( u,w_0,w_1 \colon U \to \cV \) and \( \beta \colon U \to \pow \cV \), that
\[
  U \Vdash_{\beta,w_0,u} \psi
  \quad\AND\quad
  U \Vdash_{\beta,w_1,u} \psi
  \quad\IMPLIES\quad
  w_0 = w_1.
\]
(In fact, this holds already if the second hypothesis is replaced by \( U \Vdash_{\beta, w_1, u} u \beta w \).)
Writing
\[
  P = \set{\br{\beta, w, u} \mid \psi} \tto \pow \cV \times \cV \times \cV,
\]
we thus conclude that \( i_P \pi_{02} \colon P \to \pow \cV \times \cV \) is monic.

In what follows, we decorate each occurrence of \( \cV \) with one of the variables \( x \), \( u \), or \( w \), and decorate each occurrence of \( \pow \cV \) with the variable \( \beta \).
These decorations have no substantive meaning and can be ignored; they are simply provided to help understand the argument.
Thus, for example, the two last two morphisms mentioned are \( i_P \colon P \tto \pow \cV_\beta \times \cV_w \times \cV_u \) and \( i_P \pi_{02} \colon P \to \pow \cV_\beta \times \cV_u \).

Turning to the task at hand, the underobject \( A = \set{\br{x, \beta, w} \mid \phi} \tto \cV_x \times \pow \cV_\beta \times \cV_w \) is the image under \( \pi_{012} \colon \cV_x \times \pow \cV_\beta \times \cV_w \times \cV_u \to \cV_x \times \pow \cV_\beta \times \cV_w \) of the underobject
\[
  \pi_{123}^* P \wedge \pi_{30}^* \unein_\cV \tto \cV_x \times \pow \cV_\beta \times \cV_w \times \cV_u.
\]
We thus have a commutative diagram
\[
  \begin{tikzcd}
    \pi_{123}^* P \wedge \pi_{30}^* \unein_\cV \pb \ar[r, >->] \ar[d, >->] \ar[dd, bend right=45pt, ->>, "g"'] & \pi_{30}^* \unein_\cV \pb \ar[d, >->] \ar[r] & \pi_{20}^* \unein_{\cV} \pb \ar[r] \ar[d, >->] & \unein_\cV \ar[d, "i_{\unein_\cV}", >->] \\
    \pi_{123}^* P \ar[r, >->] & \cV_x \times \pow \cV_\beta \times \cV_w \times \cV_u \pb \ar[d, "\pi_{012}"] \ar[r, "\pi_{013}"] & \cV_x \times \pow \cV_\beta \times \cV_u \pb \ar[d, "\pi_{01}"] \ar[r, "\pi_{20}"] & \cV_u \times \cV_x \ar[d, "\pi_1"] \\
    A \ar[r, "i_A"', >->] & \cV_x \times \pow \cV_\beta \times \cV_w \ar[r, "\pi_{01}"'] & \cV_x \times \pow \cV_\beta \ar[r, "\pi_0"'] & \cV_x
  \end{tikzcd}
\]
with pullback squares as indicated and \( g \) epic.

Now, the morphism \( i_{\unein_\cV} \pi_1 \colon \unein_\cV \to \cV_x \) is small by \jref{cor:small-closed-converse-univ-case}, hence its pullback \( \pi_{20}^* \unein_\cV \to \cV_x \times \pow \cV_\beta \) in the penultimate column is as well.
Next, the morphism \( \pi_{123}^* P \to \cV_x \times \pow \cV_\beta \times \cV_u \) in the second row is monic, since we have a pullback diagram
\[
  \begin{tikzcd}
    \pi_{123}^* P \pb \ar[r, >->] \ar[d] &[20pt] \cV_x \times \pow \cV_\beta \times \cV_w \times \cV_u \pb \ar[d, "\pi_{123}"] \ar[r, "\pi_{013}"] & \cV_x \times \pow \cV_\beta \times \cV_u \ar[d, "\pi_{12}"] \\
    P \ar[r, "i_P", >->] & \pow \cV_\beta \times \cV_w \times \cV_u \ar[r, "\pi_{02}"] & \pow \cV_\beta \times \cV_u
  \end{tikzcd}
\]
in which the bottom row is monic as observed above.

Hence, the morphism \( \pi_{123}^* P \wedge \pi_{30}^* \unein_\cV \to \pi_{20}^* \unein_\cV \) in the top row of the above diagram is the pullback of a monic, and hence itself monic, and so the morphism \( g i_A \pi_{01} \colon \pi_{123}^* P \wedge \pi_{30}^* \unein_\cV \to \cV_x \times \pow \cV_\beta \) is the composite of a mono and a small morphism, and hence small.
Hence, finally, since \( g \) is epic, it follows that \( i_A \pi_{01} \colon A \to \cV_x \times \pow \cV_\beta \) is small by \jref{propn:topoculus-epi-closed}, as desired.

\subsection{Proof of \alttext{\hyperref[thm:infinity]{Theorem~\ref*{thm:infinity}}}{Theorem~\ref*{thm:infinity}}: the Axiom of Infinity} \label{subsec:infinity-proof}
We prove the three parts of the theorem in order.

For part~\ref{item:infinity-nno-exists}, we must show that if a topoculus \( \fS \) on a topos \( \bC \) admits an hs-classifier \( \cV \), then \( \bC \) has an NNO.
We recall from \sref{subsec:axioms-of-mk} that in the language of MK, we use Latin letters for variables of sort \( \St \) and Greek letters for variables of sort \( \Cls \).
\begin{claim} \label{claim:zero-and-succ-exist}
  There exist (by extensionality unique) morphisms \( \emptyset \colon \tm \to \cV \) and \( (-)^+ \colon \cV \to \cV \) such that \( \bC \vDash \forall w.\ \neg w \ein \emptyset \) and \( \bC \vDash \forall v,w.\ w \ein v^+ \tot (w \ein v \vee w = v) \).
\end{claim}
\begin{proof}
  The existence of \( \emptyset \) is precisely the statement that \( \cV \) satisfies the Axiom of Empty Set, as was established in \jref{thm:model-of-mk}.
  Next, we use that \( \cV \) satisfies the Axioms of Pairing and Union, as were also established there.
  First, applying Pairing to \( \pi_0, \pi_1 \colon \cV \times \cV \to \cV \), we obtain \( p \colon \cV \times \cV \to \cV \) such that \( w \ein \br{x,y} p \) iff \( U \Vdash_{x,y,w} w = x \vee w = y \) for \( x,y,w \colon U \to \cV \); hence also \( w \ein x \Delta p \) iff \( w = x \), where \( \Delta = \br{\id, \id} \colon \cV \to \cV \times \cV \).
  Then, applying Union to \( \br{\id, \Delta p} p \colon \cV \to \cV \), we obtain \( (-)^+ \colon \cV \to \cV \) with \( w \ein x (-)^+ \) for \( w, x \colon U \to \cV \) iff \( w' \ein y \ein \br{x', x' \Delta p} p \) for some \( U \oot U' \tox{y} \cV \), iff \( w_0 \ein x_0 \) and \( w_1 \ein x_1 \Delta p \) for some cover \( U_0,U_1 \to U \), iff \( w_0 \ein x_0 \) and \( w_1 = x_1 \) for some cover \( U_0,U_1 \to U \), as desired.
\end{proof}
In the statement of the Axiom of Infinity in \sref{subsec:axioms-of-mk}, we explained that the expressions ``\( \emptyset \in \beta \)'' and ``\( v^+ \in \beta \)'' were merely shorthands.
However, we can also introduce a constant symbol \( \emptyset \) and function symbol \( v \colon \St \to \St \) into the language of MK, and interpret them in (the MK-object associated to) \( \cV \) using the morphisms from \jref{claim:zero-and-succ-exist}.
It is then clear (for example, by spelling out the KJ interpretation of each) that \( \cV \) satisfies the resulting modified Axiom of Infinity (with the shorthands replaced by the expressions involving the new function symbols) if and only if it satisfies the original Axiom of Infinity.
Thus, henceforth, we always understand ``Axiom of Infinity'' in this modified sense.

Now, form an underobject \( i_N \colon N \tto \cV \) defined by the formula appearing in the Axiom of Infinity:
\[
  N = \set{w \in \cV \mid \forall \beta.\ \pbig{\emptyset \in \beta \wedge \forall v.\ (v \in \beta \to v^+ \in \beta)} \to w \in \beta} \tto \cV.
\]
\begin{claim} \label{claim:n-has-zer-suc}
  The morphisms \( \emptyset \colon \tm \to \cV \) and \( i_N (-)^+ \colon N \to \cV \) both factor through \( N \), giving morphisms \( \zer \colon \tm \to N \) and \( \suc \colon N \to N \).
\end{claim}
\begin{proof}
  By the definition of \( N \), for \( w \colon U \to \cV \), we have that \( N w \) iff \( w' = \delta_{U'} w \in \beta \) for any \( \beta \colon U' \to \pow \cV \) satisfying \( U' \Vdash_\beta \emptyset \in \beta \wedge \forall v.\ (v \in \beta \to v^+ \in \beta) \).
  Thus, we must check this condition with \( w = \emptyset \colon \tm \to \cV \), and with \( w = i_N (-)^+ \colon N \to \cV \).

  In the case \( w = \emptyset \), the condition \( \delta_{U'} \emptyset \in \beta \) is immediate from the assumption on \( \beta \).

  In the case \( w = i_N (-)^+ \), given any such \( \beta \colon U' \to \pow \cV \), first note that \( N i_N \) tautologically and hence \( \delta_{U'} i_N \in \beta \).
  But this immediately implies \( \delta_{U'} i_n (-)^+ \), as desired, by the assumption on \( \beta \).
\end{proof}

\begin{claim}
  The object \( N \), equipped with the morphisms \( \zer \) and \( \suc \) from \jref{claim:n-has-zer-suc}, is an NNO.
\end{claim}

\begin{proof}
  By \cite[\S D5.1]{johnstone-elephant-vols}, it suffices to check that \( N \) satisfies the Peano postulates, i.e., that (i) \( \suc \) is monic, (ii) the underobjects \( \zer \colon \tm \to N \) and \( \suc \colon N \to N \) are disjoint, and (iii) \( p = \unex \top_N \) for any \( p \colon U \to \pow N \) with \( \zer \in \unex p \) and \( x\,\suc \in p' \) for all \( x \colon U' \to N \) with \( x \in p' \).

  To prove (i), it suffices to check that \( (-)^+ \colon \cV \to \cV \) is monic.
  Let \( x,y \colon U \to \cV \) with \( x (-)^+ = y (-)^+ \); we must show that \( x = y \).
  First, note that if we additionally assume \( x \ein y \), then \( U \) is initial; indeed, since \( y \ein y (-)^+ = x (-)^+ \), we have a cover \( U_0,U_1 \to U \) with \( y_0 \ein x_0 \) and \( y_1 = x_1 \), and thus we get that \( U_0 \) and \( U_1 \) are both initial by \jref{propn:no-loops}.
  But now since \( x \ein x (-)^+ = y (-)^+ \), we have a cover \( U_0,U_1 \to U \) with \( x_0 \ein y_0 \) and \( x_1 = y_1 \), and \( U_0 \) is initial by what we just showed.
  Hence \( U_1 \to U \) is epi and so \( x = y \).

  To prove (ii), it suffices to show that \( \emptyset \colon \tm \to \cV \) and \( (-)^+ \colon \cV \to \cV \) intersect trivially in \( \cV \), i.e., that given \( x \colon U \to \cV \), if \( x (-)^+ = \unex_U\, \emptyset \), then \( U \) is initial.
  But this is clear from the definition of \( \emptyset \) since \( x \ein x (-)^+ \).

  To prove (iii), given such a \( p \colon U \to \pow N \), it suffices to show that \( p \exists_{i_N} = \top_N \exists_{i_N} \colon U \to \pow \cV \), since \( \exists_{i_N} \colon \pow N \to \pow \cV \) is monic (because \( i_N \) is; this can be checked directly using \jref{propn:topos-extensionality}).
  Thus, we must show that \( x \in p \exists_{i_N} \) for all \( x \colon U \to \cV \) with \( N x \).
  This follows since \( p \exists_{i_N} \), by our assumption on \( p \), satisfies the hypothesis on ``\( \beta \)'' appearing in the definition of \( N \).
\end{proof}

We have thus proven part~\ref{item:infinity-nno-exists} of \jref{thm:infinity}.
Part~\ref{item:infinity-small-nno}---that the hs-classifier \( \cV \) satisfies the Axiom of Infinity iff the NNO in \( \bC \) is small---now follows immediately.
Indeed, as in the proofs of the other set-existence axioms in \sref{subsec:proof-of-small-converse}, we know that we have \( \alpha \colon \tm \to \pow \cV \) satisfying the defining condition of \( \N \), and the Axiom of Infinity holds iff \( \alpha \) factors through \( \kappa \colon \cV \to \pow \cV \).
But by the way we defined the NNO \( N \in \bC \) above, we see that \( \alpha \) is precisely the morphism classifying \( N \tto \cV \), and thus by \jref{thm:topos-universe-small-closed} and its converse \jref{propn:topos-universe-small-closed-converse}, \( \alpha \) factors through \( \kappa \) iff \( N \) is small.

Part~\ref{item:infinity-2cat-has-nno} of the theorem now also follows readily.
The setup here is that we have plentiful DOF classifiers \( \rS, \rS', \rS'' \in \cC \) in a groupoidant 2-category \( \cC \) with \( p \colon \rS_* \to \rS \) and \( \hs(\rS) \) both \( \rS' \)-small, and \( p' \colon \rS'_* \to \rS' \) and \( \hs(\rS') \) both \( \rS'' \)-small, and we must show that the NNO in the hs-topoculized topos \( \SOF_{\rS''}(\tm) \) is small, i.e., \( \rS' \)-small.
(And we must show the same in \( \SOF(\tm) \) when \( \cC \) satisfies axiom (UA), but the argument for that is exactly the same.)

Now, in the hs-topoculized topos \( \SOF_{\rS'}(\tm) \), we have, as just showed, an NNO \( N \) which is given as an underobject \( N \tto \hs(\rS) \) of the hs-classifier.
But the logical functor \( \SOF_{\rS'}(\tm) \hto \SOF_{\rS''}(\tm) \) preserves NNOs and monos, hence we have a mono \( N \to \hs(\rS) \) in \( \SOF_{\rS''}(\tm) \) from the NNO \( N \), where now \( \hs(\rS) \) is \( \rS' \)-small, and hence \( N \) is \( \rS' \)-small as well, as desired.

\printbibliography

@online {abellan-martini-2-topoi,
  TITLE        = {$(\infty,2)$-Topoi and descent},
  AUTHOR       = {Fernando Abellán and Louis Martini},
  YEAR         = {2024},
  EPRINT       = {2410.02014},
  ARCHIVEPREFIX= {arXiv},
  PRIMARYCLASS = {math.CT},
  URL          = {https://arxiv.org/abs/2410.02014},
}

@book {aczel-non-well-founded-sets,
  AUTHOR       = {Aczel, Peter},
  TITLE        = {Non-well-founded sets},
  SERIES       = {CSLI Lecture Notes},
  VOLUME       = {14},
  NOTE         = {With a foreword by Jon Barwise [K. Jon Barwise]},
  PUBLISHER    = {Stanford University, Center for the Study of Language and Information, Stanford, CA},
  YEAR         = {1988},
  PAGES        = {xx+137},
  ISBN         = {0-937073-22-9},
  MRCLASS      = {03E30 (03E65 03E70)},
  MRNUMBER     = {940014},
  MRREVIEWER   = {Alena\ Vencovsk\'{a}},
}

@article {awodey-butz-simpsons-streicher-relating,
  AUTHOR       = {Awodey, Steve and Butz, Carsten and Simpson, Alex and Streicher, Thomas},
  TITLE        = {Relating first-order set theories, toposes and categories of classes},
  JOURNAL      = {Ann. Pure Appl. Logic},
  FJOURNAL     = {Annals of Pure and Applied Logic},
  VOLUME       = {165},
  YEAR         = {2014},
  NUMBER       = {2},
  PAGES        = {428--502},
  ISSN         = {0168-0072,1873-2461},
  MRCLASS      = {03G30 (03C90 03E70 03F50 18B05)},
  MRNUMBER     = {3129724},
  MRREVIEWER   = {Robert\ S.\ Lubarsky},
  DOI          = {10.1016/j.apal.2013.06.004},
  URL          = {https://doi.org/10.1016/j.apal.2013.06.004},
}

@book {barwise-handbook,
  EDITOR       = {Barwise, Jon},
  TITLE        = {Handbook of mathematical logic},
  SERIES       = {Studies in Logic and the Foundations of Mathematics},
  VOLUME       = {90},
  NOTE         = {Edited by Jon Barwise, With the cooperation of H. J. Keisler, K. Kunen, Y. N.  Moschovakis and A. S. Troelstra},
  PUBLISHER    = {North-Holland Publishing Co., Amsterdam},
  YEAR         = {1977},
  PAGES        = {xi+1165},
  ISBN         = {0-7204-2285-X},
  MRCLASS      = {02-06},
  MRNUMBER     = {457132},
}

@article {bourke-accessible-aspects,
  AUTHOR       = {Bourke, John},
  TITLE        = {Accessible aspects of 2-category theory},
  JOURNAL      = {J. Pure Appl. Algebra},
  FJOURNAL     = {Journal of Pure and Applied Algebra},
  VOLUME       = {225},
  YEAR         = {2021},
  NUMBER       = {3},
  PAGES        = {Paper No. 106519, 43},
  ISSN         = {0022-4049},
  MRCLASS      = {18C35 (18N10)},
  MRNUMBER     = {4137713},
  MRREVIEWER   = {Hugo Luiz Mariano},
  DOI          = {10.1016/j.jpaa.2020.106519},
  URL          = {https://doi.org/10.1016/j.jpaa.2020.106519},
}

@inproceedings {cole-categories_of_sets,
  AUTHOR       = {Cole, J. C.},
  TITLE        = {Categories of sets and models of set theory},
  BOOKTITLE    = {The {P}roceedings of the {B}ertrand {R}ussell {M}emorial {C}onference ({U}ldum, 1971)},
  PAGES        = {351--399},
  PUBLISHER    = {Bertrand Russell Memorial Logic Conference, Leeds},
  YEAR         = {1973},
  MRCLASS      = {02K15 (18B05)},
  MRNUMBER     = {357116},
  MRREVIEWER   = {Andreas\ Blass},
}

@article {fourman-sheaf-models,
  AUTHOR       = {Fourman, Michael P.},
  TITLE        = {Sheaf models for set theory},
  JOURNAL      = {J. Pure Appl. Algebra},
  FJOURNAL     = {Journal of Pure and Applied Algebra},
  VOLUME       = {19},
  YEAR         = {1980},
  PAGES        = {91--101},
  ISSN         = {0022-4049},
  MRCLASS      = {03G30 (03C62 03E40 03E70 03F55 18B25)},
  MRNUMBER     = {593249},
  MRREVIEWER   = {Andreas Blass},
  DOI          = {10.1016/0022-4049(80)90096-1},
  URL          = {https://doi.org/10.1016/0022-4049(80)90096-1},
}

@incollection {hayashi-set-theories-in-toposes,
  AUTHOR       = {Hayashi, Susumu},
  TITLE        = {On set theories in toposes},
  BOOKTITLE    = {Logic {S}ymposia, {H}akone 1979, 1980 ({H}akone, 1979/1980)},
  SERIES       = {Lecture Notes in Math.},
  VOLUME       = {891},
  PAGES        = {23--29},
  PUBLISHER    = {Springer, Berlin-New York},
  YEAR         = {1981},
  MRCLASS      = {18B25 (03E40 03E70 03F55 03G30)},
  MRNUMBER     = {656960},
  MRREVIEWER   = {M. M. Richter},
}

@book {heijenoort-frege-to-godel,
  AUTHOR       = {van Heijenoort, Jean},
  TITLE        = {From {F}rege to {G}\"{o}del. {A} source book in mathematical logic, 1879--1931},
  PUBLISHER    = {Harvard University Press, Cambridge, MA},
  YEAR         = {1967},
  PAGES        = {xi+660 pp. (1 plate)},
  MRCLASS      = {02.00},
  MRNUMBER     = {209111},
  MRREVIEWER   = {R.\ C.\ Lyndon},
}

@article {helfer-topoi-in-topoi,
  AUTHOR       = {Helfer, Joseph},
  TITLE        = {Internal 1-topoi in 2-topoi},
  JOURNAL      = {Theory Appl. Categ.},
  FJOURNAL     = {Theory and Applications of Categories},
  VOLUME       = {45},
  YEAR         = {2026},
  PAGES        = {Paper No. 21, 779--862},
  MRCLASS      = {18B25 (03G30 18N10 18N25)},
  URL          = {http://www.tac.mta.ca/tac/volumes/35/35/35-35.pdf},
}

@book {johnstone-elephant-vols,
  AUTHOR       = {Johnstone, Peter T.},
  TITLE        = {Sketches of an elephant: a topos theory compendium. {V}ols. 1 and 2},
  SERIES       = {Oxford Logic Guides},
  VOLUME       = {43,44},
  PUBLISHER    = {The Clarendon Press, Oxford University Press, Oxford},
  YEAR         = {2002},
  PAGES        = {xxii+468+621+71},
  ISBN         = {0-19-853425-6,0-19-851598-7},
  MRCLASS      = {18B25 (03B15 03G30 18-02 54A05)},
  MRNUMBER     = {1953060,2063092},
  MRREVIEWER   = {Colin McLarty},
  DOI          = {10.1093/oso/9780198534259.001.0001},
}

@book {johnstone-topos-theory,
  AUTHOR       = {Johnstone, P. T.},
  TITLE        = {Topos theory},
  SERIES       = {London Mathematical Society Monographs, Vol. 10},
  PUBLISHER    = {Academic Press [Harcourt Brace Jovanovich, Publishers], London-New York},
  YEAR         = {1977},
  PAGES        = {xxiii+367},
  ISBN         = {0-12-387850-0},
  MRCLASS      = {18B99 (02K10 18A15 18C99)},
  MRNUMBER     = {0470019},
  MRREVIEWER   = {Marta C. Bunge},
}

@book {joyal-moerdijk-ast,
  AUTHOR       = {Joyal, A. and Moerdijk, I.},
  TITLE        = {Algebraic set theory},
  SERIES       = {London Mathematical Society Lecture Note Series},
  VOLUME       = {220},
  PUBLISHER    = {Cambridge University Press, Cambridge},
  YEAR         = {1995},
  PAGES        = {viii+123},
  ISBN         = {0-521-55830-1},
  MRCLASS      = {03E70 (03-02 03C90 03F55 03G25 04A10)},
  MRNUMBER     = {1368403},
  MRREVIEWER   = {Ioan Tofan},
  DOI          = {10.1017/CBO9780511752483},
  URL          = {https://doi.org/10.1017/CBO9780511752483},
}

@book {lambek-scott,
  AUTHOR       = {Lambek, J. and Scott, P. J.},
  TITLE        = {Introduction to higher order categorical logic},
  SERIES       = {Cambridge Studies in Advanced Mathematics},
  VOLUME       = {7},
  PUBLISHER    = {Cambridge University Press, Cambridge},
  YEAR         = {1986},
  PAGES        = {x+293},
  ISBN         = {0-521-24665-2},
  MRCLASS      = {03B40 (03-02 03B15 03G30 18-02 18B25 18D15)},
  MRNUMBER     = {856915},
  MRREVIEWER   = {Saunders Mac Lane},
}

@incollection {lawvere-cat-of-cats,
  AUTHOR       = {Lawvere, F. William},
  TITLE        = {The category of categories as a foundation for mathematics},
  BOOKTITLE    = {Proc. {C}onf. {C}ategorical {A}lgebra ({L}a {J}olla, {C}alif., 1965)},
  PAGES        = {1--20},
  PUBLISHER    = {Springer, New York},
  YEAR         = {1966},
  MRCLASS      = {02.00},
  MRNUMBER     = {0207517},
  MRREVIEWER   = {J. R. Isbell},
  DOI          = {10.1007/978-3-642-99902-4_1},
}

@article {lawvere-etcs,
  AUTHOR       = {Lawvere, F. William},
  TITLE        = {An elementary theory of the category of sets (long version) with commentary (2005 reprint)},
  NOTE         = {Reprinted and expanded from Proc. Nat. Acad. Sci. U.S.A.  {{\bf{5}}2} (1964) [MR0172807], With comments by the author and Colin McLarty},
  JOURNAL      = {Repr. Theory Appl. Categ.},
  FJOURNAL     = {Reprints in Theory and Applications of Categories},
  NUMBER       = {11},
  YEAR         = {1964},
  PAGES        = {1--35},
  MRCLASS      = {03A05 (00B60 03E99 18B05)},
  MRNUMBER     = {2177727},
  URL          = {http://www.tac.mta.ca/tac/reprints/articles/11/tr11.pdf},
}

@incollection {lawvere-variable,
  AUTHOR       = {Lawvere, F. William},
  TITLE        = {Variable quantities and variable structures in topoi},
  BOOKTITLE    = {Algebra, topology, and category theory (a collection of papers in honor of {S}amuel {E}ilenberg)},
  PAGES        = {101--131},
  PUBLISHER    = {Academic Press, New York},
  YEAR         = {1976},
  MRCLASS      = {02K10 (18A15)},
  MRNUMBER     = {0419232},
  MRREVIEWER   = {A. Preller},
}

@book {mac-lane-moerdijk,
  AUTHOR       = {Mac Lane, Saunders and Moerdijk, Ieke},
  TITLE        = {Sheaves in geometry and logic},
  SERIES       = {Universitext},
  NOTE         = {A first introduction to topos theory, Corrected reprint of the 1992 edition},
  PUBLISHER    = {Springer-Verlag, New York},
  YEAR         = {1994},
  PAGES        = {xii+629},
  ISBN         = {0-387-97710-4},
  MRCLASS      = {03G30 (18B25 54B40)},
  MRNUMBER     = {1300636},
  MRREVIEWER   = {M. Makkai},
  DOI          = {10.1007/978-1-4612-0927-0},
}

@online {makkai-abstract-sets,
  TITLE        = {The theory of abstract sets based on first-order logic with dependent types},
  AUTHOR       = {Makkai, M.},
  YEAR         = {2013},
  URL          = {https://www.math.mcgill.ca/makkai/Various/MateFest2013.pdf},
}

@article {makkai-avoiding-choice,
  AUTHOR       = {Makkai, M.},
  TITLE        = {Avoiding the axiom of choice in general category theory},
  JOURNAL      = {J. Pure Appl. Algebra},
  FJOURNAL     = {Journal of Pure and Applied Algebra},
  VOLUME       = {108},
  YEAR         = {1996},
  NUMBER       = {2},
  PAGES        = {109--173},
  ISSN         = {0022-4049},
  MRCLASS      = {18A15 (04A25)},
  MRNUMBER     = {1382246},
  MRREVIEWER   = {Andreas Blass},
  DOI          = {10.1016/0022-4049(95)00029-1},
  URL          = {https://doi.org/10.1016/0022-4049(95)00029-1},
}

@incollection {makkai-categorical-foundation,
  AUTHOR       = {Makkai, M.},
  TITLE        = {Towards a categorical foundation of mathematics},
  BOOKTITLE    = {Logic {C}olloquium '95 ({H}aifa)},
  SERIES       = {Lecture Notes Logic},
  VOLUME       = {11},
  PAGES        = {153--190},
  PUBLISHER    = {Springer, Berlin},
  YEAR         = {1998},
  MRCLASS      = {03B15 (00A30 03A05 03G30 18A15 18D05)},
  MRNUMBER     = {1678360},
  MRREVIEWER   = {Peter Johnstone},
  DOI          = {10.1017/9781316716830.014},
  URL          = {https://doi.org/10.1017/9781316716830.014},
}

@online {makkai-folds,
  AUTHOR       = {Makkai, M.},
  TITLE        = {First Order Logic with Dependent Sorts, with Applications to Category Theory.},
  YEAR         = {1995},
  URL          = {http://www.math.mcgill.ca/makkai/folds/foldsinpdf/FOLDS.pdf},
}

@online {makkai-fom-both,
  AUTHOR       = {Makkai, M.},
  TITLE        = {A new foundation for abstract mathematics},
  YEAR         = {2003},
  NOTE         = {Posting on the Foundations Of Mathematics mailing list.},
  URL          = {https://www.math.mcgill.ca/makkai/Foundations seminar/FoundationsSeminar.html},
  URLDATE      = {2026-07-06},
}

@article {makkai-generalized-sketches-all,
  AUTHOR       = {Makkai, M.},
  TITLE        = {Generalized sketches as a framework for completeness theorems.  {I}-{III}},
  JOURNAL      = {J. Pure Appl. Algebra},
  FJOURNAL     = {Journal of Pure and Applied Algebra},
  VOLUME       = {115},
  YEAR         = {1997},
  NUMBER       = {1-3},
  PAGES        = {49--79, 179--212, 241--274},
  ISSN         = {0022-4049},
  MRCLASS      = {03G30 (03C95 03F03 18C10)},
  MRNUMBER     = {1429299},
  MRREVIEWER   = {Anna Labella},
  DOI          = {10.1016/S0022-4049(96)00007-2},
  URL          = {https://doi.org/10.1016/S0022-4049(96)00007-2},
}

@article {makkai-reyes-focl,
  AUTHOR       = {Makkai, Michael and Reyes, Gonzalo E.},
  TITLE        = {First order categorical logic},
  SERIES       = {Lecture Notes in Mathematics, Vol. 611},
  NOTE         = {Model-theoretical methods in the theory of topoi and related categories},
  PUBLISHER    = {Springer-Verlag, Berlin-New York},
  YEAR         = {1977},
  PAGES        = {viii+301},
  ISBN         = {3-540-08439-8},
  MRCLASS      = {02J99 (02B25 18C10)},
  MRNUMBER     = {0505486},
  MRREVIEWER   = {C. J. Mikkelsen},
}

@online {makkai-set-theory,
  AUTHOR       = {Makkai, M.},
  TITLE        = {Notes on Set Theory},
  YEAR         = {2000},
  URL          = {https://www.math.mcgill.ca/makkai/},
  URLDATE      = {2026-07-06},
}

@incollection {makkai-ultraproducts,
  AUTHOR       = {Makkai, M.},
  TITLE        = {Ultraproducts and categorical logic},
  BOOKTITLE    = {Methods in mathematical logic ({C}aracas, 1983)},
  SERIES       = {Lecture Notes in Math.},
  VOLUME       = {1130},
  PAGES        = {222--309},
  PUBLISHER    = {Springer, Berlin},
  YEAR         = {1985},
  MRCLASS      = {03G30 (03C20)},
  MRNUMBER     = {799044},
  DOI          = {10.1007/BFb0075314},
  URL          = {https://doi.org/10.1007/BFb0075314},
}

@article {mitchell-boolean-topoi,
  AUTHOR       = {Mitchell, William},
  TITLE        = {Boolean topoi and the theory of sets},
  JOURNAL      = {J. Pure Appl. Algebra},
  FJOURNAL     = {Journal of Pure and Applied Algebra},
  VOLUME       = {2},
  YEAR         = {1972},
  PAGES        = {261--274},
  ISSN         = {0022-4049,1873-1376},
  MRCLASS      = {02K15 (18A15)},
  MRNUMBER     = {319757},
  MRREVIEWER   = {Andreas\ Blass},
  DOI          = {10.1016/0022-4049(72)90006-0},
  URL          = {https://doi.org/10.1016/0022-4049(72)90006-0},
}

@incollection {ml72,
  AUTHOR       = {Martin-L\"{o}f, Per},
  TITLE        = {An intuitionistic theory of types},
  BOOKTITLE    = {Twenty-five years of constructive type theory ({V}enice, 1995)},
  SERIES       = {Oxford Logic Guides},
  VOLUME       = {36},
  PAGES        = {127--172},
  PUBLISHER    = {Oxford Univ. Press, New York},
  YEAR         = {1998},
  MRCLASS      = {03B15 (03F55)},
  MRNUMBER     = {1686864},
}

@article {mostowski-undecidable-statement,
  AUTHOR       = {Mostowski, Andrzej},
  TITLE        = {An undecidable arithmetical statement},
  JOURNAL      = {Fund. Math.},
  FJOURNAL     = {Polska Akademia Nauk. Fundamenta Mathematicae},
  VOLUME       = {36},
  YEAR         = {1949},
  PAGES        = {143--164},
  ISSN         = {0016-2736,1730-6329},
  MRCLASS      = {02.0X},
  MRNUMBER     = {35721},
  MRREVIEWER   = {A.\ Heyting},
  DOI          = {10.4064/fm-36-1-143-164},
  URL          = {https://doi.org/10.4064/fm-36-1-143-164},
}

@phdthesis {nashaat-powell-valued-models,
  AUTHOR       = {Nashaat, Sherif Wagdy},
  YEAR         = {2024},
  TITLE        = {Powell-{V}alued {M}odels for {I}ntuitionistic {S}et {T}heory},
  NOTE         = {Thesis (Ph.D.)--McGill University (Canada)},
  PUBLISHER    = {ProQuest LLC, Ann Arbor, MI},
  PAGES        = {152},
  ISBN         = {979-8304-94692-6},
  MRCLASS      = {Thesis},
  MRNUMBER     = {4967872},
  URL          = {https://www.proquest.com/dissertations-theses/powell-valued-models-intuitionistic-set-theory/docview/3167808555/se-2},
}

@article {osius-categorical-set-theory,
  AUTHOR       = {Osius, Gerhard},
  TITLE        = {Categorical set theory: a characterization of the category of sets},
  JOURNAL      = {J. Pure Appl. Algebra},
  FJOURNAL     = {Journal of Pure and Applied Algebra},
  VOLUME       = {4},
  YEAR         = {1974},
  PAGES        = {79--119},
  ISSN         = {0022-4049},
  MRCLASS      = {18B05 (02K15)},
  MRNUMBER     = {364389},
  MRREVIEWER   = {A. Preller},
  DOI          = {10.1016/0022-4049(74)90032-2},
  URL          = {https://doi.org/10.1016/0022-4049(74)90032-2},
}

@article {powell-extending-godels-negative,
  ISSN         = {00224812},
  URL          = {http://www.jstor.org/stable/2271902},
  AUTHOR       = {William C. Powell},
  JOURNAL      = {The Journal of Symbolic Logic},
  NUMBER       = {2},
  PAGES        = {221--229},
  PUBLISHER    = {Association for Symbolic Logic},
  TITLE        = {Extending Gödel's Negative Interpretation to ZF},
  URLDATE      = {2026-07-06},
  VOLUME       = {40},
  YEAR         = {1975}
}

@article {power-robinson-pie-limits,
  AUTHOR       = {Power, John and Robinson, Edmund},
  TITLE        = {A characterization of pie limits},
  JOURNAL      = {Math. Proc. Cambridge Philos. Soc.},
  FJOURNAL     = {Mathematical Proceedings of the Cambridge Philosophical Society},
  VOLUME       = {110},
  YEAR         = {1991},
  NUMBER       = {1},
  PAGES        = {33--47},
  ISSN         = {0305-0041},
  MRCLASS      = {18D05},
  MRNUMBER     = {1104599},
  MRREVIEWER   = {R. H. Street},
  DOI          = {10.1017/S0305004100070092},
  URL          = {https://doi.org/10.1017/S0305004100070092},
}

@article {shulman-comparing,
  AUTHOR       = {Shulman, Michael},
  TITLE        = {Comparing material and structural set theories},
  JOURNAL      = {Ann. Pure Appl. Logic},
  FJOURNAL     = {Annals of Pure and Applied Logic},
  VOLUME       = {170},
  YEAR         = {2019},
  NUMBER       = {4},
  PAGES        = {465--504},
  ISSN         = {0168-0072,1873-2461},
  MRCLASS      = {18B25 (03E70 03F65 03G30 18B05)},
  MRNUMBER     = {3913162},
  MRREVIEWER   = {Andrzej\ Wi\'{s}nicki},
  DOI          = {10.1016/j.apal.2018.11.002},
  URL          = {https://doi.org/10.1016/j.apal.2018.11.002},
}

@online {shulman-nlab-page,
  AUTHOR       = {Shulman, Mike},
  YEAR         = {2012},
  TITLE        = {Personal page on the nLab},
  URLDATE      = {2026-07-06},
  URL          = {https://ncatlab.org/michaelshulman/show/HomePage},
}

@online {shulman-stack-semantics,
  TITLE        = {Stack semantics and the comparison of material and structural set theories},
  AUTHOR       = {Michael A. Shulman},
  YEAR         = {2010},
  EPRINT       = {1004.3802},
  ARCHIVEPREFIX= {arXiv},
  PRIMARYCLASS = {math.CT},
}

@article {van-den-berg-moerdisjk-aspects-1,
  AUTHOR       = {van den Berg, Benno and Moerdijk, Ieke},
  TITLE        = {Aspects of predicative algebraic set theory. {I}. {E}xact completion},
  JOURNAL      = {Ann. Pure Appl. Logic},
  FJOURNAL     = {Annals of Pure and Applied Logic},
  VOLUME       = {156},
  YEAR         = {2008},
  NUMBER       = {1},
  PAGES        = {123--159},
  ISSN         = {0168-0072,1873-2461},
  MRCLASS      = {03E70 (03G30 18B05)},
  MRNUMBER     = {2474446},
  MRREVIEWER   = {Nicola\ Gambino},
  DOI          = {10.1016/j.apal.2008.06.013},
  URL          = {https://doi.org/10.1016/j.apal.2008.06.013},
}

@article {weber-2-toposes,
  AUTHOR       = {Weber, Mark},
  TITLE        = {Yoneda structures from 2-toposes},
  JOURNAL      = {Appl. Categ. Structures},
  FJOURNAL     = {Applied Categorical Structures. A Journal Devoted to Applications of Categorical Methods in Algebra, Analysis, Order, Topology and Computer Science},
  VOLUME       = {15},
  YEAR         = {2007},
  NUMBER       = {3},
  PAGES        = {259--323},
  ISSN         = {0927-2852},
  MRCLASS      = {18A05 (18A15 18B25 18D05)},
  MRNUMBER     = {2320763},
  MRREVIEWER   = {Stephen Lack},
  DOI          = {10.1007/s10485-007-9079-2},
  URL          = {https://doi.org/10.1007/s10485-007-9079-2},
}

@article {zermelo-grundlagen-der-mengenlehre,
  AUTHOR       = {Zermelo, E.},
  TITLE        = {Untersuchungen \"{u}ber die {G}rundlagen der {M}engenlehre.  {I}},
  JOURNAL      = {Math. Ann.},
  FJOURNAL     = {Mathematische Annalen},
  VOLUME       = {65},
  YEAR         = {1908},
  NUMBER       = {2},
  PAGES        = {261--281},
  ISSN         = {0025-5831,1432-1807},
  MRCLASS      = {DML},
  MRNUMBER     = {1511466},
  DOI          = {10.1007/BF01449999},
  URL          = {https://doi.org/10.1007/BF01449999},
}
\end{document}